\documentclass[a4paper,11pt]{article}
\small
\normalsize

\usepackage[utf8]{inputenc}
\usepackage[T1]{fontenc}
\usepackage[english]{babel}
\usepackage{times}
\usepackage{ulem}

\usepackage[pdftex]{hyperref}

\usepackage{theorem}

\usepackage[all]{xy}

\usepackage{stmaryrd}
\usepackage{graphicx}
\usepackage{paralist}
\usepackage{amsfonts}
\usepackage[leqno]{amsmath}
\usepackage{amssymb,mathrsfs}
\usepackage{color}
\usepackage{times}
\usepackage{enumerate,vmargin}
\usepackage{tocvsec2}
\usepackage{verbatim}
\usepackage[labelfont=bf,font=it]{caption}
\usepackage[numbers,sort&compress]{natbib}

\setmarginsrb{2.9cm}{2.6cm}{2.5cm}{1.7cm}{0cm}{0mm}{0cm}{9mm}

\hypersetup{
    unicode      = false,     % non-Latin characters in AcrobatÃÂ¢ÃÂÃÂs bookmarks
    pdftoolbar   = true,      % show AcrobatÃÂ¢ÃÂÃÂs toolbar?
    pdfmenubar   = true,      % show AcrobatÃÂ¢ÃÂÃÂs menu?
    pdffitwindow = true,      % page fit to window when opened
    pdfnewwindow = true,      % links in new window
    colorlinks   = true,      % false: boxed links; true: colored links
    linkcolor    = blue,      % color of internal links
    citecolor    = red,      % color of links to bibliography
    filecolor    = blue,      % color of file links
    urlcolor     = green       % color of external links
}

\definecolor{vio}{rgb}{.5,.1,.5}
\definecolor{grey}{rgb}{.5,.5,.5}
\definecolor{cnic}{rgb}{.9,0.5,0}
\definecolor{mycolor}{rgb}{.6,0,.8}
\definecolor{darkred}{RGB}{165,42,42}

\def \mm#1{{\color{mycolor}#1\color{black}}}
\def \tt#1{{\color{darkred}#1\color{black}}}
\def \ms#1{{\color{mycolor}\sout{#1}\color{black}}}

\newtheorem{lem}{Lemma}[section]
\newtheorem{thm}[lem]{Theorem}
\newtheorem{prop}[lem]{Proposition}

{\theorembodyfont{\rmfamily}}
{\theorembodyfont{\rmfamily}}
{\theorembodyfont{\rmfamily}}
{\theorembodyfont{\rmfamily}\newtheorem{rem}{Remark}[section]}
{\theorembodyfont{\rmfamily}}

\newcommand{\noi}{\noindent}

\newcommand{\lgeo}{\llbracket}
\newcommand{\rgeo}{\rrbracket}

\newcommand{\un}{{\bf 1}}

\newcommand{\cC}{\mathcal C}

\newcommand{\cE}{\mathcal E}

\newcommand{\cL}{\mathcal L}
\newcommand{\cM}{\mathcal M}

\newcommand{\cP}{\mathcal P}

\newcommand{\cR}{\mathcal R}
\newcommand{\cS}{\mathcal S}

\newcommand{\cV}{\mathcal V}

\newcommand{\bE}{\mathbf E}

\newcommand{\bP}{\mathbf P}

\newcommand{\bm}{\mathbf m}
\newcommand{\bp}{\mathbf p}

\newcommand{\bnu}{\boldsymbol \nu}
\newcommand{\bmu}{\boldsymbol \mu}

\newcommand{\bbU}{\mathbb U}
\newcommand{\bbG}{\mathbb G}
\newcommand{\bbQ}{\mathbb Q}

\newcommand{\bbN}{\mathbb N}

\newcommand{\bbR}{\mathbb R}
\newcommand{\bbT}{\mathbb T}

\newcommand{\ccF}{\mathscr F}
\newcommand{\ccG}{\mathscr G}
\newcommand{\ccE}{\mathscr E}

\newcommand{\elldo}{{\ell}^{_{\, \downarrow}}}

\def\cq{$\hfill \square$}
\def\cqfd{$\hfill \blacksquare$}
\def\ino{ \! \in \! }

\def\bR{\mathbf{R}}
\def\bC{\mathbf{C}}
\def\bD{\mathbf{D}}

\def\bT{\mathbf{T}}

\def\cJ{\mathcal{J}}
\def\cV{\mathscr{V}}
\def\cH{\mathcal{H}}
\def\cK{\mathcal{K}}
\def\ccX{\mathscr{X}}

\def\ccZ{\mathscr{Z}}

\def\bv{\mathtt{v}}
\def\subo{\gamma}

\newcommand{\Ptt}{{\boldsymbol{\Pi}}} 
\newcommand{\Htt}{{\boldsymbol{\mathtt{H}}}} 
\newcommand{\Ytt}{{\boldsymbol{\mathtt{Y}}}} 
\newcommand{\bdelta}{{\boldsymbol{\delta}}} 

\newcommand{\cG}{\boldsymbol{\mathcal G}}
\newcommand{\cT}{\boldsymbol{\mathcal T}}
\newcommand{\bw}{\mathtt w}
\newcommand{\JJ}{J}
\newcommand{\Jtt}{{\mathtt J}}
\newcommand{\bq}{\mathbf{q}}

\def\epp{\varepsilon}
\usepackage[refpage,noprefix]{nomencl}                          %%%
\makeglossary
\makenomenclature                                           %%%
\graphicspath{{.}{figures/}}

\title{ \textsc{Limits of multiplicative inhomogeneous random graphs and L\'evy trees: Limit theorems
}}
\date{}
\author{Nicolas \textsc{Broutin}
\thanks{Sorbonne Universit\'e, Campus Pierre et Marie Curie, 
Case courrier 158, 4 place Jussieu,  
75252 Paris Cedex 05,  
France. Email: nicolas.broutin@upmc.fr} 
\and Thomas \textsc{Duquesne}
\thanks{Sorbonne Universit\'e, Campus Pierre et Marie Curie, 
Case courrier 158, 4 place Jussieu,  
75252 Paris Cedex 05,  
France.
Email: thomas.duquesne@upmc.fr}
\and Minmin \textsc{Wang}
\thanks{University of Sussex, Department of Mathematics, Falmer Campus,  
Brighton, BN1 9QH,  United Kingdom.  
Email: minmin.wang@sussex.ac.uk }
}

\begin{document}

\maketitle

\begin{abstract}
We consider a natural model of inhomogeneous random graphs that extends the classical Erd\H os--R\'enyi graphs and shares a close connection with the multiplicative coalescence, as pointed out by Aldous [\emph{Ann.\ Probab.}, vol.~25, pp.~812--854, 1997]. In this model, the vertices are assigned weights that govern their tendency to form edges. It is by looking at the asymptotic distributions of the masses (sum of the weights) of the connected components of these graphs that Aldous and Limic [\emph{Electron. J. Probab.}, vol.~3, pp.~1--59, 1998] have 
identified the entrance boundary of the multiplicative coalescence, which is intimately related to the excursion lengths of certain  L\'evy-type processes. 
 We, instead, look at the metric structure of these components and prove their Gromov--Hausdorff--Prokhorov convergence to a class of (random) compact measured metric spaces that have been introduced in a companion paper \cite{BDW1}. Our asymptotic regimes relate directly to the general convergence condition appearing in the work of Aldous and Limic. Our techniques provide a unified approach for this general ``critical'' regime, and relies upon two key ingredients: an encoding of the graph by some L\'evy process as well as an embedding of its connected components into Galton--Watson forests. This embedding transfers asymptotically into an embedding of the limit objects into a forest of L\'evy trees, which allows us to give an explicit construction of the limit objects from the excursions of the L\'evy-type process. %\ts{As a consequence of our construction, we give a transparent and explicit condition for the compactness of the limit objects and determine their fractal dimensions.} 
The mains results combined with the ones in the other paper allow us to extend and complement several previous results that had been obtained via model- or regime-specific proofs, for instance: the case of Erd\H os--R\'enyi random graphs obtained by Addario-Berry, Goldschmidt and B. [\emph{Probab. Theory Rel. Fields}, vol. 153, pp. 367--406, 2012], the \emph{asymptotic homogeneous} case as studied by Bhamidi, Sen and Wang [\emph{Probab Theory Rel. Fields}, vol.~169, pp.~565--641, 2017], or the \emph{power-law} case as considered by Bhamidi, Sen and van der Hofstad [\emph{Probab. Theory Rel. Fields}, vol.~170, pp.~387--474, 2018]. 
\end{abstract}

\section{Introduction}
\label{introsec}

\begin{figure}[tp]
\centering
\includegraphics[height = 6cm]{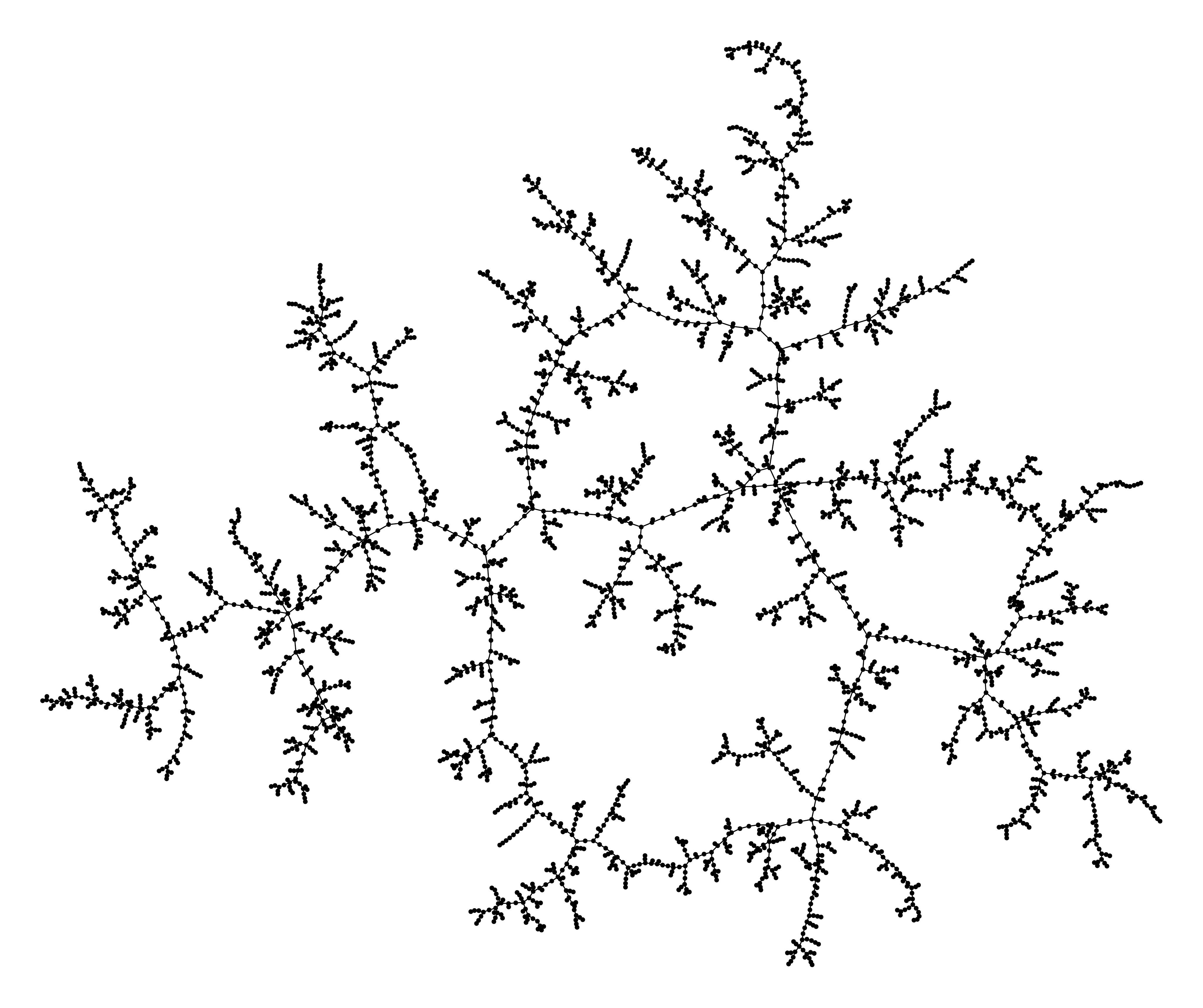}
\includegraphics[height = 6cm]{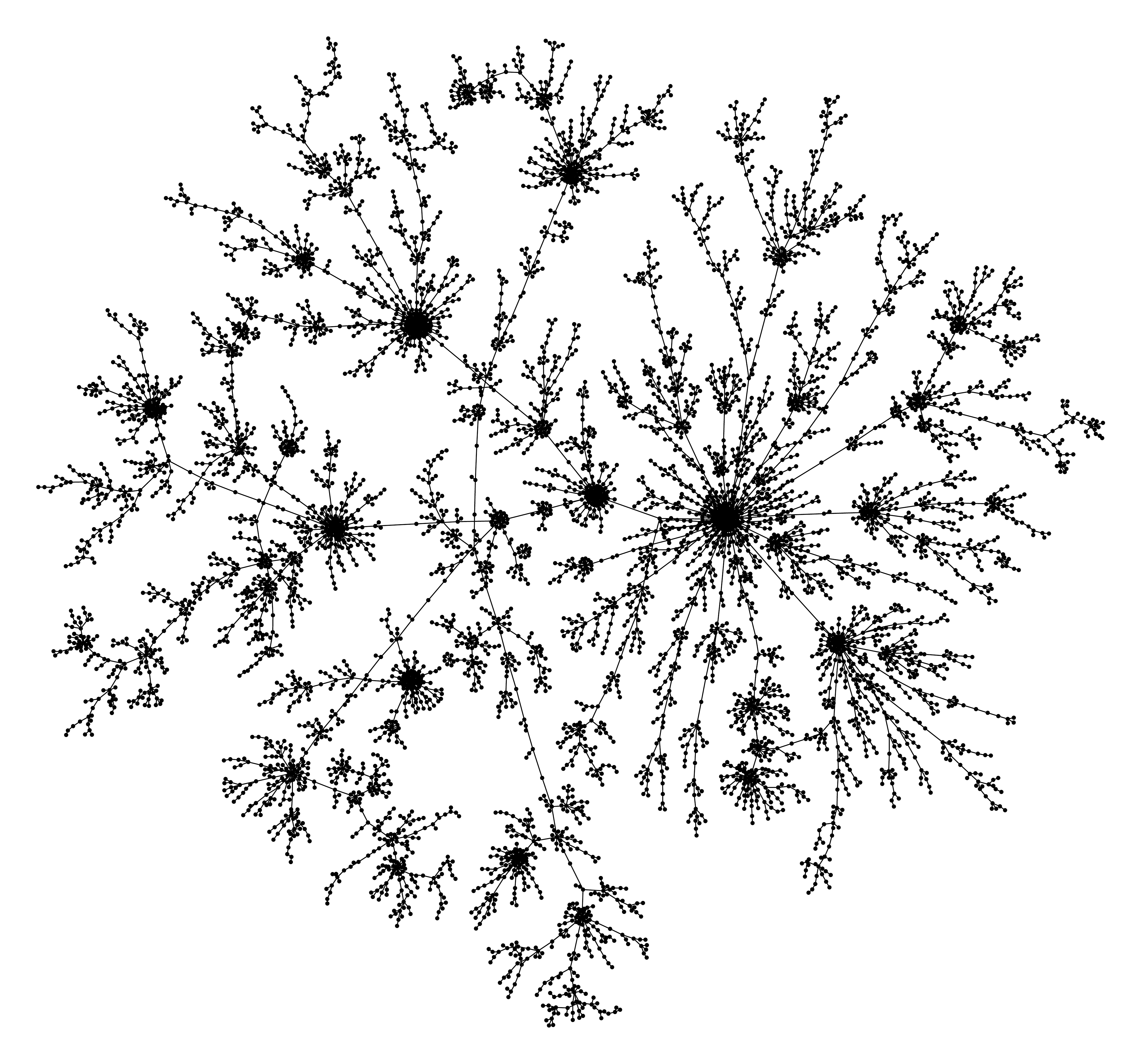}
\caption{Left: a picture of a large connected component of $G(n, p)$. Right: a picture of a large connected component of $\cG_\bw$. Observe the presence of ``hubs'' (nodes of high degrees) in the latter.}
\label{fig:simu}
\end{figure}

\noi
\textbf{Motivation and model. }\label{ssec:motivation} Random graphs have generated a large amount of literature. This is even the case for one single model: the Erd\H os--R\'enyi graph $G(n, p)$ (graph with $n$ vertices connected pairwise in an i.i.d.~way with probability $p\in [0, 1]$). 
Since its introduction by Erd\H os and R\'enyi \cite{ErRe1959} more than fifty years ago, and the discovery of a phase transition where a ``giant connected component'' gets born, the pursuit of a deeper understanding of its structure has never stopped. Many landmark results by Bollob\'as \cite{Bollobas1984}, \L uczak \cite{Luczak1990a}, Janson, Knuth, \L uczak and Pittel \cite{JaKnucPi1993a} have shaped our grasp of this phase transition. From the point of view of precise asymptotics, one of the most important papers is certainly the contribution of Aldous \cite{Al97}, who introduced a stochastic process point of view and paved the way towards the study of scaling limits of critical random graphs. In that paper, he obtained the asymptotics for the sequence of sizes of the connected components of $G(n, p)$ in the so-called critical window where the phase transition actually occurs. His work made possible the construction by Addario-Berry, Goldschmidt and B.~\cite{AdBrGo12a} of the scaling limits of these connected components, seen as metric spaces, which also confirmed the limiting fractal (Brownian) nature.

Following \cite{AdBrGo12a}, the question of identifying the scaling limits has been investigated for more general models of random graphs. Particular attention has been paid to the so-called \emph{inhomogeneous random graphs}, which exhibit heterogeneity in the node degrees and  whose behaviours are often quite different from the Erd\H os--R\'enyi graph. (See Fig.~\ref{fig:simu} for an illustration of this difference). Besides being a theoretical object with intriguing properties, these graphs are also commonly believed to offer more realistic modelling for the complex real-world networks \cite[see, e.g.][]{Newman2003a}. 

In the present work, we consider such an inhomogeneous random graph model that is defined as follows. 
Let $\bw\! =\! (w_1, w_2, \dots, w_n)$ be a sequence of $n$ positive real numbers sorted in nonincreasing order. Interpreting $w_i$ as the propensity of vertex $i$ to form edges, we define a random graph 
$\cG_{\! \bw}$ as follows: the set of its vertices is $\{1, 2, \dots, n\}$, the events $\big\{ \{ i,j\} \, \textrm{is an edge of $\cG_{\! \bw}$} \big\}$, $1 \! \leq \! i \! < \! j\! \leq \! n$, are independent and   
\begin{equation}
\label{defGwdir}
\bP \big( \{ i,j\}  \, \textrm{is an edge of $\cG_{\! \bw}$}\big)=1-\exp \big(\! -\! w_i w_j/\sigma_1(\bw)\big), \quad  \text{ where } \quad \sigma_1(\bw)=w_1+ \ldots + w_n. 
\end{equation}
The graph $\cG_\bw$ extends the classical Erd\H os--R\'enyi random graph in allowing edges to be drawn with non uniform probabilities, while keeping the independence among edges. 

The graph $\cG_{\! \bw}$ has come under different names in the literature, for instance, Poisson random graph in \cite{NoRe06, BhHoLe2010a}, the Norros--Reittu graph in \cite{BhHoLe2010a} or rank-$1$ model in \cite{BhHoLe2012b, BoJaRi07, vdHofstad_book, vdHofstad_book2, BhHoSe15}. 
%\marginpar{\small\mm{refer to ... as ...}}
Here, we will refer to it as the \emph{multiplicative graph} to emphasise its close connection with the \emph{multiplicative coalescent} as 
pointed out by Aldous in \cite{Al97}. This connection is the starting point of the work \cite{AlLi98} of Aldous \& Limic who identify the entrance boundary of multiplicative coalescent by looking at the asymptotic distributions of the \emph{sizes} of the connected components found in $\cG_{\! \bw}$. The asymptotic regime and the limiting processes found in Aldous \& Limic  \cite{AlLi98} lie at the heart of this paper. 
Namely, we extend this result to the \emph{geometry} of the connected components of $\cG_{\! \bw}$ by proving the weak convergence of these connected components as it has been done by  Addario-Berry, Goldschmidt and B.~\cite{AdBrGo12a} for the critical Erd\H os--R\'enyi graphs. 
Our approach relies on the results of a companion paper \cite{BDW1} where we provide a specific coding of $\cG_{\! \bw}$ and an embedding of $\cG_{\! \bw}$ into a Galton-Watson forest,  and where we construct the \textit{continuous multiplicative graphs }that are proven here to be the scaling limits of the discrete models.

More precisely,  
we equip $\cG_{\! \bw}$ with the graph distance $d_{\textrm{gr}}$ and we introduce the weight measure $\bm_\bw\! = \! \sum_{1\leq i\leq n} w_i \delta_{i}$ 
on $\cG_\bw$. The goal of our article can be roughly rephrased as follows: we construct 
a class of (pointed and measured) compact random metric spaces $(\mathbf{G}, d, \bm)$ such that the graphs 
$(\cG_{\bw_n}, \varepsilon_n d_{\textrm{gr}},  \varepsilon^\prime_n \bm_{w_n})$ weakly converge to $ (\mathbf{G}, d, \bm)$ along suitable subsequences $(\bw_n, \varepsilon_n,\varepsilon_n^\prime)$. 
We also prove a similar result 
where $ \bm_{\bw_n}$ is replaced by the counting measure, the limit $\mathbf{G}$ being the same. Of course, here the scaling parameters, 
$\varepsilon_n$ and $ \varepsilon^\prime_n$ go to $0$, so that $\mathbf{G}$ is not discrete. The limits we consider hold in the sense of the weak convergence
corresponding to Gromov--Hausdorff--Prokhorov topology on the space of (isometry classes of) compact metric spaces equipped with finite measures. 
%\fmm{Throughout the document, replaced Prohorov by Prokhorov: \tt{OK fine.}}
%\fmm{Mention the convergence of uniform measures as well? \tt{Done.}}
To achieve the construction of the possible limiting graphs and to prove the convergence of rescaled multiplicative graphs, we rely on two main new ideas: (1) we code multiplicative graphs by processes derived from a LIFO-queue; (2) we embed multiplicative graphs into Galton--Watson trees whose scaling limits are well-understood. 
Before discussing further the connections with previous works and in order to explain the advantages of our approach, let us give a brief but precise overview of our results and of the two 
above mentioned ideas.

\bigskip

\noi
\textbf{Overview of the results. }\label{ssec:overview} Our approach relies first on a specific coding of $\bw$-multiplicative graphs $\cG_{\! \bw}$ via a \textit{LIFO-queue} and a related stochastic process; the queue actually yields an exploration of $\cG_{\! \bw}$ and a spanning tree that encompasses almost all the metric structure of the graph. The LIFO-queue is defined as follows. 

\begin{compactenum}

\smallskip

\item[-] \textit{A single server is visited by $n$ clients labelled by $1, \ldots , n$; }

\smallskip

\item[-] \textit{Client $j$ arrives at time $E_j$ and she/he requests an amount of time of service $w_j$; }

\smallskip

\item[-] \textit{The $E_j$ are independent exponentially distributed r.v.~such that $\bE [E_j] \! =\! \sigma_1 (\bw)/ w_j$; }

\smallskip

\item[-] \textit{A LIFO (last in first out) policy applies: whenever a new client arrives, the server interrupts the service of the current client (if any) and serves the newcomer; when the latter leaves the queue, the server resumes the previous service. }
\end{compactenum}

\smallskip

\noi
As mentioned above, the LIFO-queue yields a tree $\cT_{\!\!\bw}$ whose vertices are the clients: namely, \textit{the server is the root (Client $0$) and Client $j$ is a child of Client $i$ in $\cT_{\! \! \bw}$ if and only if  Client $j$ interrupts the service of Client $i$ (or arrives when the server is idle if $i\! =\!  0$).} We introduce the following.
%\ts{Note that the LIFO-queue is coded by either of the two following processes defined for all $t\! \in \! [0, \infty)$ by: }
\begin{equation}
\label{defYcHH}
 Y^{\bw}_t= \! - t + \!\! \sum_{1\leq i\leq n} \!\! w_i \un_{\{ E_i \leq t\}}, \quad J^\bw_t\! = \! \! \inf_{s\in [0, t]} \!\!  Y^\bw_s \quad \textrm{and} \quad \cH^\bw_t = \# \Big\{ s\in [0, t] : \inf_{r\in [s, t]} Y^\bw_r > Y^\bw_{s-} \Big\} \; .
 \end{equation}
The quantity $Y^\bw_t\! -\! J^\bw_t$ is the \textit{load of the server}, i.e.~the amount of service due at time $t$. We call sometimes $Y^\bw_t$ the \textit{algebraic load} of the server. Note that the LIFO-queue is coded by $Y^\bw$. Then, observe that $\cH^\bw_t$ is the number of clients waiting in the queue at time $t$. We easily see that $\cH^\bw$ is the contour (or the depth-first exploration) of $\cT_{\!\! \bw}$; this entails that the graph-metric of $\cT_{\!\! \bw}$ is entirely encoded by $\cH^\bw$: namely, the distance between the vertices/clients served at times $s$ and $t$ in $\cT_{\!\!\bw}$ is $\cH^\bw_t+\cH^\bw_s \! - \!  2 \!  \min_{r \in [s  \wedge t , s\vee t] }   \cH^\bw_r $. 
%
%The tree $\cT_{\!\! \bw}$ contains most of the metric information of $\cG_{\!\bw}$, but not all. Surplus edges are added to $\cT_{\!\! \bw}$ to obtain $\cG_{\bw}$ as follows: conditionally on $Y^\bw$, 
%let $\sum_{\; 1\leq p\leq \bp_\bw} \!\!  \delta_{(t_p, y_p) }$ be a Poisson point measure on 
%$[0, \infty)\! \times \! [0, \infty)$ with intensity $\frac{{_1}}{{^{\sigma_1 (\bw)}}}  \un_{\{ 0< y <Y^\bw_t -\inf_{[0, t]} Y^\bw \}} \, dt\, dy$ and set 
%$$s_p\! = \! \inf \Big\{ s\ino [0, t_p] :  \inf_{u\in [s, t_p]} (Y^\bw_u \! - \! \inf_{[0, u]}Y^\bw ) > y_p   \Big\} \; .$$ 
%Then we define the set of additional edges $\cS_\bw$ as the set of the edges connecting the clients served at times $s_p$ and $t_p$, for all 
%$1\! \leq \! p \! \leq \! \bp_\bw$. Theorem \ref{loiGG} asserts that the graph obtained by removing the root $0$ from $\cT_{\!\!\bw }$ and adding the edges $\cS_{\bw}$ is distributed as $\cG_\bw$, a $\bw$-multiplicative graph. Namely, 
% $$ \cG_{\bw} \overset{\textrm{(d)}}{=}  (\cT_{\! \bw} \backslash \{ 0\}) \cup \cS_\bw \; .$$ 
%

To get to the graph from the tree $\cT_{\!\! \bw}$, we need to include some surplus edges which are sampled from a Poisson point measure. More precisely, conditional on $Y^\bw$, let 
\begin{equation}
\label{Poissurpl}
\cP_\bw\! = \!\!\! \! \!\! \sum_{\; 1\leq p\leq \bp_\bw} \!\!\!  \delta_{(t_p, y_p) }
 \; \textrm{be a Poisson pt.~meas.~on $[0, \infty)^2$ 
  with intensity $\frac{{_1}}{{^{\sigma_1 (\bw)}}}  \un_{\{ 0< y <Y^\bw_t -\JJ^\bw_t \}} \, dt\, dy$.}  
\end{equation}
Note that a.s.~$\mathbf{p}_\bw\! < \! \infty$, since $Y^\bw\! -\! \JJ^\bw$ is null eventually. We set:
\begin{equation}
\label{pin1}
 \Ptt_\bw= \big( (s_p, t_p)\big)_{1\leq p \leq \bp_\bw} \quad \textrm{where} \quad  s_p\! = \! \inf \big\{ s\ino [0, t_p] :  \inf_{u\in [s, t_p]} Y^\bw_u \! - \! \JJ^\bw_{u} > y_p   \big\}, \; 1\! \leq \! p \! \leq \! \bp_\bw \; .
\end{equation}
Next, we define the set of additional edges $\cS_\bw$ as the set of the edges connecting the clients served at times $s_p$ and $t_p$, for all 
$1\! \leq \! p \! \leq \! \bp_\bw$ and we then define the graph $\cG_{\bw}$ by %\ts{setting}
%\[
%\cV(\cG_{\bw})=\{1, 2, \dots, n\} \quad \text{ and } \quad \ccE(\cG_{\bw}) = \cS_{\bw}\cup \big\{\{i, j\}\in \ccE(\cT_{\bw}): i, j\ge 1\big\}\,. 
%\]
$$ \cG_{\bw}
% \overset{\textrm{(d)}}{=} 
:= (\cT_{\! \bw} \backslash \{ 0\}) \cup \cS_\bw \; .$$ 
Namely, $\cG_{\bw}$ is the graph obtained by removing the root $0$ from $\cT_{\!\!\bw }$ and adding the edges in $\cS_{\bw}$. The following is proved in the companion paper \cite{BDW1}. 
\begin{thm}[Theorem \tt{2.1} in \cite{BDW1}]
\label{loiGG} $\cG_{\bw}$ is distributed as a $\bw$-multiplicative random graph as specified in (\ref{defGwdir}). 
\end{thm}

From this representation of the discrete graphs, one expects that if $Y^\bw$ converges, then the graph should also converge, at least in a weak sense. 
However, since $Y^\bw$ is not Markovian, it is difficult to obtain a limit for the local-time functional $\cH^\bw$, which is the function that encodes the metric. 
To circumvent this technical difficulty, we embed the non-Markovian LIFO-queue governed by $Y^\bw$ into a Markovian one that is defined as follows. 

\smallskip

\begin{compactenum}
\item[-] \textit{A single server successively receives an infinite number of clients; }

\smallskip

\item[-] \textit{A LIFO policy applies; }

\smallskip

\item[-]\textit{Clients arrive at unit rate; }

\smallskip

\item[-] \textit{Each client has a type that is an integer ranging in $\{ 1, \ldots, n\}$; the amount of service required by a client of type $j$ is $w_j$; types are i.i.d.~with law $\nu_\bw \! = \! \frac{1}{\sigma_1 (\bw)}\sum_{1 \leq j\leq n} w_j \delta_{j} $. }
\end{compactenum}

\smallskip

\noi
Namely, let $\tau_k$ be the arrival-time of the $k$-th client and let $\Jtt_k$ be the type of the $k$-th client; then the Markovian LIFO queueing system is entirely characterised by  $\sum_{k\geq 1} \delta_{(\tau_k , \Jtt_k)}$ that is a Poisson point measure on $[0, \infty) \! \times \! \{ 1, \ldots, n\}$ with intensity $\ell  \otimes  \nu_\bw$, where $\ell $ stands for the Lebesgue measure on $[0, \infty)$.
To simplify the explanation of the main ideas, \textit{we concentrate in this Overview only
on the (sub)critical cases} where the Markovian queue is recurrent, which amounts to assume that 
$$\sigma_2 (\bw) \! \leq \! \sigma_1 (\bw)\; .$$ 
Here, for all $r \! \in \! (0, \infty)$, we use the notation 
$\sigma_{r} (\bw)\! = \! \sum_{1\leq j\leq n} w_j^r$.

The Markovian queue yields a tree $\bT_{\! \bw}$ that is defined as follows: \textit{the server is the root of $\bT_{\! \bw}$ and the $k$-th client to enter the queue is a child of the $l$-th one if the $k$-th client enters when the $l$-th client is being served. }
One easily checks that $\bT_{\! \bw} $ is a sequence of i.i.d.~Galton--Watson trees glued at their root and that 
their common offspring distribution is 
\begin{equation}
\label{muGWdef}
\mu_\bw (k) \! =\!  \sum_{1\leq j\leq n} \frac{w_j^{k+1}}{\sigma_1 (\bw) k!}e^{-w_j}, \quad k\! \in \! \bbN. 
\end{equation} 
Observe that $\sum_{k\in \bbN} k\mu_\bw (k)\! = \! \sigma_2 (\bw)/ \sigma_1 (\bw) \! \leq \! 1$, which implies that the GW-trees are finite.  
The tree $\bT_{\! \bw}$ is then coded by its contour process $(H^\bw_t)_{t\in [0, \infty)}$: namely, $H^\bw_t$ stands for the number of clients waiting in the Markovian queue at time $t$ and it is given by 
 \begin{equation}
 \label{HXdisdef}
 H^\bw_t = \# \Big\{ 
 s\in [0, t] : \inf_{r\in [s, t]} X^\bw_r > X^\bw_{s-} \Big\} \quad \textrm{where} \quad  X^\bw_{t}  =  -t + \sum_{k\geq 1} w_{\Jtt_k}\un_{[0, t]} (\tau_k), \; t \! \in \! [0, \infty), 
 \end{equation}
is the (algebraic) load of the Markovian server. These definitions make sense in the supercritical cases. Note that $X^\bw$ is a spectrally positive L\'evy process with initial value $0$; it is characterised by 
its Laplace exponent defined by $\bE[ e^{-\lambda X^\bw_t}]\! = \! e^{t\psi_\bw (\lambda)}$, for $t, \lambda\! \in \! [0, \infty)$, that is explicitly given by:  
$$\psi_\bw (\lambda)  =  \alpha_\bw \lambda + \!\!\!\!   \sum_{1\leq j\leq n}\!  \! \frac{_{w_j} }{^{\sigma_1 (\bw)}} \big( e^{-\lambda w_j}\! -\! 1\! + \! \lambda w_j \big) \quad \textrm{and} \quad \alpha_\bw \! := \! 1\! -\! \frac{_{\sigma_2 (\bw)}}{^{\sigma_1 (\bw)}} \; .$$  

From this tractable model, we derive the LIFO-queue and the tree $\cT_{\! \! \bw}$ governed by $Y^\bw$ by a time-change that ``skips'' some time intervals and that is defined as follows. We colour in \textit{blue} or \textit{red} the clients of the Markovian queue in the following recursive way:  
\begin{compactenum} 
\item[\textit{(i)}] \textit{if the type $\Jtt_k$ of the $k$-th client already 
appeared among the types of the blue clients who previously entered the queue, then the $k$-th client is red;}
\item[\textit{(ii)}] \textit{  
otherwise the $k$-th client inherits her/his colour from the colour of the client who is 
currently served when she/he arrives (and this colour is blue if there is no client served when she/he arrives: namely, we consider that the server is blue).} 
\end{compactenum}
Note that 
a client who is the first arriving of her/his type is not necessarily coloured in blue. We easily check that 
exactly $n$ clients are coloured in blue and their types are necessarily distinct. Moreover, while 
a blue client is served, note that the other clients waiting in the line (if any) are blue too. Actually, the sub-queue of blue clients corresponds to the previous LIFO queue governed by $Y^\bw$. 
More precisely, we set   
$$\mathtt{Blue}\! = \! \big\{t\ino [0, \infty)\! : \! \textrm{a blue client is served at $t$} \big\} \quad \textrm{and} \quad 
\theta^{\mathtt{b}, \bw}_t \! \! = \! \inf \Big\{ s\ino [0, \infty)\!  : \!  \!  \int_0^s \!\!\!  \un_{\mathtt{Blue}} (u)  du  >\! t  \Big\}.$$  
We refer to (\ref{thetabw}) in Section \ref{redbldissec} for a precise definition of $\theta^{\mathtt{b}, \bw}$. Then, 
$$ (Y^\bw_t, \cH^\bw_t)_{t\in [0, \infty)}   = \big( X^\bw_{\theta^{\mathtt{b}, \bw}_t} , H^\bw_{\theta^{\mathtt{b}, \bw}_t}  \big)_{t\in [0 ,\infty)}   \; .$$
We refer to Proposition \ref{Xwfrombr} and Lemma \ref{Hthetalem} in 
Section \ref{redbldissec} for a more precise statement of this equality. 
This explains how to code $\cG_{\! \bw}$ in terms of the two tractable processes $X^\bw$ and $H^\bw$ derived from the Markovian queue.

  Such Markovian queues and their coding processes $(X^\bw, H^\bw)$ 
 have analogues in the continous time and space setting. In our context, the parameters governing such processes are those identified by  Aldous \& Limic  \cite{AlLi98} 
for the eternal multiplicative coalescent. Namely: 
%\begin{equation}
%\label{parconing}
%\alpha \ino \bbR , \quad \beta \ino [0, \infty) , \quad \kappa \ino (0, \infty) , \quad \mathbf{c}\! = \! (c_j)_{j\geq 1} \ino \elldo_3 \; .
%\end{equation}
\begin{equation}
\label{parconing}
 \alpha \in \bbR, \; \beta  \in [0, \infty),  \; \kappa\in(0, \infty) \quad \textrm{and} \quad \mathbf{c}\! = \! (c_j)_{j\geq 1} \; \textrm{decreasing and such that} \;  \sum_{j\geq 1} c_j^3<\infty \, . 
 \end{equation}
The load of service of the continuous analogue of the Markovian queue is a spectrally positive L\'evy process $(X_t)_{t\in [0, \infty)}$ starting at $X_0\! = \! 0$ whose Laplace exponent $\psi$ is given by 
\begin{equation}\label{granouff}
\tfrac{1}{t}\log \big( \bE \big[ e^{-\lambda X_t}\big]\big)\! := \psi(\lambda)  \! = \!   
 \alpha \lambda +\frac{_{_1}}{^{^2}} \beta \lambda^2 +  \! \sum_{j\geq 1}   \kappa c_j   \big( e^{-\lambda c_j}\! -\! 1\! + \! \lambda c_j \big), \quad t, \lambda\ino [0, \infty).
\end{equation}  
To simplify, we restrict our explanations to the cases where $X$ does not drift to $\infty$, which is equivalent to assuming that $\alpha \! \in \! [0 , \infty)$. 
The tree corresponding to the clients of the continuous analogue of the Markovian queue that is driven by $X$, is actually the L\'evy tree yielded by $X$, which is defined through its contour process as introduced by Le Gall \& Le Jan \cite{LGLJ98}. 
To that end, we assume that $\psi$ (as defined in (\ref{granouff})) satisfies the following: 
\begin{equation}\label{Grey}
\int^\infty \frac{d\lambda}{\psi(\lambda)}<\infty, 
\end{equation}
which implies that either $\sum_j c_j^2=\infty$ or $\beta \! \neq \! 0$; therefore $X$ has infinite variation sample paths. 
Under Assumption (\ref{Grey}), Le Gall \& Le Jan \cite{LGLJ98} (see also Le Gall \& D.~\cite{DuLG02}) prove that there exists a continuous process $(H_t)_{t\in [0, \infty)}$ such that the following limit holds true for all $t\ino [0, \infty)$ in probability: 
\begin{equation}
\label{llapproHdef}
H_t = \lim_{\varepsilon \rightarrow 0} \frac{1}{\varepsilon} \!  \int_0^{t} \!  \un_{\{  X_s - \inf_{r\in [s, t]} X_r  \leq \varepsilon\}} \, ds \; . 
\end{equation}
We explain further how to make sense of this definition in the supercritical cases. The process $H$ is called the \textit{height process associated with $X$} and the processes $(X, H)$ are the continuous analogues of $(X^\bw, H^\bw)$. 

We explain in Section \ref{hautcontset} how to colour the Markovian queue driven by $X$: namely, we explain how to define a right-continuous increasing time-change 
$(\theta^{\mathtt{b}}_t)_{t\in [0, \infty)}$ that is the analogue of the discrete one $\theta^{ \mathtt{b}, \bw}$. We refer to (\ref{thetabdef}) in Section \ref{hautcontset} for a formal definition of $\theta^{\mathtt{b}}$. 
Then we define the c\`adl\`ag process 
\begin{equation}\label{Ydef1}
Y_t\! = \! X_{\theta^{\mathtt{b}}_t}, \quad t\ino [0, \infty) ,
\end{equation}
that represents the load driving the analogue of the LIFO-queue (without repetitions). As we will see in (\ref{Ydef}) Section \ref{hautcontset}, $Y$ can be written under the following form: 
\begin{equation}\label{plasouilloc}
\forall t \in [0, \infty) , \quad Y_t=  -\alpha t \! -\! \frac{_1}{^2}\kappa \beta t^2 + \sqrt{\beta} B_t +  \sum_{j\geq 1} \! c_j (\un_{\{ E_j \leq t \}} \! -\! c_j \kappa t), 
\end{equation}
where $(B_t)_{t\in [0, \infty)}$ is a standard linear Brownian motion starting at $0$ and where the $E_j$ are independent exponentially distributed r.v.~that are independent from $B$ and such that 
$\bE [E_j]\! = \! (\kappa c_j)^{-1}$. The sum in (\ref{plasouilloc}), as it is, is informal: it has to be understood in the sense of $L^2$ semimartingales (see Section \ref{hautcontset} for a precise explanation). The latter expression of $Y$ can be found in Aldous \& Limic  \cite{AlLi98} who proved that the lengths of the excursions of $Y$ above its infimum (ranked in decreasing order) are distributed as the multiplicative coalescent. We refer to Theorem \ref{Xdefthm} in Section \ref{hautcontset} for a precise statement of (\ref{Ydef1}). 

As it is proved in Theorem 2.6 in \cite{BDW1} (that is 
recalled in Theorem \ref{cHdefthm}, Section \ref{hautcontset}) , there exists a continuous process $(\cH_t)_{t\in [0, \infty)}$ 
that is an adapted functional of $Y$ such that 
\begin{equation}
\label{defcH}
 \forall t \in [0, \infty), \quad \cH_t = H_{\theta^{\mathtt{b}}_t} \; .
 \end{equation}
Here, $\cH$ is a.s.~a continuous process that is called the \textit{height process associated with $Y$} and we claim that $(Y,\cH)$ is the continuous analogue of $(Y^\bw, \cH^\bw)$, as justified by limit theorems stated further. 

As proved in  \cite{BDW1} (and recalled in Lemma \ref{AHeuer}, Section \ref{hautcontset}), the excursion intervals of $\cH$ above $0$ and the excursion intervals of $Y$ above its infimum are 
the same. Moreover, Proposition 14 in Aldous \& Limic \cite{AlLi98} (that is recalled in Proposition \ref{AldLim1}, Section \ref{hautcontset}) asserts that these excursions can be indexed in the decreasing order of their lengths. Namely, 
\begin{equation}
\label{excuHY}
\big\{  t\ino [0, \infty) : \cH_t >0  \big\}= \Big\{  t\ino [0, \infty) : Y_t > \inf_{[0, t]} Y  \Big\}= \bigcup_{k\geq 1} (l_k , r_k) 
\end{equation}
where the sequence $\zeta_k\! = \! l_k \! -\! r_k$, decreases. The continuous analogue of $\cG_{\!\bw}$ is derived from $(Y, \cH)$ as follows:
first, for all $s, t \! \in \! [0, \infty)$, we define the usual tree 
pseudometric associated with $\cH$: $ d_\cH (s, t)= \cH_s + \cH_t -2\min_{u\in [s\wedge t, s\vee t]} \cH_u $. 
Then, we set 
\begin{equation}
\label{Jdef}
\forall t\ino [0, \infty), \quad  
\JJ_t \! = \! \inf_{s\in [0, t]} Y_s \; 
\end{equation} 
and conditionally given $Y$, let 
\begin{equation}
\label{Poisurcon}
\cP\! = \! \sum_{\; p\geq 1} \!  \delta_{(t_p, y_p) }
 \; \textrm{be a Poisson pt.~meas.~on $[0, \infty)^2$ 
 with intensity $ \kappa\un_{\{ 0< y <Y_t -J_t \}} \, dt\, dy$.} 
\end{equation}
%\mm{Comment: In fact, the proof of Theorem \ref{HYcvth} on p.\pageref{thm} indicates that the intensity here should be $\kappa \un_{\{ 0< y <Y_t -J_t \}} \, dt\, dy$.}
Then, we set  
\begin{equation}
\label{pinchset}
\Ptt\! = \! \big( (s_p, t_p) \big)_{p\geq 1} \quad \textrm{where} \quad s_p\! = \! \inf \big\{ s\ino [0, t_p] :  \inf_{u\in [s, t_p]} 
Y_u\! -\! J_u > y_p   \big\}, \; \,  p\! \geq \! 1. 
\end{equation}
Here $\Ptt$ plays the role of $\Ptt_\bw$. 
%%%We claim that the processes $(Y, \cH, \Ptt)$ completely characterise the continuous version of the multiplicative graph 
%as explained in the next section. 
%
%
%
%
%Then, we introduce 
%$\sum_{ p \geq 1}\delta_{(t_p, y_p) }$ distributed (\mm{conditionally} on $Y$) as 
%a Poisson point measure on $[0, \infty)\times [0, \infty)$ with intensity $\kappa \un_{\{ 0< y <Y_t -\inf_{[0, t]} Y\}} \, dt\, dy$. 
%Next for all $p\! \geq \! 1$, we set 
%$$s_p\! = \! \inf \Big\{ s\ino [0, t_p] :  \inf_{u\in [s, t_p]} (Y_u \! -\!  \inf_{[0, u]} Y) \,  > \! y_p   \Big\} \; .$$
Fix $k\! \geq \! 1$. One can prove that if $t_p \in [l_k, r_k]$, then $s_p \! \in \! [l_k, r_k]$.
We define $\mathbf{G}_k$ as the set $[l_k, r_k]$ where we have identified points $s, t\in [l_k, r_k]$ such that either 
$d_\cH (s,t)\! = \! 0$ or $(s,t) \! \in \! \{ (s_p, t_p) ; p\! \geq \! 1 :  t_p \! \in \! [l_k, r_k] \}$. It actually yields a metric denoted by $\mathrm{d}_{k}$, on $\mathbf{G}_k$; note that $l_k$ and $r_k$ are identified 
and we denote by $\varrho_{k}$ the corresponding point in $\mathbf{G}_k$; we denote by 
$\bm_{k}$  the measure induced by the Lebesgue measure on $[l_k, r_k]$. The continuous analogue of $\cG_{\! \bw}$ is then the sequence of pointed measured compact metric spaces 
\begin{equation}
\label{contGdef}
\mathbf{G}\! = \! \big( (\mathbf{G}_k, \mathrm{d}_{k},\varrho_{k} ,\bm_{k})\big)_{k\geq 1}\; ,
\end{equation}
that is called the \textit{$(\alpha, \beta, \kappa, \mathbf{c})$-continuous multiplicative graph}. 
We refer to Section \ref{ccvmulsec} (and more specifically see (\ref{redefGmu})) for a more precise definition.

As already mentioned, the main goal of the paper is to prove that $\mathbf{G}$ is the scaling limit of sequences of rescaled discrete graphs $\cG_{\! \bw_n}$ for a suitable sequence of weights with finite support $\bw_n\! = \! (w^{_{(n)}}_{^j})_{j\geq 1}$ that are listed 
in the nonincreasing order: namely, $w^{_{(n)}}_{^j} \! \geq \! w^{_{(n)}}_{^{j+1}}$, and  $w^{_{(n)}}_{^j} \! = \! 0$ for all sufficiently large $j$.
Here, we first set 
\begin{equation}
\label{jjnndef}
\mathbf{j}_n \! := \! \sup \big\{ j\! \geq \! 1 \! : \! w^{_{(n)}}_{^j}\! >\! 0 \big\} < \infty\,.
\end{equation}
We don't require that $\mathbf{j}_n$ is equal to $n$ but we want $\lim_{n\rightarrow \infty} 
\mathbf{j}_n\! = \! \infty$. 
Our main result (Theorem \ref{HYcvth} in Section \ref{limithouf}) asserts the following.  

\smallskip

\begin{compactenum}
\item[]\textit{If the 
Markovian processes $(X^{\bw_n}, H^{\bw_n})$, properly 
rescaled in time and space, weakly converge to $(X, H)$, then $(Y^{\bw_n}, \cH^{\bw_n})$ converges weakly to $(Y, \cH)$  with the same scaling.} 
\end{compactenum}

\smallskip

\noi
More precisely, the graphs $\cG_{\! \bw_n}$, or their coding functions, are rescaled 
by two factors $a_n$ and $b_n$ tending to $\infty$; $a_n$ is a weight factor and $b_n$ is an exploration-time factor. Namely, the rescaled processes to consider are 
$\frac{_1}{^{a_n}}X^{_{\bw_n}}_{^{b_n \cdot}}$ (or $\frac{_1}{^{a_n}}Y^{_{\bw_n}}_{^{b_n \cdot}}$) and it is natural to require a priori that 
$b_n \! = \! O(a_n^2)$ by standard results on L\'evy processes. 
Moreover, if the largest weight "persists" in the limit, then $a_n \asymp w^{_{(n)}}_1$ and in general $w^{_{(n)}}_1\! = \! O (a_n)$. 
In the limit, if two large weights persist, they cannot fuse and they tend not to be connected by an edge. 
Namely, if the two largest weights persist, then 
$1\! -\! \exp ( - w^{_{(n)}}_{^1}w^{_{(n)}}_{^2}\! / \sigma_1 (\bw_n)) \! \rightarrow \! 0$ 
and since $w^{_{(n)}}_{^1}  \asymp w^{_{(n)}}_{^2}  \asymp  a_n$, it entails 
$\lim_{n\rightarrow \infty} a_n^2/ \sigma_1 (\bw_n) \! = \! 0$. 
Next, since $b_n$ is an exploration-time factor, we require that $b_n \asymp \bE [C_n]$, where $C_n$ stands for the number of clients who are served before the arrival of Client $1$ (i.e.~the client corresponding to the largerst weight $w^{_{(n)}}_{^1}$) 
in the $\bw_n$-LIFO queue coding $\cG_{\! \bw_n}$. Let us denote by 
$D_n$ the sum of the weights of the vertices explored before visiting Client 1. It is easy to see that $ \bE [ C_n ] \! = \! \sum_{j\geq 2} w^{_{(n)}}_{^j}/ (w^{_{(n)}}_{^j}+w^{_{(n)}}_{^1})  $ and that 
$ \bE [ D_n ] =  \sum_{j\geq 2} (w^{_{(n)}}_{^j})^2/(w^{_{(n)}}_{^j}+w^{_{(n)}}_{^1}) $. So, when $w_{^1}^{_{(n)}}$ persists, we get 
$ \sigma_1 (\bw_n) \asymp  a_n \bE [C_n ]$ and $ \sigma_2 (\bw_n) \asymp  a_n \bE [D_n ]$. Moreover, in the asymptotic regime that we consider, we require that the number of visited vertices has to be of the same order of magnitude as the sum of the corresponding 
weights: namely, $\bE [C_n] \asymp \bE [D_n]$, which corresponds to the criticality assumption: $\sigma_1(\bw_n) \asymp \sigma_2 (\bw_n)$ that also implies 
$a_nb_n \asymp  \sigma_1 (\bw_n) $. 
These constraints amount to assuming the following a priori estimates: 
\begin{equation}
\label{apriori}
 \lim_{n\rightarrow \infty} a_n \! = \! \lim_{n\rightarrow \infty} \frac{b_n}{a_n} \! = \! \infty, \quad  \lim_{n\rightarrow \infty}\frac{b_n}{a^2_n}\! =  :\!  \beta_0 \ino [0, \infty), \quad  w^{_{(n)}}_1\! = \! O (a_n) , \quad \lim_{n\rightarrow \infty}\frac{a_nb_n}{\sigma_1 (\bw_n)}= \kappa. 
\end{equation} 
Note here that $\beta_0\! = \! 0$ possibly. Then, a more precise statement of Theorem \ref{HYcvth} is: \textit{If $(a_n, b_n, \bw_n)$ satisfies (\ref{apriori}), }
\begin{equation}
\label{glurniglup}
\textit{and if} \quad \big( \frac{_{_1}}{^{^{a_n}}} X^{\bw_n}_{b_n \cdot }\,  , \frac{_{_{a_n}}}{^{^{b_n}}} H^{\bw_n}_{b_n \cdot } \big)\;  
\underset{n\rightarrow \infty}{-\!\!\! -\!\!\! -\!\!\! \longrightarrow} \; \big( X,  H \big) 
\end{equation}
\textit{weakly on $\bD([0, \infty), \bbR)\! \times \! \bC ([0, \infty) , \bbR)$ equipped with the product of the Skorokhod and the continuous topologies, then the following joint convergence }
 \begin{equation}
\label{hurkoko}
 \big( \frac{_{_1}}{^{^{a_n}}} X^{\bw_n}_{b_n \cdot }\,  , \frac{_{_{a_n}}}{^{^{b_n}}} H^{\bw_n}_{b_n \cdot } \, ,  
 \big( \frac{_{_1}}{^{^{b_n}}} \theta^{\mathtt{b}, \bw_n}_{b_n \cdot } , \frac{_{_1}}{^{^{a_n}}} Y^{\bw_n}_{b_n \cdot } \big),  \frac{_{_{a_n}}}{^{^{b_n}}} \cH^{\bw_n}_{b_n \cdot}  \big) 
\;  \underset{n\rightarrow \infty}{-\!\!\! -\!\!\! -\!\!\! \longrightarrow} \; \big( X, H, (\theta^\mathtt{b}, Y), \cH \big) 
\end{equation}
\textit{holds weakly on $\bD([0, \infty), \bbR)\! \times \! \bC ([0, \infty) , \bbR)\! \times \!  \bD([0, \infty), \bbR^2)\! \times \! 
\bC ([0, \infty) , \bbR) $ equipped with the product topology.} 

Necessary and sufficient conditions on the $(a_n, b_n, \bw_n)$ for (\ref{glurniglup}) to hold can be derived  from previous results due to Le Gall \& D.~\cite{DuLG02} (let us mention it is not direct: see Proposition \ref{HMarkcvprop}). Namely, suppose that $(a_n, b_n, \bw_n)$ satisfy (\ref{apriori}); then (\ref{glurniglup}) holds if and only if  the following condition are satisfied
\begin{equation}
\label{hurkiki}
 (A): \; \,  \frac{_{_1}}{^{^{a_n}}} X^{\bw_n}_{b_n } \! \overset{\textrm{(weakly)}}{\underset{n\rightarrow \infty}{-\!\!\! -\!\!\! -\!\!\! \longrightarrow}} \! X_1 \quad \textrm{and} \quad  (B): \quad 
 \exists \,\delta \in \! (0, \infty) , \quad \liminf_{n\rightarrow \infty} \bP \big( Z^{\bw_n}_{\lfloor b_n \delta /a_n \rfloor} = 0 \big) >0 \; 
 \end{equation}
where $(Z^{\bw_n}_k)_{k\in \bbN}$ stands for a Galton--Watson Markov chain with offspring distribution $\mu_{\bw_n}$ given by (\ref{muGWdef}) and with initial state $Z^{\bw_n}_0\! = \! \lfloor a_n \rfloor$. Let us mention that 
Proposition \ref{Hcritos} shows that for all $\alpha\! \in \! \bbR$, $\beta \! \in \! [0, \infty)$, $\beta_0 \! \in \! [0, \beta]$, $\kappa\in (0, \infty)$ and $\mathbf{c}$ such that $\sum_{j\geq 1} c_j^3 \! <\!  \infty$ and such that Grey's condition \eqref{Grey} is satisfied, there exists a sequence $(a_n, b_n, \bw_n)_{n\in \bbN}$ satisfying (\ref{apriori}) and (\ref{hurkiki}), so that (\ref{hurkoko}) holds. Proposition~\ref{Hcritos} 
also shows that in 
(\ref{hurkiki}), $(A)$ does not necessarily imply $(B)$. 
Moreover, Proposition~\ref{Hcritos} also provides a more tractable  condition that implies $(B)$ in (\ref{hurkiki}) and that is satisfied in all the examples that have been considered previously.  

By soft arguments (see Lemma \ref{codconGHP}), the convergence (\ref{hurkoko}) of the coding functions implies that the rescaled sequence of graphs $\cG_{\bw_n}$ converges, as random metric spaces. As already mentioned, 
the convergence holds weakly on the space $\bbG$ of (pointed and measure preserving) isometry classes of pointed measured compact metric spaces endowed with the Gromov--Hausdorff--Prokhorov distance (whose definition is recalled in (\ref{defGHP}) in Section \ref{ccvmulsec}). Actually, the convergence holds jointly for the connected components of $\cG_{\bw_n}$: namely, equip $\cG_{\bw_n}$ with the weight-measure $\bm^{\bw_n}\! = \! \sum_{j\geq 1} w^{_{(n)}}_{^j} \delta_j$; let $\mathbf{q}_{\bw_n}$ be the number of connected components of $\cG_{\bw_n}$; we index these connected component  
$(\cG_{\! k}^{\bw_n})_{1\leq k \leq \mathbf{q}_{\bw_n}}$ in the decreasing order of their $\bm^{\bw_n}$-measure: namely, 
\begin{equation}
\label{order1}
\bm^{\bw_n} (\cG_{\! 1}^{\bw_n}) \! \geq \! \ldots  \! \geq \! 
\bm^{\bw_n} (\cG_{\! \mathbf{q}_{\bw_n} }^{\bw_n}) .
\end{equation}
For the sake of convenience, we complete this finite sequence of connected components by 
point graphs with null measure to get an infinite sequence of $\bbG$-valued 
r.v.~$\big( (\cG_{\! k}^{\bw_n}, d_{k}^{\bw_n}, \varrho_k^{\bw_n}, \bm_k^{\bw_n})\big)_{k\geq 1}$, where $d_{k}^{\bw_n}$ stands for the graph-metric on $\cG^{\bw_n}_{\! k}$, where 
$\varrho_k^{\bw_n}$ is the first vertex/client of 
$\cG^{\bw_n}_{\! k}$ who enters the queue and where $\bm_k^{\bw_n}$ is the restriction of  $\bm_{\bw_n}$ to $\cG_{\! k}^{\bw_n}$. Then, Theorem \ref{graphcvth} asserts the following: \textit{ If $(a_n, b_n, \bw_n)$ satisfy (\ref{apriori}) and (\ref{glurniglup}), then} 
\begin{equation}
\label{hurkuku}
\big(\big( \cG_{\! k}^{\bw_n} ,  \frac{_{_{a_n}}}{^{^{b_n}}}d_{k}^{\bw_n} , \varrho_k^{\bw_n}, \frac{_{_{1}}}{^{^{b_n}}}\bm_k^{\bw_n}  \big) \big)_{k\geq 1}
\;  \underset{n\rightarrow \infty}{-\!\!\! -\!\!\! -\!\!\! \longrightarrow} \; \big(\big( \mathbf{G}_{k} , \mathrm{d}_{k}, \varrho_{k} , \bm_{k} \big) \big)_{k\geq 1}
\end{equation} 
\textit{holds weakly on $\bbG^{\bbN^*}$ equipped with the product topology.} Moreover, Theorem \ref{graphcvth} also asserts first that we can replace in (\ref{hurkuku}) the weight-measure $\bm^{\bw_n}$ by the counting measure $\# \! = \! \sum_{1\leq j \leq \mathbf{j}_n} \delta_j$, where 
$\mathbf{j}_n \! := \! \sup \{ j\! \geq \! 1 \! : \! w^{_{(n)}}_{^j}\! >\! 0\}$, and it also asserts that under the additional assumption $\sqrt{\mathbf{j}_n}/ b_n \! \rightarrow \! 0$, the connected components can be listed in the decreasing order of their number of vertices: namely,  
\begin{equation}
\label{order2}
\# (\cG_{\! 1}^{\bw_n}) \! \geq \! \ldots  \! \geq \! 
\#  (\cG_{\! \mathbf{q}_{\bw_n} }^{\bw_n}). 
\end{equation}

\noi
\textbf{Discussion} \label{ssec:discussion} We now briefly 
discuss connections to other works. We refer to Section \ref{prevresults} for more detailed comments on related papers.  

\medskip

\noindent \textit{A unified and exhaustive treatment of the limiting regimes:} While important progress has been made on the Gromov--Hausdorff scaling limits of the multiplicative graphs, notably in Bhamidi, Sen \& X.~Wang and Bhamidi, van der Hofstad \& Sanchayan \cite{BhHoSe15, BhSeWa14}, previous works have distinguished two seemingly orthogonal cases depending on whether the inhomogeneity is mild enough to be washed away in the limit as in  Addario-Berry, B.~\& Goldschmidt, Bhamidi, B., Sen \& X.~Wang and Bhamidi, Sen \& X.~Wang  \cite{BhSeWa14,AdBrGo12a,BhBrSeWa14}, or strong enough to persist asymptotically as  in Bhamidi, van der Hofstad \& Sanchayan and Bhamidi, van der Hofstad \& van Leeuwaarden \cite{BhHoLe12,BhHoSe15}: the so-called asymptotic (Brownian) homogeneous case and the \emph{power-law} case. In these papers the proof strategies greatly differ in these two cases. On the other hand, the remarkable work of Aldous and Limic \cite{AlLi98} about the weights of large critical connected components deals with the inhomogeneity in a transparent way. We provide here such a unified approach for the geometry, which works not only for both cases but also for graphs which can be seen as a mixture of the two cases. 

Furthermore, an easy correspondence (see \eqref{bungabunga} below) allows us to link our parameters $(\alpha, \beta,$ $\kappa, \mathbf c)$ for the limit objects to the ones parametrising all the extremal eternal multiplicative coalescents, as identified by Aldous \& Limic in \cite{AlLi98}. 
We note that our limit theorems are valid in the Gromov--Hausdorff--Prokhorov topology, which controls the distances between all pairs of points, and not just in the Gromov--Prokhorov topology where only distances between finitely many typical points are controlled. (A general result has already been proved by Bhamidi, van der Hofstad \& Sen \cite{BhHoSe15} for the Gromov--Prokhorov topology in the special case when $\beta\!  =\! 0$.) In light of this, we believe our work contains an exhaustive treatment of all the possible limits related to those multiplicative coalescents. In the mean time, we remove some technical conditions that had been imposed on the weight sequences in some of the previous works. 

\medskip
\noindent\textit{Avoiding to compute the law of connected components:}
%\fmm{Remove this paragraph?}
    The connected components of the random graphs may be described as the result of the addition of ``shortcut edges'' to a tree; this picture is useful both for the discrete models and the limit metric spaces. The work of Bhamidi, Sen \& X.\ Wang and Bhamidi, van der Hofstad \& Sen \cite{BhSeWa14,BhHoSe15} yields an explicit description of the law of the random tree to which one should add shortcuts in order to obtain connected components with the correct distribution. As in the case of classical random graphs treated in Addario-Berry, B.~\& Goldschmidt \cite{AdBrGo12a}, this law involves a change of measure from one of the ``classical'' random trees, whose behaviour is in general difficult to control asymptotically. Our connected components are described as the metric induced on a subset of a Galton--Watson tree; the bias of the law of the underlying tree is somewhat transparently handled by the procedure that extracts the relevant subset. 

% \medskip
% \noindent\textit{Coding function and embedding provide a simplified framework:} 
% On important class of objects on which rely the the constructions of Bhamidi, Sen and Wang \cite{BhSeWa14} and of Bhamidi, van der Hofstad and Sen \cite{BhHoSe15} are the birthday trees and their scaling limits: inhomogeneous continuum random trees \cite{AlPi2000a,CaPi2000}. These trees and the corresponding limiting regimes are not well understood in full generality. For instance, the existence of the local times analogue to those in \eqref{llapproHdef} with the function $Y$ instead of $X$ is not trivial. The Galton--Watson trees and their scaling limits, Levy trees, which appear in our constructions are much more regular, and the conditions for convergence are fully understood and sharp thanks to the works of Le Gall \& Le Jan \cite{LGLJ98} and D.\ \& Le Gall \cite{DuLG02}. As an exemple, the embedding makes the existence of the local times for $Y$ an easy consequence of the existence of similar local times for the Levy process $X$.

\medskip
\noi
\textit{More general models of random graphs.}
While we focus on the model of the multiplicative graphs,  the theorems of Janson \cite{Ja10} on asymptotic equivalent models (see Section \ref{prevresults}) and the expected universality of the limits confers on the results obtained here potential implications that go beyond the realm of this specific model: for instance, random graphs constructed by the celebrated configuration model where the sequence of degrees has asymptotic properties similar to the weight sequence of the present paper are believed to exhibit similar scaling limits; see Section~3.1 in \cite{BhHoSe15} for a related discussion.

% We emphasize the fact that, while these results concern multiplicative inhomogeneous random graphs, the contiguity theorems of Janson \cite{Ja10} (see Section XXX) and the expected universality of the limits confers on them potential implications that go beyond the realm of this specific model: for instance, random graphs constructed by the celebrated configuration model where the sequence of degrees has asymptotic properties similar to the weight sequence of the present paper are believed to exhibit similar scaling limits. Some of the regimes have been studied by Bhamidi, Dhara, van der Hofstad and Sen \cite{BhDhHoSe2017a}; we conjecture that all the limits obtained here may also be obtained from random graphs generated by the configuration model.

\medskip
\noi
\textit{Upcoming work. } The current version of the limit theorems consider the sequences of connected components in the product topology. The embedding of the graphs in a forest of Galton--Watson forest actually also yields a control on the tail of the sequence, which would allow to strengthen the convergence to $\ell^p$-like spaces as in \cite{AdBrGo12a} or \cite{BhSeWa14}; this will be pursued somewhere else as well.

\bigskip

\noi
\textbf{Organisation of the paper}\label{ssec:organization}
In Section \ref{sec:results}, we state in precise forms the main results of the paper and Section \ref{prevresults} is devoted to the connection with previous results. Section \ref{discrsec} provides results on the discrete model: more specifically, in Section \ref{redbldissec} a precise definition of the red and blue coding of the Markovian queue is recalled from \cite{BDW1} and new estimates are proved in Section \ref{estisec}. In Section \ref{prevcontsec}, we recall from \cite{BDW1} the precise construction of the continuous state-space coding processes $Y$, $\cH$, $\theta^{\mathtt{b}}$, etc. The proof of the main limit theorems is done in Section \ref{Limisec} and it proceeds through a sequence of lemmas. Section \ref{Markosec} is devoted to the proof of the limit theorem for the processes related to the Markovian queues. 
An appendix collects some general results (on Laplace transform, Skorokod's topology, limit theorems for random walks, L\'evy processes and branching processes) that have been tailored in specific forms to adapt to our need here. We believe this facilitates the reading process. 

%An Appendix recall the main general results (on Laplace transform, Skorokod's topology, limit theorems for random walks, L\'evy processes and branching processes) that are needed and they are stated within our notation and under the required specific form, to facilitate the reading process. 

%%%%%%
%%%%%%

%{\footnotesize 
%\tableofcontents 
%}

\section{Main results}\label{sec:results}

\noi
\textbf{Notation}. Throughout the paper, $\bbN$ stands for the set of nonnegative integers and 
$\bbN^*\! = \! \bbN \backslash \{ 0\}$. 
A sequence of \textit{weights} refers to an element of the set  
$\elldo_{\infty}\! =\!  \big\{ (w_j)_{j\geq 1} \ino [0, \infty)^{\bbN^*} \! \! \!\! : \,  w_j \! \geq \! w_{j+1} \big\}$. For all $r\ino (0, \infty)$ and 
all $\bw\! = \!  (w_j)_{j\geq 1}\ino \elldo_\infty$, we set $\sigma_r (\bw)\! = \!  \sum_{j\geq 1} w_j^r \ino [0, \infty]$. 
The following subsets of $\elldo_{\infty}$ will be of particular interest to us. 
$$ \elldo_{{r}} = \big\{ \bw \ino \elldo_\infty : \sigma_r (\bw) \! < \! \infty \big\},  
 \quad \textrm{and} \quad  \elldo_{{\! f}}= 
\big\{ \bw \ino \elldo_\infty  :  \exists j_0 \! \geq \! 1 : w_{j_0} \! = \! 0  \big\} .$$

\subsection{Convergence results for the Markovian queue.} 
We fix a sequence $\bw_n\ino \elldo_{{\! f}}$, and two sequences $a_n, b_n \ino (0, \infty)$ that satisfy the a priori assumptions (\ref{apriori}).   
As already mentioned the convergence of the graphs $\cG_{\! \bw_n}$ is obtained thanks 
to the convergence of rescaled versions of 
$Y^{\bw_n} $ and $\cH^{\bw_n}$ and the convergence of these two processes is also obtained by the convergence of the Markovian processes into which they are embedded: 
namely, 
%$X^{\bw_n}\! $ and $H^{\bw_n}$: 
the asymptotic regimes of $(Y^{\bw_n} ,\cH^{\bw_n})$ and of $(X^{\bw_n}, H^{\bw_n})$ should be the same. The purpose of this section is to state weak limit-theorems for 
$X^{\bw_n}$ and $H^{\bw_n}$. Let us mention that many results of this section rely on standard 
limit-theorems on random walks, on results due to Grimvall in \cite{Gr74} on branching processes and on results due to Le Gall \& D.~in 
%\margmm{Galton--Watson trees}
\cite{DuLG02} on the height processes of Galton--Watson trees. However, the specific form of the jumps and of the offspring distribution of the trees 
actually requires a careful analysis done in the Proof-Section \ref{Markosec}.

\medskip

Recall from (\ref{HXdisdef}) the definition of $X^{\bw_n}$; recall 
that the Markovian queueing system induced by $X^{\bw_n}$ yields a tree that is 
an i.i.d.~sequence of Galton-Watson trees with offspring distribution $\mu_{\bw_n}$ whose definition is given by (\ref{muGWdef}). 
Denote by $(Z^{\bw_n}_k)_{k\in \bbN}$ a Galton-Watson Markov chain with offspring distribution $\mu_{\bw_n}$ and with initial state 
$Z^{\bw_n}_0\! = \! \lfloor a_n \rfloor$. The following proposition is mainly based on Theorem 3.4 in Grimvall \cite{Gr74} p.1040, that proves weak convergence for Galton-Watson processes to \textit{Continuous States Branching Processes} (CSBP for short). 
Recall that a (conservative) CSBP is a $[0, \infty)$-valued Markov process obtained from 
spectrally positive L\'evy processes via Lamperti's time-change; the law of the CSBP is completely characterised by the L\'evy process and thus by its Laplace exponent that is usually called the \textit{branching mechanism} of the CSBP: we refer to Bingham \cite{Bi76} for more details on CSBP (and see Appendix Section \ref{CSBPApp} for a very brief account). 
We denote by $\bD ([0,\infty) , \bbR)$ the space of c\`adl\`ag functions from $[0, \infty)$ to $\bbR$ equipped with Skorokod's topology and we denote by $\bC ([0,\infty) , \bbR)$ the space of continuous functions from $[0, \infty)$ to $\bbR$, equipped with the topology of uniform convergence on all compact subsets. 
 
% \margmm{Need to suppose $X_0=0$ and $Z_0=1$.} 
\begin{prop}
\label{cvmarkpro} Let $a_n , b_n \ino (0, \infty) $ and $\bw_n \ino \elldo_f$, $n\ino \bbN$, satisfy (\ref{apriori}). Recall from above the definition of $X^{\bw_n}$ and $Z^{\bw_n}$. Let $(X_t)_{t\in [0, \infty)}$ and $(Z_t)_{t\in [0, \infty)}$ be two c\`adl\`ag processes such that $X_0=0$ and $Z_0=1$. 
Then, the following holds true. 
\begin{compactenum}

\smallskip

\item[$(i)$] The following convergences are equivalent.  
\begin{compactenum}

\smallskip

\item[] $(i$-a$)$ There exists $t\ino (0, \infty)$ such that $\frac{1}{a_n} X^{\bw_n}_{b_n t }\! \rightarrow \! X_t$ weakly on $\bbR$.

\smallskip

\item[]  $(i$-b$)$ $(\frac{1}{a_n} X^{\bw_n}_{b_n t })_{t\in [0, \infty)}  \! \longrightarrow \! (X_t)_{t\in [0, \infty)}$ weakly  on $\bD ([0,\infty) , \bbR)$.

\smallskip

\item[] $(i$-c$)$ $(\frac{1}{a_n} Z^{\bw_n}_{\lfloor b_n t/a_n \rfloor})_{t\in [0, \infty)}  \! \longrightarrow \! (Z_t)_{t\in [0, \infty)}$  weakly  on $\bD ([0,\infty) , \bbR)$.
\end{compactenum}

\smallskip

\noi
If any of the three convergences in $(i)$ holds true, then $X$ is a spectrally L\'evy process and $Z$ a conservative CSBP; moreover 
there exist $\alpha \ino \bbR$, $\beta \ino [\beta_0, \infty)$, $\kappa \ino (0, \infty)$ and $\mathbf{c}\! = \! (c_j)_{j\geq 1}\ino \elldo_3$ such that the branching mechanism of $Z$ and the Laplace exponent of $X$  are equal to the same function $\psi$ given by:  
%\margmm{$j\ge 1$}
\begin{equation}
\label{abkcpsi}
\forall \lambda \ino [0, \infty), \quad \; 
\psi(\lambda)  \! = \!   
 \alpha \lambda +\frac{_{_1}}{^{^2}} \beta \lambda^2 +  \! \sum_{j\geq 1}   \kappa c_j   \big( e^{-\lambda c_j}\! -\! 1\! + \! \lambda c_j \big). 
\end{equation}
%\margmm{Add a comment that \eqref{abkcpsi} implies non-explosion?}
\item[$(ii)$] Any of the three convergences in $(i)$ is equivalent to the following three conditions:  
\begin{equation}
\label{unalphcv}
%\mathbf{(C1):} \quad \frac{a_nb_n}{\sigma_1 (\bw_n)}\;  \underset{n\rightarrow \infty}{-\!\!\! -\!\!\! \longrightarrow} \; \kappa \, , \qquad 
\mathbf{(C1):} \quad \frac{b_n}{a_n} \Big( 1-\frac{\sigma_2 (\bw_n)}{\sigma_1 (\bw_n)} \Big)\;  \underset{^{n\rightarrow \infty}}{-\!\!\! -\!\! \!  \longrightarrow}  \; \alpha \qquad \mathbf{(C2):} \quad  \frac{b_n}{a^2_n}\!  \cdot \!  \frac{\sigma_3 (\bw_n)}{\sigma_1 (\bw_n)} \;  \underset{^{n\rightarrow \infty}}{-\!\!\! -\!\! \!  \longrightarrow}  \; \beta + \kappa \sigma_3 (\mathbf{c}) \, ,
\end{equation}
\begin{equation}
\label{sig3cvcj}
%\mathbf{(C3):} \quad  \frac{b_n}{a^2_n}\!  \cdot \!  \frac{\sigma_3 (\bw_n)}{\sigma_1 (\bw_n)} \;  \underset{^{n\rightarrow \infty}}{-\!\!\! -\!\! \!  \longrightarrow}  \; \beta + \kappa \sigma_3 (\mathbf{c}) \, , \qquad 
\mathbf{(C3):} \quad \forall j \in \bbN^*, \quad \frac{w^{(n)}_j}{a_n } \;  \underset{^{n\rightarrow \infty}}{-\!\!\! -\!\! \! \longrightarrow}  \;  c_j
\, .
\end{equation}
\item[$(iii)$] Any of the three convergences of $(i)$ is equivalent to $\mathbf{(C1)}$ and the following limit for all $\lambda \ino (0, \infty)$: 

\vspace{-5mm}

\begin{equation}
\label{Laplexpo}
\frac{a_n b_n}{\sigma_1 (w_n)} \sum_{j\geq 1} \frac{w_j^{(n)}}{a_n} \Big( e^{-\lambda w^{(n)}_j /a_n}-1 + \lambda\,  w^{(n)}_j \!\! /a_n \Big) \; \underset{^{n\rightarrow \infty}}{-\!\!\! -\!\! \! -\!\! \!\longrightarrow} \; \psi (\lambda) -\alpha \lambda \; , 
\end{equation}

\vspace{-2mm}

\item[$(iv)$] For all $\alpha \ino \bbR$, $\beta \ino [0, \infty)$, $\kappa \ino (0, \infty)$ and $\mathbf{c}\! = \! (c_j)_{j\geq 1}\ino \elldo_3$, there are sequences $a_n , b_n \ino (0, \infty) $, $\bw_n \ino \elldo_f$, $n\ino \bbN$, satisfying (\ref{apriori}) with $\beta_0 \ino [0, \beta]$, $\mathbf{(C1)}$, $\mathbf{(C2)}$, $\mathbf{(C3)}$ and $\sqrt{\mathbf{j}_n}/ b_n \! \rightarrow \! 0$ where we recall that $\mathbf{j}_n \! = \! \max \{ j\! \geq \! 1 \! : \! w^{_{(n)}}_{^j} \! >\! 0 \}$. 
\end{compactenum}
\end{prop}
\textbf{Proof.} See Section \ref{Markosec} (and more specifically Section \ref{pfcvmarkpro}). As already mentioned, Proposition \ref{cvmarkpro} $(i)$ strongly relies on Theorem 3.4 in Grimvall \cite{Gr74} 
p.~1040. However, $(ii)$, $(iii)$ and $(iv)$ require specific arguments. \cqfd 

\begin{rem}
\label{newasspt} The condition  $\sqrt{\mathbf{j}_n}/b_n \! \rightarrow \! 0$ is explained by the second statement in Theorem \ref{graphcvth} that provides a scaling limit for the connected components of multiplicative graphs listed in the decreasing order of their numbers of vertices. \cq  
 \end{rem}

\medskip

Recall from (\ref{HXdisdef}) the definition of $H^{\bw_n}$, 
the height process associated with $X^{\bw_n}$. Note here that we also deal with supercritical cases. 
\begin{prop}
\label{HMarkcvprop} Let $\alpha \ino \bbR$, $\beta \ino [0, \infty)$, $\kappa \ino (0, \infty)$ and $\mathbf{c}\! = \! (c_j)_{j\geq 1}\ino \elldo_3$ and let $\psi$ be given by (\ref{abkcpsi}).  We assume that $\psi$ satisfies (\ref{Grey}): namely, $\int^\infty d\lambda / \psi (\lambda ) \! < \! \infty$.  
%
%
%
%\vspace{-5mm}
%
%
%
%
%\begin{equation} 
%\label{subcriGrey}
%\alpha \geq 0 \quad \textrm{and} \quad \int^\infty\!\!\!  \frac{d\lambda}{\psi (\lambda)} < \infty \; .
%\end{equation}
Let $X$ be a spectrally positive L\'evy process with Laplace exponent $\psi$. Let $H$ 
%$(H_{t})_{t\in [0, \infty)}$ 
be its height process as defined in (\ref{llapproHdef}). Let $a_n , b_n \ino (0, \infty) $, $\bw_n \ino \elldo_f$, $n\ino \bbN$, satisfy (\ref{apriori}) with $\beta_0 \ino [0, \beta]$, $\mathbf{(C1)}$, $\mathbf{(C2)}$ and $\mathbf{(C3)}$. 
%Suppose that $\sigma_2(\bw_n)\le \sigma_1(\bw_n)$ for all $n$.
We also assume the following:  
%Then, the weak convergence in $\bC ([0, \infty) , \bbR)$ 
%$(\frac{b_n}{a_n} H^{\bw_n}_{b_n t })_{t\in [0, \infty)} \! \rightarrow H$ is equivalent to the following condition: 

\vspace{-5mm}

\begin{equation}
\label{scalheight} 
\mathbf{(C4):} \quad \exists \, \delta \! \in \! (0, \infty) , \qquad \liminf_{n\rightarrow \infty} \bP \big( Z^{\bw_n}_{\lfloor b_n \delta /a_n \rfloor} \! =\!  0 \big) >0 \; .
\end{equation}
%Moreover if  (\ref{apriori}) and $\mathbf{(C1)}$ -- $\mathbf{(C5)}$ hold true, then the joint convergence holds true: 
Then, the joint convergence holds true 
\begin{equation}
\label{jointecon}
\big( (\frac{_{_1}}{^{^{a_n}}} X^{\bw_n}_{b_n t })_{t\in [0, \infty)} , (\frac{_{_{a_n}}}{^{^{b_n}}} H^{\bw_n}_{b_n t })_{t\in [0, \infty)} \big) \underset{n\rightarrow \infty}{-\!\!\! -\!\!\! -\!\!\! \longrightarrow } (X, H)  
\end{equation}
weakly on $\bD ([0, \infty), \bbR) \times \bC ([0, \infty), \bbR)$ equipped with the product topology. We also get:  
\begin{equation}
\label{cvtpsvie} 
\forall t\ino [0, \infty), \quad  \lim_{n\rightarrow \infty} \bP \big( Z^{\bw_n}_{\lfloor b_n t /a_n \rfloor} \! = \! 0 \big) = e^{-v_\psi(t)} \quad \textrm{where} \quad  \int_{v_\psi (t)}^\infty \! \frac{d\lambda}{\psi (\lambda)}= t . 
\end{equation}
%
%Then, the following assertions holds true. 
%%Recall from above the definition of $X^{\bw_n}$, $H^{\bw_n}$ and $Z^{\bw_n}$. 
%
%\begin{compactenum}
%
%
%\smallskip
%
%
%
%\item[$(i)$] Suppose that $(\frac{1}{{a_n}} X^{\bw_n}_{b_n t })_{t\in [0, \infty)}  \! \rightarrow \! (X_t)_{t\in [0, \infty)}$ a.s.~for the Skorokod topology on $\bD ([0,\infty) , \bbR)$. 
%Then, for all $t\ino [0, \infty)$, 
%$\frac{a_n}{b_n} H^{\bw_n}_{b_n t} \rightarrow H_t$ holds in probability. 
%This entails the weak convergence of the finite dimensional marginals of $(\frac{b_n}{a_n} H^{\bw_n}_{b_n t })_{t\in [0, \infty)}$ to those of $H$. 
%
%
%\smallskip
%
%
%\item[$(ii)$] Under (\ref{apriori}) and $\mathbf{(C1)}$ -- $\mathbf{(C4)}$, the 
%convergence $(\frac{b_n}{a_n} H^{\bw_n}_{b_n t })_{t\in [0, \infty)} \! \rightarrow H$ in law on $\bC ([0, \infty) , \bbR)$ is equivalent to the following condition: 
%
%
%\vspace{-6mm}
%
%
%
%
%\begin{equation}
%\label{scalheight} 
%\mathbf{(C5):} \quad \exists \delta \ \in \! (0, \infty) , \qquad \liminf_{n\rightarrow \infty} \bP \big( Z^{\bw_n}_{\lfloor b_n \delta /a_n \rfloor} = 0 \big) >0 \; .
%\end{equation}
%
%\vspace{-2mm}
%
%\noi
%Moreover if  (\ref{apriori}) and $\mathbf{(C1)}$ -- $\mathbf{(C5)}$ hold true, then 
%
%\vspace{-6mm}
%
%
%
%
%\begin{equation}
%\label{cvtpsvie} 
%\forall t\ino [0, \infty), \quad  \lim_{n\rightarrow \infty} \bP \big( Z^{\bw_n}_{\lfloor b_n t /a_n \rfloor} = 0 \big) = e^{-v_\psi(t)} \quad \textrm{where} \quad  \int_{v_\psi (t)}^\infty \! \frac{d\lambda}{\psi (\lambda)}= t . 
%\end{equation}
%
%\vspace{-4mm}
%
%
%
%\end{compactenum}
\end{prop}
\noi
\textbf{Proof.} See Section \ref{Markosec} (and more specifically Section \ref{pfcvmarkpro}). Proposition \ref{HMarkcvprop} strongly relies on Theorem 2.3.1 in Le Gall \& D.~\cite{DuLG02}. However, its proof requires more care than expected at first glance because $H^{\bw_n}$ is not exactly the \textit{height process} 
as defined in \cite{DuLG02} (it is actually a time-changed version of the so-called \textit{contour} process as in Theorem 2.4.1 \cite{DuLG02} p.~68). \cqfd

\medskip

The following proposition provides a practical criterion to check $\mathbf{(C4)}$: in particular, 
it shows that \textit{$\mathbf{(C4)}$ is always true when $\beta_0 \! >\!  0$}; it also shows that Proposition \ref{HMarkcvprop} is never void. 
\begin{prop}
\label{Hcritos} Let $\alpha\ino \bbR, \beta  \ino [0, \infty)$, $\kappa \ino (0, \infty)$ and 
$\mathbf{c}\! = \! (c_j)_{j\geq 1}\ino \elldo_3$. Let $\psi$ be given by (\ref{abkcpsi}) and assume that $\psi$ satisfies (\ref{Grey}): namely, $\int^\infty d\lambda / \psi (\lambda ) \! < \! \infty$. Let us recall from (\ref{jjnndef}) that $\mathbf{j}_n\! = \! \max \{ j \ge 1:  w^{_{(n)}}_{^j} \! >\! 0\}$. 
%$$\mathbf{j}_n\! = \! \max \{  \mm{j \ge 1:}  w^{_{(n)}}_{^j} \! >\! 0\} \; .$$
Then, the following holds true. 
\begin{compactenum}

\smallskip

\item[$(i)$] Let $a_n , b_n \ino (0, \infty) $, $\bw_n \ino \elldo_f$, $n\ino \bbN$, satisfy (\ref{apriori}), $\mathbf{(C1)}$, $\mathbf{(C2)}$ and $\mathbf{(C3)}$. Denote by 
$\psi_n$ the Laplace exponent of $(\frac{1}{a_n} X^{\bw_n}_{b_n t })_{t\in [0, \infty)}$: namely, for all $\lambda \ino [0, \infty)$, 
\begin{equation}
\label{psindef}
\psi_n (\lambda)\! = \! \frac{b_n}{a_n} \Big( 1-\frac{\sigma_2 (\bw_n)}{\sigma_1 (\bw_n)} \Big) \lambda + \frac{a_n b_n}{\sigma_1 (w_n)} \sum_{j\geq 1} \frac{w_j^{(n)}}{a_n} \Big( e^{-\lambda w^{(n)}_j /a_n}-1 + \lambda\,  w^{(n)}_j \!\! /a_n \Big) . 
\end{equation} 
Then, $\mathbf{(C4)}$ holds true if 
\vspace{-4mm}

\begin{equation}
\label{crinullo}
\lim_{y\rightarrow \infty} \limsup_{n\rightarrow \infty} \int_y^{a_n} \!\! \frac{d\lambda}{\psi_n (\lambda)} = 0 \; .
\end{equation}

\vspace{-2mm}

In particular, if $\beta_0 \! >\! 0$  in (\ref{apriori}), then (\ref{crinullo}) is always satisfied and $\mathbf{(C4)}$ holds true.

\smallskip

\item[$(ii)$] There are sequences 
$a_n , b_n \ino (0, \infty) $, $\bw_n \ino \elldo_f$, $n\ino \bbN$, satisfying (\ref{apriori}) with $\beta_0\! = \! 0$,  $\sqrt{\mathbf{j}_n}/b_n \! \rightarrow \! 0$, $\mathbf{(C1)}$, $\mathbf{(C2)}$ and $\mathbf{(C3)}$ but not 
$\mathbf{(C4)}$. 

\smallskip

\item[$(iii)$] There exist $a_n , b_n \! \in \! (0, \infty)$, and  $\bw_n \! \in \! \elldo_f$, $n\ino \bbN$, that satisfy 
(\ref{apriori}) with any $\beta_0\ino [0, \beta]$, $\sqrt{\mathbf{j}_n}/b_n \! \rightarrow \! 0$, $\mathbf{(C1)}, \mathbf{(C2)}, \mathbf{(C3)}$ and $\mathbf{(C4)}$.

%
%\smallskip
%
%
%
%
%\item[$(ii)$] Let $a_n , b_n \ino (0, \infty) $, $\bw_n \ino \elldo_f$, $n\ino \bbN$, satisfying (\ref{apriori}) and $\mathbf{(C1)}$, $\mathbf{(C2)}$ and $\mathbf{(C3)}$. Denote by 
%$\psi_n$ the Laplace exponent of $(\frac{1}{a_n} X^{\bw_n}_{b_n t })_{t\in [0, \infty)}$. Then, $\mathbf{(C4)}$ holds true if 
%
%\vspace{-4mm}
%
%\begin{equation}
%\label{crinullo}
%\lim_{y\rightarrow \infty} \limsup_{n\rightarrow \infty} \int_y^{a_n} \!\! \frac{d\lambda}{\psi_n (\lambda)} = 0 \; .
%\end{equation}
%
%
%\vspace{-2mm}
%
%
%\item[$(iii)$] There exist $a_n , b_n \! \in \! (0, \infty)$, and  $w_n \! \in \! \elldo_f$, $n\ino \bbN$, satisfying 
%(\ref{apriori}), $\mathbf{(C1)}, \mathbf{(C2)}, \mathbf{(C3)}$ and $\mathbf{(C4)}$. 
%
%For all $y\ino (0, \infty)$, set $ s(y) \! =\!  \sum_{j\geq1} \un_{[0, y]} (c_j) \, c_j^3 $. 
%Note that $s \downarrow 0$ as $y \downarrow 0+$. Assume that 
%%$\lim_{y\rightarrow 0+}s(y)\! = \! 0$. Assume that 
%
%\vspace{-6mm}
%
%\begin{equation}
%\label{pastop}
%\textrm{either} \quad \beta >0  \qquad \textrm{or} \quad \exists q \ino (1, \infty) : \quad  \liminf_{y \rightarrow 0+} s(q y) / s(y) > 1 \; . 
%\end{equation}
%
%\vspace{-4mm}
%
%\noi
%Then, there exist $a_n , b_n \! \in \! (0, \infty)$, and  $w_n \! \in \! \elldo_f$, $n\ino \bbN$, satisfying 
%(\ref{apriori}) and $\mathbf{(C1)}$ -- $\mathbf{(C4)}$ and $\mathbf{(C5)}$. 
\end{compactenum}
\end{prop}
\textbf{Proof.} See Sections \ref{pfHcritos1}, \ref{pfHcritos2} and \ref{pfHcritos3}. \cqfd  

%\begin{rem}
%\label{qdmmtop} 
%We believe that for all $\alpha \ino \bbR$, $\beta \ino [0, \infty)$, $\kappa \ino (0, \infty)$ and $\mathbf{c}\! = \! (c_j)_{j\geq 1}\ino \elldo_3$ such that $\psi$ given by (\ref{abkcpsi}) satisfies (\ref{subcriGrey}), 
%there exist $a_n , b_n \! \in \! (0, \infty)$, $n\! \in \! \bbN$, and  $w_n \! \in \! \elldo_f$, that satisfy 
%(\ref{apriori}) and $\mathbf{(C1)}$ -- $\mathbf{(C5)}$ (namely, Proposition \ref{HMarkcvprop} $(ii)$ is never void). 
%Condition (\ref{pastop}) is not optimal but it holds true for all $\mathbf{c}$ such that $s$ is regularly varying with positive exponent: for instance 
%$c_j \! = \! j^{-1/\tau}$ with $\tau\ino [2, 3)$.   \cq 
%\end{rem}
%\begin{rem}
%\label{newasspt} The condition  $\sqrt{\mathbf{j}_n}/b_n \! \rightarrow \! 0$ is explained by the second statement in Theorem \ref{graphcvth} that provide a scaling limit for the connected components of multiplicative graphs listed in the decreasing order of their number of vertices. \cq  
% \end{rem}

\subsection{Convergence of the processes coding the multiplicative graphs.}
%\label{limithouf} 
%\subsubsection*{Limit theorems for multiplicative random graphs.}
\label{limithouf} 
%\noi
%\textit{A convention.} 
To deal with limits of sequences of pinching times, it is convenient 
to embed $([0, \infty)^2)^p$
into $(\bbR^2)^{\bbN^*}$ by extending any sequence 
$((s_i,t_i))_{1\leq i\leq p} \in ([0, \infty)^2)^p$ by setting $(s_i, t_i)\! = \! (-1, -1)$, for all $i\! >\! p$. Here, $(-1, -1)$ plays 
the role of an unspecific cemetery point. We equip $ (\bbR^2)^{\bbN^*}$ 
with the product topology. 
Then, the main theorem of paper is the following. 
\begin{thm}
\label{HYcvth} Let $\alpha\in \bbR$, $\beta  \ino [0, \infty)$, $\kappa \ino (0, \infty)$ and 
$\mathbf{c}\! = \! (c_j)_{j\geq 1}\ino \elldo_3$. Let $\psi$ be given by (\ref{abkcpsi}) and assume that $\psi$ satisfies (\ref{Grey}): namely, $\int^\infty d\lambda / \psi (\lambda ) \! < \! \infty$. 
Recall from (\ref{plasouilloc}) the definition of $Y$ and from (\ref{defcH}) the definition of $\cH$; 
%recall from Theorem \ref{Xdefthm} the definition of $X$; recall from (\ref{llapproHdef}) the definition of $H$; 
%recall from Theorem \ref{cHdefthm} the definition of $\cH$; 
recall from (\ref{pinchset}) the definition of $\Ptt$. 
Let $a_n , b_n \! \in \! (0, \infty)$, and  $\bw_n \! \in \! \elldo_f$, $n\ino \bbN$, satisfy (\ref{apriori}), $\mathbf{(C1)}\! -\! \mathbf{(C4)}$ as specified in (\ref{unalphcv}), (\ref{sig3cvcj}) and (\ref{scalheight}). 
Recall from (\ref{defYcHH}) the definition of $(Y^{\bw_n}, \cH^{\bw_n})$; 
%recall from Proposition \ref{Xwfrombr} the definition of $X^{\bw_n}$;
%recall (\ref{JHdef}) the definition of $\cH^{\bw_n}$; 
%recall from (\ref{XJHdef}) the definition of $H^{\bw_n}$; 
recall from (\ref{pin1}) the definition of $\Ptt_{\bw_n}$.
Then, the joint convergence 
\begin{equation}
\label{jtcvYHP}
\big(
% \frac{_{_1}}{^{^{a_n}}} A^{\bw_n}_{b_n \cdot }\,  , 
\frac{_{_1}}{^{^{a_n}}} Y^{\bw_n}_{b_n \cdot } ,  \,     \frac{_{_{a_n}}}{^{^{b_n}}} \cH^{\bw_n}_{b_n \cdot}  \, ,   \frac{_{_1}}{^{^{b_n}}} \Ptt_{ \bw_n} \big) 
\;  \underset{n\rightarrow \infty}{-\!\!\! -\!\!\! -\!\!\! \longrightarrow} \; \big( Y,  \cH ,  \Ptt  \big)
\end{equation} 
%
%
%
%\begin{equation}
%\label{jtcvYHP}
%\big( \frac{_{_1}}{^{^{a_n}}} X^{\bw_n}_{b_n \cdot }\,  , \frac{_{_1}}{^{^{a_n}}} Y^{\bw_n}_{b_n \cdot } ,  \,   \frac{_{_{a_n}}}{^{^{b_n}}} H^{\bw_n}_{b_n \cdot } \, ,   \frac{_{_{a_n}}}{^{^{b_n}}} \cH^{\bw_n}_{b_n \cdot}  \, ,   \frac{_{_1}}{^{^{b_n}}} \Ptt_{ \bw_n} \big) 
%\;  \underset{n\rightarrow \infty}{-\!\!\! -\!\!\! -\!\!\! \longrightarrow} \; \big( X, Y, H, \cH ,  \Ptt  \big)
%\end{equation} 
%
%
%
%\begin{equation}
%\label{jtcvYHP}
%\big( \frac{_{_1}}{^{^{a_n}}} X^{\bw_n}_{b_n \cdot }\,  , \frac{_{_{a_n}}}{^{^{b_n}}} H^{\bw_n}_{b_n \cdot } \, ,  
% \big( \frac{_{_1}}{^{^{b_n}}} \theta^{\mathtt{b}, \bw_n}_{b_n \cdot } , \frac{_{_1}}{^{^{a_n}}} Y^{\bw_n}_{b_n \cdot } \big),  \frac{_{_{a_n}}}{^{^{b_n}}} \cH^{\bw_n}_{b_n \cdot}  \, ,   \frac{_{_1}}{^{^{b_n}}} \Ptt_{ \bw_n}\big) 
%\;  \underset{n\rightarrow \infty}{-\!\!\! -\!\!\! -\!\!\! \longrightarrow} \; \big( X, H, (\theta^\mathtt{b}, Y), \cH ,  \Ptt  \big)
%\end{equation} 
holds weakly on $\bD([0, \infty), \bbR) \! \times \! \bC ([0, \infty) , \bbR)
 \! \times \! (\bbR^2)^{\bbN^*}\! $ equipped with the product topology.  
%holds weakly on $\bD([0, \infty), \bbR)Â \! \times \! \bC ([0, \infty) , \bbR)Â \! \times \!  \bD([0, \infty), \bbR^2)Â \! \times \! 
%\bC ([0, \infty) , \bbR) \! \times \! (\bbR^2)^{\bbN^*}$ equipped with the product topology.  
\end{thm}
\noi
\textbf{Proof.} See Section \ref{proofmain}. Let us mention that we actually prove a joint convergence of all the involved processes such as $X^{\bw_n}, H^{\bw_n}, \theta^{\mathtt{b}, \bw_n}$,  ... to their continuous counterparts.  \cqfd  

\medskip

Theorem \ref{HYcvth} implies 
%that rescaled versions of $(\cH^{\bw_n} , \Ptt_{\bw_n})$ converge to $(\cH, \Ptt)$. 
%This entails 
the convergence of the coding processes of the connected components of $\cG_{\! \bw_n}$, because each connected component of 
$\cG_{\! \bw_n}$ is coded by an excursion above $0$ of $\cH^{\bw_n}$ and the corresponding pinching points. 
More precisely, denote by $(l^{\bw_n}_k, r^{\bw_n}_k)$, $1\! \leq \! k \! \leq \! \bq_{\bw_n}$, the excursion intervals of $\cH^{\bw_n}$ above $0$, that are exactly the excursion intervals of $Y^{\bw_n}$ above its infimum process $J^{\bw_n}_t \! = \! \inf_{s\in [0, t]} Y^{\bw_n}_s$. Namely, 
\begin{equation}
\label{excdibor}
 \bigcup_{1\leq k\leq \bq_{\bw_n}} [l^{\bw_n}_k, r^{\bw_n}_k) = \big\{  t\ino [0, \infty) : \cH^{\bw_n}_t >0  \big\}= \big\{  t\ino [0, \infty) : Y^{\bw_n}_t > \JJ^{\bw_n}_t  \big\} \, .
\end{equation} 
Here the indexation is such that $\zeta^{\bw_n}_k \! \geq \! \zeta^{\bw_n}_{k+1}$, where we have set $\zeta^{\bw_n}_k\! = \! r^{\bw_n}_k \! -\! l^{\bw_n}_k$ 
(if $\zeta^{\bw_n}_k\! = \! \zeta^{\bw_n}_{k+1}$, then we agree on the convention that $l^{\bw_n}_k \! < \! l^{\bw_n}_{k+1}$); \
the \textit{excursions processes} are then defined as follows: 
\begin{equation}
\label{excdiscret}
 \forall k \ino \{ 1, \ldots , \bq_{\bw_n}\}, \, \forall t\ino [0, \infty), \quad  
\Htt^{\bw_n}_k (t)\! =\!  \cH^{\bw_n}_{(l^{\bw_n}_k + t)\wedge r^{\bw_n}_k} .
%\; \,  \textrm{and} \; \,   \Ytt^{\bw_n}_k (t)\! = \! Y^{\bw_n}_{(l^{\bw_n}_k + t)\wedge r^{\bw_n}_k}- \JJ^\bw_{l^{\bw_n}_k} . 
\end{equation}
We next define the sequences of \textit{pinching points of the excursions}: to that end, 
recall from (\ref{pin1}) 
the definition of $\Ptt_{\bw_n}\! = \! \big( (s_p, t_p)\big)_{1\leq p\leq \bp_{\bw_n}}\! $ that is the sequence of 
pinching times of $\cG_{\bw_n}$; observe that if $t_p \ino [l^{\bw_n}_k, r^{\bw_n}_k]$, 
then $s_p \ino [l^{\bw_n}_k, r^{\bw_n}_k]$; then, it allows to define the following for all $k\ino \{ 1, \ldots, \bq_{\bw_n}\}$: 
\begin{multline}
\label{defPiwk}
\Ptt_{k}^{\bw_n}\! = \! \big( (s^k_p, t^k_p)\big)_{1\leq p\leq \bp^\bw_k} \; \textrm{where $(t^k_p)_{1\leq p\leq \bp^{\bw_n}_k}$ increases and where} \\ \textrm{ the $(l^{\bw_n}_k +s^k_p, l^{\bw_n}_k+t^k_p)$'s are exactly the terms $(s_{p^\prime} , t_{p^\prime})$ of $\Ptt_{\bw_n}$ such that $t_{p^\prime}\ino [l^{\bw_n}_k, r^{\bw_n}_k]$.}
\end{multline}
As already specified, we trivially extend each finite sequence 
$\Ptt_{k}^{\bw_n}$ as a random element of $(\bbR^2)^{\bbN^*}$. 
We pass to the limit for rescaled versions of $((\Htt_k^{\bw_n} , l_k^{\bw_n}, r_k^{\bw_n}, \Ptt_k^{\bw_n}))_{1\leq k \leq q_{\bw_n}}$. 
Since $q_{\bw_n}$ tends to $\infty$, it is convenient to extend this sequence by taking for all $k\! >\! q_{\bw_n}$, 
$\Htt_k^{\bw_n}$ as the null function, $l_k^{\bw_n} \! = \! r_k^{\bw_n}\! = \! 0 $ and $\Ptt_k^{\bw_n}$ as the sequence constant to $(-1, -1)$.

Similarly, recall from (\ref{excuHY}) the definition of the excursion intervals of $\cH$ above $0$: $\bigcup_{k\geq 1} (l_k, r_k)\! = \! \{ t\ino [0, \infty): \cH_t \! >\! 0 \}$, where indexation is chosen in such a way that the sequence $\zeta_k:= r_k\! -\! l_k$, $k\! \geq \! 1$, decreases.  We define the \textit{excursion processes} as follows. 
\begin{equation}
\label{exccont}
 \forall k \! \geq \! 1, \; \forall t\ino [0, \infty), \qquad 
\Htt_{k}(t)= \cH_{(l_k + t)\wedge r_k} .
%\quad \textrm{and} \quad  \Ytt_{k}(t)= Y_{(l_k + t)\wedge r_k}- \JJ_{l_k} . 
\end{equation}
The \textit{pinching times} are defined as follows: recall from (\ref{Poisurcon}) and (\ref{pinchset}) 
the definition of $\Ptt\! = \! \big( (s_p, t_p)\big)_{p\geq 1}$. 
If $t_p \ino [l_k, r_k]$, then note that $s_p \ino [l_k, r_k]$, by definition of $s_p$. 
For all $k\! \geq \! 1$, we then define:  
\begin{multline}
\label{defPik}
\Ptt_{k}\! = \! \big( (s^k_p, t^k_p)\big)_{1\leq p\leq \bp_k} \; \textrm{where $(t^k_p)_{1\leq p\leq \bp_k}$ increases 
and where} \\ \textrm{ the $(l_k +s^k_p, l_k+t^k_p)$'s are exactly the terms $(s_{p^\prime} , t_{p^\prime})$ of $\Ptt$ such that $t_{p^\prime}\ino [l_k, r_k]$.}
\end{multline}
Then the following theorem holds true. 
\begin{thm}
\label{excucvth} Under the same assumptions as in Theorem \ref{HYcvth}, the following convergence 
\begin{equation}
\label{Hkbwn}
\big(\big(  (\frac{_{_{a_n}}}{^{^{b_n}}} \Htt_k^{\bw_n} (b_n t))_{t\in [0, \infty)} , \, \frac{_{_{1}}}{^{^{b_n}}}l_k^{\bw_n} , \, \frac{_{_{1}}}{^{^{b_n}}}r_k^{\bw_n} , \,  \frac{_{_{1}}}{^{^{b_n}}}\Ptt_k^{\bw_n} \big) \big)_{k\geq 1}
\;  \underset{n\rightarrow \infty}{-\!\!\! -\!\!\! -\!\!\! \longrightarrow} \; \big(\big(\Htt_k, l_k, r_k, \Ptt_k\big) \big)_{k\geq 1}
\end{equation} 
holds weakly on $((\bC([0, \infty), \bbR)\! \times \! [0, \infty)^2 \! \times \!  (\bbR^2)^{\bbN^*})^{\bbN^*}$ equipped with the product topology. 
\end{thm}
\noi
\textbf{Proof.} See Section \ref{excucvpf}. \cqfd  

\medskip

\subsection{Convergence of the multiplicative graphs.} 
\label{ccvmulsec}
We recall here a generic procedure described in \cite{BDW1} which allows us to extract the $\bw$-graph $\cG_{\bw}$ from the coding processes $(Y^{\bw}, \cH^{\bw}, \Ptt_{\bw})$ and the continuous multiplicative graph from $(Y, \cH, \Ptt)$. We begin with the coding of trees by real-valued functions. 

\smallskip
\noi
\textbf{Coding trees.}
\label{codtreesec} Let $h\! :\!  [0, \infty) \! \rightarrow \! [0, \infty)$ be a c\`adl\`ag function such that 
\begin{equation}
\label{zetadef}
%h(0)= 0 \quad \textrm{and} \quad 
\zeta_h\! = \! \sup \{ t\ino [0,  \infty) \! : \! h(t) \! >\! 0 \}< \infty \; .
\end{equation}
We further assume that one of the following conditions is satisfied:
\begin{align}
\label{codecond}
\text{either} \quad \text{(a)} \quad h \text{ takes finitely many values or } 
\quad \text{(b)} \quad h \text{ is continuous.}
\end{align}
Note that the (discrete) height  process $\cH^{\bw}$ as defined in (\ref{defYcHH}) verifies Condition (a), while in the continuous setting, the process $\cH$ defined in (\ref{defcH}) verifies Condition (b), as asserted by Theorem \ref{cHdefthm}. 
For all $s, t\ino [0, \zeta_h)$, we set 
\begin{equation}
\label{pseudometric}
b_h(s,t)=\!\!\!\!\!\!\!\!\! \inf_{\quad r\in[s\wedge t,s\vee t]} \!\!\!\!\!\!\!\!\! h(r) \qquad {\rm and} \qquad d_h(s,t)=h(s)+h(t)-2b_h(s,t).
\end{equation}
We readily check that $d_{h}$ satisfies the {\it four-point inequality}: for all $s_1, s_2, s_3, s_4 $ belonging to $[0, \zeta_h)$, 
$d_h(s_1,s_2) + d_h(s_3, s_4) \! \leq \! \big(d_h(s_1, s_3) + d_h(s_2, s_4)\big) \! \vee \!  \big( d_h(s_1, s_4) + d_h(s_2, s_3)  \big) $. It follows that $d_h$ is a pseudometric on $[0, \zeta_h)$. We denote by $s\! \sim_h \! t$ the equivalence relation
 $d_h(s,t) \! =\! 0$ and we set 
\begin{equation}
\label{codef}
T_h= [0, \zeta_h) / \! \sim_h \; .  
\end{equation}
%\margmm{$[0, \zeta_h)$. Also, I suggest  mentioning $p_h$ is not necessarily continuous.}
Then, $d_h$ induces a true metric on the quotient set $T_h$ that we keep denoting by $d_h$ and we denote by $p_h\! :\!  [0, \zeta_h) \!  \rightarrow \! T_h $ 
the \textit{canonical projection}. Note that  if $h$ is continuous, then $p_h$ is a continuous map. It follows that in that case the metric space $T_h$ is a {\it compact real tree}, namely a compact metric space where any pair of points is joined by a unique injective path that turns out to be a geodesic (see Evans \cite{Ev08} for more references on this topic). If, on the other hand, $h$ satisfies Condition (a) in \eqref{codecond}, then $T_{h}$ is compact but not connected. It is still tree-like, as  $d_{h}$ satisfies the four-point inequality. 

%\vspace{-2mm}
%
%\begin{rem}
%\label{pataques}
%The metric space $(T_h, d_h)$ is tree-like but in general it is not necessarily connected or compact. However, we shall consider the following cases.  
%
%\smallskip
%
%\begin{compactenum}
%\item[$(a)$]  $h$ takes finitely many values. 
%% jumps and for all $t\ino [0, \infty)$ there is $\epp_t  \ino  (0, \infty)$ such that $h$ is constant on $s\ino[t, t+\epp_t]$. 
%
%\smallskip
%
%
%\item[$(b)$] $h$ is continuous. 
%\end{compactenum}
%
%\smallskip
%
%%\margmm{Suggest: $T_h$ is in fact formed by a finite number of points. In particular, $\mathcal H^\bw$ is in this case.}
%\noi
%\textit{In Case $(a)$}, $T_h$ is not connected but it is compact; $T_h$ is in fact formed by a finite number of points. In particular, $\mathcal H^\bw$ is in this case: by (\ref{gloupii}),  the exploration tree $\cT_{\!\! \bw}$ as defined in (\ref{explotree}) is actually isometric to $T_{\cH^\bw}$, that is the tree coded by the height process $\cH^{\bw}$ that is derived from $Y^\bw$ by (\ref{JHdef}). 
%
%\smallskip
%
%\noi
%\textit{In Case $(b)$}, $T_h$ is compact and connected; the metric $d_h$ satisfies the four-points condition: it is therefore a compact real tree, namely a compact metric space such that any pair of points is joined by a unique injective path that turns out to be a geodesic (see Evans \cite{Ev08} for more references on this topic). \cq 
%\end{rem}

The metric space $(T_{h}, d_{h})$ also inherits from $h$ the following features: a distinguished point $\rho_h\! = \! p_h (0)$, called the \textit{root} of $T_h$, and the 
%\margmm{$[0, \zeta_h)$}
\textit{mass measure} $m_h$, which satisfies that for any Borel measurable function $f \! : \! T_h \! \rightarrow \!  [0, \infty) $, we have $\int_{T_h}  f(\sigma ) \, m_h (d\sigma) = \! \int_{[0, \zeta_h]} f(p_h(t)) \, dt.$

\medskip

\noi
\textbf{Pinched metric spaces.}
\label{pinmetpar} 
%We next briefly explain how the metric of a graph is modified when several vertices are added. 
%Since we deal with rescaled version of graphs, and continuous 
%limits of such graphs, we discuss this point in a general setting. 
%Our construction of the graphs both in the discrete and the continuous settings consists in creating cycles in  the spanning trees of the graphs. In the discrete setting, we do this  by adding surplus edges; in the continuous setting, we identify pairs of points in a real tree. In what follows, we provide  a unified way to deal with both operations on the metrics. 
Let $(E, d)$ be a metric space and let $\Ptt\! = \! ((x_i, y_i))_{1\leq i\leq p}$ where the elements $(x_i, y_i)\ino E^2$, $1\! \leq \! i\! \leq \! p$, are referred to as \textit{pinching points}. Let $\epp \ino [0, \infty)$, that is interpreted as the length of the edges that are added to $E$ (if $\epp \! = \! 0$, then each $x_i$ is identified with $y_i$). 
Set $A_E \! =\! \{ (x,y) : x,y\ino E\}$ and for all $e\! = \! (x,y)\ino A_E$, 
set $\underline{e}\! = \! x$ and $\overline{e}\! = \! y$. A path $\gamma$ joining $x$ to $y$ is a sequence of $e_1, \ldots , e_q \ino A_E$ such that $\underline{e}_1\! = \! x$, $\overline{e}_q\! = \! y$ and $\overline{e}_i\! = \! \underline{e}_{i+1}$,  for all $1\! \leq \! i \! < \! q$. 
For all  $e\! = \! (x, y)\ino A_E$, we then define its length by $l_e \! =\!  \epp \! \wedge \! d(x_i, y_i)$ if $(x,y)$ or $(y,x)$ is equal to some $(x_i,y_i)\ino \Ptt$; otherwise we set $l_e \! = \! d(x, y)$. The length of a path $\gamma\! = \! (e_1, \ldots , e_q)$ is given by 
$l(\gamma)\! = \! \sum_{1\leq i\leq q} l_{e_i}$, and we set:  
\begin{equation}
\label{pinchdist}
 \forall x, y\ino E, \quad d_{\Ptt, \epp}  (x,y)= \inf \big\{ l (\gamma) :\;  \textrm{$\gamma$ is a path joining $x$ to $y$} \big\} \; .
\end{equation}
We set $A_{\Ptt}= \{ (x_i, y_i), (y_i, x_i); 1\! \leq \! i \! \leq \! p \}$ and we easily check that 
\begin{multline}
\label{vnaglurns}
d_{\Ptt , \epp} (x,y)\! = \! d(x,y) \wedge \min \big\{ \, l (\gamma)\; : \;  \gamma \!  = \! (e_0, e^\prime_0, \ldots ,e_{r-1}, e^\prime_{r-1}, e_r), \\ \textrm{a path joining $x$ to $y$ such that}\;  e_0^\prime, \ldots e^\prime_{r-1}\ino A_{\Ptt} \; \textrm{and} \; r\! \leq \!  p  
\big\} .  
\end{multline}

%We refer to Section \ref{Pinfrac} for more details. 
%\margmm{Suggest including a reference to Section \ref{Pinfrac}.}
Clearly, $d_{\Ptt, \epp}$ is a pseudo-metric
% (namely, it is symmetric and it satisfies the triangle inequality). 
and we denote the equivalence relation $d_{\Ptt, \epp}  (x,y)\! = \! 0$ by $x \equiv_{\Ptt, \epp} y$; 
the \textit{$(\Ptt, \epp)$-pinched metric space associated with $(E, d)$} is then the quotient space $E/\!\! \equiv_{\Ptt, \epp}$ equipped with $d_{\Ptt, \epp}$. First note that if $(E, d)$ is compact or connected, so is the associated $(\Ptt, \epp)$-pinched metric space since the canonical projection $\varpi_{\Ptt, \epp}\! : \! E \! \rightarrow \!  E/\!\! \equiv_{\Ptt, \epp}$ is $1$-Lipschitz. Of course when $\epp \! >\! 0$, $d_{\Ptt , \epp}$ on $E$ is a true metric, $E\! = \!  E/\!\! \equiv_{\Ptt, \epp}$ and $\varpi_{\Ptt, \epp}$ is the identity map on $E$. 

\medskip

\noi
\textbf{Coding pinched trees.} 
\label{codpintre} Let $h\! :\!  [0, \infty) \! \rightarrow \! [0, \infty)$ be a c\`adl\`ag function that satisfies (\ref{zetadef}) and (\ref{codecond}); let $\Pi\! = \! ((s_i, t_i))_{1\leq i \leq p}$ where $0 \leq s_i \! \leq \! t_i \! < \! \zeta_h$, for all $1\! \leq \! i \! \leq \! p$ and let $\epp \ino [0, \infty) $. 
Then, the \textit{compact measured metric space coded by $h$ and the pinching setup $(\Pi, \epp)$} is the $(\Ptt, \epp)$-pinched metric space associated with $(T_h, d_h)$ and the pinching points %\margmm{pinching}
$\Ptt\! = \!  ((p_h (s_i), p_h (t_i)))_{1\leq i\leq p}$, where $p_h \! : \! [0, \zeta_h ) \! \rightarrow \! T_h$ stands for the canonical projection. 
%\marginpar{\footnotesize\mm{Slight permutation in the presentation}}
%We shall use the following notation:  
%Recall that $p_h \! : \! [0, \zeta_h ) \! \rightarrow \! T_h$ is the canonical projection. We then take 
%$\Ptt \! = \! ((p_h (s_i), p_h (t_i)))_{1\leq i\leq p}$ as the set of pinching points of $T_h$ and 
%the \textit{compact measured metric space coded by $h$ and the pinching setup ($\Pi, \epp$)} is then denoted by 
%\begin{equation}
%\label{defgrgraf}
%G (h, \Pi, \epp)= \big( G_{h, \Pi, \epp} , d_{h, \Pi, \epp}, \varrho_{h, \Pi, \epp} , m_{h, \Pi, \epp} \big).
%\end{equation}
We shall denote by $p_{h, \Pi, \epp}$ the composition of the canonical projections 
$\varpi_{\Ptt, \epp} \circ  p_h \! :  \! [0, \zeta_h) \! \rightarrow \! G_{h, \Pi, \epp}$;  then 
$\varrho_{h, \Pi, \epp}\! = \! p_{h, \Pi, \epp} (0)$ and 
$m_{h, \Pi, \epp}$ stands for the pushforward measure of the Lebesgue on $[0, \zeta_h )$ via $p_{h, \Pi, \epp}$. We shall use the following notation:  
%Recall that $p_h \! : \! [0, \zeta_h ) \! \rightarrow \! T_h$ is the canonical projection. We then take 
%$\Ptt \! = \! ((p_h (s_i), p_h (t_i)))_{1\leq i\leq p}$ as the set of pinching points of $T_h$ and 
%the \textit{compact measured metric space coded by $h$ and the pinching setup ($\Pi, \epp$)} is then denoted by 
\begin{equation}
\label{defgrgraf}
G (h, \Pi, \epp)= \big( G_{h, \Pi, \epp} , d_{h, \Pi, \epp}, \varrho_{h, \Pi, \epp} , m_{h, \Pi, \epp} \big).
\end{equation}

%
%\smallskip
%
%\begin{compactenum}
%%\item[$-$] $(G_{h, \Pi, \epp} , d_{h, \Pi, \epp})$ is the $(\Ptt, \epp)$-pinched metric space associated with $(T_h, d_h)$; 
%
%
%\smallskip
%
%
%
%\item[$-$] Denote by $p_{h, \Ptt, \epp}$ the composition of the canonical projections 
%$\varpi_{\Ptt_h, \epp} \circ  p_h \! :  \! [0, \zeta_h) \! \rightarrow \! G_{h, \Ptt, \epp}$;  then 
%$\varrho_{h, \Ptt, \epp}\! = \! p_{h, \Ptt, \epp} (0)$ and 
%$m_{h, \Ptt, \epp}$ stands for the pushforward of the Lebesgue on $[0, \zeta_h )$ via $p_{h, \Ptt, \epp}$. 
%\end{compactenum}
%
%   
 
\medskip 
 
\noi 
\textbf{Convergence of metric spaces.}
%\subsubsection{Convergence of metric spaces.}
%\label{codetspa} 
%\label{GHPdist}
Let $(G_1, d_1, \rho_1, m_1)$ and $(G_2, d_2, \rho_2, m_2)$ be two pointed compact measured metric spaces. The pointed \textit{Gromov-Hausdorff-Prokhorov distance} of $G_1$ and $G_2$ is then defined by 
\begin{multline}
\label{defGHP}
\bdelta_{\mathrm{GHP}} (G_1, G_2)\! = \! \inf \!  \Big\{ d_{E}^{\textrm{Haus}} \! \big(\phi_1 (G_1), \phi_2 (G_2)\big)  \\ 
+ 
d_E (\phi_1 (\rho_1), \phi_2 (\rho_2))  + d_{E}^{\textrm{Proh}} \big(m_1 \! \circ \! \phi^{-1}_1\! \! , m_2 \! \circ \! \phi^{-1}_2\big)  \Big\}. 
\end{multline}
\textit{Here}, the infimum is taken over all Polish spaces $(E, d_E)$ and all isometric embeddings $\phi_i: G_i
\hookrightarrow E$, $i\ino \{ 1, 2\}$; 
$d_{E}^{\textrm{Haus}}$ stands for the Hausdorff distance on the space of compact subsets of $E$, $d_{E}^{\textrm{Proh}}$ stands for the Prokhorov distance on the space of finite Borel measures on $E$ and for all 
$i\ino \{ 1, 2\}$, $m_i \! \circ \! \phi^{-1}_i$ stands for the pushforward measure of $m_i$ via $\phi_i$. 

  We recall from Theorem 2.5 in Abraham, Delmas \& Hoscheit \cite{AbDeHo13} 
the following assertions: $\bdelta_{\mathrm{GHP}} $ is symmetric and it satisfies the triangle inequality; 
$\bdelta_{\mathrm{GHP}} (G_1, G_2)\! = \! 0$ if and only if  $G_1$ and $G_2$ are \textit{isometric}, namely if and only if  there exists a bijective isometry $\phi \! : \! G_1 \! \rightarrow \! G_2$ such that $\phi(\rho_1)\! = \! \rho_2$ and such that $m_2 \! = \!  m_1 \circ  \phi^{-1}$. Denote by $\bbG$ the isometry classes of pointed compact measured metric spaces. Then, we recall the following result. 
\begin{thm}[Theorem 2.5  \cite{AbDeHo13} ]
\label{rappGHP} $(\bbG, \bdelta_{\mathrm{GHP}})$ is a complete and separable metric space. 
\end{thm}
Actually in our paper, weak-limits are proved for coding functions, which entail $\bdelta_{\mathrm{GHP}}$-limits 
as asserted by the following lemma. 
\begin{lem}
\label{codconGHP} Let $h, h^\prime\! :\!  [0, \infty) \! \rightarrow \! [0, \infty)$ be two c\`adl\`ag functions such that $\zeta_h$ and $\zeta_{h^\prime}$ are finite and that satisfy \eqref{codecond}. 
 Let $\Pi\! = \! ((s_i, t_i))_{1\leq i \leq p}$ and $\Pi^\prime \! = \! ((s^\prime_i, t^\prime_i))_{1\leq i \leq p}$ be two sequences such that 
$0 \! \leq \! s_i \! \leq \! t_i\! < \! \zeta_h $ and $0 \! \leq \! s^\prime_i \! \leq \! t^\prime_i\! < \! \zeta_{h^\prime}$, and let $\delta \ino (0, \infty)$ be such that 
\begin{equation}
\label{hypodelt}
 \forall i \ino \{ 1, \ldots , p\}, \quad \lvert s_i \! -\! s^\prime_i \rvert \! \leq \! \delta   \quad \textrm{and} \quad   \lvert t_i \! -\! t^\prime_i \rvert \! \leq \! \delta \; .
\end{equation} 
Let 
$\epp, \epp^\prime \ino [0, \infty)$ and recall from (\ref{defgrgraf}) the definition of the pointed compact measured metric spaces $G\! := \! G(h, \Pi, \epp)$ and 
$G^\prime\! := \! G(h^\prime, \Pi^\prime, \epp^\prime)$. Then, we get: 
\begin{equation}
\label{contGHP}
\bdelta_{\mathrm{GHP}} (G, G^\prime) \leq 6 (p+1) \big( \lVert h \! -\! h^\prime \rVert_{\infty} +  \omega_{\delta} (h)  \big) + 
3p (\epp \!  \vee \! \epp^\prime) + \lvert \zeta_h \! -\! \zeta_{h^\prime} \rvert \; , 
\end{equation}
where $\omega_\delta (h)\! = \! \max \big\{  \lvert h(s) \! -\! h(t) \rvert :\; s, t \ino [0, \infty)\!  : |s\! -\! t| \! \leq \! \delta \big\}$ and where $ \lVert \cdot \rVert_{\infty}$ stands for the uniform norm on $[0, \infty)$.  
\end{lem}
\textbf{Proof.} See  Appendix Section \ref{cvGHPpf}. The proof is 
partly adapted from Theorem 2.1 in Le Gall \& D.~\cite{DuLG05},  Proposition 2.4 Abraham, Delmas \& Hoscheit \cite{AbDeHo13} and Lemma 21 in Addario-Berry, Goldschmidt \& B.~in \cite{AdBrGo12a}. \cqfd

\medskip
   
\noi
\textbf{Limit theorems for multiplicative graphs.} 
 Recall from \eqref{defGwdir} the $\bw_n$-multiplicative graph $\cG_{\!\bw_n}$. We equip its vertex set with a measure $\bm_{\bw_n} \! = \! \sum_{1\leq j\leq \mathbf{j}_n} 
 w^{_{(n)}}_{^j} \delta_{j}$. Recall from (\ref{excdibor}) the definition of the excursion intervals $[l^{\bw_n}_k, r^{\bw_n}_k)$, $1\! \leq\! k\! \leq \! \mathbf{q}_{\bw_n}$, of $\cH^{\bw_n}$ above $0$, recall from  (\ref{excdiscret}) the definition of the corresponding excursions $\mathtt{H}^{\bw_n}_k (\cdot)$ of $\cH^{\bw_n}$ above $0$ and recall from (\ref{defPiwk}) the corresponding sets of pinching times $\Ptt^{\bw_n}_k$. We recall that each excursion $\mathtt{H}^{\bw_n}_k (\cdot)$ corresponds to a connected component $\cG_{\! k}^{\bw_n}$ of $\cG_{\!\bw_n}$ and we have $\bm_{\bw_n} \big( \cG_{\! k}^{\bw_n} \big)\! = \! \zeta^{\bw_n}_k\! = \! r^{\bw_n}_k \! -\! l^{\bw_n}_k$. Thus, we get 
  \begin{equation}
\label{order1bis}
\bm^{\bw_n} (\cG_{\! 1}^{\bw_n}) \! \geq \! \ldots  \! \geq \! 
\bm^{\bw_n} (\cG_{\! \mathbf{q}_{\bw_n} }^{\bw_n}) .
\end{equation}
Then, $\cG_{\! k}^{\bw_n}$ is the pinched (measured pointed) metric space coded by $(\Htt^{\bw_n}_k,  \Ptt_k^{ \bw_n})$. Namely, 
\begin{equation}
\label{disctry}
G ( \Htt^{\bw_n}_k,  \Ptt_k^{ \bw_n}, 1) \; \, \textrm{is isometric to} \; \, (\cG_{\! k}^{\bw_n} , d^{\bw_n}_k, \varrho_k^{\bw_n}, \bm^{\bw_n}_k) \; ,
\end{equation}
Thus, they 
define the same random element in the space $\bbG$ of the isometry classes of pointed compact measured metric 
spaces equipped with the Gromov-Hausdorff-Prokhorov distance $\bdelta_{\mathrm{GHP}}$ defined in (\ref{defGHP}). Here, we have denoted 
by $ d^{\bw_n}_k$ the graph-distance, by $\varrho_k^{\bw_n}$ the first vertex explored via the LIFO coding and $\bm^{\bw_n}_k$ stands for the restriction to $\cG_{\! k}^{\bw_n} $ of the weight measure $\bm_{\bw_n}$. Since $\mathbf{q}_{\bw_n}$ tends to $\infty$, it is convenient to extend the sequence $(\cG_{\! k}^{\bw_n})_{1\leq k\leq \mathbf{q}_{\bw_n}}$ by taking $\cG_{\! k}^{\bw_n}$ equal to the point space equipped with the null measure for all $k \! >\!\mathbf{q}_{\bw_n}$.    

Similarly, recall from (\ref{excuHY}) the definition of the excursion intervals $(l_k, r_k)$, $k\! \geq \! 1$, of 
$\cH$ above $0$. Recall from (\ref{exccont}) the definition of the excursion $\mathtt{H}_k (\cdot)$ of $\cH$ above $0$ and recall from (\ref{defPik}) 
the definition of the set of pinching times $\Ptt_k$. We recall from (\ref{contGdef}) the definition of the continuous $(\alpha, \beta, \mathbf{c}, \kappa)$-multiplicative graph $\mathbf{G}\! = \! ((\mathbf{G}_k, d_k, \varrho_k, \bm_k ))_{k\geq 1}$, where for all $k\! \geq \! 1$, $\mathbf{G}_k$ is the pinched (measured pointed) metric space coded by $(\mathtt{H}_k, \Ptt_k, 0)$. Namely, 
\begin{equation}
\label{redefGmu}
\mathbf{G}_{ k}:= G ( \Htt_k,  \Ptt_k, 0)
\end{equation}
Thus, they 
define the same random element in the space $\bbG$.  
Then, Theorem \ref{excucvth} and Lemma \ref{codconGHP} entail the following theorem. 
\begin{thm}
\label{graphcvth} Under the same assumptions as in Theorem \ref{HYcvth}, the following convergence 
\begin{equation}
\label{graphbwn}
\big(\big( \cG_{\! k}^{\bw_n} ,  \frac{_{_{a_n}}}{^{^{b_n}}}d_{k}^{\bw_n} , \varrho_k^{\bw_n}, \frac{_{_{1}}}{^{^{b_n}}}\bm_k^{\bw_n}  \big) \big)_{k\geq 1}
\;  \underset{n\rightarrow \infty}{-\!\!\! -\!\!\! -\!\!\! \longrightarrow} \; \big(\big( \mathbf{G}_{k} , \mathrm{d}_{k}, \varrho_{k} , \bm_{k} \big) \big)_{k\geq 1}
\end{equation} 
holds weakly on $\bbG^{\bbN^*}$ equipped with the product topology. Denote by $\bmu^{\bw_n}_k\! = \! \sum_{j\in  \cG_{\! k}^{\bw_n}} \delta_j$ the counting measure on $ \cG_{\! k}^{\bw_n}$. Then, the following convergence  
\begin{equation}
\label{globabab}
\big(\big( \cG_{\! k}^{\bw_n} ,  \frac{_{_{a_n}}}{^{^{b_n}}}d_{k}^{\bw_n} , \varrho_k^{\bw_n}, \frac{_{_{1}}}{^{^{b_n}}}\bmu_k^{\bw_n}  \big) \big)_{k\geq 1}
\;  \underset{n\rightarrow \infty}{-\!\!\! -\!\!\! -\!\!\! \longrightarrow} \; \big(\big( \mathbf{G}_{k} , \mathrm{d}_{k}, \varrho_{k} , \bm_{k} \big) \big)_{k\geq 1}
\end{equation} 
holds weakly on $\bbG^{\bbN^*}$ equipped with the product topology. 

Recall notation 
$\mathbf{j}_n\! = \! \max \{ j\! \geq \!1:  w^{_{(n)}}_{^j} \! >\! 0\}$.  
If we furthermore assume that $\sqrt{\mathbf{j}_n}/b_n \! \rightarrow \! 0$, then (\ref{globabab}) holds when the connected components are listed in the decreasing order of their numbers of vertices: namely, when $ \bmu_1^{\bw_n} \big( \cG_{\! 1}^{\bw_n}\big)  \geq \ldots \geq \bmu_{q_{\bw_n}}^{\bw_n} \big( \cG_{\! q_{\bw_n}}^{\bw_n}\big)$.  
\end{thm}
\noi
\textbf{Proof.} See Section \ref{graphcvpf}. \cqfd  

\begin{rem}
\label{hypsuppl} The assumption $\sqrt{\mathbf{j}_n}/b_n \! \rightarrow \! 0$ may not be optimal for  (\ref{globabab}) to hold when the connected components are listed in the decreasing order of their 
numbers of vertices. However, for all 
$\alpha\in \bbR$, $\beta  \ino [0, \infty)$, $\kappa \ino (0, \infty)$ and 
$\mathbf{c}\! = \! (c_j)_{j\geq 1}\ino \elldo_3$ satisfying (\ref{contH}), this statement is never void since the examples of $(a_n, b_n, \bw_n)$ provided in 
Proposition \ref{Hcritos} $(iii)$ satisfy $\sqrt{\mathbf{j}_n}/b_n \! \rightarrow \! 0$. Moreover, let us mention that all the cases that have been considered previously by other authors satisfy this assumption, as it is pointed out in the next Section \ref{prevresults}. \cq  
\end{rem}
%
%\medskip
%
%A similar result holds for the graphs $\cG^{\bw_n}$ equipped with the empirical measure on the set of  
%vertices.  
%\begin{thm}
%\label{graphcvth} We make the same assumption as in Theorem \ref{HYcvth}. 
%
% the following convergence 
%\begin{equation}
%\label{graphbwn}
%\big(\big( \cG_{\! k}^{\bw_n} ,  \frac{_{_{a_n}}}{^{^{b_n}}}d_{k}^{\bw_n} , \varrho_k^{\bw_n}, \frac{_{_{1}}}{^{^{b_n}}}\bm_k^{\bw_n}  \big) \big)_{k\geq 1}
%\;  \underset{n\rightarrow \infty}{-\!\!\! -\!\!\! -\!\!\! \longrightarrow} \; \big(\big( \mathbf{G}_{k} , \mathrm{d}_{k}, \varrho_{k} , \bm_{k} \big) \big)_{k\geq 1}
%\end{equation} 
%holds weakly on $\bbG^{\bbN^*}$ equipped with the product topology. 
%\end{thm}
%
%
 
\subsection{Connections with previous results.}
\label{prevresults}  
\noindent 
\textbf{Entrance boundary of the multiplicative coalescent.} The model of $\bw$-multiplicative random graphs appears in the work of Aldous \cite{Al97} as an extension of Erd\H{o}s--R\'enyi random graphs that have close connections with multiplicative coalescent processes. 
Relying upon this connection, Aldous and Limic determine in \cite{AlLi98} the extremal eternal versions of the multiplicative coalescent in terms of the excursion lengths of L\'evy-type processes $Y$ (up to rescaling, as explained below); to that end, they 
consider in Proposition 7 \cite{AlLi98} asymptotics of the masses of the connected components of sequences of multiplicative random graphs. The asymptotic regime in Proposition 7 \cite{AlLi98} is very close to Assumptions (\ref{apriori}) and ($\mathbf{C1}$) -- ($\mathbf{C3}$) in our Theorem \ref{graphcvth}. 

Let us briefly recall Proposition 7 in \cite{AlLi98} since it is used in the proof of Theorem \ref{graphcvth}. 
Aldous \& Limic fix a sequence of weights $\mathtt{x}_n\ino \elldo_{\! f}$, $n\ino \bbN$, and their notation for 
multiplicative graphs is the following: let $(\xi_{i,j})_{j>i\geq 1}$ be an array of independent %\margmm{$j>i\ge 1$}
and exponentially distributed r.v.~with unit mean; let $N(\mathtt{x}_n)\! = \! \max \{ j\! \geq \! 1: x^{_{(n)}}_{^j} \! >\! 0\}$; then for all $q\ino [0, \infty)$, Aldous \& Limic consider the random graph $G_n (q)$ whose set of vertices is $\cV (G_n (q))\! = \! 
\{ 1, \ldots , N(\mathtt{x}_n)\}$ and whose set of edges $\ccE (G_n (q))$ is such that 
$\{ i,j\} \ino \ccE (G_n (q))$ if and only if  $\xi_{i,j} \! \leq \! q x^{_{(n)}}_{^i}\!  x^{_{(n)}}_{^j}\! $; the multiplicative  graph $G_n(q)$ is equipped with the measure $m_{n}\! = \! \sum_{j\geq 1}  x^{_{(n)}}_{^j} \delta_j$; let 
$\zeta_1 (\mathtt{x}_n , q)\! \geq \! \ldots  \! \geq \zeta_k (\mathtt{x}_n, q)\! \geq \ldots $ stand for the (eventually null) sequence of the $m_n$-masses of the connected components of $G_n (q)$. Then, it is easy to check that $\mathbf{X}_n \! : \! q\! \mapsto\! (
\zeta_k (\mathtt{x}_n, q))_{k\geq 1}$ is a multiplicative coalescent process with finite support. 
Aldous \& Limic describe the limit of the processes $\mathbf{X}_n$ in terms of the excursion-lengths of a process $(W^{\kappa_{\textrm{AL}}, -\tau_{\textrm{AL}}, \mathbf{c}_{\textrm{AL}}  }_s)_{s\in [0, \infty)}$ whose law is characterized by 
three parameters: $\kappa_{\textrm{AL}}\ino [0, \infty)$, $\tau_{\textrm{AL}} \ino \bbR$ and $ \mathbf{c}_{\textrm{AL}}\ino \elldo_3$; this process is 
connected to the $(\alpha, \beta , \kappa, \mathbf{c})$-process $Y$ defined in (\ref{plasouilloc}) as follows: 
\begin{equation}
\label{bungabunga}
\forall s\ino [0, \infty), \quad W^{\kappa_{\textrm{AL}}, -\tau_{\textrm{AL}}, \mathbf{c}_{\textrm{AL}}  }_s  = Y_{s/ \kappa}, \; \textrm{where} \quad \kappa_{\textrm{AL}} = \frac{\beta}{\kappa}, \quad \tau_{\textrm{AL}}= \frac{\alpha}{\kappa} \quad \textrm{and} \quad \mathbf{c}_{\textrm{AL}}=\mathbf{c}. 
\end{equation}
Proposition 7 \cite{AlLi98} assumes the following: 
\begin{equation}
\label{assAL}   
\lim_{n\rightarrow \infty} \frac{\sigma_3 (\mathtt{x}_n)}{(\sigma_2 (\mathtt{x}_n))^3}\! =\!  \kappa_{\textrm{AL}}+ 
\sigma_3 (\mathbf{c}_{\textrm{AL}}), \quad \forall j\ino \bbN^*\! , \;  \lim_{n\rightarrow \infty} \frac{x^{_{(n)}}_{^j}}{\sigma_2 (\mathtt{x}_n)}
\! = \! c^{\textrm{AL}}_j \quad \textrm{and} \quad  \lim_{n\rightarrow \infty}\sigma_2 (\mathtt{x}_n)\! = \! 0 , 
\end{equation}
%\margmm{Since $W^{\kappa_{\textrm{AL}}, -\tau_{\textrm{AL}}, \mathbf{c}_{\textrm{AL}}  }_s  = Y_{s/ \kappa}$, we have $\zeta^{\textrm{AL}}_k=\kappa \zeta_k$.}
and it asserts that for all $\tau_{\textrm{AL}}\ino \bbR$, 
$\mathbf{X}_n ( \sigma_2 (\mathtt{x}_n)^{-1} \! -\!  \tau_{\textrm{AL}})\rightarrow (\zeta_k)_{k\geq 1}$, weakly in $\elldo_2$, where $(\zeta_k)_{k\geq 1}$ are the excursion-lengths 
of $ W^{\kappa_{\textrm{AL}}, -\tau_{\textrm{AL}}, \mathbf{c}_{\textrm{AL}}  }$ 
above its infimum, listed in the decreasing order. 

  Assumptions (\ref{assAL}) are close to ($\mathbf{C2}$) and 
  ($\mathbf{C3}$). More precisely, let $(\alpha, \beta, \kappa, \mathbf{c})$ be connected with $\kappa_{\textrm{AL}}$, $\tau_{\textrm{AL}}$ 
and $\mathbf{c}_{\textrm{AL}}$ as in (\ref{bungabunga}); let 
$a_n , b_n \ino (0, \infty)$ and $\bw_n \ino \elldo_f$ satisfy (\ref{apriori}) and ($\mathbf{C1}$) -- ($\mathbf{C3}$); then, set:  
$$\forall j \ino \bbN^*, \quad  x^{_{(n)}}_{^j}= \frac{\kappa w^{_{(n)}}_{^j} }{{b_n}}  \quad \textrm{and} \quad \tau^n_{\textrm{AL}}= \frac{b^2_n }{\kappa^2 \sigma_2 (\bw_n)} \Big(1\! -\! \frac{\sigma_2 (\bw_n)}{\sigma_1 (\bw_n)} \Big)
\underset{n\rightarrow \infty}{-\!\!\!-\!\!\! -\!\!\! \longrightarrow}\;  \frac{\alpha}{\kappa}\! = \! \tau_{\textrm{AL}} .  $$
Recall from (\ref{defGwdir}) the definition of $\cG_{\! \bw_n}$, the $\bw_n$-multiplicative graph. Recall that $\bm_{\bw_n}\! = \! \sum_{j\geq 1} 
w^{_{(n)}}_{^j} \delta_j$. Recall from Section \ref{ccvmulsec} that the $\cG^{\bw_n}_{\! k}$ stand for the connected components of $\cG_{\! \bw_n}$ listed in the nonincreasing order of their $\bm_{\bw_n}$-mass. Then, it is easy to check the following. 
\begin{equation}  
\label{macarena}
G_n \big(\sigma_2 (\mathtt{x}_n)^{-1} \! -\!  \tau^n_{\textrm{AL}} \big) = \cG_{\! \bw_n} \quad \textrm{and} \quad 
\zeta_k \big(\mathtt{x}_n\, , \, \sigma_2 (\mathtt{x}_n)^{-1} \! -\!  \tau^n_{\textrm{AL}} \big)\! = \! \frac{\kappa}{b_n} \bm_{\bw_n} \big(\cG^{\bw_n}_{\! k} \big)=:  \kappa\zeta_k^{n} . 
\end{equation}
%\margmm{See my previous comment. } 
Note that the $\zeta^{n}_{k}$  are the excursion-lengths of 
$(\frac{1}{{a_n}} Y^{\bw_n}_{ b_n t  })_{t\in [0, \infty)}$ above its infimum. Since $\tau^n_{\textrm{AL}}\! \rightarrow 
\alpha / \kappa$ and since multiplicative coalescent processes have no fixed time-discontinuity, Proposition 7 in \cite{AlLi98} immediately entails the following proposition that is used in Section \ref{excucvpf} to prove Theorems \ref{excucvth} and \ref{graphcvth}. 
\begin{prop}[Proposition 7 \cite{AlLi98}]
\label{AldLim98}
Let $a_n , b_n \ino (0, \infty)$ and $\bw_n \ino \elldo_f$ satisfy (\ref{apriori}) and $(\mathbf{C1})$--$(\mathbf{C3})$, with $\alpha\ino \bbR$, $\beta  \ino [0, \infty)$, $\kappa \! \in \! (0, \infty)$ and $\mathbf{c}\ino \elldo_3$. 
Recall from (\ref{defYcHH}) (resp.~from (\ref{plasouilloc})) the definition of $Y^{\bw_n}$ (resp.~of $Y$). 
Let $(\zeta^n_k)_{1\leq k \leq \mathbf{q}_{\bw_n}}$ (resp.~$(\zeta_k)_{k\geq 1}$) be the excursion-lengths of $(\frac{1}{{a_n}} Y^{\bw_n}_{ b_n t })_{t\in [0, \infty)}$ (resp.~of $Y$) above its infimum.  
Then,  
\begin{equation}
\label{AldLiml2}
\big( \zeta^n_k\big)_{1\leq k \leq \mathbf{q}_{\bw_n}} \overset{\textrm{weakly in $\elldo_2$}}{\underset{n\rightarrow \infty}{-\!\!\! -\!\!\! -\!\!\! \longrightarrow}}\; (\zeta_k)_{k\geq 1} .
\end{equation}
%weakly in $\elldo_2$. 
\end{prop}

\medskip

\medskip

\noindent 
\textbf{Limits of Erd\H{o}s--R\'enyi graphs in the critical window.} The first result proving metric convergence in a strong Hausdorff sense of rescaled Erd\H os-R\'enyi graphs and their inhomogeneous extensions is due to 
Addario-Berry, Goldschmidt \& B.~in \cite{AdBrGo12a}. 
In this paper, they study the scaling limits of the largest components of Erd\H{o}s--R\'enyi random graph $G(n, p_n)$ in the critical window $p_n\! =\! n^{-1}\! -\! \alpha n^{-4/3}$,
%$ \frac{1}{n}- \frac{\alpha}{n^{4/3}}$
with $\alpha\ino \bbR$, which corresponds to the graph $\cG_{\bw_n}$ where
% the elements of $\bw_n\! = \! (w^{_{(n)}}_j)_{j\geq 1} $ are given by 
$w^{_{(n)}}_j \! =\! \un_{\{ j\leq n \}} n\log(\frac{1}{1 - p_n})$, $j\! \geq \! 1$.  
%and $w^{_{(n)}}_j\! = \! 0$, for all $j\! >\! n$. 
Taking, $a_n \! =\! n^{1/3}$ and $b_n\! =\! n^{2/3}$, we immediately see that $a_n , b_n$ and $\bw_n$ satisfy (\ref{apriori}) with $\kappa\!  =\! \beta_0 \! = \! 1$, 
$(\mathbf{C1})$, $(\mathbf{C2})$, $(\mathbf{C3})$ and $\sqrt{\mathbf{j_n}}/b_n\! = \! n^{-1/6}\! \rightarrow \! 0$, with the parameters $\alpha \in\bbR$, $\beta\! = \! 1$ 
and $\mathbf{c}\! = \! 0$. Namely, the branching mechanism is $\psi (\lambda) \! = \! \alpha \lambda + \frac{_1}{^2} \lambda^2$. Since $\beta_0\! >\! 0$, Proposition \ref{Hcritos} $(i)$ implies that 
$(\mathbf{C4})$ is automatically satisfied and Theorem \ref{graphcvth} applies: 
in this case, Theorem \ref{graphcvth} is a weaker version of Theorem 2 in Addario-Berry, Goldschmidt \& B.~\cite{AdBrGo12a}, p.~369: the result in \cite{AdBrGo12a} actually provides precise estimates on the size of metric components.  
Let us mention that  \cite{AdBrGo12a} also contains tail-estimates on the diameters of the small components. Such estimates seem difficult to obtain in the case of general $\bw_n$. 

\medskip

%\noindent\mm{Add literature on uniform measures--IN PROGRESS.}

\noindent 
\textbf{Multiplicative graphs in the same basin of attraction as Erd\H os-R\'enyi graphs.}
%\margmm{new}
Bhamidi, van der Hofstad \& van Leeuwaarden in \cite{BhHoLe10} prove the scaling limit of the component sizes (number of vertices) for examples of  multiplicative graphs which behave asymptotically like the Erd\H os-R\'enyi graphs. Bhamidi, Sen \& X.~Wang in \cite{BhSeWa14} and Bhamidi, Sen, X.~Wang \& B.~in \cite{BhBrSeWa14}
investigate instead the scaling limits of these graphs seen as measured metric spaces. The conditions under which these authors prove their limit theorems slightly differ. We give here a detailed account of these conditions so as to make a connection with our results.  In all the cases covered by \cite{BhHoLe10, BhSeWa14, BhBrSeWa14}, the scalings appear to be $a_n \! =\! n^{1/3}$, $b_n\! =\! n^{2/3}$ and $\bw_n$ is a sequence of length $n$ having the following asymptotic behaviour: 
%In \cite{BhHoLe10}, the authors assume the following asymptotic behaviour on $(\bw_n)_{n\ge 1}$. Let $W: \Omega\to [0, \infty)$ be a r.v.~such that $\mathbf E[W]=\mathbf E[W^2]$ and $\mathbf E[W^{3}]<\infty$. Suppose that  
%\begin{align*}
%& \frac{w^{(n)}_1}{n^{1/3}}\to 0 \quad \text{and} \quad \tfrac{1}{n} \sum_{i}\mathbf 1_{\{w^{(n)}_{i}\le x\}} \to \mathbf P(W\le x) \quad \text{for all } x\ge 0,\\
%& \sigma_i (\bw_n)\! = \! n \mathbf E[W^{i}] \! +\!  o(n^{2/3}), \; i\ino \{ 1, 2\} \; \textrm{and} \;  
%\sigma_3 (\bw_n)\! = \! n \mathbf E[W^{3}]  \!+ \!  o(n).  
%\end{align*}
%In \cite{BhSeWa14}, the authors consider the following (more general) conditions:}
%More precisely,  
%they consider the cases where 
%Bhamidi, Sen \& X.~Wang in \cite{BhSeWa2014} 
%prove the Gromov--Hausdorff convergence of the largest components of $\mathcal{G}_{\bw_n}$ in 
%the following regime: 
%$a_n \! =\! n^{1/3}$, $b_n\! =\! n^{2/3}$ and where $\bw_n \ino \elldo_{\! f}$ is such that $w^{_{(n)}}_j \! =\! 0$ for all $j\! >\! n$ and 
\begin{equation}
\label{bassBro}
\frac{w^{(n)}_1}{n^{1/3}}\to 0, \quad \exists \, \sigma,  \sigma^\prime \ino (0, \infty)\! : \, \sigma_i (\bw_n)\! = \! n \sigma \! +\!  o(n^{2/3}), \; i\ino \{ 1, 2\} \; \textrm{and} \;  
\sigma_3 (\bw_n)\! = \! n \sigma^\prime  \!+ \!  o(n).  
\end{equation}
%
%
%\begin{compactenum}
%%\item[$(i)$] $w^{_{(n)}}_j \! =\! 0$ if $j\! >\! n$; 
%
%\smallskip
%
%
%\item[$(i)$] $\exists \, \sigma,  \sigma^\prime \ino (0, \infty): \; \sigma_i (\bw_n)\! = \! n \sigma + o(n^{2/3})$, $i\ino \{ 1, 2\}$, 
%and  $\sigma_3 (\bw_n)\! = \! n \sigma^\prime  + o(n)$.  
%
%\smallskip
%
%\item[$(ii)$] $\exists \, \eta_0 \ino (0, 1/6): \; \max_{j\geq 1} w^{_{(n)}}_j \! =\! o(n^{1/6-\eta_0})$ and $\exists \, \gamma_0\ino (0, \infty) : \; 1/ \min_{j\geq 1} w^{_{(n)}}_j \! =\! o(n^{\gamma_0})$. 
%
%\smallskip
%
%\end{compactenum}
%\margmm{not really relevant here}
%\ms{Observe that $\sigma^\prime \! \geq \! \sigma$ since $\sigma_2 (\bw_n) \! \leq \! \sqrt{\sigma_3 (\bw_n)} \sqrt{\sigma_1 (\bw_n)}$ by Cauchy-Schwarz. 
%By Cauchy-Schwarz again, $\sigma_1(\bw_n)\le \sqrt{n\sigma_2(\bw_n)}$, which implies $ \sigma \! \leq \! 1$.}
For all $\alpha \in\bbR$, set 
$$ \bw_n (\alpha)\! =\!  \big(1\! -\! \alpha n^{-\frac{1}{3}} \big) \bw_n \! =\!  
\big( \big(1\! -\! \alpha n^{-\frac{1}{3}} \big) w^{_{(n)}}_j  \big)_{j\geq 1} \; .$$
This is a situation covered by Theorem \ref{graphcvth}. Indeed, (\ref{bassBro}) easily implies that $a_n , b_n , \bw_n(\alpha)$ satisfy (\ref{apriori}), $(\mathbf{C1})$, $(\mathbf{C2})$, $(\mathbf{C3})$, $\sqrt{\mathbf{j_n}}/b_n\! = \! n^{-1/6}\! \rightarrow \! 0$, with the parameters $\alpha \in\bbR$, $\beta_0\! = \! 1$, $\beta \! = \! \sigma^\prime / \sigma $, $\kappa \! = 1/ \sigma $ and $\mathbf{c}\! = \! 0$. 
Thus, the branching mechanism is $\psi (\lambda)\! = \! \alpha \lambda + \frac{1}{2} \frac{\sigma^\prime}{\sigma} \lambda^2$. 
%\margmm{new; see footnote}
Since $\beta_0\! = \! 1$, Proposition \ref{Hcritos} $(i)$ implies $(\mathbf{C4})$. Then, Theorem \ref{graphcvth} applies in this setting, which allows us to extend 
\begin{compactenum}[-]
\item[$-$] Theorem 1.1 in \cite{BhHoLe10}, which has been proved under the supplementary assumption that there exists a r.v.~$W: \Omega\to [0, \infty)$ such that 
\[
\tfrac{1}{n} \sum_{i}\mathbf 1_{\{w^{(n)}_{i}\le x\}} \to \mathbf P(W\le x) \quad \text{for all } x\ge 0, \text{ and } \sigma=\mathbf E[W]=\mathbf E[W^{2}], \sigma'=\mathbf E[W^{3}].
\]
(Assumption (b) in \cite{BhHoLe10}. ) 
\item[$-$] Theorem 3.3 in Bhamidi, Sen \& X.~Wang in \cite{BhSeWa14} that has been proved by quite different methods and 
under two additional technical assumptions (Assumptions 3.1 (c) and (d)).
\end{compactenum}
%\fmm{remark on the $\ell^{2}$-convergence for the sizes of the components? }
%as well as 
%Theorem 3.3 in Bhamidi, Sen \& X.~Wang in \cite{BhSeWa14} that have been proved by quite different methods and 
%under two additional technical assumptions (Assumptions 3.1 (c) and (d)). 
Turova in \cite{Turova} also proved a result similar to Theorem 1.1 of \cite{BhHoLe10} for i.i.d.~random weight sequences. 
Let us mention that the convergence under the sole assumptions (\ref{bassBro}), that we proved, has been conjectured in \cite{BhSeWa14}, Section 5, part (c). 

\medskip

\noindent
\textbf{Gromov--Prokhorov convergence of multiplicative graphs without Brownian component.} 
In light of the above mentioned result of Aldous \& Limic \cite{AlLi98} on the convergence of the component masses of the multiplicative graph in   the asymptotic regime \eqref{assAL}, it is natural to expect that the graph itself should also converge in some sense. 
The first affirmation in this direction is due to Bhamidi, van der Hofstad and Sen who prove  the following in \cite{BhHoSe15}:  Denote by $\mathscr C_i(q)$ the $i$-largest (in $m_n$-mass) connected component of $G_n(q)$, that is,  $m_n(\mathscr C_i(q))=\zeta_i(\mathtt{x}_n, q)$. Equip each component $\mathscr C_i(-\tau_{\textrm{AL}}+\sigma_2(\mathtt{x}_n)^{-1})$ with its graph distance rescaled by $\sigma_2(\mathtt{x}_n)$ and with the mass measure $m_n$, 
%\margmm{corrected}
they prove that under \eqref{assAL} with $\kappa_{\mathrm{AL}}=0$, the collection of rescaled metric spaces converge in the sense of Gromov--Prokhorov topology to a collection of measured metric spaces, which are not necessarily compact. They also give an explicit construction of the limiting spaces based upon a model of continuum random tree called ICRT.  
The Gromov--Prokhorov convergence is equivalent to the convergence of mutual distance of an i.i.d.~sequence with law $m_n$, which is weaker than the Gromov--Hausdorff--Prokhorov that we 
obtain in Theorem \ref{graphcvth} under the compactness assumption $\int^\infty d\lambda / \psi (\lambda)$. Let us mention that our approach via coding processes is quite distinct from that of 
Bhamidi, van der Hofstad \& Sen in \cite{BhHoSe15}.

\medskip

\noindent 
\textbf{Power-law cases.} We extend the power-law cases investigated in Bhamidi, van der Hofstad \& van 
Leeuwaarden \cite{BhHoLe12} and 
Bhamidi, van der Hofstad \& Sen \cite{BhHoSe15}. Let $W\! :\!  \Omega \! \rightarrow [0, \infty)$ be a r.v.~such that 
\begin{equation}
\label{Wassu}
r= \bE [W] = \bE [ W^2] < \infty \quad \textrm{and} \quad \bP (W \geq x) = x^{-\rho} L(x), 
\end{equation}
where $\rho\ino (2, 3)$ (in the notations of \cite{BhHoSe15}, $\tau \! = \! \rho +1\ino (3, 4)$) and where $L$ is slowly varying at $\infty$.  We then 
set for all $y\ino [0, \infty)$, 
\begin{equation}
\label{GGdef}
G(y) = \sup \big\{ x\ino [0, \infty) : \bP (W\! \geq \! x) \geq 1\! \wedge \! y \big\} . 
\end{equation}
Note that $G(y)\! = \! 0$ for all $y\ino [1, \infty)$ and that $G(y) \! =\!  y^{-1/\rho} \, \ell (y)$, where $\ell$ is slowly varying at $0$. 
We assume:  
\begin{equation}
\label{chitithypo}
\forall n\ino \bbN^*, \quad \bP \big( W\! = \!  G(1/n) \big)= 0 \; .  
\end{equation}
Let $\kappa , q \ino (0, \infty)$ and let $a_n, b_n , \bw_n$ be such that 
\begin{equation}
\label{pwrlwabw}
 a_n \underset{{n\rightarrow \infty}}{\sim} q^{-1} G(1/n), \quad  \forall\,  j\! \geq \! 1, \; \, w_j^{(n)}\! \! =\!  G(j/n), \quad b_n
 \underset{{n\rightarrow \infty}} {\sim} \kappa \sigma_1 (\bw_n) /a_n \; .
\end{equation} 
 Then, the following lemma holds true. 
\begin{lem}
\label{powerlaw} We keep the notations from above and we assume (\ref{chitithypo}).  Then $a_n \sim q^{-1}n^{\frac{1}{\rho}} \ell (1/n)$, $b_n \sim 
q\kappa \, n^{1-\frac{1}{^\rho}}\! / \ell (1/n)$ and $a_n , b_n $ and $\bw_n$ satisfy (\ref{apriori}) with $\beta_0 \! =\! 0$ and $\sqrt{\mathbf{j}_n}/b_n\!  \sim \! \frac{\ell (1/n)}{q\kappa} 
 n^{\frac{1}{^\rho}-\frac{1}{2}} \! \rightarrow \! 0$.  
Next, for all integers $j\! \geq \! 1$ and for all $\alpha \ino\bbR$, set:  
\begin{equation}
\label{glutamate}
 w^{(n)}_j (\alpha)  \! = \! \big(1 \! -\! \frac{_{a_n}}{^{b_n}} (\alpha - \alpha_0)   \big) w^{(n)}_j , \; \textrm{where} \quad 
 \alpha_0 \! =\!  2\kappa q^2    \Big(   \! \int_0^1 \!\!\!   y \{ y^{-\rho}\} \, dy  +\!  \frac{_1}{^{\rho  - 2}} \Big) 
\end{equation} 
and where $\{ \cdot \} $ stands for the fractional part function. Then, $a_n , b_n , \bw_n(\alpha)$ satisfy (\ref{apriori}), $(\mathbf{C1})$--$(\mathbf{C4})$ and $\sqrt{\mathbf{j}_n}/b_n \! \rightarrow \! 0$, with the parameters $\alpha \ino \bbR$, $\kappa \ino (0, \infty)$, $\beta\! = \! \beta_0\! = \! 0$ and $c_j = q \,  j^{-\frac{1}{\rho}}$, for all $j\! \geq \! 1$. 
\end{lem}
%\margmm{here and below: $c_j = q \,  j^{-\frac{1}{\rho}}$}
\noi
\textbf{Proof.} See section \ref{powerlawpf}. \cqfd 

\medskip

Lemma \ref{powerlaw} implies that Theorem \ref{graphcvth} applies to $a_n, b_n$ and $\bw_n (\alpha)$ as defined above. This extends 
%\margmm{new}
Theorem 1.1 in Bhamidi, van der Hofstad \& van 
Leeuwaarden \cite{BhHoLe12} that proves the convergence of the component sizes under the more restrictive assumption that $L(x)\! = \! x^\rho \bP (W\! \geq\!  x) \! \rightarrow \! c_F \ino (0, \infty)$ as $x\! \rightarrow \! \infty$ (see (1.6) in \cite{BhHoLe12}) as well as 
Theorem 1.2 in Bhamidi, van der Hofstad \& Sen \cite{BhHoSe15} (Section 1.1.3) that asserts the convergence of the components as measured metric spaces  
under the supplementary assumptions that 
%$L(x)\! = \! x^\rho \bP (W\! \geq\!  x) \! \rightarrow \! c_F \ino (0, \infty)$ as $x\! \rightarrow \! \infty$ (see (1.3) in  \cite{BhHoSe15}, Section 1.1.1) and that 
$\bP (W\ino dx) \! = \! f(x) dx$, where $f$ is a continuous function 
whose support is of the form $[\epp , \infty)$ with $\epp \! >\! 0$, and such that $x\ino [\epp , \infty) \mapsto xf(x)$ is nonincreasing (see Assumption 1.1 in \cite{BhHoSe15}, Section 1.1.3). Let us mention that both proofs in \cite{BhHoLe12} and in \cite{BhHoSe15} are quite different from ours.

Let us also mention that a solution to the Conjecture 1.3 on fractal dimensions of the components of $\mathbf{G}$ right after Theorem 1.2 in \cite{BhHoSe15} is given in the companion paper \cite{BDW1}, Proposition 2.7.

\medskip
\noindent 
 \textbf{General inhomogeneous Erd\H{o}s--R\'enyi graphs that are close to be multiplicative.} In \cite{Ja10}, Janson investigates strong asymptotic equivalence of general inhomogeneous 
Erd\H{o}s--R\'enyi graphs that are defined as follows: denote by $P$ the set of  arrays 
$\mathbf{p}\! = \! (p_{{i,j}})_{j>i\geq 1}$ of real numbers in $[0, 1]$ such that 
$N_{\mathbf{p}}\! = \! \sup \{ j\! \geq 2\! : \sum_{1\leq i < j} p_{i,j} \! >\! 0 \} \! < \! \infty$; 
the $\mathbf{p}$-inhomogeneous Erd\H{o}s--R\'enyi graph $G(\mathbf{p})$ is the random graph whose set of vertices is $\{ 1, \ldots , N(\mathbf{p})\}$ and whose random set of edges 
 $\ccE (G(\mathbf{p}))$ is such that the r.v.~$(\un_{\{ \{ i,j\} \in \ccE(G(\mathbf{p}))\}})_{1 \leq i < j \leq N(\mathbf{p})}$ are independent and such that 
 $\bP (\{ i,j\} \ino \ccE (G(\mathbf{p})))\! = \! 
 p_{i,j}$. The asymptotic equivalence is measured through the following function $\rho$ that is defined  
for all $p, q\ino [0, 1]$, by $\rho (p, q)\! = \! (\sqrt{p} \! -\! \sqrt{q})^2 \! + \!  (\sqrt{1\! -\! p} \! -\! \sqrt{1\! -\! q})^2$. More precisely, let $\mathbf{p}_n, \mathbf{q}_n\ino P$, $n\ino \bbN$; 
% be such that $\lim_{n\rightarrow \infty}  
%N_{\mathbf{p}_n} \! \wedge\! N_{\mathbf{q}_n}\! =\! 0$;  
then Theorem 2.2 in Janson \cite{Ja10} implies that 
there are couplings of $G(\mathbf{p}_n)$ and $G(\mathbf{q}_n)$ such that $\lim_{n\rightarrow \infty} \bP (G(\mathbf{p}_n)\! \neq \! G(\mathbf{q}_n))\! = \! 0$ if and only if   
\begin{equation}
\label{contigu}
\lim_{n\rightarrow \infty}  \sum_{j>i\geq 1} \rho (p^{_{(n)}}_{^{i,j}} ,  q^{_{(n)}}_{^{i,j}})  = 0 \; .
\end{equation}
We then apply this result as follows: let $a_n , b_n\ino (0, \infty)$ and $\bw_n\ino \elldo_{\! f}$, $n \ino \bbN$, satisfy the assumptions of 
Theorem \ref{graphcvth}; we set 
%\margmm{what is $q_n$ here?}
\begin{equation}
\label{defpniju}
 \forall j\! >\! i\!  \geq\!  1, \quad  p^{_{(n)}}_{^{i,j}} = \frac{ w^{_{(n)}}_{^{i}} w^{_{(n)}}_{^{j}}}{\sigma_1 (\bw_n)} \quad \textrm{and} \quad u^{_{(n)}}_{^{i,j}}= \left\{
\begin{array}{ll} 
\frac{q^{_{(n)}}_{^{i,j}} }{p^{_{(n)}}_{^{i,j}} } \! -\!  1\, , & \textrm{if} \; p^{_{(n)}}_{^{i,j}} \! >\! 0 \\
0\, ,  & \textrm{if} \; p^{_{(n)}}_{^{i,j}} \! =\! 0 .
\end{array}\right.
 \end{equation} 
% \margmm{$p^{_{(n)}}_{^{i,j}} \! =\! 0$}
% 
%\begin{equation}
%\label{defqnj} 
%q_{^j}^{_{(n)}}=\left\{
%\begin{array}{ll}
%c_j & \text{if} \quad  j \! \in \!  \big\{ 1, \ldots ,    n \!  \big\} 
%%1\! \leq \! j\! \leq \! \lfloor  \rho_n 2^n \! \rfloor 
%,  \\
%(\beta/\kappa)^{\frac{1}{3}}n^{-1} & \text{if} \quad  j \! \in \!  \big\{n  +1 , \ldots , n   + n^3  \big\}, \\
%n^{-3} & \text{if} \quad  j \! \in \!  \big\{ n   +n^3 +1 , \ldots , n   + n^3+n^8 \big\}, \\
%% \lfloor \rho_n 2^n\! \rfloor +1 \! \leq \! j\! \leq \! \lfloor \rho_n 2^n\! \rfloor + n2^n , \\
%\; \; \; 0 & \text{if} \quad  j >  n + n^3 + n^8.
%\end{array}\right.
%\end{equation} 
% 
First note that 
$\max _{j>i\geq 1}   p^{_{(n)}}_{^{i,j}} \! = \! O (( w^{_{(n)}}_{^{1}}\! / a_n)^2 a_n / b_n) \! \rightarrow 0$ by (\ref{apriori}); next, as proved in Janson \cite{Ja10} (2.5) p.~30, if $p\! \leq \! 0.9$, then $\rho (p, q) \asymp |p\! -\! q| \big(1 \wedge |q/p \! -\! 1| \big)$. Thus, (\ref{contigu}) is equivalent to 
\begin{equation}
\label{contiguu}
\lim_{n\rightarrow \infty}  \sum_{j>i\geq 1} p^{_{(n)}}_{^{i,j}} |u^{_{(n)}}_{^{i,j}}| \big(1 \wedge  |u^{_{(n)}}_{^{i,j}}|  \big) = 0 ,\quad \textrm{with the convention} \;  p^{_{(n)}}_{^{i,j}} |u^{_{(n)}}_{^{i,j}}|\! = \!  q^{_{(n)}}_{^{i,j}} \; \textrm{if} \;  p^{_{(n)}}_{^{i,j}}= 0. 
\end{equation}
In particular, let $h: [0, 1] \rightarrow [0, 1]$ be such that $h(x) \! = \! x + O (x^2)$. If we set $q^{_{(n)}}_{^{i,j}}\! =\!  
h( p^{_{(n)}}_{^{i,j}})$, then there exists $C\ino (0, \infty)$ such that $|u^{_{(n)}}_{^{i,j}}| \leq C  p^{_{(n)}}_{^{i,j}}$. In this case, for all sufficiently large $n$, 
$$  \sum_{j>i\geq 1} p^{_{(n)}}_{^{i,j}} |u^{_{(n)}}_{^{i,j}}| \big(1 \wedge  |u^{_{(n)}}_{^{i,j}}|  \big) \leq C^2  \sum_{j>i\geq 1}
 (p^{_{(n)}}_{^{i,j}})^3\leq  C^2 \frac{\sigma_3 (\bw_n)^2}{\sigma_1 (\bw_n)^3} \sim C^\prime (a_n/b_n)^3 \longrightarrow 0$$ 
by ($\mathbf{C2}$) and (\ref{apriori}). 
%Thus 
%there are couplings of $G(\mathbf{p}_n)$ and $G(\mathbf{q}_n)$ such that 
%$\lim_{n\rightarrow \infty} \bP (G(\mathbf{p}_n)\! \neq \! G(\mathbf{q}_n))\! = \! 0$. 
%\margmm{In  \cite{EsHoHoo08}, didn't they also take $h(x)=x\wedge 1$?}
Cases where $h(x)\! = \! 1\! \wedge \! x$ have been studied by Chung \& Lu \cite{ChLu02} and van der Esker, van der Hofstad \& Hooghiemstra  \cite{EsHoHoo08};  the cases 
where $h(x)\! = \! 1\! -\! e^{-x}$, was first studied by Aldous \cite{Al97} and Aldous \& Limic \cite{AlLi98} 
and the previous cited papers \cite{AdBrGo12a, BhHoLe10, BhSeWa14, BhBrSeWa14, BhHoLe12, BhHoSe15}, 
including this paper; cases where $h(x)\! = \! x/ (1+ x)$ have been investigated by 
Britton, Deijfen \& Martin-L\"of \cite{BrDeMaLo06}. To summarise, Janson's Theorem 2.2 \cite{Ja10}, p.~31 combined with Theorem \ref{graphcvth} imply the following result. 
\begin{thm}[Theorem 2.2 in Janson \cite{Ja10}]
\label{Jansoncsq} Assume that $a_n , b_n, \bw_n$ satisfy the same assumptions as in Theorem \ref{HYcvth} 
(and thus as in Theorem \ref{graphcvth}). 
We furthermore assume that $\sqrt{\mathbf{j}_n}/b_n \! \rightarrow \! 0$. 
We define $\mathbf{p}_n$ by (\ref{defpniju}). Let $\mathbf{q}_n \ino P$. We define 
$(u^{_{(n)}}_{^{i,j}})_{j>i\geq 1}$ by (\ref{defpniju}) and we suppose (\ref{contiguu}). Then, there exist couplings of 
$G(\mathbf{q}_n)$ and $\cG_{\bw_n}$ such that   
\begin{equation}
\label{stroassequ}
\lim_{n\rightarrow \infty} \bP (\cG_{\bw_n}\! \neq \! G(\mathbf{q}_n))\! = \! 0
\end{equation}
and the weak limit (\ref{globabab}) in Theorem \ref{graphcvth} holds true in the same scaling for the connected components of $G(\mathbf{q}_n)$ that are listed in the decreasing order of their numbers of vertices and that are equipped with the graph distance and the the counting measure. In particular, its holds true when $u^{_{(n)}}_{^{i,j}}\! = \! h(p^{_{(n)}}_{^{i,j}})$, $j\! >\! i\! \geq \! 1$, for all functions $h\! : \! [0, 1] \! 
\rightarrow \! [0, 1]$ such that $h(x)\! = \! x + O (x^2)$. %\margmm{$u^{_{(n)}}_{^{i,j}}$--> $q^{_{(n)}}_{^{i,j}}$}
\end{thm}

\section{Preliminary results on the discrete model.}
\label{discrsec}

\subsection{Height and contour processes of Galton-Watson trees.}
\label{HeightApp}
Let us briefly recall basic notation about the coding of trees. 
We first denote by $\bbU\! = \! \bigcup_{n\in \bbN} (\bbN^*)^n$ the set of finite words written with positive integers; here, $(\bbN^*)^0$ is taken as $\{ \varnothing \}$. The set 
$\bbU$ is totally ordered by the \textit{lexicographical order} $\leq_{\mathtt{lex}}$ (the strict order is denoted by $<_\mathtt{lex}$).  

  Let $u\! = \! [i_1, \ldots , i_n]\ino \bbU$ be distinct from $\varnothing$. We set $|u|\! = \! n$ that is 
the \textit{length} or the \textit{height} of $u$, with the convention that $|\varnothing |\! = \! 0$.  
We next set $\overleftarrow{u}\! = \! [i_1, \ldots, i_{n-1}]$ that is interpreted as the \textit{parent of $u$} (and if $n\! = 1$, then $\overleftarrow{u}$ is taken as $\varnothing$). 
More generally, for all $p\ino \{ 1, \ldots, n\}$, we set 
$u_{| p}\! = \! [i_1, \ldots , i_p]$, with the convention: $u_{| 0}\! = \! \varnothing$. Note that $\overleftarrow{u}\! = \! u_{|n-1}$. 
We will also use the following notation:  
$\lgeo\varnothing  , u\rgeo \! = \! \{ \varnothing, u_{|1}, \ldots , u_{| n-1} , u\}$, 
$\rgeo\varnothing  , u\rgeo \! = \! \lgeo\varnothing  , u\rgeo  \backslash \{ \varnothing\}$, $\lgeo\varnothing  , u\lgeo  \,  = \! \lgeo\varnothing  , u\rgeo  \backslash \{ u\}$ and 
$\rgeo\varnothing  , u\lgeo  \,  = \! \lgeo\varnothing  , u\rgeo  \backslash \{\varnothing, u \}$.
For all $v\! = \! [j_1, \ldots, j_m]\ino \bbU$, we also set 
$u\ast v\! = \! [i_1, \ldots , i_n, j_1, \ldots,  j_m]$ that is 
the \textit{concatenation} of $u$ with $v$, with the convention that 
$\varnothing \ast u\! = \! u \ast  \varnothing\! = \! u$. 
%Then for all $u, v\ino \bbU$, we denote by $u\wedge v$ their most recent common ancestor, namely $u\wedge v\! = \! u_{|l}$

A \textit{rooted ordered tree} can be viewed as a subset $t\! \subset \! \bbU$ such that the following holds true. 
\begin{compactenum}

\smallskip

\item[$(a):$] $\varnothing \ino t$. 

\smallskip

\item[$(b):$] If $u\ino t\backslash \{ \varnothing\}$, then $\overleftarrow{u}\ino t$. 

\smallskip

\item[$(c):$] For all $u\ino t$, there exists $k_u (t)\ino \bbN\cup \{ \infty\}$ such that $u \ast [i]\ino t$ 
if and only if  $1\! \leq \! i \! \leq \! k_u (t)$.  

\smallskip

\end{compactenum}

\noi
Here, $k_u (t)$ is interpreted as the \textit{number of children of $u$} and if 
$1\! \leq \! i \! \leq \! k_u (t)$, then $u\ast [i]$ is the\textit{ $i$-th child of $u$}; $k_u (t)+1$ is the degree of the vertex $u$ in the graph $t$ 
%\margmm{rewritten}
when $u$ is distinct from the root. 
Implicitly, if $k_u (t)\! = \! 0$, then there is no child stemming from $u$ and assertion $(c)$ is trivially satisfied. %empty  
%\ts{We next set $\theta_u t\!  = \! \{ v\ino \bbU: u\ast v\ino t \}$ that is also a rooted ordered tree: it is viewed  as the subtree} \ms{of the descendants} \ts{stemming from $u$}.  
Note that the subtree stemming from $u$ that is $\theta_u t\!  = \! \{ v\ino \bbU: u\ast v\ino t \}$ is also a rooted ordered tree.

Let $\bbT$ be the set of rooted ordered trees that is equipped with the sigma-field $\ccF(\bbT)$ generated by the sets $\{ t\ino \bbT\! : \! u\ino t\}$, $u\ino \bbU$. Then, a \textit{Galton-Watson tree with offspring distribution $\mu$} (a \textit{GW($\mu$)-tree}, for short) is a $(\ccF, \ccF(\bbT))$-measurable r.v.~$\tau\! :\! \Omega \! \rightarrow \! \bbT$ that satisfies the following. 
\begin{compactenum}

\smallskip

\item[$(a^\prime):$] $k_\varnothing (\tau)$ has law $\mu$. 

\smallskip

\item[$(b^\prime):$] For all $k\! \geq\!  1$ such that $\mu(k)\! >\! 0$, the subtrees $\theta_{[1]} \tau, \ldots , \theta_{[k]} \tau$ under $\bP (\, \cdot \, | k_\varnothing (\tau)\! = \! k)$ are independent with the same law as $\tau$ under $\bP$. 
\end{compactenum}

\smallskip

%\margmm{new}
\noi
Assume that $\mu(1)\! <\! 1$. Recall that $\tau$ is a.s.~finite if and only if  $\mu$ is critical or subcritical: namely, if and only if  $\sum_{k\geq 1} k\mu (k)\! \leq \! 1$. 

A \textit{Galton-Watson forest with offspring distribution $\mu$} (a \textit{GW($\mu$)-forest}, for short) is a random tree $\bT$ such that $k_\varnothing (\bT)\! = \! \infty$ and such that the subtrees $(\theta_{[k]} 
\bT )_{k\geq 1}$ stemming from $\varnothing$ are i.i.d.~GW($\mu$)-trees. 
We next recall how to encode a GW($\mu$)-forest $\bT$ thanks to three processes: its \textit{Lukasiewicz path}, its \textit{height process} and its \textit{contour process}. We denote by $(u_l)_{l\in \bbN}$ the sequence of vertices of $\bT$ such that $u_0\! = \! \varnothing$ and such that for all $l$, $u_{l+1}$ is the smallest vertex of $\bT$ with respect to the lexicographical order that is larger than $u_l$. If $\mu$ is critical or subcritical, then $(u_l)_{l\in \bbN}$ exhausts all the vertices of $\bT$; however, if $\mu$ is supercritical (namely if $\sum_{k\geq 1} k\mu (k)\! >\! 1$), then $(u_l)_{l\in \bbN}$ exhausts the vertices of $\bT$ that are situated before (or on) the first infinite line of descent. %\ts{in the lexicographical order}. 
We first set: 
\begin{equation}
\label{coding}
V^{\bT}_0\! = \! 0,  \quad \forall l\! \geq \! 0, \quad V_{l+1}^\bT \! = \! V_{l}^\bT  + k_{u_{l+1}} (\bT)\! -\! 1 \quad \textrm{and} \quad \mathtt{Hght}^\bT_l\! =\!  |u_{l+1}|\! -\! 1. 
\end{equation}
The process $(V_l^\bT)_{l\in \bbN}$ is the \textit{Lukasiewicz path associated with $\bT$} and $(\mathtt{Hght}^\bT_ l)_{l\in \bbN}$ is the \textit{height process associated with $\bT$}. We recall from 
Le Gall \& Le Jan \cite{LGLJ98} the following results. 

\begin{compactenum}
\item[$(i)$]  $V^\bT$ is distributed as a random walk starting from $0$ and 
with jump-law $\nu (k)\! = \! \mu (k+1)$, $k\ino \bbN \cup \{ -1\}$. 

\smallskip

\item[$(ii)$] We set $\underline{V}^\bT_{\, 0}\! = \! 0$ and for all $l\! \geq \! 1$, $\underline{V}^\bT_{\, l}\! = \! \inf_{0\leq k \leq l-1} V^\bT_k -1$. Note that $u_l\ino \theta_{[p]} \bT$ if and only if  $(u_l)_{| 1}\! =\!  p$. Then, we get:  
\begin{equation}
\label{faure}
- \underline{V}^\bT_{\, l} \! = \! (u_l)_{| 1} \quad \textrm{and} \quad V^\bT_l \! -\! \underline{V}^\bT_{\, l}\! = \! \# \big\{ v \ino \bT: u_{l}<_\mathtt{lex} v  \; \textrm{and} \;  \overleftarrow{v}\! \in\,  \rgeo \varnothing , u_l \rgeo  \big\} 
\end{equation}
\item[$(iii)$]The height process $\mathtt{Hght}^\bT$ is derived from 
$V^\bT$ by setting $\mathtt{Hght}^\bT_0\! = \! 0$ and 
\begin{equation}
\label{codheight}
\forall l \! \geq \! 1, \quad \mathtt{Hght}^\bT_l= \# \big\{ m\ino \{ 0, \ldots , l\! -\! 1\}: V^\bT_m \! = \! \inf_{m\leq j\leq l} V^\bT_j  \big\}. 
\end{equation}
\end{compactenum}

  The \textit{contour process of $\bT$}
is informally defined as follows: suppose that $\bT$ is embedded in the oriented half plane in such a way that edges have length one and that orientation reflects lexicographical order of visit; 
we think of a particle starting at time $0$ from $\varnothing$ and exploring the tree from the left to the right, backtracking as less as possible and 
moving continuously along the edges at unit speed. 
In (sub)critical cases, the particle crosses each edge twice (upwards first and then downwards). 
In supercritical cases, the particle only explores the edges that 
are situated before (or on) the first infinite line of descent in the lexicographical order: the edge on the infinite line of descent are visited once (upwards only) and the edge strictly before the infinite line of descent are visited twice (upwards first and then downwards). 
For all $s\ino [0, \infty)$, we define $C^\bT_s $ 
as the distance at time $s$ of the particle from the root $\varnothing$. %\ms{The contour process is close to the height process.} 
The associated distance $d_{C^\bT}$ 
as defined in (\ref{pseudometric}) is the graph distance of $\bT$ in the (sub)critical cases. We refer to Le Gall \& D.~\cite{DuLG02} (Section 2.4, Chapter 2, pp.~61-62) for a formal definition and a formula relating the contour process to the height process.

\subsection{Coding processes related to the Markovian queueing system.}
\label{connecMqu}
%\tt{SLIGHTLY RE-ORGANISED SECTION FROM HERE.} 
We fix $\bw\! = \! (w_1, \ldots , w_n, 0, 0 , \ldots ) \ino \elldo_f$ and we briefly recall the definition of the Markovian queue as in Introduction: 
a single server is visited by infinitely many clients; clients arrive according to a Poisson process with unit rate; 
each client has a \textit{type} that is a positive integer ranging in $\{ 1, \ldots, n\}$; the amount of time of service required by a client of type $j$ is $w_j$; types are i.i.d.~with law 
\begin{equation}
\label{defnuw}
\nu_\bw \! = \! \frac{1}{\sigma_1 (\bw)}\sum_{j=1}^n w_j \delta_{j} \; .
\end{equation}
Let $\tau_l$ stand for the time of arrival of the $l$-th client in the queue and let $\Jtt_l$ stand for her/his type; then, the queueing system is entirely characterised by the point measure 
%$\ccX_\bw  \! = \sum_{l\geq 1} \delta_{(\tau_l , \Jtt_l)}$
\begin{equation}
\label{rXPoisdef}
\ccX_\bw  \! = \sum_{k\geq 1} \delta_{(\tau_k , \Jtt_k)},
\end{equation} 
that is distributed as a Poisson point measure on $[0, \infty) \! \times \! \{ 1, \ldots, n\}$ whith intensity 
$\ell \! \otimes \! \nu_\bw$, where $\ell $ stands for the Lebesgue measure on $[0, \infty)$. 
We next introduce the following:  
\begin{equation}
\label{rXwdef}
 \forall t\ino [0, \infty), \quad X^\bw_t  =  -t + \sum_{l\geq 1} w_{\Jtt_l}\un_{[0, t]} (\tau_l) \;  \quad \textrm{and} \quad I^\bw_t \! = \! \inf_{s\in [0, t]} X^\bw_s . 
\end{equation} 
%Note that when the server is idle, the load is counted negatively, hence the name \textit{algebraic}. The actual load of the server at time $t$ is then 
%$X^\bw_t \! -\! \inf_{s\in [0, t]} X^\bw_s $. 
Then, $X^\bw_t \! -\! I^\bw_t$ is interpreted as the load of the Markovian queueing system at time $t$ and $X^{\bw}_t$ is the algebraic load of the queue.  
Note that $X^\bw$ is a spectrally positive L\'evy process whose law is determined by 
its Laplace exponent 
$\psi_\bw\! : \! [0, \infty) \! \rightarrow \! \bbR$ in the following way: 
\begin{multline}
\label{psibwdef}
 \bE \big[ e^{-\lambda X^\bw_t}\big]\! = \! e^{t\psi_\bw (\lambda)} \quad \textrm{where} \\
 \psi_\bw (\lambda)  =  \alpha_\bw \lambda + \!\!\!\!   \sum_{1\leq j\leq n}\!  \! \frac{_{w_j} }{^{\sigma_1 (\bw)}} \big( e^{-\lambda w_j}\! -\! 1\! + \! \lambda w_j \big) \quad \textrm{and} \quad \alpha_\bw \! := \! 1\! -\! \frac{_{\sigma_2 (\bw)}}{^{\sigma_1 (\bw)}} \; .
\end{multline}
Here, recall that $\sigma_r (\bw)\! = \! w_1^r + \ldots + w_n^r$, $r\ge 0$.  
We call the queueing system {\it recurrent} if a.s.~$\liminf_{t\rightarrow \infty}$ $X^\bw_t \!= \! -\infty$, which means that all the clients will eventually depart. Let us observe that the system is recurrent if and only if    
$\sigma_2 (\bw) / \sigma_1 (\bw)\! \leq \!  1$. If, on the other hand, $\sigma_2 (\bw) / \sigma_1 (\bw)\! >\!  1$, then $\alpha_\bw \! <\! 0$ and a.s.~$\lim_{t\rightarrow \infty} X^\bw_t \!= \! \infty$ (the queue will see an accumulation of infinitely many clients). 
In the sequel, we shall refer to the following cases: 
\begin{equation}
\label{soususcri}
\textrm{supercritical:} \; \sigma_2 (\bw)\! >\! \sigma_1 (\bw), \quad   \textrm{critical:} \;  \sigma_2 (\bw)\! =\! \sigma_1 (\bw), \quad   \textrm{subcritical:} \; \sigma_2 (\bw)\! < \! \sigma_1 (\bw). 
\end{equation}

The LIFO queueing system governed by $\ccX_\bw$ generates a tree 
that can be informally defined as follows: 
\textit{the clients are the vertices and the server is the root (or the ancestor); the $j$-th client to enter the queue is a child of the $i$-th one if the $j$-th client enters when the $i$-th client is served; among siblings, clients are ordered according to their time of arrival.} In critical or subcritical cases, it fully defines a Galton-Watson forest; however in supercritical cases, it only defines the part of a Galton-Watson forest situated before the first infinite line of descent. To circumvent this problem, we actually define the tree first and then we couple it with the queueing system as follows.  

In what follows, what we mean by a Poisson random subset $\Pi$ on $[0, \infty)$ with unit rate is the set of atoms of a unit rate Poisson random measure: namely, it is the random subset $\{ \mathtt{e}_1+ \ldots + \mathtt{e}_n ; n\! \geq \! 1\}$, where the $\mathtt{e}_n$ are i.i.d.~exponentially distributed r.v.~with unit mean.  
For all $u\ino \bbU\backslash \{ \varnothing\}$, let $J(u)$ and $\Pi_u$ be independent r.v.~whose laws are given as follows: $J(u)$ has law $\nu_\bw$ as defined in (\ref{defnuw}) and $\Pi_u$ is a Poisson random subset of $[0, \infty)$ with unit rate. We next define $\Pi_\varnothing$ as a Poisson random subset of $[0, \infty)$ with unit rate that is assumed to be independent of $(J(u), \Pi_u)_{u\in \bbU\backslash \{ \varnothing\}}$ and by convenience, we set 
 %\margmm{rewritten} 
  $J(\varnothing) \! = \! 0$. For all $u\ino \bbU$, we index the points of $\Pi_u$ using the children of $u$. Formally, we define a map 
  $\sigma: \{ u \! \ast \! [p]\, ;\; p\! \geq \! 1 \} \to \Pi_{u}$ as follows:  
 \begin{equation}
 \label{notaMu}
 \Pi_u \! = \! \big\{ \sigma (u \! \ast \! [p]) \, ; \; p\! \geq \! 1\big\}, \,\;  \textrm{where}\; \, 
 \sigma (u \! \ast \! [p])  <  \sigma (u \! \ast \! [p+1]), \; \,  p\! \geq \! 1.
 \end{equation} 
Note that it defines a collection $(\sigma (u))_{u\in \bbU \backslash \{ \varnothing \}}$ of r.v. 
It is easy to check that here is a unique random tree $\bT_\bw\! : \! \Omega \! \rightarrow \! \bbT$ such that
%\margmm{see footnote}
% \fmm{there is some problem with the formal definition of $\bT_{\bw}$: since $u$ is not all in $bT_{\bw}$, $k_{u}(\bT_{\bw})$ is not always well defined. \tt{I think it is better now with the correction. -tt}}
  \begin{equation}
\label{treeTrts}
\forall u \ino \bT_\bw \backslash \{ \varnothing \}, \quad k_u (\bT_\bw)\! = \! \# \big( \Pi_u \cap [ 0, w_{J(u)} ] \big)  \quad \textrm{and} \quad k_\varnothing (\bT_\bw) \! = \! \infty . 
\end{equation}
Clearly $\bT_\bw$ is distributed as a GW($\mu_\bw$)-forest where $\mu_\bw$ is given by 
\begin{equation}
\label{rmupoissw}
\forall k \ino \bbN, \quad \mu_\bw (k) \! =\!  \sum_{j\geq 1} \frac{{w_j^{k+1} e^{-w_j} }}{{ \sigma_1 (\bw) \,  k!}}  \; .
\end{equation} 
Namely, $k_\varnothing (\bT_\bw)\! = \! \infty$ and the subtrees $(\theta_{[k]} 
\bT_\bw )_{k\geq 1}$ stemming from $\varnothing$ are i.i.d.~GW($\mu_\bw$)-trees. Note that 
$\sum_{k\geq 0} k\mu_\bw (k)\! = \! \sigma_2 (\bw) / \sigma_1 (\bw)$. %\tt{TO HERE.}

  Then we define the point process $\ccX_\bw$ governing the Markovian queueing system as follows: denote by $(u_l)_{l\in \bbN}$ the sequence of vertices of $\bT_\bw$ such that $u_0\! = \! \varnothing$ and such that for all $l$, $u_{l+1}$ is the smallest vertex of $\bT_\bw$ (with respect to the lexicographical order) that is larger than $u_l$. Then we set 
 \begin{equation}
 \label{treecoupl}
J_l = J(u_l) \, ,  \quad \tau_l \! = \!\! \! \!\! \!\!\! \!
\sum_{\quad \substack{ v\in \bT_\bw : v <_{\mathtt{lex}} u_l \\  \textrm{and} \; v\notin \lgeo \varnothing , u_l \rgeo } } \!\! \!\!\! \!\!\! \!\!\! \! w_{J(v)}  +\sum_{v\in \rgeo \varnothing , u_l \rgeo} \sigma (v) \quad \textrm{and} \quad \ccX_\bw= \sum_{l\geq 1} \delta_{(\tau_l, J_l)}. 
\end{equation}
We also set 
\begin{equation}
\label{NNwwdef}
\forall t\ino [0, \infty), \quad N^\bw (t)= \sum_{l\geq 1} \un_{[0, t]} ( \tau_l) \; . 
\end{equation}
Recall from (\ref{coding}) that $(V^{\bT_\bw}_l)_{l\geq 0}$ stands for the Lukasiewicz path associated with $\bT_\bw$; we also recall the notation $\underline{V}^{\bT_\bw}_{\, l} $ for 
the quantity $\inf_{0\leq k\leq l-1} V^{\bT_\bw}_k -1$. 
\begin{lem}
\label{LukaXw} We keep the notation from above. Then $\ccX_\bw$ as defined by (\ref{treecoupl}) is a 
Poisson point measure on $[0, \infty) \! \times \! \{ 1, \ldots, n\}$ with intensity 
$\ell \! \otimes \! \nu_\bw$ and therefore $N^\bw$ as defined by (\ref{NNwwdef}) is a Poisson process on $[0, \infty)$ with unit rate. Let $X^\bw$ and $I^\bw$ be derived from $\ccX_\bw$ by (\ref{rXwdef}). 
For all $t\ino [0, \infty)$, then the following holds true:  

\smallskip

\begin{compactenum}

\item[$(i)$] Conditionally given $X^\bw_t \! -\! I^\bw_t$, $V^{\bT_\bw}_{N^\bw(t)} \!- \!  
\underline{V}^{\bT_\bw}_{N^\bw(t)}$ is distributed as a Poisson r.v.~with mean $X^\bw_t \! -\! I^\bw_t $. 

\smallskip

\item[$(ii)$] $\bP$-almost surely: $-\underline{V}^{\bT_\bw}_{N^\bw(t)} \! =\! \# (\Pi_\varnothing \cap [0, -I^\bw_t] ) $. 
\end{compactenum}

\smallskip

\noi
Then, for all $a , x\ino (0, \infty)$, we get 
\begin{equation}
\label{LuapproXw}
\bP \big( \big| V^{\bT_{\! \bw}}_{N^\bw (t)} \! -\! X^\bw_t \big| > 2a\big) \leq 1 \! \wedge \! (4x/a^2) + \bP \big( \! -\! I^\bw_t \! >\!  x) + \bE \big[ 1 \wedge\!  \big( (X^\bw_t \! -\! I^\bw_t)/a^2 \big)\big] . 
\end{equation}
\end{lem}
\noi
\textbf{Proof.}  
%\fmm{This proof is a bit hard to follow. If the reader has the discrete exploration process in mind, it is ok to see where the proof is going; otherwise he might get lost. \tt{I agree and I think it is parttly written somewhere ... -tt}}
We first explain how $(\tau_{l+1}, 
X^\bw_{\tau_{l+1}})$ is derived from $(\tau_l, X^\bw_{\tau_l})$ in terms of the r.v.~$(J(u), \Pi_u)$, $u\in \bbU$. To that end, we need notation: fix $u\ino \bbU\backslash \{ \varnothing\}$; then for all $0\! \leq \! p \! < \! |u|$, we set: 
$$R^u_p\! = \! \big( w_{J(u_{|p})} \! -\! \sigma (u_{| p+1}) \big)_+ \quad \textrm{and} \quad Q^u_p\! = \! \big\{ \sigma (v)\! -\! \sigma (u_{| p+1})\, ; \; v\ino \bbU: u_{|p+1} \!\!  <_{\mathtt{lex}} \!  v \; \textrm{and} \; \overleftarrow{v}\! = \! u_{| p} \big\}  .$$
Note that $J(u_{| 0})\! = \! J (\varnothing)\! = \! 0$ and that $w_0 \! = \! \infty$ (by convention); thus, $R_0^u\! = \! \infty$. We also set $R^u_{|u|}\! = \! w_{J(u)} $, $Q^u_{|u|}\! = \! \Pi_u$, 
$R(u)\! = \! (R^u_1, \ldots, R^u_{|u|})$ and $Q(u)\! = \! (Q^u_0, \ldots, Q^u_{|u|})$. By convention, we finally set $R^\varnothing_0\! = \! \infty$, $R(\varnothing)\! = \! \varnothing$ and $Q(\varnothing)\! = \! (\Pi_\varnothing)$.  

We next denote by $\ccG (u)$ the sigma-field generated by the r.v.~$(\sigma(v), J(v), \Pi_v \cap [0, \sigma (v)))_{ v  \in \, \rgeo \varnothing , u \rgeo }$ and $(J (v), \Pi_v)_{ v  <_\mathtt{lex}  u \; \textrm{and} \; v \notin \lgeo \varnothing , u \rgeo }$. Elementary properties of Poisson point processes imply that conditionnally given $\ccG(u)$, the 
%\margmm{changed}
$Q^u_p$, $0\! \leq \! p\! < \! |u|$ are independent Poisson random subsets of $[0, \infty)$ with unit rate: they are therefore independent of $\ccG(u)$; by construction they are 
also independent from the r.v.~$(J (v), \Pi_v)$, $u \! <_\mathtt{lex} v$.     

 For all $u\ino \bbU\backslash \{ \varnothing\}$, we next define $s(u)\ino \bbU$ and $e(u)\ino [0, \infty)$  that satisfy $s(u_{l})\! = \! u_{l+1}$ and $\tau_{l+1}\! = \! e(u_l)+ \tau_l$. To that end, we first set $\mathbf{q}\! = \! \sup \big\{ p \ino \{ 0, \ldots , |u| \} : \# (Q^u_p \cap [0, R^u_p])\! \geq \! 1 \big\}$ that is well-defined since $R^u_0\! = \! \infty$. 

\smallskip

\begin{compactenum}

\item[$\bullet$] If $\mathbf{q}\! = \! |u|$, then we set $s(u)\! = \! u\! * \! [1]$ and $e(u)\! = \! \sigma (u \! \ast \! [1])$. 

\smallskip

\item[$\bullet$]  If  $\mathbf{q}\! <\! |u|$, then $|u|\! \geq \! 1$ and we set $s(u)\! = \! [i_1, \ldots , i_{\mathbf{q}} ,  i_{\mathbf{q}+1} \! +\! 1]$ (namely, $s(u)\! = \![i_1\! +\! 1]$ if $\mathbf{q}\! = \! 0$), 
where  $u \! = \! [i_1, \ldots , i_{|u|} ]$. We also set: 
\begin{equation}
\label{s-chanff}
e(u)\! = \! \sigma (s(u))\! -\! \sigma (u_{|\mathbf{q} +1}) + \sum_{\mathbf{q}< p \leq |u|} R^u_p \; .
\end{equation}
\end{compactenum}
%\margmm{see footnote}
Elementary properties on Poisson point processes imply that $e(u)$, $J(s(u))$ and $\ccG (u)$ are independent, that $e(u)$ is exponentially distributed with unit mean and that $J(u)$ has law $\nu_\bw$. 
%\fmm{Are these really obvious to see?}
Then, we easily derive from (\ref{treeTrts}) that for all $l\ino \bbN$, 
$u_{l+1}\! = \! s(u_l)$, as already mentioned. It is also easy to deduce from (\ref{treecoupl}) that $\tau_{l+1}\! = \! e(u_l)+ \tau_l$. Thus, $\tau_{l+1}\! -\! \tau_l$,  $J(u_{l+1})$ are independent and they are also independent of $\ccG(u_l)$ and therefore of the r.v.~
%\margmm{changed}
$((\tau_k, J(u_k)))_{1\leq k\leq l}$. 
%moreover, $\tau_{l+1}\! -\! \tau_l$  is exponentially distributed with unit mean and $J(u_l)$ has law $\nu_\bw$. 
It implies that $\ccX_\bw$ is a Poisson point measure on $[0, \infty) \! \times \! \{ 1, \ldots, n\}$ with intensity $\ell \! \otimes \! \nu_\bw$.   

  We next prove inductively that for all $l\! \geq \! 1$, 
\begin{equation}
\label{Xinfblop}
Z_l\! := \!\!\! \sum_{1\leq p\leq |u_l|} \!\!\!  R^{u_l}_p =  X^\bw_{\tau_l} -I^\bw_{\tau_l} \quad \textrm{and} \quad \sigma ((u_l)_{| 1}) \! = \!  -I^\bw_{\tau_l}. 
\end{equation}
\textit{Proof.} Clearly, (\ref{Xinfblop}) holds for $l\! = \! 1$. Assume it holds true for $l$. 
Set $k\! = \! (u_l)_{| 1}$; namely $u_l\ino \theta_{[k]} \bT_\bw$. Since $u_{l+1}\! = \! s(u_l)$, $u_{l+1} \ino \theta_{[k]} \bT_\bw$ if and only if   $\mathbf{q}\! = \! \sup \big\{  p \ino \{ 0, \ldots,   |u_l|\} : \# (Q^{u_l}_p \cap [0, R^{u_l}_p])\! \geq \! 1 \big\} \! \geq \!  1$. We first suppose that $\mathbf{q}\! \geq \! 1$. 
By comparing (\ref{s-chanff}) and (\ref{Xinfblop}), we see that $e(u_l) \! <\! Z_l$.  
Since $\tau_{l+1} \! - \tau_l\! = \! e(u_l)$ and since $X^\bw$ does not jump on $[\tau_l, \tau_{l+1})$ (by the definition (\ref{rXwdef})), we get:    
\begin{equation}
\label{saens}
\inf_{s\in [\tau_l, \tau_{l+1}] } X^\bw_s \! = \! X^{\bw}_{\tau_{l+1}-}\! = \! X^\bw_{\tau_l} \! -\! (\tau_{l+1} \! -\! \tau_l) \! = \! X^\bw_{\tau_l} -e(u_l)\! = \! Z_l -e(u_l) + I^\bw_{\tau_l} 
\end{equation}
and thus $-I^\bw_{\tau_{l+1}}\! = \! -I^\bw_{\tau_l}$. Since 
$u_{l+1} \ino \theta_{[k]} \bT_\bw$, $k\! = \! (u_{l+1})_{| 1}\! = \!   (u_{l})_{| 1}$ and thus $\sigma ( (u_{l+1})_{| 1})\! = \! \sigma ((u_{l})_{| 1})$. Then (\ref{Xinfblop}) entails $-I^\bw_{\tau_{l+1}}\! = \! \sigma ( (u_{l+1})_{| 1})$. We also check easily that $Z_{l+1} \! -\! Z_l\! = \! 
w_{J(u_{l+1}) } \! -\! e(u_l)\! = \! X^{\bw}_{\tau_{l+1}} \! -\! X^\bw_{\tau_l}$, which easily entails (\ref{Xinfblop}) for $l+1$.  

  Suppose now that $\mathbf{q} \! = \! 0$, which is equivalent to $u_{l+1}\! = \! [k+1]$. Thus, $R(u_{l+1})\! =\!  ( w_{J(u_{l+1})})$ and $Z_{l+1} \! = \! w_{J(u_{l+1})}$.  
Since $\mathbf{q} \! = \! 0$, $e(u_l)\! = \! Z_l+ \sigma ([k+1])\! -\! \sigma ([k])$ by (\ref{s-chanff}). 
As in (\ref{saens}), we get $ \inf_{s\in [\tau_l, \tau_{l+1}] } X^\bw_s \! = \! Z_l -e(u_l) + I^\bw_{\tau_l} = I^\bw_{\tau_l} \! -\! \sigma ([k+1])+ \sigma ([k]) \! = \!  - \sigma ([k+1])$, the last equality being a consequence of (\ref{Xinfblop}) \mm{for $l$}. It implies that $-I^\bw_{\tau_{l+1}}\! = \! \sigma ([k+1])\! >\! \sigma ([k])\! = \! -I^\bw_{\tau_{l}}$. Therefore, $X^\bw_{\tau_{l+1}} \! -\! I^\bw_{\tau_{l+1}}\! = \! \Delta X^\bw_{\tau_{l+1}}\! = \! w_{J(u_{l+1})}\! = \! Z_{l+1}$. This proves that (\ref{Xinfblop}) holds for $l+1$. It also completes the proof of (\ref{Xinfblop}) by induction. \cq

\smallskip

Next, it is easy to check that 
$$ \# \big\{ v \ino \bT_\bw: u_{l}<_\mathtt{lex} v  \; \textrm{and} \;  \overleftarrow{v}\! \in\,  \rgeo \varnothing , u_l \rgeo  \big\} = \sum_{1 \leq p \leq |u_l |} \# \big( Q^{u_l}_p\cap [0, R^{u_l}_p ]   \big) , $$
and by (\ref{faure}) we get $\sum_{1 \leq p \leq |u_l|} \# \big( Q^{u_l}_p\cap [0, R^{u_l}_p ]   \big) \! = \!  V^{\bT_\bw}_l \! -\! \underline{V}^{\bT_\bw}_{\, l}$. By (\ref{Xinfblop}) and  elementary 
properties of Poisson point processes, it shows that conditionally given $X^\bw_{\tau_l} -I^\bw_{\tau_l}$, $V^{\bT_\bw}_l \! -\! \underline{V}^{\bT_\bw}_{\, l}$ is distributed as a Poisson r.v.~with mean 
$X^\bw_{\tau_l} -I^\bw_{\tau_l}$ and elementary arguments entail $(i)$. 

Next recall from (\ref{faure}) that $-\underline{V}^{\bT_\bw}_{\, l}\! = \! (u_l)_{| 1}$ and recall from  (\ref{Xinfblop}) that $\sigma ((u_l)_{| 1}) \! = \!  -I^\bw_{\tau_l}$. Namely, $-\underline{V}^{\bT_\bw}_{\, l}\! = \! \# \big(\Pi_\varnothing \cap [0,   -I^\bw_{\tau_l} ] \big)$ and elementary arguments entail 
$-\underline{V}^{\bT_\bw}_{\, l}\! = \!  \# \big(\Pi_\varnothing \cap [0,   -I^\bw_{t} ] \big)$ for all $t \ino [\tau_l, \tau_{l+1})$. This easily proves $(ii)$ for all $t\ino [\tau_1, \infty)$. For all $t\ino [0, \tau_1)$, observe that $N^\bw_{t}\! = \! 0$ and $-I^\bw_t\! = t\! < \! \tau_1$. Since $\underline{V}^{\bT_\bw}_{\, 0}\! = \! 0$, it entails $(ii)$ for all $t\ino [0, \tau_1)$, which completes the proof of $(ii)$.   

We next prove (\ref{LuapproXw}). We fix $t\ino [0, \infty)$ and to simplify we set 
$$D\! = \! V^{\bT_\bw}_{N^\bw(t)} \!- \!  
\underline{V}^{\bT_\bw}_{N^\bw(t)}, \; Z\! = \! X^\bw_t \! -\! I^{\bw}_t\, , \; D^\prime\! = \! -\underline{V}^{\bT_\bw}_{N^\bw(t)} \quad \textrm{and} \quad Z^\prime\! = \! -I^\bw_t \; .$$
By $(i)$, $\bE \big[ (D\! -\! Z)^2| Z \big]\! = \! Z$; thus 
$\bP (|D\! -\! Z| \! >\! a) \! \leq \! \bE [1 \wedge (Z/a^2)]$. By $(ii)$, $D^\prime\! = \! \# (\Pi_\varnothing \cap [0, Z^\prime])$; then, for all $x\ino (0, \infty)$, we get 
$ \bP (|D^\prime\! -\! Z^\prime| \! >\! a) \! \leq \! 
\bP (\sup_{z\in [0, x]}|\# (\Pi_\varnothing \cap [0, z])\! -\! z| \! > \! a )+ \bP (Z^\prime \! > \! x) \! \leq\!  1\! \wedge \! (4x/a^2) + \bP (Z^\prime \! > \! x)$
by Doob $L^2$-inequality for the martingale $z\! \mapsto \! \# (\Pi_\varnothing \cap [0, z])\! -\! z$. It implies (\ref{LuapproXw}), which completes the proof of Lemma \ref{LukaXw}. \cqfd

\bigskip   
   
\noi
\textbf{The contour of $\bT_{\! \bw}$: estimates.}
\label{contourTw}   
Recall from (\ref{HXdisdef}) that $H^\bw_t$ stands for the number of clients waiting in the line right after time $t$. More precisely, for all $s, t \ino [0, \infty)$ such that $s \! \leq \! t$, we get 
\begin{equation}
\label{XJHdef}
H^\bw _t \! = \# \cK_t, \;   \textrm{where} \;    \cK_t \! =\!  \big\{ s\ino [0, t]\!  : \! I^{\bw, s-}_{t} \! <\! I^{\bw, s}_{t} \big\}   \;    \textrm{and where} \;     \forall s\ino [0, t], \;  I^{\bw, s}_t= \inf_{r\in [s, t]}X^\bw_r . 
\end{equation}
% $I^{\bw, s}_t\! = \! \inf_{r\in [s, t]}X^{\bw}_r  $ and $H^{\bw}_t \! = \# \big\{ s\ino [0, t]\, : \;I^{\bw, s-}_{t} \! <\! I^{\bw, s}_{t} \big\}$. 
%\begin{equation}
%\label{rXJHdef}
%H^{\bw}_t \! = \# \big\{ s\ino [0, t]\, : \;I^{\bw, s-}_{t} \! <\! I^{\bw, s}_{t} \big\}. 
%\end{equation}
The process $H^\bw$ is called the \textit{height process associated with $X^\bw$} by analogy with (\ref{codheight}), but $H^\bw$ is actually closer to the contour process of 
$\bT_{\! \bw}$. 

To see this, recall that $(u_l)_{l\in \bbN }$ stands for the sequence of vertices of $\bT_{\! \bw}$ listed in the lexicographical order; we identify $u_l$ with the $l$-th client to enter the queueing system. For all $t\ino [0, \infty)$, we denote by $\mathbf{u} (t)$ the client currently served right after time $t$: namely $\mathbf{u} (t) \! = \! u_l$ where 
$l\! = \! \sup \{ k\ino \bbN\! :  \tau_k \! \leq \! t  \;  \textrm{and} \;  X^\bw_{\tau_k-} \! < \! \inf_{s\in [\tau_k, t]}X^\bw_s  \}$. Then, the length of the word $\mathbf{u} (t)$ is the number of clients waiting in the line right after time $t$: $| \mathbf{u} (t)| \! = \! H^{\bw}_t $.   

  We next denote by $(\xi_m)_{m\geq 1}$ the sequence of jump-times of $H^{\bw}$: namely, $\xi_{m+1} \! = \! \inf \{ s\! >\! \xi_m \! : H^\bw_s \! \neq \! H^{\bw}_{\xi_m} \}$, for all $m\ino \bbN$, 
with the convention $\xi_0 \! = \! 0$. 
We then set: 
 \begin{equation}
\label{Mbwddef}
\forall t\ino [0, \infty), \quad M^\bw_t  =\sum_{m\geq 1} \un_{[0, t]} (\xi_m)  \; .
\end{equation}
Note that $(\xi_m)_{m\geq 1}$ is also the sequence of jump-times of $\mathbf{u}$ and that for all $m\! \geq \! 1$, $(\mathbf{u} (\xi_{m-1}), \mathbf{u} (\xi_m))$ is necessarily an oriented edge of $\bT_{\! \bw}$. We then set $\bT_{\! \bw} (t)\! = \! \{\mathbf{u} (s); s\ino [0, t]\}$, that represents the set of the clients who entered the queue before time $t$ (and the server $\varnothing$); $\bT_{\! \bw} (t)$ has $N^\bw (t)+1$ vertices (including the server represented by $\varnothing$); therefore, $\bT_{\! \bw} (t)$ has 
$2N^\bw (t)$ oriented edges.  
Among the $2N^\bw (t)$ oriented edges of $\bT_{\! \bw} (t)$, the $|\mathbf{u} (t)|$ edges going down from $\mathbf{u} (t)$ to $\varnothing$ does not belong to the subset 
$\{  (\mathbf{u} (\xi_{m-1}), \mathbf{u} (\xi_m)); m\! \geq \! 1 \! : \xi_m \! \leq \! t \}$. Thus, we get   
 \begin{equation}
\label{crossedg}
\forall t\ino [0, \infty), \quad M^\bw_t = 2N^\bw (t) \! -\! H^\bw_t \; .
\end{equation}
Recall from Section \ref{HeightApp} the definition of the contour and the height processes of $\bT_{\! \bw}$, denoted resp.~by $(C^{\bT_{\! \bw}}_t)$ and $(\mathtt{Hght}^{\bT_{\! \bw}}_k)$. Then, observe that 
\begin{equation}
\label{couCHght}
\forall t\ino [0, \infty), \quad C^{\bT_{\! \bw}}_{M^\bw (t)} \! = \! H^\bw_t \quad \textrm{and} \quad \sup_{s\in [0, t]} H^\bw_s  \leq 1+ \sup_{s\in [0, t]} \mathtt{Hght}^{\bT_{\! \bw}}_{N^\bw_s} . 
\end{equation}
%\label{page_ht}\mm{Comment: I feel there is a small problem. In \eqref{codheight}, $\mathtt{Hght}^{\bT_{\! \bw}}$ is actually the height process of the forest $\bT_{\! \bw}\setminus\{\varnothing\}$, so that the first client corresponds to $\mathtt{Hght}^{\bT_{\! \bw}}_0$ and it has height $0$. But here, it has  height $1$.}
Since $N^\bw$ is a homogeneous Poisson process with unit rate, Doob's $L^2$-inequality combined with (\ref{crossedg}) and (\ref{couCHght}) imply the following inequality: 
\begin{equation}
\label{estMbw}
\forall t, a\ino (0, \infty), \quad \bP \big(\! \sup_{s\in [0, t]}\!   |M^\bw_s \! -\! 2s |> 2a\big) \leq 1\!  \wedge\!  (16t/a^2) + \bP \big(1+ \!\! \sup_{s\in [0, t]} \! \mathtt{Hght}^{\bT_{\! \bw}}_{N^\bw_s} > a \big) . 
\end{equation}

\subsection{Red and blue processes.}
\label{redbldissec}
%\tt{"NEW "SUBSECTION MADE WITH MATERIAL FROM THE LONG INTRODUCTION FROM THE "BIG" PAPER.}
This section contains no new result and we recall here more precisely the embedding  
of the LIFO queue without repetition coding the multiplicative graph $\cG_{\! \bw}$ into the Markovian queue considered in the previous Section \ref{connecMqu}. This embedding has been introduced in \cite{BDW1} (and it is informally recalled in Introduction).  
This embedding uses two auxiliary processes, the so-called {\it blue} and {\it red} processes, that are defined as follows. 
First, we introduce two independent random point measures on $[0, \infty) \! \times \! \{ 1, \ldots , n\} $: 
\begin{equation}
\label{ccXbrwdef}\ccX_\bw^{\mathtt{b}}\! = \! \sum_{k\geq 1} \delta_{(\tau^{\mathtt{b}}_k , \Jtt^{\mathtt{b}}_k)} \quad \textrm{and} \quad  \ccX_\bw^{\mathtt{r}}\! = \! \sum_{k\geq 1} \delta_{(\tau^{\mathtt{r}}_k , \Jtt^{\mathtt{r}}_k)}, 
\end{equation}
that are Poisson point measures with intensity $\ell \! \otimes \! \nu_\bw$, where we recall that $\ell $ stands for the Lebesgue measure and that $\nu_\bw \! = \! \frac{1}{\sigma_1 (\bw)}\sum_{1 \leq j\leq n} w_j \delta_{j} $. The blue process $X^{\mathtt{b}, \bw}$ and the red process $X^{\mathtt{r}, \bw}$ are defined respectively by
 \begin{equation}
\label{Xbrwdef}
X^{\mathtt{b}, \bw}_t  =  -t + \sum_{k\geq 1} w_{\Jtt^{\mathtt{b}}_k}\un_{[0, t]} (\tau^{\mathtt{b}}_k) \quad  \text{and} \quad  X^{\mathtt{r}, \bw}_t  =  -t + \sum_{k\geq 1} w_{\Jtt^{\mathtt{r}}_k}\un_{[0, t]} (\tau^{\mathtt{r}}_k). 
\end{equation}
Note that $X^{\mathtt{b}, \bw}$ and $X^{\mathtt{r}, \bw}$ are two independent spectrally positive L\'evy processes with Laplace exponent $\psi_\bw$ given by (\ref{psibwdef}). 
For all $j\ino \{ 1, \ldots , n\}$ and all $t\in[0, \infty)$, we next set:  
\begin{equation}
\label{Njbdef}
 N^{\bw}_j (t)\! = \! \ccX_\bw^{\mathtt{b}} \big( [0, t] \! \times \! \{ j\} \big) \quad \textrm{and} \quad E^\bw_j = \inf \big\{ t\ino [0, \infty) : \ccX^{\mathtt{b}}_\bw ([0, t] \! \times \! \{ j\})\! = \! 1 \big\}. 
 \end{equation} 
Then the $N^\bw_j$ are independent homogeneous Poisson processes with jump-rate $w_j/ \sigma_1 (\bw)$ and the r.v.~$(\frac{w_j}{\sigma_1 (\bw)} E^\bw_j)_{1\leq j\leq n}$ are i.i.d.~exponentially distributed r.v.~with unit mean. %Note that $X^{\mathtt{b}, \bw}_t \! =\! -t + \sum_{1\leq j\leq n} w_j N^{\bw}_j (t)$. 
%Recall from \eqref{defYcHH} the $\bw$-LIFO queue and its encoding process $Y^{\bw}$. 
%We can now construct a version of $Y^{\bw}$ using the r.v.~contained in $\ccX_\bw^{\mathtt{b}}$: Set
We next set 
\begin{equation}
\label{YwAwSigw}
 Y^\bw_t \! = \! -t  \, +\!\!\!   \sum_{1\leq j \leq n} \!\!\!  w_j \un_{\{ E^\bw_j \leq t \}}
 \quad \textrm{and} \quad A^\bw_t =   X^{\mathtt{b}, \bw}_t -Y^\bw_t  =  
  \! \sum_{1\leq j\leq n} \!\! w_j (N^\bw_j (t)\! -\! 1)_+ \; .
\end{equation}
Here $Y^\bw$ is the algebraic load of the following queue without repetition that codes the multiplicative graph $\cG_{\! \bw}$ (as explained in Introduction):  
\textit{a single server is visited by $n$ clients labelled by $1, \ldots , n$; Client $j$ arrives at time $E^\bw_j$ and she/he requests an amount of time of service $w_j$; 
%\item[-] the $E_j$ are independent exponentially distributed r.v.~such that $\bE [E_j] \! =\! \sigma_1 (\bw)/ w_j$; 
a LIFO (last in first out) policy applies: whenever a new client arrives, the server interrupts the service of the current client (if any) and serves the newcomer; when the latter leaves the queue, the server resumes the previous service.}

We embed this queue without repetition into a Markovian one that is obtained from $(Y^\bw, A^\bw)$ and $X^{\mathtt{r}, \bw}$ as follows. 
We first introduce the following  time-change process that will play a prominent role: 
\begin{equation}
\label{thetabw}
 \theta^{\mathtt{b}, \bw}_t \!\! = t + \subo^{\mathtt{r}, \bw}_{A^\bw_t}, \;   \textrm{where for all $x\ino [0, \infty)$, we have set:} \;\;  \subo^{\mathtt{r}, \bw}_x \!  = \! \inf \big\{ t\ino [0, \infty ) \! :  X^{\mathtt{r} , \bw}_t \!\!\!  < \! - x \big\}, 
\end{equation}
with the convention that $\inf \emptyset \! = \! \infty$. 
We next recall various properties of $\theta^{\mathtt{b}, \bw}$ that are used in the sequel. 
To that end, let us first 
note that standard results on L\'evy processes (see e.g.~Bertoin's book \cite{Be96} Chapter VII) assert that 
$(\subo^{\mathtt{r}, \bw}_x)_{x\in [0, \infty)}$ is a (possibly killed) subordinator whose Laplace exponent is given by: 
\begin{equation}
\label{gamwexpo}
\forall\,\lambda\ino [0, \infty): \bE \big[ e^{-\lambda \subo^{\mathtt{r}, \bw}_x } \big]= e^{-x\psi^{-1}_\bw (\lambda) } \quad  \textrm{where} \quad  \psi^{-1}_\bw (\lambda)\! = \! \inf \big\{ u \ino [0, \infty ) : \psi_\bw (u) \! >\! \lambda \big\}.  
\end{equation}
Set $\varrho_{\bw}\! = \!  \psi^{-1}_\bw(0)$, the largest root of $\psi_\bw$. Then, $\varrho_{\bw}\! = \! 0$ in the subcritical or critical cases, while $\varrho_{\bw}\! > \! 0$ in the supercritical case. Moreover, in the latter case, $-I^{\mathtt{r}, \bw}_\infty:=-\inf_{t\in [0, \infty )}X^{\mathtt{r}, \bw}_{t}$ is exponentially distributed with parameter $\varrho_\bw$ and $\gamma^{\texttt{r}, \bw}_x \! <\! \infty$ if and only if  $x\! <\! -I^{\mathtt{r}, \bw}_\infty$. It follows that the explosion time for $\theta^{\mathtt{b}, \bw}$ is given by
\begin{equation}
\label{T*wdef} 
T^*_\bw\! = \! \sup \{ t\ino [0, \infty)\! :  \theta^{\mathtt{b}, \bw}_t \! < \infty \}=   \sup \{ t\ino [0, \infty)\! : A^\bw_t \! < \! - I^{\mathtt{r}, \bw}_\infty \} \; ,
\end{equation}
which is infinite in the critical and subcritical cases and which is a.s.~finite in the supercritical cases. Note that $ \theta^{\mathtt{b}, \bw} (T^*_\bw-) \! < \! \infty$ in the supercritical cases. 
%We next recall from (\ref{Lambdefi}) the definition of the time-changes $\Lambda^{\mathtt{b}, \bw}$ and $\Lambda^{\mathtt{r}, \bw}$: 
%\begin{equation}
%\label{gagaoutch}
%\Lambda^{\mathtt{b}, \bw}_t \!  = \! \int_0^t \!\! \un_{\mathtt{Blue}} (s) \, ds =  \inf \big\{ s\ino [0, \infty) \! :  \theta^{\mathtt{b}, \bw}_s \! >\! t \}  \quad \textrm{and} \quad 
% \Lambda^{\mathtt{r}, \bw}_t \!\! = t \! -\!  \Lambda^{\mathtt{b}, \bw}_t \! \! = \! \! \int_0^t \!\! \un_{\mathtt{Red}} (s) \, ds. 
% \end{equation}

We shall also introduce the following processes:  
\begin{align}
\label{Lambdefi}
\Lambda^{\mathtt{b}, \bw}_t \!  =   \inf \big\{ s\ino [0, \infty) \! :  \theta^{\mathtt{b}, \bw}_s \! >\! t \}  \quad \textrm{and} \quad 
 \Lambda^{\mathtt{r}, \bw}_t \!\! = t \! -\!  \Lambda^{\mathtt{b}, \bw}_t \!.
 \end{align}
Both processes $\Lambda^{\mathtt{b}, \bw}$ and $\Lambda^{\mathtt{r}, \bw}$ are continuous and nondecreasing. Moreover, a.s.~$\lim_{t\rightarrow \infty}  \Lambda^{\mathtt{r}, \bw}_t \! = \! \infty$. In critical and subcritical cases, we also get a.s.~$\lim_{t\rightarrow \infty}  \Lambda^{\mathtt{b}, \bw}_t \! = \! \infty$ and $\Lambda^{\mathtt{b}, \bw} (\theta^{\mathtt{b}, \bw}_t)\! = \! t$ for all $t\ino [0, \infty)$.  However, in supercritical cases, $ \Lambda^{\mathtt{b}, \bw}_t \! = \! T^*_\bw$ for $t\ino [\theta^{\mathtt{b}, \bw} (T^*_\bw-), \infty)$ and a.s.~for all $t\ino [0, T^*_\bw)$, 
$\Lambda^{\mathtt{b}, \bw} (\theta^{\mathtt{b}, \bw}_t)\! = \! t$. 
The following proposition was proved in \cite{BDW1}.   
\begin{prop}
%[Proposition 2.2 \cite{BDW1}]
\label{Xwfrombr} 
We keep the previous notation and 
we define the process $X^\bw$ by:   
\begin{equation}
\label{redblumix}
\forall t\ino [0, \infty) , \qquad X^{\bw}_t =  X^{\mathtt{b}, \bw}_{ \Lambda^{\mathtt{b}, \bw}_t } + X^{\mathtt{r}, \bw}_{ \Lambda^{\mathtt{r}, \bw}_t } \; .
\end{equation}
Then, $X^{\bw}$ has the same law as $X^{\mathtt{b}, \bw}$ and $X^{\mathtt{r}, \bw}$: namely, it is a spectrally positive L\'evy process with Laplace exponent $\psi_\bw$ as defined in (\ref{psibwdef}). 
Furthermore, 
we have  
\begin{equation}
\label{YXtheta}
\textrm{a.s.}\; \forall\,t\ino [0, T^*_\bw) , \quad Y^\bw_t \! = \! X^\bw_{\theta^{\mathtt{b}, \bw}_t } . 
\end{equation}
\end{prop}
\noi
\textbf{Proof.} See Proposition 2.2 in \cite{BDW1}. \cqfd 

\medskip

\begin{figure}[tb]
\centering
\includegraphics[scale=.7]{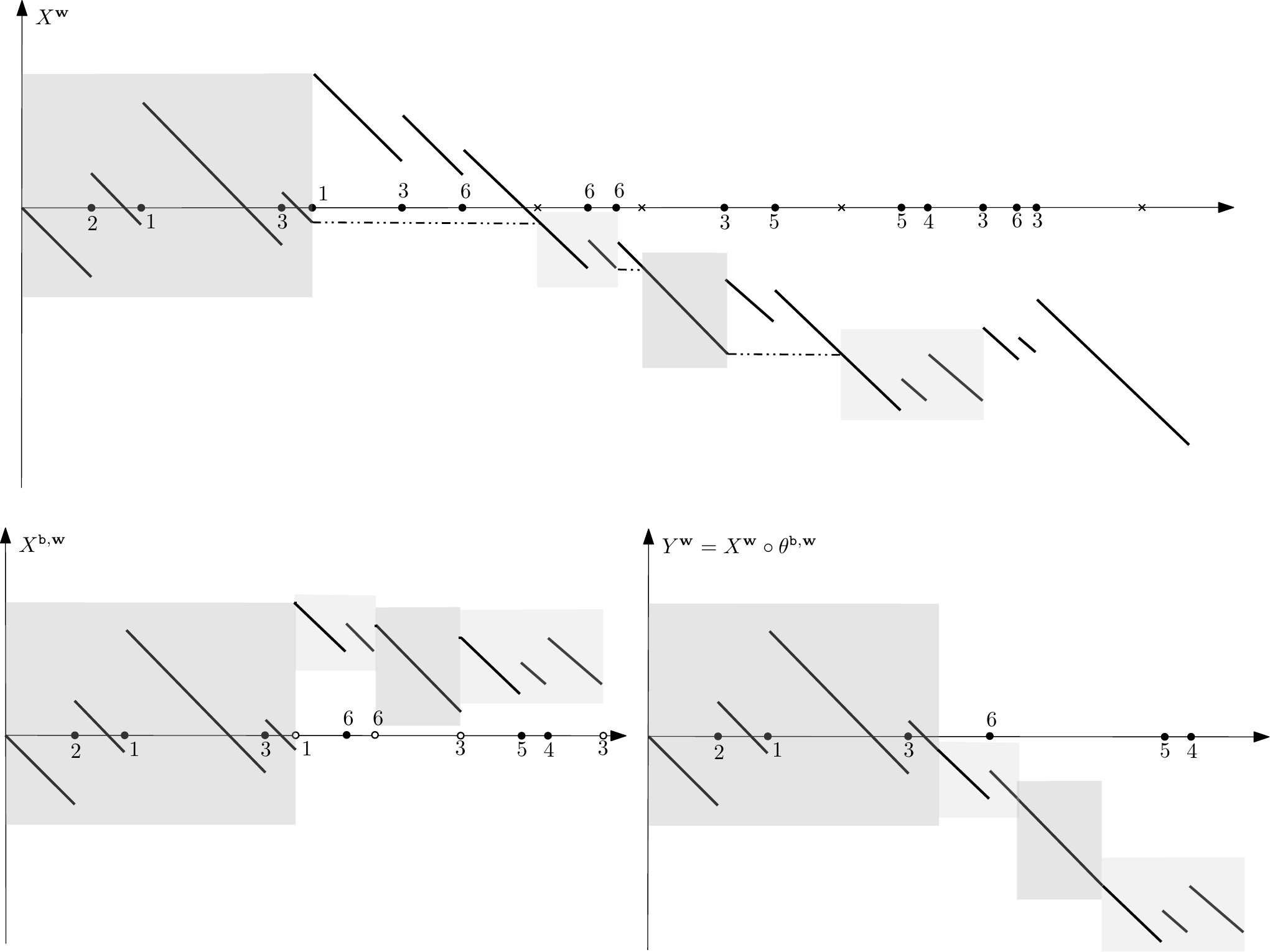}
\caption{{\small Decomposition of $X^\bw$ into $X^{\mathtt b, \mathbf w}$ and $X^{\mathtt r, \mathbf w}$. Above, the process $X^\bw$: clients are in bijection with its jumps; their types are the numbers next to the jumps. The grey blocks correspond to the set $\mathtt{Blue}$. 
Concatenating these blocks yields the blue process $X^{\mathtt b, \mathbf w}$. The remaining pieces give rise to the red process $X^{\mathtt r, \mathbf w}$. Concatenating the grey blocks but \textbf{without} the final jump of each block yields $Y^\bw$. Alternatively, we can obtain $Y^\bw$ by removing the temporal gaps between the grey blocks in $X^{\bw}$: this is the graphic representation of  $Y^\bw=X^\bw\circ\theta^{\mathtt b, \mathbf w}$.
Observe also that each connected component of $\mathtt{Red}$ begins with the arrival of a client whose type is a repeat among the types of the previous blue ones,  and ends with the departure of this red client, marked by  ${\scriptstyle\times}$ on the abscissa. \cq}}
\label{fig:proc_color} 
\end{figure}

Recall that $\mathtt{Blue}$ and 
$\mathtt{Red}$ are the sets of times during which respectively blue and red clients are served (the server is considered as a blue client). Then formally these sets are given by: 
\begin{equation}
\label{gramidos}
\mathtt{Red} = \! \!\!\!\!\!\!\!\! \!\!\!\!\! \bigcup_{\qquad t\in [0, T^*_\bw]: \Delta \theta^{\mathtt{b}, \bw}_t >0}  \!\!\!\!\!\!\!\! \!\!\!\!\!  \!\! \big[  \theta^{\mathtt{b}, \bw}_{t-}, 
\theta^{\mathtt{b}, \bw}_t \big) \quad \textrm{and} \quad \mathtt{Blue} \! = \! [0, \infty) \backslash \mathtt{Red} . 
\end{equation}
%Here we interpret $ \theta^{\mathtt{b}, \bw}_t$ as the time in the scale of the two-colours Markovian LIFO queue that  corresponds to $t$ in the "blue" scale of time of the 
%$\bw$-LIFO queue. Namely, thanks to $\theta^{\mathtt{b}, \bw}$, we define the set of times 
%$\mathtt{Blue}$ (resp.~$\mathtt{Red}$) that are the times 
%during which the blue clients (resp.~the red clients) are served (and 
%we recall that the server is considered as a blue client). More precisely, we set:  
%\begin{equation} 
%\label{blureddef}
%\mathtt{Red}\! = \! \bigcup_{t\in [0, \infty)} \big[  \theta^{\mathtt{b}, \bw}_{t-}, 
%\theta^{\mathtt{b}, \bw}_t \big) \quad \textrm{and} \quad \mathtt{Blue} \! = \! [0, \infty) \backslash \mathtt{Red} . 
%\end{equation}
Note that the union defining $\mathtt{Red}$ is countably infinite in critical and subcritical cases and that it is a finite union  
in supercritical cases where 
$\big[  \theta^{\mathtt{b}, \bw} (T^*_\bw-),  \theta^{\mathtt{b}, \bw} (T^*_\bw))\! = \! [\theta^{\mathtt{b}, \bw} (T^*_\bw-), \infty)$. We next recall from (\ref{Lambdefi}) the definition of the time-changes $\Lambda^{\mathtt{b}, \bw}$ and $\Lambda^{\mathtt{r}, \bw}$; then, we easily check that 
\begin{equation}
\label{gagaoutch}
\Lambda^{\mathtt{b}, \bw}_t \!  = \! \int_0^t \!\! \un_{\mathtt{Blue}} (s) \, ds 
%=  \inf \big\{ s\ino [0, \infty) \! :  \theta^{\mathtt{b}, \bw}_s \! >\! t \} 
 \quad \textrm{and} \quad 
 \Lambda^{\mathtt{r}, \bw}_t \!\! = t \! -\!  \Lambda^{\mathtt{b}, \bw}_t \! \! = \! \! \int_0^t \!\! \un_{\mathtt{Red}} (s) \, ds. 
 \end{equation}
We have the following properties of $X^{\bw}$, $\theta^{\mathtt{b}, \bw}$, etc.~that are recalled from \cite{BDW1} (see Figure \ref{fig:proc_color}).
\begin{lem}
%[Lemma 4.1 in \cite{BDW1}]
\label{gumodis} A.s.~for all $b\ino [0, T^*_\bw]$ such that $\theta^{\mathtt{b}, \bw}_{b-}<\theta^{\mathtt{b}, \bw}_{b}$, we get the following  for all $s\ino [\theta^{\mathtt{b}, \bw}_{b-} , \theta^{\mathtt{b}, \bw}_{b} )$:
\begin{equation}
\label{flaccidos}
 %\forall s\ino [\theta^{\mathtt{b}, \bw}_{t-} , \theta^{\mathtt{b}, \bw}_{t} ), \quad 
 X^{\bw}_s \! > \! X^\bw_{(\theta^{\mathtt{b}, \bw}_{b-})-}\!\!\!  = \! Y^{\bw}_b , \; \, \, \Delta X^\bw_{\theta^{\mathtt{b}, \bw}_{b-}}\! = \! \Delta A^{\bw}_b \; \,\,  \textrm{and} \; \, X^\bw_{(\theta^{\mathtt{b}, \bw}_{b-})-}\! \! \! = \! X^\bw_{\theta^{\mathtt{b}, \bw}_{b}}=X^\bw_{(\theta^{\mathtt{b}, \bw}_{b})-}  \; \, \textrm{if}  \; \, \theta^{\mathtt{b}, \bw}_{b}\!  < \infty. 
\end{equation}
Thus, a.s.~for all $s\ino [0, \infty)$, $X^\bw_s \! \geq \! Y^\bw(\Lambda^{\mathtt{b}, \bw}_s)$. Moreover, a.s.~for all $s_1, s_2 \ino [0, \infty)$ such that $\Lambda^{\mathtt{b}, \bw}_{s_1} \! < \!  \Lambda^{\mathtt{b}, \bw}_{s_2}$, then 
\begin{equation}
\label{flaccid}
\inf_{b\in [\Lambda^{\mathtt{b}, \bw}_{s_1} ,  \Lambda^{\mathtt{b}, \bw}_{s_2}]} Y^\bw_b= \inf_{s\in [s_1, s_2]} X^\bw_s\; .
\end{equation}
We next introduce the red time-change: 
\begin{equation}
\label{thetarw}
\theta^{\mathtt{r}, \bw}_t = \inf \big\{s \ino [0, \infty): \Lambda^{\mathtt{r}, \bw}_s >t \,  \big\} 
\end{equation}
Then, for all $s, t\ino [0, \infty)$, $\theta^{\mathtt{r}, \bw}_{s+t} \! -\! \theta^{\mathtt{r}, \bw}_t \! \geq \! s$ and if $\Delta \theta^{\mathtt{r}, \bw}_t \! >\! 0$, then $\Delta 
X^{\mathtt{r}, \bw}_t \! = \! 0$. 
\end{lem}
\noi
\textbf{Proof.} See Lemma 4.1 in \cite{BDW1}. \cqfd 

\bigskip

\noi
\textbf{Embedding of the tree.}  
The previous embedding of the LIFO queue without repetition governed by $Y^\bw$ into the Markovian queue governed by $X^\bw$ 
yields a related embedding of the trees associated with these queues.
More precisely, consider first the queue governed by $Y^\bw$: the LIFO rule implies that Client $i$ arriving at time $E_i$ will leave the queue at the moment $\inf\{t\ge E_i: Y^\bw_t<Y^\bw_{E_i-}\}$, namely the first moment when the service load falls back to the level right before her/his arrival. It follows that the number of clients waiting in queue at time $t$ is given by 
\begin{equation}
\label{JHdef}
%\forall t\ino [0, \infty) , \; \,  
\cH^\bw_t \! = \# \cJ_t, \;   \textrm{where} \;     \cJ_t\! = \! \big\{ s\ino [0, t] \!  : \!  J^{\bw,s -}_{t} \! <\! J^{\bw,s }_{t} \big\}   \; \textrm{and where} \;  \forall s\ino [0, t], \; J^{\bw,s }_t\! = \! \inf_{r\in [s, t]} \! Y^\bw_r . 
\end{equation}
Recall that we denote by $\cT_{\!\! \bw}$ the tree formed by the clients in the queue governed by $Y^\bw$.  
The process $\cH^{\bw}$ is actually the contour (or the depth-first exploration) of $\cT_{\!\! \bw}$ and the graph-metric $d_{\cT_{\!\! \bw} }$ of $\cT_{\!\! \bw}$ is encoded by $\cH^\bw$ in the following way: if we denote by $V_{t}\in \{0, 1, \dots, n\}$ the label of the client served at time $t$ (with the understanding that $V_{t}=0$ if the server is idle), then
\begin{equation}
\label{gloupii}
\forall s, t\ino [0, \infty), \quad d_{\cT_{\!\! \bw} } (V_s, V_t) = \cH^\bw_t+\cH^\bw_s - \, 2 \!\!\!\!\!\!  \min_{\quad r \in [s  \wedge t , s\vee t] }  \!\!\!\!\! \!    \cH^\bw_r  \; .
\end{equation}

Similarly for the Markovian queue governed by the process $X^{\bw}$ given in Proposition \ref{Xwfrombr}, 
we define its  associated \textit{height process} $H^\bw$ by setting $H^\bw _t$ to be the number of the clients waiting at time $t$, namely, 
\begin{equation}
\label{XJHdef}
H^\bw _t \! = \# \cK_t, \;   \textrm{where} \;    \cK_t \! =\!  \big\{ s\ino [0, t]\!  : \! I^{\bw, s-}_{t} \! <\! I^{\bw, s}_{t} \big\}   \;    \textrm{and where} \;     \forall s\ino [0, t], \;  I^{\bw, s}_t= \inf_{r\in [s, t]}X^\bw_r . 
\end{equation}
Then $H^{\bw}$ is the contour process of the i.i.d.~Galton--Watson forest $\bT_{\bw}$ with offspring distribution $\mu_{\bw}$ characterized by \eqref{treeTrts}. 
Note that in (sub)critical cases, $H^\bw$ fully explores the whole tree $\bT_{\bw}$. 
However in supercritical cases, the exploration of $H^\bw$ does not go beyond the first infinite line of descent. 
We shall use the following form of the previsouly mentioned embedding of $\cT_{\!\! \bw}$ into $\bT_{\bw}$ that is recalled from \cite{BDW1}. 
\begin{lem}
%[Lemma 2.3 \cite{BDW1}]
\label{Hthetalem} Following the previous notation, we have
\begin{equation}
\label{HYenHX}
\textrm{a.s.} \; \forall t\ino [0, T^*_\bw), \quad \cH^\bw_t = H^\bw_{\theta^{\mathtt{b} , \bw}_t}  \; . 
\end{equation}
\end{lem}  
\noi
\textbf{Proof.} See Lemma 2.3 in \cite{BDW1}. \cqfd

\subsection{Estimates on the coloured processes.}
\label{estisec}
We keep notation from the Section \ref{redbldissec}. 
In this section, we now provide estimates for $A^\bw$, $X^{\mathtt{b}, \bw}_{\Lambda^{\mathtt{b}, \bw}}$ and 
$X^{\mathtt{r}, \bw}_{\Lambda^{\mathtt{r}, \bw}}$ that are used in the proof of Theorem \ref{HYcvth}.

Recall from (\ref{YwAwSigw}) that $A^\bw_t \! = \! \sum_{j\geq 1} \!\! w_j (N^\bw_j (t)\! -\! 1)_+$ where the 
$N^\bw_j (\cdot)$ 
are independent homogeneous Poisson processes with respective jump-rate $w_j/ \sigma_1 (\bw)$. 
Let $(\ccF_t)_{t\in [0, \infty)}$ be a filtration such that for all $j\! \geq \! 1$, $N_j^{\bw}$ is a $(\ccF_t)$-homogeneous Poisson process. Namely, 
\begin{compactenum}
\item[$\bullet$] $N^\bw_j$ is $(\ccF_t)$-adapted; 
\item[$\bullet$] for all a.s.~finite $(\ccF_t)$-stopping time $T$, 
set $N^{_{\bw , T}}_{^j} (t)\! = \! N^\bw_j (T+ t) \! -\! N^\bw_j (T)$. Then, the sequence of processes 
$(N^{_{\bw , T}}_{^j})_{j\geq 1}$ is independent of $\ccF_T$ and distributed as $(N^\bw_j)_{j\geq 1}$.
\end{compactenum}
 Thus, the process $A^{\bw, T} \! = \! \sum_{j\geq 1} \!\! w_j (
N^{\bw, T}_j (\cdot)\! -\! 1)_+$ is independent of $\ccF_T$ and distributed as $A^\bw$. 
We easily get 
\begin{equation}
\label{almindAw}
A^{\bw}_{T+t}  - A^\bw_T =  A^{\bw , T}_t + \sum_{j\geq 1}w_j  \un_{\{ E^\bw_j \leq T \}} \un_{\{ N^{\bw , T}_j (t) \geq 1 \}} , 
\end{equation}
where recall from (\ref{Njbdef}) that $E^\bw_j$ stands for the first jump-time of $N^\bw_j$;  $E^\bw_j$ is therefore exponentially distributed with mean $\sigma_1 (\bw)/ w_j$. 
Elementary calculations combined with (\ref{almindAw}) immediately entail the following lemma. 
\begin{lem}
\label{EstDifAw} We keep the notation from above. For all $(\ccF_t)$-stopping time $T$ and all 
$a, t_0 , t \ino (0, \infty)$, 
\begin{equation}
\label{Difldcrioui}
a \, \bP \big( T\leq t_0 \, ; \, A^\bw_{T+t}\! -\! A^\bw_T \geq a \big)  \leq  \bE [A^\bw_t] + \sum_{j\geq 1} w_j 
\bP (E^\bw_j \leq t_0 ) \bP (N^{\bw }_j (t) \geq 1). 
\end{equation}
Note that $\bE [A^\bw_t]\! = \! \sum_{j\geq 1} w_j (e^{-tw_j / \sigma_1 (\bw)} \! -\! 1 +  \frac{ tw_j}{\sigma_1 (\bw)} )$. Thus,  
\begin{equation}
\label{goretex}
a \, \bP \big( T\leq t_0 \, ; \, A^\bw_{T+t}\! - \! A^\bw_T \geq a \big)  \leq     t \big( t_0 + \frac{_{_1}}{^{^2}} t \big)  \frac{_{\sigma_3 (\bw)}}{^{\sigma_1 (\bw)^2}}  . 
 \end{equation} 
\end{lem}
%
%\subsubsection{Oscillations of $X^{\mathtt{b}, \bw}_{\Lambda^{\mathtt{b}, \bw}}$ and 
%$X^{\mathtt{r}, \bw}_{\Lambda^{\mathtt{r}, \bw}}$.}
%\label{oscblrd} 

\medskip

We next discuss the oscillations of $X^{\mathtt{b}, \bw}_{\Lambda^{\mathtt{b}, \bw}}$ and of 
$X^{\mathtt{r}, \bw}_{\Lambda^{\mathtt{r}, \bw}}$. To that end, let us recall that $\bD ([0, \infty) , \bbR)$ stands for the space of $\bbR$-valued c\`adl\`ag functions equipped with Skorokhod's topology. 
For all $y \ino \bD ([0, \infty) , \bbR)$ and for all intervals $I$ of $[0, \infty)$, we set 
\begin{equation}
\label{oscilldef}
\mathtt{osc} (y, I)= \sup \big\{ |y(s)\! -\! y(t) | ; \, s, t\ino I  \big\}\,, 
\end{equation} 
%\fmm{This was some inconsistence in the notation $\mathtt{osc}$. All changed to $\mathtt{osc}$ now. }
that is the \textit{oscillation} of $y$ on $I$. It is easy to check that for all $ a \! <\! b  \! <\! c$, 
\begin{equation}
\label{oscbouint}
\mathtt{osc} (y, [a, c) \, ) \leq  \mathtt{osc} (y, [a, b]) + \mathtt{osc} (y, [b, c))  \leq   \mathtt{osc} (y, [a, b) \, ) + |\Delta y(b)|  + \mathtt{osc} (y, [b, c))\; , 
\end{equation}
where we recall that $\Delta y (b)\! = \! y(b)\! -\! y(b-)$. We also recall the definition of the \textit{c\`adl\`ag modulus of continuity} of $y$: let $z, \eta \ino (0, \infty)$; then, we set 
\begin{equation}
\label{modudu}
w_z(y,\eta)=  \inf \big\{\!  \max_{1\leq i \leq r } \mathtt{osc} (y, [t_{i-1} , t_i ) \, )\;  ; \; \,   0\! = \! t_0 \! < \! \ldots  \! < \! t_r\! = \! z \;  :  \!\! \!\!\!\! \min_{\quad 1\leq i \leq r-1} \!\!\!\! (t_i\! -\! t_{i-1}) \geq \eta  \;    \big\} , 
\end{equation}
\noi
Here the infimum is taken on the set of all subdivisions $(t_i)_{0\leq i\leq r}$, of $[0, z]$, %\margmm{corrected}
$r$ being a positive integer; note that we do not require $t_r\! -\! t_{r-1} \! \geq \! \eta$. We refer to Jacod \& Shiryaev \cite{JaSh02} Chapter VI for a general introduction on Skorokod's topology. 
Recall from (\ref{T*wdef}) the definition of $T^*_\bw$ and from (\ref{Lambdefi}) the definition of $ \Lambda^{\mathtt{b}, \bw}$ and  $\Lambda^{\mathtt{r}, \bw}$. Recall from (\ref{redblumix}) in Lemma \ref{Xwfrombr} that $X^{\bw}\! = \! X^{\mathtt{b}, \bw}_{\Lambda^{\mathtt{b}, \bw}}+ X^{\mathtt{r}, \bw}_{\Lambda^{\mathtt{r}, \bw}}$. 
%%%%
%%%%To state the result on the oscillations of $X^{\mathtt{b}, \bw} \! \circ  \Lambda^{\mathtt{b}, \bw}$ and $X^{\mathtt{r}, \bw}\! \circ  \Lambda^{\mathtt{r}, \bw}$, we need to recall the following. 
%%%%By (\ref{gagaoutch}), $\Lambda^{\mathtt{b}, \bw}$ and $\Lambda^{\mathtt{r}, \bw}$ are continuous and nondecreasing and 
%%%%a.s.~$\lim_{t\rightarrow \infty}  \Lambda^{\mathtt{r}, \bw}_t \! = \! \infty$. Recall from (\ref{grumide}) the definition of $T^*_\bw$. 
%%%%In critical and subcritical cases (where $T^*_\bw\!= \! \infty$ a.s.) we get a.s.~$\lim_{t\rightarrow \infty}  \Lambda^{\mathtt{b}, \bw}_t \! = \! \infty$ and $\Lambda^{\mathtt{b}, \bw} (\theta^{\mathtt{b}, \bw}_t)\! = \! t$ for all $t\ino [0, \infty)$; in supercritical cases (where $T^*_\bw\! < \! \infty$), we a.s.~get $\Lambda^{\mathtt{b}, \bw} (\theta^{\mathtt{b}, \bw}_t)\! = \! t$, for all $t\ino [0, T^*_\bw)$ and 
%%%%$\Lambda^{\mathtt{b}, \bw}$ is constant to $T^*_\bw$ on $[ \theta^{\mathtt{b}, \bw} (T^*_\bw-) , \infty)$. 
%%%%
The following lemma is a key argument in the proof of Theorem \ref{HYcvth}.
\begin{lem}
\label{modconXL} We keep the notation from above. 
Let $ \eta \ino (0, \infty)$. Then, the following holds true. 
\begin{compactenum}

%\smallskip

\item[$(i)$] Almost surely, for all $z_0, z_1, z\ino [0, \infty)$, if $ z_1 \! \leq  \! \theta^{\mathtt{b}, \bw}_{z_0} \! \leq \! z$, 
%the event $\{ z_1 \! \leq  \! \theta^{\mathtt{b}, \bw}_{z_0} \! \leq \! z \}$, 
then we get 
\begin{equation} 
\label{Xlambcont}
%\textrm{a.s.~on $\{\Lambda^{\mathtt{b}, \bw}_{z_0} \!  \!> \! z  \! >\! \Lambda^{\mathtt{b}, \bw}_{z_1} \}$}, 
%\quad \qquad 
w_{z_1} \big( X^{\mathtt{b}, \bw}_{\Lambda^{\mathtt{b}, \bw}}\,  , \eta \big) \leq w_{z+\eta} \big(X^\bw\!  , \eta \big) + w_{z_0} \big(X^{\mathtt{b}, \bw}\! , \eta \big) . 
\end{equation}
\item[$(ii)$] Assume that we are in the supercritical cases (namely, $\alpha_\bw\! = \! 1\! -\! \frac{_{\sigma_2(\bw)}}{^{\sigma_1(\bw)}} \! < \! 0$) where 
a.s.~$T^*_{\bw}\!< \! \infty$ and $\theta^{\mathtt{b}, \bw} (T^*_\bw\! -)\! < \!\infty$. 
Then a.s.~for all $z_0, z_1, z\ino [0, \infty)$ if $z \! >\!  \theta^{\mathtt{b}, \bw} (T^*_\bw\! -)$ and $z_0 \! >\! T^*_{\bw} \! >\! 2\eta $, we get 
\begin{equation} 
\label{Xlambcontt}
%\textrm{a.s.~on $\{\Lambda^{\mathtt{b}, \bw}_{z_0} \!  \!> \! z  \! >\! \Lambda^{\mathtt{b}, \bw}_{z_1} \}$}, 
%\quad \qquad 
w_{z_1} \big( X^{\mathtt{b}, \bw}_{\Lambda^{\mathtt{b}, \bw}}\,  , \eta \big)
%w_{z+\eta} \big(X^\bw\!  , \eta \big) + w_{z_0} 
%\big(X^{\mathtt{b}, \bw}\! , \eta \big) +  2w_{z_0} \big(X^{\mathtt{b}, \bw}\! , 
%2\eta \big) 
\! \leq \! w_{z+\eta} \big(X^\bw\!  , \eta \big) +3w_{z_0} \big(X^{\mathtt{b}, \bw}\! , 
2\eta \big).
\end{equation}
\item[$(iii)$] Almost surely  on the event $\{ z  \! >\! \Lambda^{\mathtt{r}, \bw}_{z_1} \}$, we get $w_{z_1} \big( X^{\mathtt{r}, \bw}_{\Lambda^{\mathtt{r}, \bw}}\,  , \eta \big) \leq  w_{z} 
 \big(X^{\mathtt{r}, \bw}\! , \eta \big)$. 
%\begin{equation} 
%\label{Xlamrcont}
%\textrm{a.s.~on $\{ z  \! >\! \Lambda^{\mathtt{r}, \bw}_{z_1} \}$}, 
%\quad \qquad w_{z_1} \big( X^{\mathtt{r}, \bw}_{\Lambda^{\mathtt{r}, \bw}}\,  , \eta \big) \leq  w_{z} 
% \big(X^{\mathtt{r}, \bw}\! , \eta \big) . 
% \end{equation}
\end{compactenum}
\end{lem}
\noi
\textbf{Proof.} First note that for all intervals $I$, we get: 
$$ \mathtt{osc} \big( X^{\mathtt{b}, \bw}_{\Lambda^{\mathtt{b}, \bw}} , I \big) = \sup \big\{ \big|X^{\mathtt{b}, \bw}_{\Lambda^{\mathtt{b}, \bw}_t}\! -\! X^{\mathtt{b}, \bw}_{\Lambda^{\mathtt{b}, \bw}_s}   \big|  ; s, t \ino I  \big\}=  
\sup \big\{ \big|X^{\mathtt{b}, \bw}_{t}\! -\! X^{\mathtt{b}, \bw}_{s}   \big|  ; s, t \ino  \big\{\Lambda^{\mathtt{b}, \bw}_u ; u\ino I  \big\}   \big\}. $$
We fix $\eta, a, b\ino [0, T^*_\bw)$ such that $b\! - \! a \! \geq  \! \eta$. By the definition (\ref{thetabw}) of $\theta^{\mathtt{b}, \bw}$, we get 
$\theta^{\mathtt{b}, \bw}_{b-} \! \! -\!  \theta^{\mathtt{b}, \bw}_{a} \! \geq \! b\! -\! a \! \geq \! \eta$. Since $\Lambda^{\mathtt{b}, \bw}$ is non-decreasing and continuous, and since $\theta^{\mathtt{b}, \bw}$ is strictly increasing, 
we get $\{\Lambda^{\mathtt{b}, \bw}_t ; t\ino [\theta^{\mathtt{b}, \bw}_{a} , \theta^{\mathtt{b}, \bw}_{b-} )  \}\! = \! [a, b)$ and 
\begin{equation}
\label{gloubinett}  
\mathtt{osc} \big( X^{\mathtt{b}, \bw}_{\Lambda^{\mathtt{b}, \bw}} ,[\theta^{\mathtt{b}, \bw}_{a} , \theta^{\mathtt{b}, \bw}_{b-}) \, \big)  \! = \!      \mathtt{osc} \big( X^{\mathtt{b}, \bw} , [a, b) \, \big) . 
\end{equation}

We next suppose that $\Delta  \theta^{\mathtt{b}, \bw}_{b} \! \! >\! 0$. Then, $\{\Lambda^{\mathtt{b}, \bw}_t ; t\ino [\theta^{\mathtt{b}, \bw}_{a} , \theta^{\mathtt{b}, \bw}_{b} )  \}\! = \! [a, b]$ and by (\ref{oscbouint}), we get 
\begin{equation}
\label{modjump}
\mathtt{osc} \big( X^{\mathtt{b}, \bw}_{\Lambda^{\mathtt{b}, \bw}} ,[\theta^{\mathtt{b}, \bw}_{a} , \theta^{\mathtt{b}, \bw}_{b}) \, \big)  \! = \!      \mathtt{osc} \big( X^{\mathtt{b}, \bw} , [a, b] \, \big) \leq      \mathtt{osc} \big( X^{\mathtt{b}, \bw} , [a, b) \, \big) + 
|\Delta  X^{\mathtt{b}, \bw}_b| . 
\end{equation} 
Since the process $ X^{\mathtt{b}, \bw}_{\Lambda^{\mathtt{b}, \bw}} $ is constant on $[\theta^{\mathtt{b}, \bw}_{b-} , \theta^{\mathtt{b}, \bw}_{b})$, we get $\mathtt{osc} \big( X^{\mathtt{b}, \bw}_{\Lambda^{\mathtt{b}, \bw}} ,[\theta^{\mathtt{b}, \bw}_{b-} , \theta^{\mathtt{b}, \bw}_{b}) \big)\! = \! 0$ and thus 
\begin{equation}
\label{tranqmod}
 \max \Big(  \mathtt{osc} \big( X^{\mathtt{b}, \bw}_{\Lambda^{\mathtt{b}, \bw}} ,[\theta^{\mathtt{b}, \bw}_{a} , \theta^{\mathtt{b}, \bw}_{b-})\big) , \mathtt{osc} \big( X^{\mathtt{b}, \bw}_{\Lambda^{\mathtt{b}, \bw}} ,[\theta^{\mathtt{b}, \bw}_{b-} , \theta^{\mathtt{b}, \bw}_{b}) \big)\Big)=    \mathtt{osc} \big( X^{\mathtt{b}, \bw} , [a, b) \, \big) . 
 \end{equation}
 
   We next assume that $\Delta  \theta^{\mathtt{b}, \bw}_{b} \ino (0, \eta)$. We 
want to control $|\Delta  X^{\mathtt{b}, \bw}_b| $ in terms of the c\`adl\`ag $\eta$-modulus of continuity 
of $X^\bw$. 
%
%To that end, first observe that $\Delta  \theta^{\mathtt{b}, \bw}_{b}\! >\! 0$ implies that $[ \theta^{\mathtt{b}, \bw}_{b-}, \theta^{\mathtt{b}, \bw}_{b} )$ is a connected component of $\mathtt{Red}$; by (\ref{thenRed}), there exists $p\! \geq \! 1$ such that \margmm{Pb: cahnge of notation}\mm{$\theta^{\mathtt{b}, \bw}_{b-} \! = \! \frs_{j(p)}$: namely, $\theta^{\mathtt{b}, \bw}_{b-} $ is the time of arrival of a red client who interrupts a blue one.  
%%By the definition (\ref{defZp}), by (\ref{Saucisred}) and by $(a)$ in the recursive definition of the $\frs_l$, 
%By (\ref{rbrebond}), for all $t\ino [\frs_{j(p)} , \frs_{j(p)+1})$, $X^\bw_t \! \geq \! X^\bw_{\frs_{j(p)}-}\! \!\! = \! Y^\bw (\Lambda^{\! \mathtt{b}, \bw}_{\frs_{j(p)}})$}; moreover 
%$\Delta X^{\mathtt{b}, \bw}_b \! = \! \Delta X^{\bw}_{\theta^{\mathtt{b} , \bw}_{b-}}$. To summarize: 
%\begin{equation}
%\label{ptclefcont}
%\forall t\ino [\theta^{\mathtt{b} , \bw}_{b-} , \theta^{\mathtt{b} , \bw}_{b} ), \quad X^\bw_t \! >\! X^{\bw}_{(\theta^{\mathtt{b}, \bw}_{b-}) -}\!\! = \! X^{\bw}_{\theta^{\mathtt{b}, \bw}_{b}} 
%\quad \textrm{and} \quad 
%\Delta X^{\bw}_{\theta^{\mathtt{b} , \bw}_{b-}} \! = \! X^{\bw}_{\theta^{\mathtt{b}, \bw}_{b-}}\! -\! X^{\bw}_{(\theta^{\mathtt{b}, \bw}_{b-}) -} \! =\!  \Delta X^{\mathtt{b}, \bw}_b
%\end{equation}
To that end, let us introduce $z\ino (0, \infty)$ such that $ \theta^{\mathtt{b}, \bw}_{b-}\! \leq  \! z$ and 
$0\! = \! t_0 \! < \! \ldots  \! < \! t_r\! = \! z+\eta$ such that $\min_{1\leq i \leq r-1} (t_i\! -\! t_{i-1}) \! \geq \! \eta$. 
Then, there exists $i\ino \{ 1, \ldots , r\}$ such that 
$t_{i-1} \! \leq \! \theta^{\mathtt{b}, \bw}_{b-} \! < \! t_{i} $ and necessarily $i$ satisfies $t_i\! -\! t_{i-1}\! \geq \! \eta$: \textit{indeed}, it is clear if $i\! <\! r$ and if $i\! =\! r$, then $t_{r-1} \! \leq \! \theta^{\mathtt{b}, \bw}_{b-} \! \leq z \!< \!z+ \eta \! =\! t_{r} $. There are two cases to consider: 

\begin{compactenum}

\smallskip

\item[--] If $t_{i-1} \! < \! \theta^{\mathtt{b}, \bw}_{b-}$, then $\mathtt{osc} (X^\bw\! , [t_{i-1}, t_{i} )) \! \geq \! |\Delta X^\bw ( \theta^{\mathtt{b}, \bw}_{b-})|$. 
Since $\theta^{\mathtt{b}, \bw}_{b}\! < \! \infty$, (\ref{flaccidos}) in Lemma \ref{gumodis} implies that 
$|\Delta X^\bw ( \theta^{\mathtt{b}, \bw}_{b-})|\! = \! |\Delta X^{\mathtt{b}, \bw}_b|$. Thus, $\mathtt{osc} (X^\bw\! , [t_{i-1}, t_{i} )) \! \geq \! |\Delta X^{\mathtt{b}, \bw}_b|$.

\smallskip

\item[--] If $t_{i-1}\! = \! \theta^{\mathtt{b}, \bw}_{b-}$, since $\Delta  \theta^{\mathtt{b}, \bw}_{b} \ino (0, \eta)$ and since $t_i - t_{i-1}\! \geq \! \eta$, we get $ \theta^{\mathtt{b}, \bw}_{b} \! < \! t_{i} $. Then 
$\mathtt{osc} (X^\bw, [t_{i-1}, t_{i} ))$ $ \geq \! 
| X^{\bw} (\theta^{\mathtt{b}, \bw}_{b-}) \! -\!   X^{\bw} (\theta^{\mathtt{b}, \bw}_{b})|$. 
Since $\theta^{\mathtt{b}, \bw}_{b}\! < \! \infty$, (\ref{flaccidos}) in Lemma \ref{gumodis} entails  
$X^{\bw} ((\theta^{\mathtt{b}, \bw}_{b-}) -)\! = \! 
X^{\bw} (\theta^{\mathtt{b}, \bw}_{b}) $ and 
$| X^{\bw} (\theta^{\mathtt{b}, \bw}_{b-}) \! -\!   X^{\bw} (\theta^{\mathtt{b}, \bw}_{b})| \! = \!  |\Delta X^\bw ( \theta^{\mathtt{b}, \bw}_{b-})|\! = \! |\Delta X^{\mathtt{b}, \bw}_b| $. Consequently, 
$\mathtt{osc} (X^\bw, [t_{i-1}, t_{i} )) \! \geq \!  |\Delta X^{\mathtt{b}, \bw}_b|$. 

\end{compactenum}

\smallskip

\noi
We have proved that if $\Delta  \theta^{\mathtt{b}, \bw}_{b} \! \ino (0, \eta)$ and if $ \theta^{\mathtt{b}, \bw}_{b-}\! \leq  \! z$, then $|\Delta X^{\mathtt{b}, \bw}_b| \! \leq \! \max_{1\leq i\leq r} 
\mathtt{osc} (X^\bw, [t_{i-1}, t_{i} ) ) $; since it holds true for all subdivisions of $[0, z+\eta ]$ 
satisfying the conditions as above, we get 
\begin{equation}
\label{ptitcontr}
 \textrm{a.s.~on $\{ \theta^{\mathtt{b}, \bw}_{b-}\! \leq  \! z \, ; \, \Delta  \theta^{\mathtt{b}, \bw}_{b} \! \ino (0, \eta) \}$},  \quad |\Delta X^{\mathtt{b}, \bw}_b| \! \leq \! w_{z+\eta} \big( X^\bw , \eta\big). 
\end{equation}  

We are now ready to prove (\ref{Xlambcont}). 
%\margmm{corrected}
Let us fix $z_0, z \ino (0, \infty)$ and let $0\! = \! t_0 \! < \! \ldots  \! < \! t_r\! =  z_0$ be such that $\min_{1\leq i \leq r-1} (t_i\! -\! t_{i-1}) \! \geq \! \eta$. We assume that $\theta^{\mathtt{b} , \bw}_{z_0} \! \leq \! z$. 
For all $i\ino \{ 1, \ldots , r\}$, we set $S_i \! = \! \{   \theta^{\mathtt{b} , \bw}_{t_i} \}$ if $\Delta \theta^{\mathtt{b} , \bw}_{t_i} \! < \! \eta$ and we set $S_i \! = \! \{   \theta^{\mathtt{b} , \bw}_{t_i-} ,   \theta^{\mathtt{b} , \bw}_{t_i} \}$ if $\Delta \theta^{\mathtt{b} , \bw}_{t_i} \! \geq  \! \eta$; we then define $S\! =\!  \{ s_0\! = \! 0 < \ldots \! < \! s_{r^\prime}\! = \!  \theta^{\mathtt{b} , \bw}_{z_0}\} = \{ 0\} \cup S_1 \cup \ldots  \cup S_r $ that is a subdivision of $[0,  \theta^{\mathtt{b} , \bw}_{z_0}]$ such that 
$\min_{1\leq i \leq r^\prime-1} (s_i\! -\! s_{i-1}) \! \geq \! \eta$ (\textit{indeed}, recall that $\theta^{\mathtt{b} , \bw}_{t_i-}\! -\! \theta^{\mathtt{b} , \bw}_{t_{i-1}} \! \geq \! t_i \! -\! t_{i-1}$).
By (\ref{tranqmod}) (if $S_i$ has two points) and by 
(\ref{modjump}) and (\ref{ptitcontr}) (if $S_i$ reduces to a single point), we get 
$$ w_{\theta^{\mathtt{b} , \bw}_{z_0}} \big(X^{\mathtt{b}, \bw}_{\Lambda^{\mathtt{b}, \bw}} , \eta)  \leq  \max_{1\leq i\leq r^\prime} \Big( \mathtt{osc}   \big( X^{\mathtt{b}, \bw}_{\Lambda^{\mathtt{b}, \bw}} , [s_{i-1} , s_i)\big)   \Big)  
\leq  w_{z+\eta} 
\big( X^\bw\! , \eta\big) +  \max_{1\leq i\leq r} \Big( \mathtt{osc}   \big( X^{\mathtt{b}, \bw}\! ,\,  [t_{i-1} , t_i)\big)   \Big).   $$
Since it holds true for all subdivisions $(t_i)$ and since $z^\prime \! \mapsto \! w_{z^\prime} ( y(\cdot), \eta)$ 
is nondecreasing, it easily entails (\ref{Xlambcont}) if $z_1\!  \leq \!  \theta^{\mathtt{b} , \bw}_{z_0} \! \leq \! z$, which completes the proof of $(i)$. 

Let us prove $(ii)$. We assume that we are in the supercritical cases.  
The control of the c\`adl\`ag modulus of continuity of 
$X^{\mathtt{b}, \bw} \circ \Lambda^{\mathtt{b}, \bw}$ is more complicated because this process becomes eventually constant after a last jump at time $\theta^{\mathtt{b}, \bw} (T^*_\bw\! -)$. 
To simplify notation we set $\tau\! = \!  \theta^{\mathtt{b}, \bw} (T^*_\bw\! -)$. We suppose that $z \! >\! \tau $ and 
$z_0 \! >\! T^*_{\bw} \! >\! 2\eta $. We fix $z_1\ino (0, \infty)$. There are several cases to consider.

\smallskip

\noi
$\bullet$ We first assume that $z_1 \! \leq \! 
\tau$. If $z_{1}<\tau$, then there is $z^\prime_0\ino [0, T^*_{\bw})$ 
 such that $z_1 \! \leq  \! \theta^{\mathtt{b}, \bw}_{z_0^\prime}$; next, note that $ \theta^{\mathtt{b}, \bw}_{z_0^\prime}\! \leq  \! z $ and $z^\prime_0 \! \leq \! z_0$. Thus, applying (\ref{Xlambcont}) to $(z_{0}^{\prime}, z_{1}, z)$, we get that $ w_{z_1} ( X^{\mathtt{b}, \bw}_{\Lambda^{\mathtt{b}, \bw}}\,  , \eta ) \! \leq \! w_{z+\eta} (X^\bw\!  , \eta) 
 + w_{z^\prime_0} (X^{\mathtt{b}, \bw}\! , \eta ) \leq  w_{z+\eta} (X^\bw\!  , \eta ) + w_{z_0} (X^{\mathtt{b}, \bw}\! , \eta )$ for all $z_{1}<\tau$. We then extend this to $z_{1}\le \tau$ by using a basic property of the c\`adl\`ag modulus of continuity:  
$\lim_{z_1 \rightarrow \tau-} w_{z_1} ( X^{\mathtt{b}, \bw}_{\Lambda^{\mathtt{b}, \bw}}\,  , \eta )\! = \! w_{\tau} ( X^{\mathtt{b}, \bw}_{\Lambda^{\mathtt{b}, \bw}}\,  , \eta )$.
Thus we have proved for all $z_1 \ino [0, \tau]$, 
\begin{equation}
\label{gabardine}
w_{z_1} \big( X^{\mathtt{b}, \bw}_{\Lambda^{\mathtt{b}, \bw}}\,  , \eta \big)  \leq  
 w_{z+\eta} \big(X^\bw\!  , \eta \big) + w_{z_0} \big(X^{\mathtt{b}, \bw}\! , \eta \big) \; , 
 \end{equation} 
which implies (\ref{Xlambcontt}) when $z_1\l \leq \! \tau$.

\smallskip

\noi
$\bullet$ We next assume that $z_1\! >\! \tau$. Observe that 
$\Delta (X^{\mathtt{b}, \bw}\!  \circ \! \Lambda^{\mathtt{b}, \bw})(\tau)\! = \! \Delta X^{\mathtt{b} , \bw} (T^*_{\bw})$. There are two subcases to consider. 

\smallskip

\begin{compactenum}
\item[$\circ$] We first assume that $z_1\! >\! \tau$ and that $\Delta X^{\mathtt{b}, \bw} (T^*_{\bw}) \! \leq \!  w_{z_0} (X^{\mathtt{b}, \bw}\! , 2\eta )$. 
As an easy consequence of (\ref{oscbouint}) and of the definition (\ref{modudu}) 
of the c\`adl\`ag modulus of continuity, we get 
$w_{z_1} \big( X^{\mathtt{b}, \bw}_{\Lambda^{\mathtt{b}, \bw}}\,  , \eta \big)\! \leq \!  w_{\tau}  \big( X^{\mathtt{b}, \bw}_{\Lambda^{\mathtt{b}, \bw}}\,  , \eta \big)+ \Delta (X^{\mathtt{b}, \bw}\!  \circ \! \Lambda^{\mathtt{b}, \bw})(\tau) +\mathtt{osc} \big(X^{\mathtt{b}, \bw}_{\Lambda^{\mathtt{b}, \bw}}, [\tau, z_1) \big)$. Since $X^{\mathtt{b}, \bw}  \circ  \Lambda^{\mathtt{b}, \bw}$ is constant on $[\tau, \infty)$, we get $w_{z_1} \big( X^{\mathtt{b}, \bw}_{\Lambda^{\mathtt{b}, \bw}}\,  , \eta \big)\! \leq \!  w_{\tau}  \big( X^{\mathtt{b}, \bw}_{\Lambda^{\mathtt{b}, \bw}}\,  , \eta \big) + w_{z_0} (X^{\mathtt{b}, \bw}\! , 2\eta )$, which implies (\ref{Xlambcontt}) thanks to (\ref{gabardine}) and since $w_{z_0} (X^{\mathtt{b}, \bw}\! , \eta )\! \leq \! w_{z_0} (X^{\mathtt{b}, \bw}\! , 2\eta )$. 

\smallskip

\item[$\circ$] Now assume that $\Delta X^{\mathtt{b}, \bw} (T^*_{\bw}) \! > \!  w_{z_0} (X^{\mathtt{b}, \bw}\! , 2\eta )$. Then there exists a subdivision $t_0\! = \! 0 \! < \! t_1 \! < \! \ldots \! < \! t_r\! = \! z_0$ such that $\min_{1\leq i\leq r-1} (t_i-t_{i-1}) \! \geq \! 2\eta$ and such that 
$$\max_{1\leq i\leq r} \mathtt{osc} (X^{\mathtt{b}, \bw}, [t_{i-1}, t_i)) \! < \! (2 w_{z_0} (X^{\mathtt{b}, \bw}\! , 2\eta )) \wedge \Delta X^{\mathtt{b}, \bw} (T^*_{\bw}),$$ 
which, combined with the assumption  $T^*_{\bw}\! >\! 2\eta$, implies that there exists $i\ino \{ 1, \ldots , r\! -\! 1\}$ such that $t_i\! = \! T^*_{\bw}$. Thus, $ \mathtt{osc} (X^{\mathtt{b}, \bw}, [t_{i-1}, T^*_{\bw})) \! <\! 2 w_{z_0} (X^{\mathtt{b}, \bw}\! , 2\eta )$. By (\ref{gloubinett}) applied to $a\! = \! t_{i-1}$ and all $b\! < \! T_{\bw}^*$, we get $\mathtt{osc} (X^{\mathtt{b}, \bw}_{\Lambda^{\mathtt{b}, \bw}}, [\theta^{\mathtt{b}, \bw}_{t_{i-1}} ,\tau) )\! <\! 2 w_{z_0} (X^{\mathtt{b}, \bw}\! , 2\eta )$. Recall that 
$\tau\! -\!  \theta^{\mathtt{b}, \bw}_{t_{i-1}}\! \geq \! T^*_{\bw} \! -\! t_{i-1} \! \geq \! 2\eta$. Consequently, there is $z_1^\prime \ino (\theta^{\mathtt{b}, \bw}_{t_{i-1}} ,\tau \! -\! \eta)$ such that
%\ts{such that $\Delta (X^{\mathtt{b}, \bw}\!  \circ \! \Lambda^{\mathtt{b}, \bw})(z_1^\prime)\! = \! 0$. Thus, we proved the existence of 
%$z^\prime_1$ such that }
\begin{equation}
\label{gluontrou}
  \Delta (X^{\mathtt{b}, \bw}\!  \circ \! \Lambda^{\mathtt{b}, \bw})(z_1^\prime)\! = \! 0 \quad \textrm{and} \quad  \mathtt{osc} \big( X^{\mathtt{b}, \bw}_{\Lambda^{\mathtt{b}, \bw}}, [z_1^\prime ,\tau) \big) \! <\! 2 w_{z_0} (X^{\mathtt{b}, \bw}\! , 2\eta ) \; .
 \end{equation}
Let $s_0\! = 0  \! <\! s_1 \! < \! \ldots \!  < s_r\! =\! z^\prime_1 $ be such that $\min_{1\leq i \leq r-1} (s_i \! -\! s_{i-1})\!  \geq \! \eta$. We define the subdivision $(s^\prime_i)_{0\leq i\leq r+1}$ of $[0, z_1]$ by setting 
$s^\prime_i\! = \! s_i$ for all $i\ino \{ 0, \ldots , r\! -\! 1\}$ and $s_r^\prime\! = \! \tau$, $s^\prime_{r+1}\! = \! z_1$. Clearly, $\min_{1\leq i \leq r} (s^\prime_i \! -\! s^\prime_{i-1}) \! \geq \! \eta$ since $\tau\! -\! z^\prime_1 \! >\! \eta$. Note that $\mathtt{osc} 
(  X^{\mathtt{b}, \bw}_{\Lambda^{\mathtt{b}, \bw}} , [\tau , z_1)  )\! = \! 0$. 
%\margmm{rewritten}
On the other hand, since $ \Delta (X^{\mathtt{b}, \bw}\!  \circ \! \Lambda^{\mathtt{b}, \bw})(z_1^\prime)\! = \! 0$, (\ref{oscbouint}) and (\ref{gluontrou}) imply that 
\begin{eqnarray*}
\mathtt{osc} \big(  X^{\mathtt{b}, \bw}_{\Lambda^{\mathtt{b}, \bw}} , [s^\prime_{r-1} , \tau) \big) 
& \leq&  \mathtt{osc} \big(  X^{\mathtt{b}, \bw}_{\Lambda^{\mathtt{b}, \bw}} , [s_{r-1} , z^\prime_1)   \big)+ \mathtt{osc} \big(  X^{\mathtt{b}, \bw}_{\Lambda^{\mathtt{b} , \bw}}, [z^\prime_1, \tau) \big) \\ 
&< & \mathtt{osc} \big(  X^{\mathtt{b}, \bw}_{\Lambda^{\mathtt{b} , \bw}}, [s_{r-1} , z^\prime_1)   \big) + 
2 w_{z_0} (X^{\mathtt{b}, \bw}\! , 2\eta )\; . 
\end{eqnarray*}
Putting all these together, we then obtain that  
$$w_{z_1} (  X^{\mathtt{b}, \bw}_{\Lambda^{\mathtt{b} , \bw}}, \eta) \! \leq \! \max_{1\leq i \leq r+1} \mathtt{osc} \big(  X^{\mathtt{b}, \bw}_{\Lambda^{\mathtt{b}, \bw}} , [s^\prime_{i-1} , s^\prime_i)   \big) \leq  
\max_{1\leq i \leq r} \mathtt{osc} \big(  X^{\mathtt{b}, \bw}_{\Lambda^{\mathtt{b}, \bw}} , [s_{i-1} , s_i)   \big) + 2 w_{z_0} (X^{\mathtt{b}, \bw}\! , 2\eta ).$$
Since $(s_i)$ is arbitrary, we get $ w_{z_1} (  X^{\mathtt{b}, \bw}_{\Lambda^{\mathtt{b} , \bw}}, \eta) \! \leq \! w_{z^\prime_1} (  X^{\mathtt{b}, \bw}_{\Lambda^{\mathtt{b} , \bw}}, \eta)  + 2 w_{z_0} (X^{\mathtt{b}, \bw}\! , 2\eta ) $ and we obtain (\ref{Xlambcontt}) thanks to (\ref{gabardine}) and the fact that  
$w_{z_0} (X^{\mathtt{b}, \bw}\! , \eta )\! \leq \! w_{z_0} (X^{\mathtt{b}, \bw}\! , 2\eta )$. This 
completes the proof of $(ii)$. 
\end{compactenum}

\smallskip

The proof of $(iii)$ is similar and simpler. Recall from (\ref{thetarw}) that 
$ \theta^{\mathtt{r} , \bw}_{t}\! = \! \inf \{ s\ino [0, \infty) \! : \!  \Lambda^{\! \mathtt{r} , \bw}_{s}\! >\! t \}$. Let $b\! >\! a$. Recall from Lemma \ref{gumodis} that 
$\theta^{\mathtt{r} , \bw}_{b-}\! -\!  \theta^{\mathtt{r} , \bw}_{a} \! \geq \! b\! -\! a $ and observe that 
$\{\Lambda^{\mathtt{r}, \bw}_t ; t\ino [\theta^{\mathtt{r}, \bw}_{a} , \theta^{\mathtt{r}, \bw}_{b-} )  \}\! = \! [a, b)$.
Thus, $\mathtt{osc} \big( X^{\mathtt{r}, \bw}_{\Lambda^{\mathtt{r}, \bw}} ,[\theta^{\mathtt{r}, \bw}_{a} , \theta^{\mathtt{r}, \bw}_{b-}) \, \big)  \! = \!      \mathtt{osc} \big( X^{\mathtt{r}, \bw} , [a, b) \, \big) $. 
Suppose next that $\Delta  \theta^{\mathtt{r}, \bw}_{b} \! \! >\! 0$. Then, $\{\Lambda^{\mathtt{r}, \bw}_t ; t\ino [\theta^{\mathtt{r}, \bw}_{a} , 
\theta^{\mathtt{r}, \bw}_{b} )  \}\! = \! [a, b]$ but since $ |\Delta  X^{\mathtt{r}, \bw}_b|\! = \! 0 $ by Lemma \ref{gumodis}, we get 
$\mathtt{osc} \big( X^{\mathtt{r}, \bw}_{\Lambda^{\mathtt{r}, \bw}} ,[\theta^{\mathtt{r}, \bw}_{a} , \theta^{\mathtt{r}, \bw}_{b}) \, \big)  \! = \!      \mathtt{osc} \big( X^{\mathtt{r}, \bw} , [a, b] \, \big) =    
\mathtt{osc} \big( X^{\mathtt{r}, \bw} , [a, b) \, \big) $. Thus, we have proved for all $b\! >\! a$, 
$$\mathtt{osc} \big( X^{\mathtt{r}, \bw}_{\Lambda^{\mathtt{r}, \bw}} ,[\theta^{\mathtt{r}, \bw}_{a} , \theta^{\mathtt{r}, \bw}_{b-}) \, \big)  \! = \!  
\mathtt{osc} \big( X^{\mathtt{r}, \bw}_{\Lambda^{\mathtt{r}, \bw}} ,[\theta^{\mathtt{r}, \bw}_{a} , \theta^{\mathtt{r}, \bw}_{b}) \, \big)  \! = \!      \mathtt{osc} \big( X^{\mathtt{r}, \bw} , [a, b) \, \big) \; .$$
To complete the proof of $(iii)$ we then argue as in the proof of (\ref{Xlambcont}). \cqfd 

\section{Previous results on the continuous setting.}
\label{prevcontsec}

\subsection{Preliminary results on spectrally L\'evy processes and their height process}
\label{infimsec}
In this section we briefly recall the known results that we need on the analogues $(X, H)$ in the continuous setting of the processes $(X^\bw, H^\bw)$ 
coding the Markovian queue. More precisley, we fix 
$\alpha\ino \bbR$, $\beta \ino [0, \infty) $, $\kappa \ino (0, \infty) $, $\mathbf{c}\! = \! (c_j)_{j\geq 1} \ino \elldo_3$ and we set  
\begin{equation}
\label{repsidefi}
\forall \lambda \ino [0, \infty) , \quad \psi(\lambda)  \! = \!   
 \alpha \lambda +\frac{_{_1}}{^{^2}} \beta \lambda^2 +  \!  \sum_{j\geq 1}   \kappa c_j   \big( e^{-\lambda c_j}\! -\! 1\! + \! \lambda c_j \big). 
 \end{equation}
Let $(X_t)_{t\in [0, \infty)}$ be a spectrally positive L\'evy process with initial state $X_0 \! = \! 0$ and with Laplace exponent $\psi$: namely, 
$ \log \bE [ \exp ( - \lambda X_t )] \! = \!   t\psi (\lambda) $, for all $t, \lambda \ino [0, \infty)$. The L\'evy measure of $X$ is 
$\pi \! = \! \sum_{j\geq 1} \kappa c_j \delta_{c_j}$, $\beta$ is its Brownian parameter and $\alpha$ is its drift. 
%Since $\alpha \! \geq \! 0$, $X$ does not drift to $\infty$: namely a.s.~$\liminf_{t\rightarrow \infty} X_t \! = \! - \infty$. 

First, note that these cases 
include the discrete processes $X^{\bw}$ by taking $\mathbf{c}\! = \! \bw\ino \elldo_{^{\! f}}$, 
$\kappa \! = \! 1/\sigma_1 (\bw)$, $\beta\! = \! 0$ and $\alpha\! = \!  1\! -\! 
\frac{_{\sigma_2 (\bw)}}{^{\sigma_1 (\bw)}}$. However, in the sequel we shall focus on the cases where $X$ has infinite variation sample paths, which is equivalent 
to the following conditions: $\beta\! > \! 0$ or 
 $\sigma_2 (\mathbf{c})\! = \! \int_{(0, \infty)}\!  r\, \pi (dr) \! = \! \infty$, by standard results on L\'evy processes. 
%Next, recall that $X$ has infinite variation sample paths 
%if and only if  either $\beta\! > \! 0$ or 
% $\sigma_2 (\mathbf{c})\! = \! \int_{(0, \infty)}\!  r\, \pi (dr) \! = \! \infty$, which is implied by the following condition: 
% \begin{equation}
%\label{recontH}
%\alpha \geq 0 \quad \textrm{and} \quad
% \int^\infty \frac{d\lambda}{\psi (\lambda)} <\infty \; , 
%\end{equation}
%that is assumed in various places in the paper.   
If $\alpha \! \geq 0$, then a.s.~$\liminf_{t\rightarrow \infty} X_t \! = \! -\infty$ and if $\alpha \! < \! 0$, then a.s.~$\lim_{t\rightarrow \infty} X_t \! = \! \infty$. By analogy with the discrete setting, we refer to the following cases as  
\begin{equation}
\label{mmbjmboo}
\textrm{the supercritical cases if} \; \alpha \! < \! 0, \; \,   \textrm{the critical cases if} \;  \alpha \! = \! 0,  \; \,    \textrm{the subcritical cases if} \; 
 \alpha \! >\! 0. 
\end{equation}
%
%
%If $\alpha \! >\! 0$ (resp.~$\alpha \! < \! 0$), then a.s.~$\lim_{t\rightarrow \infty} X_t \! = \! -\infty$ (resp.~$\infty$) and we refer to these cases as to the \textit{subcritical} cases (resp.~as to the \textit{supercritical} cases).
%If $\alpha\! = \! 0$, then a.s.~$\liminf_{t\rightarrow \infty} X_t\! = \! -\infty$ and $\limsup_{t\rightarrow \infty} X_t \! = \! \infty$; 
%we refer to these cases as to the \textit{critical} ones.    
We next introduce the following  process. 
\begin{equation}
\label{gammadef}
 \forall x \ino [0, \infty), \quad \subo_x = \inf \{ s \ino [0, \infty):    X_s \! < \! -x \} \; .
\end{equation}
with the convention: $\inf \emptyset \! = \! \infty$. For all $t\ino [0, \infty)$, we also set 
\begin{equation}
\label{infimdef}
\forall t\ino [0, \infty), \quad I_t \! = \! \inf_{s \in [0, t]} X_s \quad \textrm{and} \quad 
 I_\infty\! = \! \lim_{t\rightarrow \infty} \! I_t .
 \end{equation}
Note that $I_\infty$ is a.s.~finite in supercritical cases and a.s.~infinite in critical or subcritical cases. 
Observe that 
$\subo_x\! < \! \infty$ if and only if $x\! < \! \! -I_\infty $.  Standard results on spectrally positive L\'evy processes (see e.g.~Bertoin's book \cite{Be96} Ch.~VII) assert that 
$(\subo_x)_{x\in [0, \infty)}$ is a subordinator (a killed subordinator in supercritical cases) whose Laplace exponent is given for all 
$\lambda \ino [0, \infty)$ by: 
\begin{equation}
\label{gamexpobs}
\bE \big[ e^{-\lambda \subo_x } \big]= e^{-x\psi^{-1}(\lambda) } \quad  \textrm{where} \quad  \psi^{-1} (\lambda)\! = \! \inf \big\{ u \ino [0, \infty ) : \psi (u) \! >\! \lambda \big\}.  
\end{equation}
We set $\varrho \! = \!  \psi^{-1}(0)$ that is the largest root of $\psi$. Note that $\varrho \! > \! 0 $ if and only if $\alpha \! < \! 0$. 
%
%As recalled in (\ref{supsub}) in Remark \ref{supercr}, supercritical spectrally positive 
%L\'evy processes are absolutely continuous 
%\margmm{rewritten}
%(\mm{up to finite times } \ms{within 
%finite horizon}) with respect to subcritical ones. More precisely, 
%we recall from  Bertoin's book \cite{Be96} Chapter VII, the following.  
%
%: 
%for all $t\ino [0, \infty)$ and for all nonnegative measurable functional $F \! : \! \bD ([0, \infty), \bbR) \! \rightarrow \! \bbR$, 
%\begin{equation}
%\label{supsubis}
%\bE \big[ F( X_{\cdot \wedge t} ) ] \! =\! \bE \big[ \exp (\varrho \overline{X}_t ) \,  F( \overline{X}_{\cdot \wedge t} )\big]  , 
%\end{equation}
%where $\overline{X}$ stands for a \textit{subcritical} spectrally L\'evy process with Laplace exponent $\psi (\varrho + \cdot)$ (see e.g.~in Bertoin's book \cite{Be96} Chapter VII for more details). We shall use this result under the following form. 
The following elementary lemma gather basic properties of $X$ that are used further in the proofs. 
\begin{lem}
\label{minipoub} Let $X$ be a spectrally positive L\'evy process with Laplace exponent $\psi$ given by (\ref{repsidefi}) and with initial value $X_0\! = \! 0$. 
%\margmm{see footnote}
We assume that there is $\lambda \ino (0, \infty)$ such that $\psi (\lambda ) \! >\! 0$. 
Let $\psi^{-1}$ be defined by (\ref{gamexpobs}) and recall that $\varrho \! = \! \psi^{-1} (0)$ that is the largest root of $\psi$. Let $\overline{X}$ stand for a spectrally L\'evy process with Laplace exponent $\psi (\varrho + \cdot)$ and with initial value $0$. 
Then, the following holds true. 
\begin{compactenum}

\smallskip

\item[$(i)$] A.s.~$\liminf_{t\rightarrow \infty} \overline{X}_t \! = \! -\infty$. Moreover, 
for all $t\ino [0, \infty)$ and for all nonnegative measurable functional $F \! : \! \bD ([0, \infty), \bbR) \! \rightarrow \! \bbR$, 
\begin{equation}
\label{supsubis}
\bE \big[ F( X_{\cdot \wedge t} ) ] \! =\! \bE \big[ \exp (\varrho \overline{X}_t ) \,  F( \overline{X}_{\cdot \wedge t} )\big] . 
\end{equation}
\item[$(ii)$] The c\`adl\`ag process 
$x\ino [0, \infty)\! \mapsto \! \gamma_x (\overline{X})\! := \! \inf \{ s \ino [0, \infty)\! : \!    \overline{X}_s \! < \! -x \}$ is a (conservative) subordinator with Laplace exponent $\psi^{-1} (\cdot) \! -\! \varrho $. 
\item[$(iii)$] For all $x\ino [0, \infty)$, we set 
\begin{equation}
\label{schnitzel2}
\overline{\gamma}_{x} = \gamma_x \quad \textrm{if $x\! <\! -I_\infty$} \quad \textrm{and} \quad \overline{\gamma}_{x} = \gamma  ((-I_\infty) -) \quad  \textrm{if $x\! \geq \! -I_\infty$.} 
\end{equation}
Let $\cE$ be an exponentially distributed r.v.~with parameter $\varrho $ that is independent from $\overline{X}$ (with the convention that a.s.~$\cE\! = \! \infty$ if $\varrho \! = \! 0$). Then,  
\begin{equation}
\label{schnitzel1}
\big( (\overline{\gamma}_{x})_{x\in [0, \infty)} \, , -I_\infty\big) \overset{\textrm{(law)}}{=}  \big( (
\gamma_{x\wedge \cE}  (\overline{X}))_{x\in [0, \infty)}\,  , \cE\big) \; .
\end{equation}
%where we have set 
%\begin{equation}
%\label{schnitzel2}
%\overline{\gamma}_{x} = \gamma_x \quad \textrm{if $x\! <\! -I_\infty$} \quad \textrm{and} \quad \overline{\gamma}_{x} = \gamma  ((-I_\infty) -) \quad  \textrm{if $x\! \geq \! -I_\infty$.} 
%\end{equation}
\item[$(iv)$] Let $(\ccG_{ x})_{x\in [0, \infty)}$ be a right-continuous filtration such that for all $x,y\ino [0, \infty)$, $\gamma_x$ is $\ccG_{ x}$-measurable and $\gamma_{x+y}\! -\! \gamma_x$ is independent of $\ccG_{ x}$. Let $T$ be a $(\ccG_{x})$-stopping time. Then, for all $x, \epp \ino (0, \infty)$, 
\begin{equation}
\label{schnitzel3}
\bP \big( \overline{\gamma}_{x+T} - \overline{\gamma}_{T} > \epp \, ; \,  T \! < \! \infty \big)
\leq \bP \big( \gamma_x >\epp\big) \leq \frac{1-e^{-x\psi^{-1} (1/\epp)}}{1-e^{-1}} \;  .
\end{equation}
\end{compactenum}
\end{lem}
\noi
\textbf{Proof.} The assertions in $(i)$, $(ii)$ and $(iii)$ are (easy consequences of) standard results that can be found e.g.~in Bertoin's book \cite{Be96} Chapter VII. We only need to prove $(iv)$. To that end, 
first note that the second inequality in (\ref{schnitzel3}) is a consequence of a standard inequality combined with (\ref{gamexpobs}). Then, note that in the critical or subcritical cases where $\overline{\gamma}\! = \! \gamma$, the first inequality in (\ref{schnitzel3}) is a straightforward consequence of the fact that $\gamma$ is a subordinator. 
Therefore, we now assume that $\varrho\! >\! 0$. 
Let $\gamma^* $ be a copy of $\gamma$ that is independent of $\ccG_{ \infty}$. 
Then, we set $\gamma^\prime\! = \! \gamma_{\, \cdot + T}  - \gamma_T$ if $T\! <\! \infty$ and $\gamma_T \! < \! \infty$,  
and $\gamma^\prime\! = \! \gamma^*$ otherwise. Then, $\gamma^\prime$ is independent of $\ccG_T$ and it is distributed as $\gamma$. We next set $\cE^\prime\! = \! \sup \{ x \ino (0, \infty) \! : \! \gamma^\prime_x \! < \! \infty \}$; we also define $\overline{\gamma}^\prime$ by setting 
$\overline{\gamma}^\prime_x \! = \! \gamma^\prime_x$ if $x\! < \! \cE^\prime$ and $\overline{\gamma}^\prime_x\! = \! \gamma^\prime (\cE^\prime-)$ if $x\! \geq  \! \cE^\prime$. Thus, 
$$\bP   (\overline{\gamma}_{x+T}\! -\! \overline{\gamma}_{T} > \epp \, ; \,  T \! < \! \infty )\! = \! \bP ( \overline{\gamma}^\prime_{x} > \epp  ;  \gamma_T \! < \! \infty;  T \! < \! \infty)\! = \! \bP ( \overline{\gamma}^\prime_{x} > \epp ) \bP   (\gamma_T \! < \! \infty;  T \! < \! \infty) \; . $$ 
Then observe that $\bP ( \overline{\gamma}^\prime_{x} \! >\!  \epp ) \! \leq \! \bP (\gamma^\prime_x\! >\! \epp) \! = \! \bP (\gamma_x \! <\! \epp)$, which completes the proof of (\ref{schnitzel3}). \cqfd 
%
%A simple extension of (\ref{supsubis}) implies that for all $x\ino (0, \infty)$, 
%$X_{\cdot \wedge \gamma_x} $ under $\bP (\,  \cdot \, | \gamma_x \! < \! \infty) $ is distributed as 
%$\overline{X}_{\gamma_x (\overline{X}) \cdot \wedge}$, where $\gamma_x (\overline{X})\! = \! \inf \{ s \ino [0, \infty):    \overline{X}_s \! < \! -x \}$. Let $\overline{X}$ stands for a \textit{subcritical} 
%spectrally L\'evy process with Laplace exponent $\psi (\varrho + \cdot)$. 
%
%
%
%  

%we assume that 
%\begin{equation}
%\label{recontH}
%%\alpha \geq 0 \quad \textrm{and} \quad
% \int^\infty \frac{d\lambda}{\psi (\lambda)} <\infty \; .
%\end{equation}

\bigskip

\noi
\textbf{Height process of $X$.} We next define the analogue of $H^\bw$. To that end, we assume that the function $\psi$ (as defined in (\ref{repsidefi})) 
satisfies 
\begin{equation}
\label{contH}
 \int^\infty \!\! \frac{d\lambda}{\psi (\lambda)}\,  <\infty . 
\end{equation}
In particular, note that \eqref{contH} implies that either $\beta \! >\! 0$ or $\sigma_2 (\mathbf{c})\! = \! \infty$: namely, (\ref{contH}) entails that $X$ has infinite variation sample paths. 
Le Gall \& Le Jan \cite{LGLJ98} (see also Le Gall \& D.~\cite{DuLG02}) prove that there exists a \textit{continuous} process $H\! = \! (H_t)_{t\in [0, \infty)}$ such that the following limit holds true for all $t\ino [0, \infty)$ in probability :
\begin{equation}
\label{approHdef}
H_t = \lim_{\varepsilon \rightarrow 0} \frac{1}{\varepsilon} \!  \int_0^{t} \!  \un_{\{  X_s - \inf_{r\in [s, t]} X_r  \leq \varepsilon\}} \, ds \; . 
\end{equation}
Note that (\ref{approHdef}) is a local time version of (\ref{XJHdef}). We refer to $H$ as to 
the \textit{height process of $X$}. 
\begin{rem}
\label{supercr} Let us mention that in Le Gall \& Le Jan \cite{LGLJ98} and 
Le Gall \& D.~\cite{DuLG02}, the height process $H$ is introduced only for critical and subcritical spectrally positive processes. However, it easily extends to supercritical cases thanks to (\ref{supsubis}). \cq 
%fact that is proved e.g.~in Bertoin's book \cite{Be96} Ch.~VII:  denote by $\bD ([0, \infty), \bbR)$ the space of c\`adl\`ag function equipped with Skorokhod's topology; then for all $t\ino [0, \infty)$ and for all nonnegative measurable functional $F \! : \! \bD ([0, \infty), \bbR) \! \rightarrow \! \bbR$, 
%\begin{equation}
%\label{supsub}
%\bE \big[ F( X_{\cdot \wedge t} ) ] \! =\! \bE \big[ \exp (\varrho \overline{X}_t ) \,  F( \overline{X}_{\cdot \wedge t} )\big]  , 
%\end{equation}
%where $\overline{X}$ stands for a \textit{subcritical} spectrally L\'evy process with Laplace exponent $\psi (\varrho + \cdot)$. \cq 
\end{rem}
We next recall here that the excursions of $X$ above its infimum process $I$ are the same as the excursions of $H$ above $0$. More specifically, $X-I$ and $H$ have the same set of zeros: 
\begin{equation}
\label{excusame}
\ccZ := \{ t\! \in \! \bbR_+ : H_t \! =\! 0 \}= \{ t \! \in \! \bbR_+ : X_t \! =\! I_t  \} 
\end{equation}
(see Le Gall \& D.~\cite{DuLG02} Chapter 1). 
We also recall that since $-I$ is a local time for $X\! -\! I$ at $0$, the topological support of the Stieltjes measure $d(-I)$ is $\ccZ$. Namely, 
\begin{equation}  
\label{zero}
\textrm{$\bP$-a.s.~for all $s, t\ino [0, \infty)$ such that $s\! < \! t$,} \quad \Big( (s, t) \cap \ccZ \neq \emptyset \Big) \Longleftrightarrow \Big( I_s \! > \! I_t \Big) 
\end{equation} 

We shall also recall here the following result: 
\begin{equation}
\label{Rayxtinct}
\forall x, a\ino (0, \infty), \quad \bP \Big( \! \! \! \!  \!\!  \sup_{\quad t\in [0, \subo_x ]} \! \! \! \! \! \! H_t\leq a \Big) = e^{-x v(a)} \quad \textrm{where} \quad \int_{v(a)}^\infty \! 
\frac{d\lambda}{\psi (\lambda)}= a \; .
\end{equation}
Here, $\gamma_x$ is given by (\ref{gammadef}) and we see that the integral equation completely determines the function $v \! :\!  (0, \infty) \! \rightarrow \! (\varrho, \infty)$ that is bijective,  decreasing and $C^\infty$. In the critical and subcritical cases, this result is a consequence from the excursion theory for $H$ and from Corollary 1.4.2 in 
 Le Gall \& D.~\cite{DuLG02}, p.~41. This result remains true in the supercritical cases thanks to (\ref{supsubis}): we leave the details to the readers.

\subsection{The red and blue processes in the continuous setting.} 
\label{hautcontset}

%\tt{THIS SECTION HAS BEEN CREATED/ RE-ORGANIZED FROM THE BIG PAPER.}
In this section we recall from \cite{BDW1} the definition of the analogues in the continuous setting of the processes $X^{\mathtt{b}, \bw}, X^{\mathtt{r}, \bw}, Y^{ \bw}, A^{ \bw}, \theta^{\mathtt{b}, \bw}$, etc. Let us start with some notation and some convention.  
%The various stochastic processes appearing in the previous construction have their natural counterparts in the continuous setting. We show in \cite{BDW1} that the analogous construction in the continuous setting yields a family of measured metric spaces, called continuous multiplicative graphs. We will see in Section \ref{sec: limit} that these continuous graphs appear in the scaling limits of  the $\bw$-multiplicative graphs. In this section, we focus on the construction of these continuous graphs. 
%\medskip

Let  $(\ccF_t)_{t\in [0, \infty)}$ be a filtration on $(\Omega, \ccF)$ that is specified further. 
A process $(Z_t)_{ t\in [0, \infty)}$ is said to be a 
$(\ccF_t)$-L\'evy process with initial value $0$ if 
a.s.~$Z$ is c\`adl\`ag, $Z_0 \! = \! 0$ and if for all a.s.~finite $(\ccF_t)$-stopping time $T$, the process 
$Z_{T+ \, \cdot}\! -\! Z_{T}$ is independent of 
$\ccF_{T}$ and has the same law as $Z$. 

Let $(M_j(\cdot))_{j\geq 1}$ 
%j\! \geq \! 1$, 
%\margmm{I think here you do not need $M_j$ to be centered.}
be a sequence of c\`adl\`ag  $(\ccF_t)$-martingales that are $L^2$-summable and orthogonal: namely, 
for all $t\ino [0, \infty)$, $\sum_{j\geq 1} \bE \big[ M_j(t)^2\big] \! < \! \infty$ and $\bE [M_j (t)M_k(t)]\! = \! 0$ if $k\! >\! j$. Then $\sum_{j\geq 1}^{_{\perp}} M_j $ stands for the (unique up to indistinguishability) c\`adl\`ag $(\ccF_t)$-martingale $M (\cdot)$ such that  for all $j\! \geq \! 1$ and all $t\ino [0, \infty)$, 
$\bE \big[ \sup_{s\in [0, t]} \big| M (s)  \! -\! \sum_{1\leq k\leq j} M_k (s) \big|^2 \big] \! \leq \! 4 \sum_{l>j} \bE [M_l (t)^2] $, by Doob's inequality.  
%this result being an elementary consequence of Doob's $L^2$ inequality for martingales. 
%\ts{Sometimes, we simply write $\sum^{_\perp}_{j\geq 1} M_j (t)$ instead of $M(t)$. }

\medskip

\noi
\textbf{Blue processes.} 
\label{blupro} We fix the parameters 
$\alpha\ino \bbR$, $\beta \ino [0, \infty) $, $\kappa \ino (0, \infty) $, $\mathbf{c}\! = \! (c_j)_{j\geq 1} \ino \elldo_3$. 
%\begin{equation}
%\label{parconing}
%\alpha \ino \bbR , \quad \beta \ino [0, \infty) , \quad \kappa \ino (0, \infty) , \quad \mathbf{c}\! = \! (c_j)_{j\geq 1} \ino \elldo_3 \; .
%\end{equation}
%These quantites are the \textit{parameters of the continuous multiplicative graph}: $\mathbf{c}$ plays the same role as $\bw$ in the discrete setting, $\alpha$ is a drift coefficient similar to $\alpha_\bw$, $\beta$ is a Brownian coefficient 
%and the interpretation of $\kappa$ is explained later. 
Let $(B_t)_{t \in [0, \infty)}$, $(N_j (t))_{t\in [0, \infty)}$, $j\! \geq \! 1$ be processes that satisfy the following. 
\begin{compactenum}

\smallskip

\item[]  $(b_1)$ \textit{$B$ is a $(\ccF_t)$-real valued standard Brownian motion}.

\smallskip

\item[] $(b_2)$ \textit{For all $j\! \geq \! 1$, $N_j $ is a $(\ccF_t)$-homogeneous Poisson process with jump-rate $\kappa c_j$}. 

\smallskip

\item[] $(b_3)$ \textit{The processes $B$, $N_j$, $j\! \geq \! 1$, are independent}.

\smallskip 
\end{compactenum}

\noi
The \textit{blue} L\'evy process is then defined by 
\begin{equation}
\label{Xblue}
\forall t\ino [0, \infty) , \quad X^{\mathtt{b}}_t = -\alpha t + \sqrt{\beta} B_t + \sum_{j\geq 1}  \!\!\! \,^{\perp}\,  c_j \big( N_j (t) \! -\! c_j \kappa t \big) . 
\end{equation}
Clearly $X^{\mathtt{b}}$ is a $(\ccF_t)$-spectrally positive L\'evy process with initial value $0$ 
with Laplace exponent $\psi$ as defined in (\ref{repsidefi}). 
We next introduce the analogues of the processes $A^\bw$ and $Y^\bw\! $ defined 
in (\ref{YwAwSigw}). To that end, note that $\bE [ c_j (N_j (t) \! -\! 1)_+ ]\! = \!   c_j   \big( e^{- c_j \kappa  t}\! -\! 1\! + \!  c_j \kappa t \big) \! \leq \! \frac{_1}{^2} (\kappa t)^2 c^3_j $. So it makes sense to define the following: 
 \begin{equation}
\label{AetYdef}
\forall t\ino [0, \infty) , \quad A_t = \frac{_{_1}}{^{^2}} \kappa \beta t^2 +  \sum_{j\geq 1}  c_j \big( N_j (t) \! -\! 1 \big)_+ \quad \textrm{and} \quad Y_t = X^{\mathtt{b}}_t \! -\! A_t . 
\end{equation}
%Throughout the paper, we will assume that either $\beta>0$ or $\sigma_2 (\mathbf{c})\! = \! \infty$, which is equivalent to assuming the sample paths of $X^{\mathtt{b}}$ to have infinite variation. Then under this assumption, a.s.~the process $A$ is strictly increasing and the process $Y$ has infinite variation sample paths (see Lemma ?? in \cite{BDW1}). 
%\begin{rem}
%\label{Yrepres}
To view $Y$ as in (\ref{plasouilloc}), set $E_j \! = \! \inf \{ t\ino [0, \infty)\! :\!  N_j (t)\! = \! 1 \}$; 
%(thus, $(\kappa c_jE_j)_{j\geq 1}$ are i.i.d.~exponential r.v.~with unit mean). 
note that $ c_j  (N_j (t) \! -\! c_j \kappa t)\! -\!  c_j (N_j(t)\! -\! 1)_+\! = \!  c_j (\un_{\{ E_j \leq t \}} \! -\! c_j \kappa t)$ and 
check that $ c_j (\un_{\{ E_j \leq t \}} \! -\! c_j \kappa t)\! = \! M^\prime_j (t) \! -\!    \kappa c^2_j (t\! -\! E_j)_+$, where $M^\prime_j$ is a centered $(\ccF_t)$-martingale such that 
$\bE [M^\prime_j (t)^2]\! = \! c_j^2 (1\! -\! e^{-c_j \kappa t}) \! \leq \!    \kappa t c_j^3 $. Since 
%\margmm{$E_j$ has rate $\kappa c_j$, not $c_j$.}
$\bE [ \kappa c^2_j (t\! -\! E_j)_+] \leq \kappa t c_j^2 (1\! -\! e^{-\kappa c_j t})\! \leq \!  \kappa^2 t^{2} c_j^3$, it 
makes sense to write for all $t\ino [0, \infty)$: 
\begin{eqnarray}
\label{Ydef}
 Y_t  &=  &  - \alpha t \! -\! \frac{_1}{^2}\kappa  \beta t^2 + \sqrt{\beta} B_t + \sum_{j\geq 1}  \!\!\! \,^{\perp} c_j \big( \un_{\{ E_j \leq t \}} \! -\! \kappa c_j (t\! \wedge E_j)\big)
 - \!\! \sum_{j\geq 1}  \kappa c^2_j (t\! -\! E_j)_+ \nonumber \\
 &  \overset{\textrm{(informal)}}{=} &
  -\alpha t \! -\! \frac{_1}{^2}\kappa \beta t^2 + \sqrt{\beta} B_t +  \sum_{j\geq 1} \! c_j (\un_{\{ E_j \leq t \}} \! -\! c_j \kappa t).
 \end{eqnarray}
Namely the jump-times of $Y$ are the $E_j$ and $\Delta Y_{E_j}\! = \! c_j$. 
%\end{rem}
%Under that (\ref{varinfinie}) holds. Recall from (\ref{AetYdef}) the definiton of $A$ and $Y$.  
%Then, a.s.~the process $A$ is strictly increasing and the process $Y$ has infinite variation sample paths. 

\medskip

\noi
\textbf{Red and bi-coloured processes.}  We next introduce the \textit{red process} $X^{\mathtt{r}}$ that satisfies the following. 
\begin{compactenum}

\smallskip

\item[] $(r_1)$ $X^{\mathtt{r}}$ is a $(\ccF_t)$-spectrally positive L\'evy process starting at $0$ and whose Laplace exponent is $\psi$ as in  (\ref{repsidefi}). 

\smallskip 

\item[] $(r_2)$ $X^{\mathtt{r}}$ is independent of the processes $B$ and $(N_j)_{j\geq 1}$.   

\smallskip 

\end{compactenum}

\noi
%To keep the filtration $(\ccF_t)$ minimal, we may assume that $\ccF_t$ is the completed sigma-field generated by 
%$B_s$, $(N_j (s))_{j\geq 1}$ and $X^{\mathtt{r}}_s$, $s\ino [0, t]$. 
We next introduce the following processes: 
\begin{equation}
\label{thetabdef}
 \forall x, t\ino [0, \infty), \quad \subo^{\mathtt{r}}_x = \inf \{ s \ino [0, \infty):    X^{\mathtt{r}}_s \! < \! -x \} \quad \textrm{and} \quad \theta_t^{\mathtt{b}}= t + \subo^{\mathtt{r}}_{A_t} ,
\end{equation}
 with the convention: $\inf \emptyset \! = \! \infty$. For all $t\ino [0, \infty)$, we set $I^\mathtt{r}_t \! = \! \inf_{s \in [0, t]} X^\mathtt{r}_s $ and 
 $I^{\mathtt{r}}_\infty\! = \! \lim_{t\rightarrow \infty} \! I^{\mathtt{r}}_t$ 
 that is a.s.~finite in supercritical cases and that is a.s.~infinite in critical or subcritical cases. 
 Note that $\subo^{\mathtt{r}}_x\! < \! \infty$ if and only if  $x\! < \! \! -I^{\mathtt{r}}_\infty $.  
%Standard results on spectrally positive L\'evy processes (see e.g.~Bertoin's book \cite{Be96} Ch.~VII.) assert that 
%$(\subo^{\mathtt{r}}_x)_{x\in [0, \infty)}$ is a subordinator (that is defective in supercritical cases) whose Laplace exponent is given for all 
%$\lambda \ino [0, \infty)$ by: 
%\begin{equation}
%\label{gamexpo}
%\bE \big[ e^{-\lambda \subo^{\mathtt{r}}_x } \big]= e^{-x\psi^{-1}(\lambda) } \quad  \textrm{where} \quad  \psi^{-1} (\lambda)\! = \! \inf \big\{ u \ino [0, \infty ) : \psi (u) \! >\! \lambda \big\}.  
%\end{equation}
Recall that $\varrho$ stands for the largest root of $\psi$: in supercritical cases, $\varrho\! >\! 0$ and $-I^{\mathtt{r}}_\infty$ is exponentially distributed with parameter $\varrho$, as recalled in Lemma \ref{minipoub} $(iii)$. We next set: 
\begin{equation}
\label{T*def} 
T^*\! = \! \sup \{ t\ino [0, \infty)\! :  \theta^{\mathtt{b}}_t \! < \infty \}=   \sup \{ t\ino [0, \infty)\! : A_t \! < \! - I^{\mathtt{r}}_\infty \} \; . 
\end{equation}
In critical and subcritical cases, $T^*\! = \! \infty$ and $ \theta^{\mathtt{b}}$ only takes 
finite values. In supercritical cases, a.s.~$T^* \! < \! \infty$ and we check that $ \theta^{\mathtt{b}} (T^*-) \! < \! \infty$. We next introduce the following. 
\begin{equation}
\label{2lambda}
\forall t \ino [0, \infty), \quad \Lambda^{\mathtt{b}}_t =  \inf \{ s \ino [0, \infty):    \theta^{\mathtt{b}}_s \! > \! t \} \quad \textrm{and} \quad \Lambda^{\mathtt{r}}_t \! = \! t-    \Lambda^{\mathtt{b}}_t . 
\end{equation} 
Both processes $\Lambda^{\mathtt{b}}$ and  $\Lambda^{\mathtt{r}}$ are continuous and nondecreasing. 
In critical and subcritical cases, we also get a.s.~$\lim_{t\rightarrow \infty}  \Lambda^{\mathtt{b}}_t \! = \! \infty$ and $\Lambda^{\mathtt{b}} (\theta^{\mathtt{b}}_t)\! = \! t$ for all $t\ino [0, \infty)$. 
However, in supercritical cases, a.s.~$ \Lambda^{\mathtt{b}}_t \! = \! T^*$ for all $t\ino[ \theta^{\mathtt{b}} (T^*-), \infty)$ and a.s.~for all $t\ino [0, T^*)$, 
$\Lambda^{\mathtt{b}} (\theta^{\mathtt{b}}_t )\! = \! t$. 
In the following theorem we recall from \cite{BDW1} the results about the previous processes that we need; in particular, it contains 
the analogue of Proposition \ref{Xwfrombr}. 
\begin{thm}
\label{Xdefthm} Let $(\alpha, \beta, \kappa, \mathbf{c})$ be as in (\ref{parconing}). 
Assume that either $\beta>0$ or $\sigma_2 (\mathbf{c})\! = \! \infty$. We keep the previous definition for $X^{\mathtt{b}}$, $A$, $Y$, $X^{\mathtt{r}}$, $\theta^{\mathtt{b}}$, $T^*$, 
$\Lambda^{\mathtt{b}}$ and $\Lambda^{\mathtt{r}}$. 
\begin{compactenum}

\smallskip

\item[$(i)$]  A.s.~the process $A$ is strictly increasing and the process $Y$ has infinite variation sample paths.

\smallskip

\item[$(ii)$] The process $\Lambda^{\mathtt{r}}$ is continuous, nondecreasing and a.s.~$\lim_{t\rightarrow \infty}  \Lambda^{\mathtt{r}}_t \! = \! \infty$.

\smallskip

\item[$(iii)$] For all $t\ino [0, \infty)$, we set 
\begin{equation}
\label{Xdef}
\forall t \ino [0, \infty) , \quad X_t = X^{\mathtt{b}}_{\Lambda^{\mathtt{b}}_t } + X^{\mathtt{r}}_{\Lambda^{\mathtt{r}}_t }\, ;  
\end{equation} 
$X$, $X^\mathrm{b}$ and  $X^\mathrm{r}$ have the same law: namely, $X$ is a spectrally positive L\'evy process with initial value $0$ and Laplace exponent $\psi$ as in (\ref{repsidefi}). Moreover, 
\begin{equation}
\label{YfromX}
\textrm{a.s.}\; \forall\,t\ino [0, T^*) , \quad Y_t \! = \! X_{\theta^{\mathtt{b}}_t } . 
\end{equation}    
%\item[$(v)$] $\bP$-a.s.~for all $a\ino [0, T^*)$, if $\Delta \theta^{\mathtt{b}}_a \! =\! 0$, then $t\! = \! \theta^\mathtt{b}_a$ is the unique $t \ino [0, \infty)$ such that $\Lambda^{\mathtt{b}}_t \! = \! a$. 
%\smallskip
%\item[$(vi)$]  $\bP$-a.s.~for all $a\ino [0, T^*]$, if $\Delta \theta^{\mathtt{b}}_a \! >\! 0$, then $\Delta X \! (\theta^{\mathtt{b}}_{a-}) \! = \! \Delta A_a $ and 
%$\Delta Y_a\! = \! 0$. Moreover, 
%$$ \forall t\ino \big( \theta^{\mathtt{b}}_{a-}, \theta^{\mathtt{b}}_{a} \big), \quad X_t  \! \geq \! X_{t-} \!  > \! 
%X_{(\theta^{\mathtt{b}}_{a-})-}\!\! =\! Y_a  \quad \textrm{and if $a<T^*\! $, then }\quad  X_{(\theta^{\mathtt{b}}_{a-})-}\!\! = \! X_{\theta^{\mathtt{b}}_a} .$$ 
%\item[$(vii)$] $\bP$-a.s.~if $(\Delta X^\mathtt{r} )(\Lambda^\mathtt{r}_t) \! >\! 0$, then there exists $a\ino [0, T^*]$ such that 
%$\theta^\mathtt{b}_{a-} \! < \! t  \! < \! \theta^\mathtt{b}_{a} $. 
%\smallskip
%\item[$(viii)$] $\bP$-a.s.~for all $b\ino [0, \infty)$ such that $\Delta X^\mathtt{r}_b\! >\! 0$, there is a unique $t\ino [0, \infty)$ such that $\Lambda^\mathtt{r}_t \! = \! b$.  
%\smallskip
%\item[$(ix)$] For all $t\ino [0, \infty)$, set $Q^\mathtt{b}_t \! = \! X^{\mathtt{b}}_{\Lambda^{\mathtt{b}}_t }$ and  $Q^\mathtt{r}_t \! = \! X^{\mathtt{r}}_{\Lambda^{\mathtt{r}}_t }$. Then, a.s.~for all $t\ino [0, \infty)$, 
%$\Delta Q^\mathtt{b}_t \Delta Q^\mathtt{r}_t\! = \! 0$.   
%\smallskip
\end{compactenum}
% 
%The process $\Lambda^{\mathtt{r}}$ is continuous, nondecreasing and a.s.~$\lim_{t\rightarrow \infty}  \Lambda^{\mathtt{r}}_t \! = \! \infty$. Furthermore, if we set 
%\begin{equation}
%\label{Xdef}
%\forall t \ino [0, \infty) , \quad X_t = X^{\mathtt{b}}_{\Lambda^{\mathtt{b}}_t } + X^{\mathtt{r}}_{\Lambda^{\mathtt{r}}_t }\, ,  
%\end{equation} 
%then, $X$ is a spectrally positive L\'evy process with initial value $0$ and Laplace exponent $\psi$ as in (\ref{psidefi}). 
%Namely, $X$, $X^\mathrm{b}$ and  $X^\mathrm{r}$ have the same law. Moreover, 
%\begin{equation}
%\label{YfromX}
%\textrm{a.s.}\; \forall\,t\ino [0, T^*) , \quad Y_t \! = \! X_{\theta^{\mathtt{b}}_t } . 
%\end{equation}    
\end{thm}
\noi
\textbf{Proof.} For $(i)$, see Lemma 2.4 in \cite{BDW1}; for $(ii)$ and $(iii)$, see Theorem 2.5 in \cite{BDW1}. \cqfd 

\medskip

%The set of intervals $\bigcup_{a\in [0, T^*]}  \big[ \theta^{\mathtt{b}}_{a-}, \theta^{\mathtt{b}}_{a} \big)$ is the set of times where red clients are served 
The red and blue processes behave quite similarly as in the discrete setting as in Lemma \ref{gumodis}. More precisely,  we recall from \cite{BDW1} the various properties concerning the red and blue processes that are used in the proof. 
\begin{lem}
\label{trajprop} We keep the assumption of Theorem \ref{Xdefthm}. 
%Let $\alpha\ino \bbR$, $\beta\ino [0, \infty)$,  $\kappa \ino (0, \infty)$ and $\mathbf{c}\ino \elldo_3$ satisfy (\ref{reassume}). 
Then, the following statements hold true. 
\begin{compactenum}
\smallskip
\item[$(i)$] $\bP$-a.s.~for all $a\ino [0, T^*)$, if $\Delta \theta^{\mathtt{b}}_a \! =\! 0$, then $t\! = \! \theta^\mathtt{b}_a$ is the unique $t \ino [0, \infty)$ such that $\Lambda^{\mathtt{b}}_t \! = \! a$. 
\smallskip
\item[$(ii)$]  $\bP$-a.s.~for all $a\ino [0, T^*]$, if $\Delta \theta^{\mathtt{b}}_a \! >\! 0$, then $\Delta X \! (\theta^{\mathtt{b}}_{a-}) \! = \! \Delta A_a $ and 
$\Delta Y_a\! = \! 0$. Moreover, 
$$ \forall t\ino \big( \theta^{\mathtt{b}}_{a-}, \theta^{\mathtt{b}}_{a} \big), \quad X_t  \! \geq \! X_{t-} \!  > \! 
X_{(\theta^{\mathtt{b}}_{a-})-}\!\! =\! Y_a  \quad \textrm{and if $a<T^*\! $, then }\quad  X_{(\theta^{\mathtt{b}}_{a-})-}\!\! = \! X_{\theta^{\mathtt{b}}_a} .$$ 
\item[$(iii)$] $\bP$-a.s.~if $(\Delta X^\mathtt{r} )(\Lambda^\mathtt{r}_t) \! >\! 0$, then there exists $a\ino [0, T^*]$ such that 
$\theta^\mathtt{b}_{a-} \! < \! t  \! < \! \theta^\mathtt{b}_{a} $. 
\smallskip
\item[$(iv)$] $\bP$-a.s.~for all $b\ino [0, \infty)$ such that $\Delta X^\mathtt{r}_b\! >\! 0$, there is a unique $t\ino [0, \infty)$ such that $\Lambda^\mathtt{r}_t \! = \! b$.  
\smallskip
\item[$(v)$] For all $t\ino [0, \infty)$, set $Q^\mathtt{b}_t \! = \! X^{\mathtt{b}}_{\Lambda^{\mathtt{b}}_t }$ and  $Q^\mathtt{r}_t \! = \! X^{\mathtt{r}}_{\Lambda^{\mathtt{r}}_t }$. Then, a.s.~for all $t\ino [0, \infty)$, 
$\Delta Q^\mathtt{b}_t \Delta Q^\mathtt{r}_t\! = \! 0$.   
\end{compactenum}
\end{lem}
\noi
\textbf{Proof.} For $(i)$ and $(ii)$, see Lemma 5.4 in \cite{BDW1}; for $(iii)$, $(iv)$ and $(v)$, 
see Lemma 5.5 in \cite{BDW1}. \cqfd 

\bigskip

\noi
\textbf{The excursions of $Y$ above its infimum.} Let $X$ be derived from $X^{\mathtt{b}}$ and $X^{\mathtt{r}}$ by (\ref{Xdef}) and recall from (\ref{infimdef}) the notation  
$I_t \! = \! \inf_{s\in [0, t]} X_s$, for the infimum process of $X$. Recall from (\ref{zero}) that $-I$ is a local-time for the set of zeros 
$\mathscr{Z}\! = \! \{ t\ino [0, \infty) : X_t \! = \! I_t \}$.  
Let $Y$ be defined by (\ref{AetYdef}) and recall from (\ref{Jdef}) the notation $J_t \! = \! \inf_{s\in [0, t]} Y_s$. The following lemma (that is recalled from \cite{BDW1}) asserts that$-J$ is a local-time for the set $\mathscr{Z}^{\mathtt{b}}\! = \! \{ t\ino [0, \infty) : Y_t \! = \! J_t \}$ (more precisely, it shows that $\mathscr{Z}^{\mathtt{b}}$ is bijectively sent to $\mathscr{Z}$ via $\Lambda^{\mathtt{b}}$). 
\begin{lem}
\label{XYexcu} We keep the assumptions of Theorem \ref{Xdefthm}. Then, the following holds true. 
\begin{compactenum}

\smallskip

\item[$(i)$] A.s.~for all $t\ino [0, \infty)$, $X_t \! \geq \! Y( \Lambda^\mathtt{b}_t)$. Then, a.s.~for all $t_1, t_2 \ino [0, \infty)$ such that $\Lambda^\mathtt{b}_{t_1} \! < \!  \Lambda^\mathtt{b}_{t_2}$, $\inf_{s\in [t_1, t_2]} X_s \! = \!  
\inf_{a\in [\Lambda^\mathtt{b} (t_1),  \Lambda^\mathtt{b} (t_2)]} Y_a$. It implies that a.s.~for all $t\ino [0, \infty)$, $I_t = J (\Lambda^{\! \mathtt{b}}_t)$.

\smallskip

\item[$(ii)$] A.s.~$\big\{ t \ino [0, \infty) : X_t \! >\! I_t \big\} =  \big\{ t \ino [0, \infty) : Y({\Lambda^{\! \mathtt{b}}_t}) \! >\! 
J({\Lambda^{\! \mathtt{b}}_t}) \big\} $.

\smallskip

\item[$(iii)$] A.s.~the set $\ccE\! = \! \big\{ a \ino [0, \infty) : Y_a \! >\! J_a \big\}$ is open. Moreover, 
if $(l,r)$ is a connected component of $\ccE$, then $Y_l\! = \! Y_r\! = \! J_l \! = \! J_r$ and 
for all $a\ino (l, r)$, we get $J_a\! = \! J_l$ and 
$Y_{a-}\! \wedge \! Y_a \! >\! J_l $.

\smallskip

\item[$(iv)$] Set $ \ccZ^\mathtt{b}\! = \! \{ a\ino [0, \infty) \! : \! Y_a\! = \! J_a \}$. Then, $\bP$-a.s.
\begin{equation}
\label{Jtpsloczer}
\textrm{$\forall a, z\ino [0, \infty)$ such that $a\! < \! z$,} \quad \Big( \ccZ^\mathtt{b}\!  \cap (a, z) \neq \emptyset  \Big) \Longleftrightarrow   \Big(  J_z \! < \!  J_a \Big).  
\end{equation} 

%%\smallskip
%
%\item[$(v)$] \mm{A.s.~$\cH$ is continuous and a.s.~$\big\{ a \ino [0, \infty) : Y_a \! >\! J_a \big\} \! =\!  \big\{ a \ino [0, \infty) : \cH_a \! >\! 0 \big\} $. }\margmm{Only in (sub)critical cases. Otherwise adapt.}
\end{compactenum}
\end{lem}
\noi
\textbf{Proof.} See Lemma 5.7 in \cite{BDW1}. \cqfd 

\medskip

We next recall the following result due to Aldous \& Limic \cite{AlLi98} (Proposition 14, p.~20) that is used in our proofs. 
\begin{prop}[Proposition 14 \cite{AlLi98}]  
\label{AldLim1} We keep the assumptions of Theorem \ref{Xdefthm} and the previous notation. Then, the following holds true. 
\begin{compactenum}

\smallskip

\item[$(i)$] For all $a\ino [0, \infty)$, $\bP (Y_a \! = \! J_a) \! = \! 0$. 

\smallskip

\item[$(ii)$] $\bP$-a.s.~the set $\{ a\ino [0, \infty)\! : \! Y_a\!  =\! J_a \}$ contains no isolated points. 

\smallskip

\item[$(iii)$] Set $M_a\! = \! \max \{ r\! -\! l\, ; \; r\! \geq \! l\! \geq \! a \! : \! \textrm{$(l,r)$ 
is an excursion interval of $Y\! -\! J$ above $0$}\}$. Then, $M_a \! \rightarrow \! 0$ in probability as $a\! \rightarrow \! \infty$.  
\end{compactenum}
\end{prop} 
\noi
\textbf{Proof.} The process $(Y_{s/\kappa})_{s\in [0, \infty)}$ is  the process $W^{\kappa^\prime, -\tau, \mathbf{c}}$ 
in \cite{AlLi98}, where $\kappa^\prime \! = \! \beta / \kappa$ and $\tau\! =\! \alpha/ \kappa$ (note that the letter 
$\kappa$ plays another role in \cite{AlLi98}). 
Then $(i)$ (resp.~$(ii)$ and $(iii)$) is Proposition 14 \cite{AlLi98} $(b)$ (resp.~$(d)$ and $(c)$). \cqfd 

\bigskip

\noi
Thanks to Proposition \ref{AldLim1} $(iii)$, the excursion intervals of $Y\! -\! J$ above $0$ can be listed as follows
% (or of $\cH$ above $0$ by Lemma \ref{HYconexcu} $(iii)$).   
\begin{equation}
\label{excYre}
\{ a\ino [0, \infty)\! :  Y_a\! >  \! J_a \}= \bigcup_{k\geq 1} (l_k, r_k) \; .
\end{equation}
where $\zeta_k \! = \! r_k \! -\! l_k $, $k\! \geq \! 1$, is decreasing. 
Then, as a consequence of Theorem 2 in Aldous \& Limic \cite{AlLi98}, p.~4, we recall the following. 
\begin{prop}[Theorem 2 \cite{AlLi98}]  
\label{AldLim2} We keep the assumptions of Theorem \ref{Xdefthm} and the previous notation. 
Then, $(\zeta_k)_{k\geq 1}$, that is the ordered sequence of lengths of the excursions of $Y\! -\! J$ above $0$, is distributed as the 
$(\beta / \kappa, \alpha/ \kappa, \mathbf{c})$-multiplicative coalescent (as defined in \cite{AlLi98}) taken at time $0$. In particular, we get a.s.~$\sum_{k\geq 1} \zeta_k^2 \! < \! \infty$. 
\end{prop}

\noi
\textbf{Height process of $Y$.}  We define the analogue of $\cH^\bw$ in the continuous setting thanks to the following theorem 
that is recalled from various results in \cite{BDW1}. 

\begin{thm}   
\label{cHdefthm} Let $(\alpha, \beta, \kappa, \mathbf{c})$ be as in (\ref{parconing}) and 
assume that (\ref{contH}) holds: namely, $\int^\infty d\lambda / \psi (\lambda) \! <\! \infty$, which implies the assumptions of Theorem \ref{Xdefthm}. Let $X$ be derived from $X^{\mathtt{b}}$ and $X^{\mathtt{r}}$ by (\ref{Xdef}). Let $H$ be the height process associated with $X$ as defined by (\ref{approHdef}) (and by Remark \ref{supercr} in the supercritical cases). Then, there exists a continuous process $(\cH_t)_{t\in [0, \infty)}$ such that for all $t\ino [0, \infty)$, $\cH_t$ is a.s.~equal to a measurable functional of $(Y_{\cdot \wedge t}, A_{\cdot \wedge t}) $ and such that 
\begin{equation}
\label{cHfromX}
\textrm{a.s.}\; \forall\,t\ino [0, T^*) , \quad \cH_t \! = \! H_{\theta^{\mathtt{b}}_t } . 
\end{equation}    
We refer to $\cH$ as the {\rm height process associated with $Y$}. 
\end{thm}
\textbf{Proof.} See Theorem 2.6 in \cite{BDW1}. \cqfd 
%
%\begin{enumerate}
%\item
%(Le Gall \& Le Jan \cite{LGLJ98}, Le Gall \& D.~\cite{DuLG02}) There exists a \textit{continuous} process $H\! = \! (H_t)_{t\in [0, \infty)}$ such that the following limit holds true for all $t\ino [0, \infty)$ in probability :
%\begin{equation}
%\label{approHdefbis}
%H_t = \lim_{\varepsilon \rightarrow 0} \frac{1}{\varepsilon} \!  \int_0^{t} \!  \un_{\{  X_s - \inf_{r\in [s, t]} X_r  \leq \varepsilon\}} \, ds \; . 
%\end{equation}
%We refer to $H$ as the  \textit{height process of $X$}.
%\item
%(Theorem ??, \cite{BDW1}) 
%There exists a continuous process $(\cH_t)_{t\in [0, \infty)}$ such that for all $t\ino [0, \infty)$, $\cH_t$ is a.s.~equal to a measurable functional of $(Y_{\cdot \wedge t}, A_{\cdot \wedge t}) $ and such that 
%\begin{equation}
%\label{cHfromX}
%\textrm{a.s.}\; \forall\,t\ino [0, T^*) , \quad \cH_t \! = \! H_{\theta^{\mathtt{b}}_t } . 
%\end{equation}    
%We refer to $\cH$ as the {\rm height process associated with $Y$}. 
%\end{thm}
%

\medskip

\noi
As for $H$ and $X\! -\! I$, the following lemma (that is recalled from \cite{BDW1}) asserts that the excursion intervals of $\cH$ and $Y\! -\! J$ 
above $0$ are the same.  
\begin{lem}
\label{AHeuer} We keep the same assumptions as in Theorem \ref{cHdefthm}. Then, the following holds true. 

\smallskip

\begin{compactenum}
\item[$(i)$]  Almost surely for all $t\ino [0, \infty)$, $H_t \! \geq \! \cH( \Lambda^\mathtt{b}_t)$ and a.s.~for all $t_1, t_2 \ino [0, \infty)$ such that $\Lambda^\mathtt{b}_{t_1} \! < \!  \Lambda^\mathtt{b}_{t_2}$, 
$\inf_{s\in [t_1, t_2]} H_s \! = \!  
\inf_{a\in [\Lambda^\mathtt{b} (t_1),  \Lambda^\mathtt{b} (t_2)]} \cH_a$.

\smallskip

\item[$(ii)$] Almost surely 
$\big\{ a \ino [0, \infty) : Y_a \! >\! J_a \big\} \! =\!  \big\{ a \ino [0, \infty) : \cH_a \! >\! 0 \big\}$. 
\end{compactenum}

\end{lem}
\noi
\textbf{Proof.} See Lemma 5.11 in \cite{BDW1}. \cqfd

\section{Limit theorems.}
\label{Limisec}

\subsection{Proof of Theorem \ref{HYcvth}.}
\label{proofmain}

In this section, \textit{we admit Proposition \ref{cvmarkpro} and Proposition \ref{HMarkcvprop}} whose proofs are later given in Section \ref{pfcvmarkpro}. 
%Recall that we admit Proposition \ref{cvmarkpro} and Proposition \ref{HMarkcvprop} in this section. 
As a consequence of Proposition \ref{cvmarkpro}, we get the following lemma. 
\begin{lem} 
\label{expcovcv} 
Let $(\alpha, \beta, \kappa, \mathbf{c})$ be as in (\ref{parconing}). Recall from (\ref{repsidefi}) the definition of $\psi$ and assume that (\ref{contH}) holds: namely, $\int^\infty d\lambda / \psi (\lambda) \! <\! \infty$.  
%Let $\alpha \ino \bbR$, $\beta\ino [0, \infty)$,  $\kappa \ino (0, \infty)$ and $\mathbf{c}\ino \elldo_3$ satisfy (\ref{reassume}).  
Let $a_n , b_n \! \in \! (0, \infty)$ and  $\bw_n \! \in \! \elldo_f$, $n\ino \bbN$, 
satisfy (\ref{apriori}) and $\mathbf{(C1)}$--$\mathbf{(C3)}$ as in (\ref{unalphcv}) and in (\ref{sig3cvcj}). 
%Recall from (\ref{reassume}) the definition of $\psi$ and recall from (\ref{Laplwn}) the definition of $\psi_{\bw_n}$, 
For all $\lambda \in [0, \infty)$, we set 
\begin{equation}
\label{Laplwn}
\psi_{\bw_n} (\lambda)  =  \alpha_{\bw_n}  \lambda + \!\! \sum_{1\leq j\leq n}\!  \! \frac{_{w^{(n)}_j} }{^{\sigma_1 (\bw_{n})}} \big( e^{-\lambda w^{(n)}_j}\! -\! 1\! + \! \lambda w^{(n)}_j \big) \quad \textrm{and} \quad 
 \alpha_{\bw_n}  \! := \! 1\! -\! \frac{_{\sigma_2 (\bw_{n})}}{^{\sigma_1 (\bw_{n})}} \; .
\end{equation}
For all $\lambda \ino [0, \infty)$, set $\psi^{-1}(\lambda)\! = \! \inf \{ r\ino [0, \infty): \psi (r) \! >\! \lambda \}$ and $\psi^{-1}_{\bw_n}(\lambda)\! = \! \inf \{ r\ino [0, \infty): \psi_{\bw_n} (r) \! >\! \lambda \}$ that are the inverses of resp.~$\psi$ and $\psi_{\bw_n}$. 
Recall that $\varrho\! =\! \psi^{-1} (0)$ and $\varrho_{\bw_n}\! = \!   \psi^{-1}_{\bw_n} (0)$ 
are the largest roots of the convex functions $\psi$ and $\psi_{\bw_n}$.Then, for all $\lambda\ino [0, \infty)$, 
\begin{equation}
\label{eierschwa}
\lim_{n\rightarrow \infty} b_n \psi_{\bw_n} (\lambda/ a_n)\! = \! \psi (\lambda) \quad \textrm{and} \quad \lim_{n\rightarrow \infty} a_n \psi^{-1}_{\bw_n} (\lambda/ b_n)\! = \! \psi^{-1} (\lambda)
\end{equation}
and thus, $\lim_{n\rightarrow \infty} a_n \varrho_{\bw_n} \! = \! \varrho$. 
\end{lem}
\noi
\textbf{Proof.} The limit of $ b_n \psi_{\bw_n} (\lambda/ a_n)$ is a direct 
consequence of the equivalence $(ii) \Leftrightarrow (iii)$ asserted in Proposition \ref{cvmarkpro}. Since the $\psi_{\bw_n}$ are convex functions, the convergence is uniform in $\lambda$ on all compact subsets of $[0, \infty)$, which easily entails the convergence of the inverses. \cqfd

\medskip

We will use several times the following result from Ethier \& Kurtz \cite{EtKu86}.  
 \begin{lem}[Lemma 3.8.2 \cite{EtKu86}]
 \label{deltime} 
 For all $n\ino \bbN$, let $(\mathbf{s}^n_k)_{k\in \bbN}$ be a nondecreasing $[0, \infty]$-valued sequence of r.v.~such that 
 $\mathbf{s}^n_0\! = \! 0$, a.s.~$\lim_{k\rightarrow \infty} \mathbf{s}^n_k\! = \! \infty$ 
 and $\mathbf{s}^n_k\! < \! \mathbf{s}^n_{k+1}$ for all 
 $k\ino\bbN$ such that $\mathbf{s}^n_k\! < \! \infty$. 
 Fix $z\ino (0, \infty)$ and set $\mathbf{k}_n\! = \! \max \{ k \ino \bbN: \mathbf{s}^n_k \! < \! z \}$. Then 
$$ \lim_{\eta \rightarrow 0+} \, \sup_{n\in \bbN} \, \bP \Big( \!\!\!\!\!\! \min_{\quad 0\leq k\leq \mathbf{k}_n}\!\!\!\!  \mathbf{s}^n_{k+1}\! - \!  \mathbf{s}^n_{k} \! 
 < \eta \Big) \! =\!  0 \quad 
\Longleftrightarrow \quad \lim_{\eta \rightarrow 0+} \, \sup_{n\in \bbN} \, \sup_{k\in \bbN} \, \bP \big(  \mathbf{s}^n_k \! < \! z  \, ; \,  \mathbf{s}^n_{k+1} \! -\!  \mathbf{s}^n_{k} \! < \! \eta \big)\! = \! 0 . $$
%\margmm{Needs double checking \tt{OK}.}
\end{lem}
\noi
\textbf{Proof.} 
%\fmm{Citation from E \& K's book actually requires three digits to identify the result: the first refer to the chapter number, the second one the section number.}
See Lemma 3.8.2 in Ethier \& Kurtz \cite{EtKu86} (p.~134). Note that Lemma 3.8.2 in \cite{EtKu86} only deals with sequences that take finite values but the proof extends immediately to our case.  \cqfd 
%
%
% \begin{lem}[Lemma 8.2 \cite{EtKu86}]\label{deltime} For all $n\ino \bbN$, let $0\! = \! \mathbf{s}^n_0\! <\!  \mathbf{s}^n_1 \! < \! \mathbf{s}^n_2 \! < \! \ldots $ be a sequence of r.v.~such that a.s.~$\lim_{k\rightarrow \infty} \mathbf{s}^n_k\! = \! \infty$. Fix $z\ino (0, \infty)$ and set $\mathbf{k}_n\! = \! \max \{ k \ino \bbN: \mathbf{s}^n_k \! < \! z \}$. Then 
%$$ \lim_{\eta \rightarrow 0+} \, \sup_{n\in \bbN} \, \bP \Big( \!\!\!\!\!\! \min_{\quad 1\leq k\leq \mathbf{k}_n}\!\!\!\!  \mathbf{s}^n_k\! - \!  \mathbf{s}^n_{k-1} \! 
% < \eta \Big) \! =\!  0 \quad 
%\Longleftrightarrow \quad \lim_{\eta \rightarrow 0+} \, \sup_{n\in \bbN} \, \sup_{k\in \bbN} \, \bP \big(  \mathbf{s}^n_k \! < \! z  \, ; \,  \mathbf{s}^n_{k+1} \! -\!  \mathbf{s}^n_{k} \! < \! \eta \big)\! = \! 0 . $$
%\end{lem}
%\noi
%\textbf{Proof:} see Lemma 8.2 in Ethier \& Kurtz \cite{EtKu86} (Chapter 3, p.~134). \cqfd 

\medskip

Let $y\ino  \bD ([0, \infty), \bbR)$ that is the space of c\`adl\`ag functions equipped with Skorokod's topology, and let $z, \eta \ino (0, \infty)$.
Recall from (\ref{modudu}) the notation $w_z (y, \eta) $ for the 
c\`adl\`ag modulus of continuity of $y\ino  \bD ([0, \infty), \bbR)$.
%%%
%%%: for all $z, \eta \ino (0, \infty)$:  
%%%\begin{equation}
%%%\label{remodudu}
%%%w_z(y,\eta)=  \inf \big\{\!  \max_{1\leq i \leq r } \mathtt{osc} (y, [t_{i-1} , t_i ) \, )\;  ; \; \,   0\! = \! t_0 \! < \! \ldots  \! < \! t_r\! = \! z \;  :  \!\! \!\!\!\! \min_{\quad 1\leq i \leq r-1} \!\!\!\! (t_i\! -\! t_{i-1}) \geq \eta  \;    \big\} , 
%%%\end{equation}
%%%where for all interval $I$, $\mathtt{osc} (y, I) \! = \! \sup \{ |y(s)\! -\! y(t) | ; \, s, t\ino I  \}$.  
We shall use a tightness result for \textit{nondecreasing processes} that is a consequence of Proposition 3.8.3 in Ethier \& Kurtz \cite{EtKu86}. To recall this statement, we need to introduce the following notation: assume that $y(\cdot) $ is 
nonnegative and nondecreasing; 
then for all $\epp\ino (0, \infty)$, we inductively define times $(\tau^\epp_k (y))_{k\in \bbN}$ by setting 
%\margmm{corrected}
\begin{equation}
\label{taueppdef}
\tau^{\epp}_0 (y)\! = \! 0 \quad \textrm{and} \quad \tau_{k+1}^\epp (y) = \inf \big\{ t >\tau^\epp_k (y) \!  : \, y(t)\! -\! y \big( \tau^\epp_k (y) \big) > \epp \big\} , 
\end{equation}
with the convention that $\inf \emptyset \! = \! \infty$. Observe that if $z\! >\! \eta$ and if 
$w_z (y, \eta) \! >\! \epp$, then there exists $k \! \geq \! 1$ such that $\tau^\epp_k (y) \! \leq  \! z$ and $ \tau_{k}^\epp (y) \! - \!  
\tau_{k-1}^\epp (y)  \! < \! \eta$. \textit{Indeed}, set $r\! = \! 1+ \max \{ k \ino \bbN\! : \! \tau^\epp_k (y) \! < z \}$. Note that $z\! >\! \eta$ and $w_z (y, \eta) \! >\! \epp$ imply that $r\! \geq \! 2$; then for all 
$i\ino \{ 0, \ldots , r\! -\! 1\}$, set $t_i\! = \! \tau^\epp_i(y)$ and $t_r\! = \! z$. By definition of the 
$\tau^\epp_i(y)$, we get $\max_{1\leq i \leq r } \mathtt{osc} (y, [t_{i-1} , t_i ) \, ) \! \leq \epp$. Since $w_z (y, \eta) \! >\! \epp$, we necessarily get $\min_{1\leq i \leq r-1} (t_i\! -\! t_{i-1}) \! < \! \eta$, which is the desired result.  
This observation combined with Lemma 3.8.2 of \cite{EtKu86} (recalled above as Lemma \ref{deltime}) immediately entails the following. 
\begin{lem} 
\label{tightincrea} For all $n\ino \bbN$, let $(R^n_t)_{t\in [0, \infty)}$ be a c\`adl\`ag nonnegative and nondecreasing process.  
%such that 
%a.s. $\lim_{t\rightarrow \infty} R_n (t)$ $\! = \! \infty$. 
Then, the laws of the $R^n$ are tight in $\bD ([0, \infty), \bbR)$ if for all $t\ino [0, \infty)$ 
the laws of the $R_n (t)$, $n\ino \bbN$ are tight on $\bbR$ and if 
\begin{equation}
\label{croiscrit}
\forall z, \epp \ino (0, \infty), \qquad  \lim_{\eta \rightarrow 0+} \, \limsup_{n\in \bbN} \, \sup_{k\in \bbN} \, \bP \big(  \tau^\epp_k (R^n) \! < \! z  \, ; \,   \tau^\epp_{k+1} (R^n)  \! -\!   \tau^\epp_k (R^n)  \! < \! \eta \big)\! = \! 0 .
\end{equation}
\end{lem}
\noi
\textbf{Proof.} See the previous arguments or Lemma 3.8.1 and Proposition 3.8.3 in Ethier \& Kurtz \cite{EtKu86} (pp.~134-135). \cqfd

\bigskip

We immediately apply Lemma \ref{tightincrea} in combination with the estimates in Lemma \ref{EstDifAw} to prove the tightness of a rescaled version of $A^{\bw_n}$. 
\begin{lem} 
\label{Atight} Let $(\alpha, \beta, \kappa, \mathbf{c})$ be as in (\ref{parconing}). Recall from (\ref{repsidefi}) the definition of $\psi$ and assume that (\ref{contH}) holds: namely, $\int^\infty d\lambda / \psi (\lambda) \! <\! \infty$.  
Let $a_n , b_n \! \in \! (0, \infty)$ and  $\bw_n \! \in \! \elldo_f$, $n\ino \bbN$, 
satisfy (\ref{apriori}) and $\mathbf{(C1)}$--$\mathbf{(C3)}$ as in (\ref{unalphcv}) and in (\ref{sig3cvcj}). 
%Let $\alpha \ino \bbR$, $\beta\ino [0, \infty)$,  $\kappa \ino (0, \infty)$ and $\mathbf{c}\ino \elldo_3$ satisfy (\ref{reassume}).  
%Let $a_n , b_n \! \in \! (0, \infty)$ and  $\bw_n \! \in \! \elldo_f$, $n\ino \bbN$, 
%satisfy (\ref{apriori}) and $\mathbf{(C1)}$--$\mathbf{(C3)}$. 
Recall from (\ref{YwAwSigw}) the definition of $A^{\bw_n}$. 
Then, the laws of $(\frac{_1}{^{a_n}}A^{\bw_n}_{b_n t})_{t\in [0,\infty)}$ are tight on $\bD ([0,\infty), \bbR)$. 
\end{lem}
\noi
\textbf{Proof.} We repeatedly use the following estimates on Poisson r.v.~$N$ with mean $r\ino (0, \infty)$: 
\begin{equation}
\label{estPoipoi}
\bE \big[ (N\! -\! 1)_+\big]\! =\!  e^{-r} \! -\! 1 + r \quad \textrm{and} \quad \mathbf{var} \big( (N\! -\! 1)_+\big)\! =\! r^2 \! -\!  (e^{-r} \! -\! 1 + r)(e^{-r}  + r) \leq r^2 .   
\end{equation}
By the definition (\ref{YwAwSigw}), we get $ \bE [A^{\bw_n}_t ]\! = \! \sum_{j\geq 1} w^{_{(n)}}_{^j} (e^{-tw^{_{(n)}}_{^j} / \sigma_1 (\bw_n)} \! -\! 1 +  \frac{ tw^{_{(n)}}_{^j}}{\sigma_1 (\bw_n)} ) \leq \frac{ t^2 \sigma_3 (\bw_n)}{2\sigma_1 (\bw_n)^2}$. Thus, by $\mathbf{(C1)}$--$\mathbf{(C3)}$ and the Markov inequality, we get 
$$ \limsup_{n\rightarrow \infty} \bP \Big(\frac{_1}{^{a_n}}A^{\bw_n}_{b_n t} \geq x \Big) \leq \frac{_1}{^2} x^{-1}t^2 \kappa \big( \kappa \sigma_3 (\mathbf{c}) + \beta \big) \underset{x\rightarrow \infty}{-\!\!\! -\!\!\! -\!\!\! \longrightarrow} 0. $$
This shows that for any $t\ino [0, \infty)$, the laws of the $\frac{_1}{^{a_n}}A^{\bw_n}_{b_n t}$ are tight on $\bbR$.  

   We next prove (\ref{croiscrit}) with $R^n_t \! = \!  \frac{_1}{^{a_n}}A^{\bw_n}_{b_n t}$, $t\ino [0, \infty)$. To that end, we fix $z, \epp \ino (0, \infty)$ and $k\ino \bbN$, and we set $T_n\! :=\! \tau^\epp_k (R^n) $. Then, (\ref{goretex}) in Lemma \ref{EstDifAw} with $a\! = \! a_n \epp$, $T \! =\!  b_nT_n$, $t\! = \! b_n \eta$ and $t_0\! = \! b_n z$ implies the following: 
\begin{eqnarray*}
\bP \big(  \tau^\epp_k (R^n) \! < \! z  \, ; \,   \tau^\epp_{k+1} (R^n)  \! -\!   \tau^\epp_k (R^n)  \! < \! \eta \big) \!\!\! & = &   \!\!\!
\bP \big( b_nT_n \! < \! b_n z  \, ; \,     A^{\bw_n}_{b_n T_n+ b_n \eta}\! -\!  A^{\bw_n}_{b_n T_n}  \! > \! a_n \epp \big) \\
 \!\!\! & \leq Â &  \!\!\! (a_n\epp)^{-1} b_n \eta \big( b_n z + \frac{_1}{^2} b_n \eta \big) \frac{\sigma_3 (\bw_n)}{\sigma_1 (\bw_n)^2} \\
  \!\!\! & \leq Â &  \!\!\! \epp^{-1} \eta (z + \eta) \frac{a_n b_n}{\sigma_1 (\bw_n)}  \frac{b_n \sigma_3 (\bw_n)}{a_n^2\sigma_1 (\bw_n)} . 
\end{eqnarray*}
%The inequality is a consequence of (\ref{goretex}) in Lemma \ref{EstDifAw} with $a\! = \! a_n \epp$, $T \! =\!  b_nT_n$, $t\! = \! b_n \eta$ and $t_0\! = \! b_n z$. 
Then $\mathbf{(C1)}$--$\mathbf{(C3)}$ entails (\ref{croiscrit}) and Lemma \ref{tightincrea} completes the proof. \cqfd  
%\medskip
%
%To simplify notation in the proof, we introduce the following for all $n\ino \bbN$, all $t\ino [0, \infty)$ and all $j\! \geq \! 1$: 
%\begin{equation}
%\label{simplinot}
%\XXX=  \frac{_{_1}}{^{^{a_n}}}X^{\mathtt{b},\bw_n}_{b_n t}, \quad \AAA=  \frac{_{_1}}{^{^{a_n}}}A^{\bw_n}_{b_n t} , \quad \NNNj  =  N^{\bw_n}_j(b_n t), \quad v^{(n)}_j= w^{(n)}_j / a_n \; .
%\end{equation}

\medskip

Recall from (\ref{Xbrwdef}) the definition of $X^{\mathtt{b}, \bw}$ and recall from (\ref{Njbdef}) the definition of the Poisson processes $N^{\bw}_j (\cdot)$, $j\! \geq \! 1$. 
Recall from (\ref{Xblue}) the definition of $X^{\mathtt{b}}$ and that of the Poisson processes $N_j (\cdot)$, $j\! \geq \! 1$. 
\begin{lem}
\label{cvjtNj} 
Let $(\alpha, \beta, \kappa, \mathbf{c})$ be as in (\ref{parconing}). Recall from (\ref{repsidefi}) the definition of $\psi$ and assume that (\ref{contH}) holds: namely, $\int^\infty d\lambda / \psi (\lambda) \! <\! \infty$.  
Let $a_n , b_n \! \in \! (0, \infty)$ and  $\bw_n \! \in \! \elldo_f$, $n\ino \bbN$, 
satisfy (\ref{apriori}) and $\mathbf{(C1)}$--$\mathbf{(C3)}$ as in (\ref{unalphcv}) and in (\ref{sig3cvcj}). Then, the following convergence
%Let $\alpha\ino \bbR$, $\beta\ino [0, \infty)$,  $\kappa \ino (0, \infty)$ and $\mathbf{c}\ino \elldo_3$ satisfy (\ref{reassume}). Let $a_n , b_n \! \in \! (0, \infty)$ and  $\bw_n \! \in \! \elldo_f$, $n\ino \bbN$, 
%satisfy (\ref{apriori}) and $\mathbf{(C1)}$--$\mathbf{(C3)}$. Then 
\begin{equation}
\label{megajtcv}
 \big(\big( \frac{_{_1}}{^{^{a_n}}}X^{\mathtt{b},\bw_n}_{b_n t} \big)_{t\in [0, \infty)} , (N^{\bw_n}_j(b_n t))_{t\in [0, \infty)} ; \, j\geq 1\, \big)  \underset{n\rightarrow \infty}{-\!\!\! -\!\!\! -\!\!\! \longrightarrow} \big( X^\mathtt{b}, N_j ; \, j\geq 1  \big) 
\end{equation} 
holds weakly on $(\bD ([0, \infty), \bbR))^\bbN$ equipped with the product topology.  
\end{lem}
\textbf{Proof.} Let $u \ino\bbR$. Note that 
$$\bE \big[ \! \exp (i u N^{\bw_n}_j(b_n t))\big]\! = \! 
\exp (- tb_n w^{_{(n)}}_{^j}\! (1\! -\! e^{iu }) / \sigma_1 (\bw_n) ) \! \longrightarrow \exp (- t\kappa c_j  (1\! -\! e^{iu }) )$$ 
by (\ref{apriori}) 
and $\mathbf{(C3)}$. Thus, for all $t\ino [0, \infty)$, $N^{\bw_n}_j(b_n t)\! \rightarrow \! N_j (t)$ in law. 
Next, fix $k\! \geq \! 1$ and set: 
$$\forall t\ino [0, \infty), \quad Q^n_t \! = \! \frac{_1}{^{a_n}}X^{\mathtt{b}, \bw_n}_{b_n t} \! -\! \sum_{1 \leq j\leq k} a_n^{-1}w^{_{(n)}}_{^j} N^{\bw_n}_j (b_nt) \quad \textrm{and} \quad 
Q_t \! = \! X^{\mathtt{b}}_t \! -\! \sum_{1\leq j\leq k} c_j  N_j (t) \; .$$ 
Since we assume that Proposition \ref{cvmarkpro} holds true,  $\frac{_1}{^{a_n}}X^{\mathtt{b}, \bw_n}_{b_n t} \! \rightarrow \! X^\mathtt{b}_t$ weakly on $\bbR$. Since $Q^n_t$ (resp.~$Q_t$) is independent of $(N^{\bw_n}_j)_{1\leq j\leq k}$ (resp.~independent of $(N_j)_{1\leq j\leq k}$), we easily check $$ \bE \big[ e^{iu Q^n_t }\big]\! = \! \bE \big[ e^{iu X^{\mathtt{b}, \bw_n}_{b_n t}/ a_n}\big] \big/ \!\! \!\! 
\prod_{1\leq j \leq k} \!\!\!\!  \bE \big[ e^{-i u \frac{w^{_{(n)}}_{^j}}{a_n} N^{\bw_n}_j(b_n t)} \big] \underset{n\rightarrow \infty}{ -\!\!\! -\!\!\! \longrightarrow} \bE \big[ e^{iu X^{\mathtt{b}}_{t}}\big] \big/ \!\!\!\! 
\prod_{1\leq j \leq k} \!\!\!\!  \bE \big[e^{-i u c_j N_j(t)} \big] \! =\!  \bE \big[ e^{iuQ_t}\big]. $$
Thus, $Q^n_t \! \rightarrow \! Q_t$ weakly on $\bbR$. Since L\'evy processes weakly converge in $\bD([0, \infty), \bbR)$ if an only if  unidimensional marginals weakly converge on $\bbR$ (see 
Lemma \ref{cvLevy} in Appendix Section \ref{SkoAppsec}, with precise references), we get 
$Q^n \! \rightarrow \! Q$ and for all $j\! \geq \! 1$, $N^{\bw_n}_j(b_n  \cdot )\! \rightarrow \! N_j$, weakly on $\bD([0, \infty), \bbR)$. 

  Since $Q^n, N^{\bw_n}_1, \ldots , N^{\bw_n}_k$ are independent L\'evy processes, they have a.s.~no common jump-times and 
Lemma \ref{jointos} (in Appendix, Section \ref{SkoAppsec}) asserts that 
$$(Q^n_t, N^{\bw_n}_1(b_n t), \ldots, N^{\bw_n}_k(b_n t))_{t\in [0, \infty)} \! \longrightarrow \! 
(Q, N_1, \ldots, N_k) \; \,  \textrm{weakly on $\bD ([0, \infty), \bbR^{k+1} )$.}$$ 
Since $X^{\mathtt{b}, \bw_n}$ is a linear combination of $Q^n$ and the $(N^{\bw_n}_j)_{1\leq j\leq k}$, we get: 
$$\big( (\frac{_{_1}}{^{^{a_n}}}X^{\mathtt{b}, \bw_n}_{b_n t}, N^{\bw_n}_1(b_n t), \ldots, N^{\bw_n}_k(b_n t) \big)_{t\in [0, \infty)} \! \longrightarrow \! 
(X^\mathtt{b}, N_1, \ldots, N_k) \; \,  \textrm{weakly on $\bD ([0, \infty), \bbR^{k+1} )$,}$$ which 
implies the weaker statement: 
$(\frac{_1}{^{a_n}}X^{\mathtt{b}, \bw_n}_{b_n  \cdot }, N^{\bw_n}_1(b_n \cdot), \ldots, N^{\bw_n}_k(b_n \cdot)) \! \longrightarrow \! 
(X^\mathtt{b}, N_1, \ldots, N_k)$, weakly on $(\bD ([0, \infty), \bbR ))^{k+1}$ equipped with the product topology. Since it holds 
true for all $k$, an elementary result (see Lemma \ref{obvious} in Appendix, 
Section \ref{SkoAppsec}) 
entails (\ref{megajtcv}). \cqfd

\medskip

Recall from (\ref{YwAwSigw}) the definition of $A^\bw$ and recall from (\ref{AetYdef}) the definition of $A$.
\begin{lem}
\label{cvjntAXb} 
Let $(\alpha, \beta, \kappa, \mathbf{c})$ be as in (\ref{parconing}). Recall from (\ref{repsidefi}) the definition of $\psi$ and assume that (\ref{contH}) holds: namely, $\int^\infty d\lambda / \psi (\lambda) \! <\! \infty$.  
Let $a_n , b_n \! \in \! (0, \infty)$ and  $\bw_n \! \in \! \elldo_f$, $n\ino \bbN$, 
satisfy (\ref{apriori}) and $\mathbf{(C1)}$--$\mathbf{(C3)}$ as in (\ref{unalphcv}) and in (\ref{sig3cvcj}). Then,
%%Let $\alpha \ino  \bbR$, $\beta\ino [0, \infty)$,  $\kappa \ino (0, \infty)$ and $\mathbf{c}\ino \elldo_3$ satisfy (\ref{reassume}). Let $a_n , b_n \! \in \! (0, \infty)$ and  $\bw_n \! \in \! \elldo_f$, $n\ino \bbN$, 
%%satisfy (\ref{apriori}) and $\mathbf{(C1)}$--$\mathbf{(C3)}$. Then 
\begin{equation}
\label{weakerXA}
 \big(\big( \frac{_{_1}}{^{^{a_n}}}X^{\mathtt{b},\bw_n}_{b_n t} \big)_{t\in [0, \infty)},
 \big( \frac{_{_1}}{^{^{a_n}}} A^{\bw_n}_{b_n t}\big)_{t\in [0, \infty)} \big)  \underset{n\rightarrow \infty}{-\!\!\! -\!\!\! -\!\!\! \longrightarrow} \big( X^\mathtt{b} , A \big) \quad \textrm{weakly on $(\bD ([0, \infty), \bbR))^2$.}
\end{equation} 
%$$ \big(\big( \frac{_1}{^{a_n}}X^{\mathtt{b},\bw_n}_{b_n t}, \frac{_1}{^{a_n}}A^{\bw_n}_{b_n t}\big)\big)_{t\in [0, \infty)}  \underset{n\rightarrow \infty}{-\!\!\! -\!\!\! -\!\!\! \longrightarrow} \big( (X^\mathtt{b}_t , A_t) \big)_{t\in [0, \infty)} \quad \textrm{weakly on $\bD ([0, \infty), \bbR^2)$. }$$
%$\big(\big( \frac{_1}{^{a_n}}X^{\mathtt{b},\bw_n}_{b_n t}, \frac{_1}{^{a_n}}A^{\bw_n}_{b_n t}\big)\big)_{t\in [0, \infty)} $ tend to 
%$((X^\mathtt{b}_t , A_t))_{t\in [0, \infty)}$ 
\end{lem}
\textbf{Proof.}
%we first prove the weaker statement: 
%\begin{equation}
%\label{weakerXA}
% \big(\big( \frac{_{_1}}{^{^{a_n}}}X^{\mathtt{b},\bw_n}_{b_n t} \big)_{t\in [0, \infty)},
% \big( \frac{_{_1}}{^{^{a_n}}} A^{\bw_n}_{b_n t}\big)_{t\in [0, \infty)} \big)  \underset{n\rightarrow \infty}{-\!\!\! -\!\!\! -\!\!\! \longrightarrow} \big( X^\mathtt{b} , A \big) \quad \textrm{weakly on $\bD([0, \infty), \bbR)^2$}
%\end{equation} 
%equiped with the product topology. 
%To that end, note that 
Lemma \ref{Atight} and Lemma \ref{cvjtNj} imply that the laws of  
$(\frac{_{1}}{^{{a_n}}}A^{\bw_n}_{b_n \cdot}, \frac{_{1}}{^{{a_n}}}X^{\mathtt{b},\bw_n}_{b_n \cdot} , N^{\bw_n}_j (b_n \cdot) ; j\! \geq \! 1)$ are tight on $(\bD ([0, \infty), \bbR))^\bbN$ equipped with the product topology. We want to prove that there is a unique limiting law: let $(n(p))_{p\in \bbN}$ be an increasing sequence of integers such that 
\begin{equation}
\label{grosscv}
 (\frac{_{1}}{^{{a_{n(p)}}}}A^{\bw_{n(p)}}_{b_{n(p)} \cdot}, \frac{_{1}}{^{{a_{n(p)}}}}X^{\mathtt{b},\bw_{n(p)}}_{b_{n(p)} \cdot} , N^{\bw_{n(p)}}_j (b_{n(p)} \cdot) ; j\! \geq \! 1) \underset{p\rightarrow \infty}{-\!\!\! -\!\!\! -\!\!\! \longrightarrow} \big( A^\prime , X^\mathtt{b} , N_j ; j\geq 1\big) , 
\end{equation} 
holds  weakly on $(\bD ([0, \infty), \bbR))^\bbN$. Since $(\bD ([0, \infty), \bbR))^\bbN$ equipped with the product topology is a Polish space,  Skorokod's representation theorem applies and 
 without loss of generality (but with a slight abuse of notation), we can assume that (\ref{grosscv}) holds true $\bP$-almost surely on $(\bD ([0, \infty), \bbR))^\bbN$. 
  
Recall that $A_t \! = \! \frac{_{_1}}{^{^2}} \kappa \beta t^2 +  \sum_{j\geq 1} c_j \big( N_j (t) \! -\! 1 \big)_+$, $t\ino [0, \infty)$. Then, 
to prove (\ref{weakerXA}), we claim that it is sufficient to prove that for all $t\ino [0, \infty)$, 
\begin{equation}
\label{proAA}
\frac{_{1}}{^{{a_{n(p)}}}}A^{\bw_{n(p)}}_{b_{n(p)} t} \! \longrightarrow \! A_t \quad \textrm{in probability.}
\end{equation}
%$\frac{{1}}{a_{n(p)}}A^{\bw_{n(p)}}_{b_{n(p)} t} \! \rightarrow \! A_t$ in probability 
\textit{Indeed}, let $t$ be such that $\Delta A^\prime_t\! = \! \Delta A_t\! = \! 0$ and let $q, q^\prime$ be  rational numbers 
 such that $q\! < \! t \! < \! q^\prime$; thus, $A^{_{\bw_{n(p)}}}_{^{b{(n(p))} q}} \! \leq \!A^{_{\bw_{n(p)}}}_{^{b{(n(p))} t}} \! \leq \! A^{_{\bw_{n(p)}}}_{^{b{(n(p))} q^\prime}}$; since 
 $\Delta A^\prime_t\! = \! 0$, we get a.s.~$A^{_{\bw_{n(p)}}}_{^{b_{n(p)} t}}/ a_{n(p)} \! \rightarrow \! A^\prime_t$; the convergence in probability entails that $A_q\! \leq \! A^\prime_t \! \leq \! A_{q^\prime}$; since it holds true 
 for all rational numbers $q,q^\prime$ such that 
 $q\! < \! t \! < \! q^\prime$, we get $A_{t-} \! \leq \! A^\prime_t \! \leq \! A_t$ which implies $A_t \!  = \! A^\prime_t$ since $\Delta A_t\! = \! 0$. Thus, a.s.~$A$ and $A^\prime$ coincide on the dense subset $\{ t\ino [0, \infty): \Delta A^\prime_t\! = \! \Delta A_t\! = \! 0\}$: it entails that a.s.~$A\! = \! A^\prime$ and 
the law of 
$(A, X^\mathtt{b},  N_j ; j\geq 1)$ 
is the unique weak limit of the laws of 
$(\frac{_{1}}{^{{a_n}}}A^{\bw_n}_{b_n \cdot}, \frac{_{1}}{^{{a_n}}}X^{\mathtt{b},\bw_n}_{b_n \cdot} , N^{\bw_n}_j (b_n \cdot) ; j\! \geq \! 1)$. 

Let us prove (\ref{proAA}). To simplify notation we define $\bv_n \ino \elldo_f$ by 
\begin{equation}
\label{vnjndef}
\forall j \ino \bbN^*, \quad v_j^{(n)}= w_j^{(n)}/ a_n \; .
%\quad \textrm{and} \quad j_n = \sup 
%\Big\{ j\geq 1 : \sum_{1\leq i\leq j} \! (v_{^i}^{_{(n)}})^3 \leq \sigma_3 (\mathbf{c} ) \Big\} , 
\end{equation}
%with the convention that $\sup \emptyset \! = \! 0$. 
By $\mathbf{(C3)}$, $v_j^{(n)} \! \rightarrow \! c_j$; by (\ref{apriori}) and 
$\mathbf{(C2)}$, $b_n/ \sigma_1 (\bv_n) \! \rightarrow \! \kappa$ and  $\sigma_3 (\bv_n)\! \rightarrow \! \sigma_3 (\mathbf{c}) + \beta/ \kappa$. 
We next claim that there exists $j_n \rightarrow \infty$ such that 
%This easily implies that 
\begin{equation}
\label{jncvsig3}
\lim_{n\rightarrow \infty}\!  v^{(n)}_{j_n}\! = \! 0, \quad 
 \lim_{ n\rightarrow \infty}\!\! \!\!\!\! \sum_{\quad 1\leq j \leq j_n} \!\!\!\! (v_{^j}^{_{(n)}})^3 \! = \! \sigma_3 (\mathbf{c})  
\quad  \!\! \textrm{and} \quad  \!\! \lim_{n\rightarrow \infty} \sum_{j > j_n} (v_{^j}^{_{(n)}})^3\!  =\! 
 \beta/ \kappa . 
\end{equation}
\textit{Proof of (\ref{jncvsig3}).} Indeed, suppose first that $\sup \{ j\! \geq \! 1\! : \! c_j \! >\! 0 \} \! = \! \infty$ and 
set $j_n \! =\!  \sup \big\{  j\geq 1 \! :\!  v_{^i}^{_{(n)}} \! >\! 0 \; \textrm{and} \;  \sum_{1\leq i\leq j}  (v_{^i}^{_{(n)}})^3 \! \leq \! \sigma_3 (\mathbf{c} )\big\}$, with the convention that $\sup \emptyset \! = \! 0$. Here $j_n \! \to \! \infty$, and it is easy to check that it satisfies (\ref{jncvsig3}).

Next suppose that $j_*\! = \! \sup \{ j\! \geq \! 1\! : \! c_j \! >\! 0 \} \! < \! \infty$. 
Clearly $\sum_{1\leq j\leq j_*} (v^{_{(n)}}_{^j})^3\! \to \! \sigma_3 (\mathbf{c})$ and 
$\sum_{j> j_*} (v^{_{(n)}}_{^j})^3\! \to \!  \beta/\kappa $. 
Since for all $j\! >\! j^*$, $v^{_{(n)}}_{^j}\! \rightarrow \! 0$ it is possible to 
find a sequence $(j_n)$ that tends to $\infty$ sufficiently slowly to get $\sum_{ j_*< j\leq j_n} \! (v^{_{(n)}}_{^j})^3 \!  \rightarrow \! 0$, which implies (\ref{jncvsig3}). \cq

\smallskip

Next, we use (\ref{jncvsig3}) to prove (\ref{proAA}). To that end, we fix $t\ino [0, \infty)$ 
and we fix $k\ino \bbN$  that is specified further; since $j_n \! \rightarrow \! \infty$, we can assume $p$ is such that $k\! < \! j_{n(p)}$. To simplify, we set $\xi^{p}_j \! = \!  v^{(n(p))}_j \big(N^{\bw_{n(p)}}_j (b_{n(p)} t)\! -\! 1 \big)_+$ and $\xi_j\! = \!  c_j \big(N_j (t)\! -\! 1 \big)_+$ and 
$$ D^{k,p}_t \! = \! \!\! \sum_{1\leq j\leq k} \!\!\! \xi_j^p \! -\! \xi_j, \; \,  R^{k,p}_t \! = \!  \!\! \! \!\! \sum_{\; k<  j \leq j_{n(p)}}  \!\!\! \!\! \!\!  \xi^p_j  - \sum_{j>k} \xi_j , \; \, C^{p}_t \! = \!  \!\! \! \!\!  \!\!  \sum_{\quad j>j_{n(p)}} \!\!  \!\!  \!\!   \xi_j^p \! -\! \bE [\xi^p_j] \quad   \textrm{and} \; \,  d_p(t)\! =\! 
 \frac{_1}{^2}\kappa  \beta t^2 \! -\!\! \!\! \!\! \!\!  \sum_{\quad j>j_{n(p)}}\!\! \!\! \!\!  \bE [\xi^p_j]. $$
Thus, $A^{\bw_{n(p)}} (b_{n(p)} t )/ a_{n(p)} - A_t\! = \! D^{k,p}_t + R^{k,p}_t + C^{p}_t  - d_p (t)$ and we 
prove that each term in the right-hand side goes to $0$ in probability. 

We first show that $d_p (t) \! \rightarrow \! 0$. 
Since $N^{_{\bw_{n(p)}}}_{^j} (b_{n(p)} t)$ is a Poisson r.v.~with mean 
$r_{p,j}$ that is equal to $v^{_{(n(p))}}_{^j} b_{n(p)} t/\sigma_1 (\bv_{n(p)})$,  
by (\ref{estPoipoi}) we get $\bE \big[ \xi_j^p\big] \! = \! v^{_{(n(p))}}_{^j} \big( e^{-r_{p,j}}\! -\! 1+ r_{p,j})$. 
We next use the following elementary inequality: 
\begin{equation}
\label{ordeexp}
\forall y\ino [0, \infty), \quad  0 \leq \frac{_1}{^2} y^2 - \big( e^{-y} \! -\! 1 + y\big)  \leq  \frac{_1}{^2}y^2 (1\! -\! e^{-y}) \leq  \frac{_1}{^2}  \, y^2 \! \wedge \! y^3 \; , 
\end{equation}
that holds true since $y^{-2} (e^{-y} \! -\! 1 + y)\! = \! \int_0^1 \! dv \int_0^v \! dw \, e^{-wy}$. Thus 
$$0\!  \leq \! \!\!  \!\!  \!\!  \sum_{\quad j>j_{n(p)}}  \!\!  \!\!  \!\! \frac{_1}{^2}v^{_{(n(p))}}_{^j}\!  r_{p,j}^2 \! \! -\!  \bE \big[ \xi_j^p\big]  
 \leq   \! \!\!  \!\!  \!\!  \sum_{\quad j>j_{n(p)}}  \!\!  \!\!  \!\!  \frac{_1}{^2} v^{_{(n(p))}}_{^j} r_{p,j}^3 \leq \frac{_1}{^2} v^{_{(n(p))}}_{^{j_{n(p)}}} \frac{_{(b_{n(p)}t)^3}}{^{\sigma_1 (\bv_{n(p)})^3}}  \! \!\!  \!\!  \!\!  \sum_{\quad j>j_{n(p)}}  \!\!  \!\!  \!\!   (v^{_{(n(p))}}_{^j} )^3 \longrightarrow 0, $$
by (\ref{jncvsig3}). Next, note that $ \sum_{ j>j_{n(p)}}\! v^{_{(n(p))}}_{^j} r_{p,j}^2 \! =\!  (b_{n(p)}t/ \sigma_1(\bv_{n(p)})^2  \sum_{ j>j_{n(p)}}  \!  (v^{_{(n(p))}}_{^j} )^3\!  \longrightarrow \!  \kappa\beta  t^2$, which implies that $d_p (t) \! \rightarrow 0$ as $p\! \rightarrow \! \infty$. 

We next consider $C^p_t$: by (\ref{estPoipoi}), $\mathbf{var} (\xi^p_j) \! \leq \! (v^{_{(n(p))}}_{^j})^2 r_{p,j}^2$. 
Since the $\xi^p_j$ are independent, we get 
$$ \bE \big[ (C^p_t)^2\big]\! =\!  \!\! \! \!\!  \!\!  \sum_{\quad j>j_{n(p)}} \!\!  \!\!  \!\!  \mathbf{var} (\xi^p_j) \leq v^{_{(n(p))}}_{^{j_{n(p)}}} \frac{_{(b_{n(p)}t)^2}}{^{\sigma_1 (\bv_{n(p)})^2}}  \! \!\!  \!\!  \!\!  \sum_{\quad j>j_{n(p)}}  \!\!  \!\!  \!\!   (v^{_{(n(p))}}_{^j} )^3 \longrightarrow 0 $$
by (\ref{jncvsig3}), which proves that $C^p_t \! \rightarrow \! 0$ in probability when $p\! \rightarrow \! \infty$. 

We next deal with $R^{k,p}_t$. By (\ref{estPoipoi}), (\ref{jncvsig3}) and (\ref{ordeexp}), we first get: 
\begin{equation}
\label{grugnouc}
  0\!  \leq \! \!\!  \!\!  \!\!  \sum_{\quad k< j\leq j_{n(p)}}  \!\!  \!\!  \!\!  \bE \big[ \xi_j^p\big]  \! \leq  \!   \! \!\!  \!\!  \!\!  \sum_{\quad k< j\leq j_{n(p)}}  \!\!  \!\!  \!\!  \frac{_1}{^2} v^{_{(n(p))}}_{^j} r_{p,j}^2 \! = \,  \frac{_1}{^2} \frac{_{(b_{n(p)}t)^2}}{^{\sigma_1 (\bv_{n(p)})^2}}  \! \!\!  \!\!  \!\! \!\!\!  \sum_{\quad k< j\leq j_{n(p)}}  \!\!  \!\!  \!\!   (v^{_{(n(p))}}_{^j} )^3\underset{p\rightarrow \infty}{-\!\!\! -\!\!\! -\!\!\! \longrightarrow} \frac{_1}{^2} (\kappa t)^2 
  \sum_{ j >k}     c_j^3.  
\end{equation}
Similarly, observe that $\bE [\xi_j]\! = \! c_j \big(e^{-\kappa tc_j} \! -\! 1 + \kappa t c_j \big)\leq \frac{_1}{^2} (\kappa t)^2 c_j^3$. This inequality combined with (\ref{grugnouc}) entails: 
\begin{equation}
\label{youkiflop}
\limsup_{p\rightarrow \infty} \bE \big[ |R^{k,p}_t|\big] \leq (\kappa t)^2  \sum_{ j >k}      c_j^3 \;  \underset{k\rightarrow \infty}{-\!\!\! -\!\!\! -\!\!\! \longrightarrow}\; 0 . 
\end{equation}

  Finally, we consider $D^{k,p}$. Since a.s.~$t$ is not a jump-time of $N_j$, a.s.~$v^{_{n(p)}}_{^j}(N^{\bw_{n(p)}}_j (b_{n(p)} t) \! -\! 1)_+ 
\! \! \rightarrow\!  c_j (N_j (t) \! -\! 1)_+ $. Thus, for all $k\ino \bbN$, a.s.~$D^{k,p}_t \! \rightarrow \! 0$. 
These limits combined with (\ref{youkiflop}) (and with the convergence to $0$ in probability of $C^p_t$ and $d_p(t)$) easily imply   
(\ref{proAA}), which completes the proof of the lemma. \cqfd  

\medskip

Recall from (\ref{YwAwSigw}) the definition of $Y^\bw$ and recall from (\ref{AetYdef}) the definition of $Y$.
\begin{lem}
\label{cvjntYAXb} 
Let $(\alpha, \beta, \kappa, \mathbf{c})$ be as in (\ref{parconing}). Recall from (\ref{repsidefi}) the definition of $\psi$ and assume that (\ref{contH}) holds: namely, $\int^\infty d\lambda / \psi (\lambda) \! <\! \infty$.  
Let $a_n , b_n \! \in \! (0, \infty)$ and  $\bw_n \! \in \! \elldo_f$, $n\ino \bbN$, 
satisfy (\ref{apriori}) and $\mathbf{(C1)}$--$\mathbf{(C3)}$ as in (\ref{unalphcv}) and in (\ref{sig3cvcj}). 
%%%Let $\alpha \ino \bbR$, $\beta\ino [0, \infty)$,  $\kappa \ino (0, \infty)$ and $\mathbf{c}\ino \elldo_3$ satisfy (\ref{reassume}). Let $a_n , b_n \! \in \! (0, \infty)$ and  $\bw_n \! \in \! \elldo_f$, $n\ino \bbN$, 
%%%satisfy (\ref{apriori}) and $\mathbf{(C1)}$--$\mathbf{(C3)}$. Then,  
%\begin{equation}
%\label{weakerXA}
% \big(\big( \frac{_{_1}}{^{^{a_n}}}X^{\mathtt{b},\bw_n}_{b_n t} \big)_{t\in [0, \infty)},
% \big( \frac{_{_1}}{^{^{a_n}}} A^{\bw_n}_{b_n t}\big)_{t\in [0, \infty)} , \big( \frac{_{_1}}{^{^{a_n}}}X^{\mathtt{b},\bw_n}_{b_n t} \big)_{t\in [0, \infty)} \big)  \underset{n\rightarrow \infty}{-\!\!\! -\!\!\! -\!\!\! \longrightarrow} \big( X^\mathtt{b} , A, Y \big) \quad \textrm{weakly on $\bD([0, \infty), \bbR)^3$}
%\end{equation} 
\begin{equation}
\label{tresjtYAX}
 \big(\big( \frac{_1}{^{a_n}}X^{\mathtt{b},\bw_n}_{b_n t}\! , \frac{_1}{^{a_n}}A^{\bw_n}_{b_n t} ,  \frac{_1}{^{a_n}}Y^{\bw_n}_{b_n t}\big)\big)_{t\in [0, \infty)} \!  \overset{_{\textrm{weakly}}}{\underset{n\rightarrow \infty}{-\!\!\! -\!\!\! -\!\!\! \longrightarrow}} \big( (X^\mathtt{b}_t , A_t, Y_t) \big)_{t\in [0, \infty)} \; \,  \textrm{in $\bD ([0, \infty), \bbR^3)$. }
\end{equation}  
%$\big(\big( \frac{_1}{^{a_n}}X^{\mathtt{b},\bw_n}_{b_n t}, \frac{_1}{^{a_n}}A^{\bw_n}_{b_n t}\big)\big)_{t\in [0, \infty)} $ tend to 
%$((X^\mathtt{b}_t , A_t))_{t\in [0, \infty)}$ 
\end{lem}
\textbf{Proof.} without loss of generality (but with a slight abuse of notation), by Skorokod's representation theorem we can assume that the convergence in (\ref{weakerXA}) holds true $\bP$-almost surely. We first prove that 
$( (\frac{_1}{^{a_n}}X^{\mathtt{b},\bw_n}_{b_n \cdot}\! , \frac{_1}{^{a_n}}A^{\bw_n}_{b_n \cdot}))\! \rightarrow \! ((X^\mathtt{b}, A))$ a.s.~in $\bD([0, \infty), \bbR^2)$ thanks to Lemma \ref{jtsko} $(iii)$ (a standard result recalled in Appendix, Section \ref{SkoAppsec}). 
To that end, first recall that by definition, the jumps of $A$ (resp.~of $A^{\bw_n}$) are jumps of $X^{\mathtt{b}}$ (resp.~of $X^{\mathtt{b}, \bw_n}$): namely if $\Delta  A_t \! >\! 0$, then $\Delta X^{\mathtt{b}}_t \! = \! \Delta  A_t $. The same holds true for $A^{\bw_n}$ and $X^{\mathtt{b}, \bw_n}$. 
 
  Let $t\ino (0, \infty)$. First suppose that $\Delta A_t \! >\! 0$. Thus, $\Delta X^\mathtt{b}_t \! = \! \Delta A_t$. By  
Lemma \ref{jtsko} $(i)$, there exists a sequence of times $t_n \! \rightarrow \! t$ such that $\frac{_1}{^{a_n}} 
\Delta A^{\bw_n}_{b_n t_n }  \! \rightarrow \! \Delta A_t$. Thus, for all sufficiently large $n$, 
$\frac{_1}{^{a_n}} \Delta A^{\bw_n}_{b_n t_n }\!> \! 0$, which entails $\frac{_1}{^{a_n}} \Delta A^{\bw_n}_{b_n t_n }\! =\! 
\frac{_1}{^{a_n}} \Delta X^{\mathtt{b}, \bw_n}_{b_n t_n }$ and we get $ \frac{_1}{^{a_n}} \Delta X^{\mathtt{b}, \bw_n}_{b_n t_n }\! \rightarrow \! \Delta A_t\! = \! \Delta X^\mathtt{b}_t$. 
Suppose next that $\Delta A_t \! =\! 0$; by  Lemma \ref{jtsko} $(i)$, there exists a sequence of times 
$t_n \! \rightarrow \! t$ such that $ \frac{_1}{^{a_n}} \Delta X^{\mathtt{b}, \bw_n}_{b_n t_n }\! \rightarrow \! 
 \Delta X^\mathtt{b}_t$. Since $\Delta A_t \! =\! 0$, Lemma \ref{jtsko} $(ii)$ entails that $\frac{_1}{^{a_n}} 
\Delta A^{\bw_n}_{b_n t_n }  \! \rightarrow \! \Delta A_t \! =\! 0$. In both cases, we have proved that for all $t\ino (0, \infty)$, there exists 
 a sequence of times $t_n \! \rightarrow \! t$ such that $ \frac{_1}{^{a_n}} \Delta X^{\mathtt{b}, \bw_n}_{b_n t_n }\! \rightarrow \! 
 \Delta X^\mathtt{b}_t$ and $\frac{_1}{^{a_n}} 
\Delta A^{\bw_n}_{b_n t_n }  \! \rightarrow \! \Delta A_t$: by Lemma \ref{jtsko} $(iii)$, it implies that 
$( (\frac{_1}{^{a_n}}X^{\mathtt{b},\bw_n}_{b_n \cdot}\! , \frac{_1}{^{a_n}}A^{\bw_n}_{b_n \cdot}))\! \rightarrow \! ((X^\mathtt{b}, A))$ a.s.~in $\bD([0, \infty), \bbR^2)$. This entails (\ref{tresjtYAX}), since 
the function $(x, a)\ino \bbR^2\!  \mapsto\!  (x, a, x\! -\! a) \ino \bbR^3$ is Lipschitz 
and since $X^{\mathtt{b}, \bw_n}\! -\! A^{\bw_n}\! = \! Y^{\bw_n}$ and $X^{\mathtt{b}}\! -\! A\! = \! Y$. \cqfd 

\medskip

Recall that $X^{\mathtt{r}, \bw}$ (resp.~$X^\mathtt{r}$) is an independent copy of $X^{\mathtt{b}, \bw}$ (resp.~of 
$X^\mathtt{b}$). Recall from (\ref{thetabw}) (resp.~from (\ref{gammadef})) the definition of 
$\subo^{\mathtt{r}, \bw}$ (resp.~of $\subo^{\mathtt{r}}$). Recall that $I_t^{\mathtt{r}, \bw}\! = \! \inf_{s\in [0, t] } X_s^{\mathtt{r}, \bw}$ and recall the notation $I_\infty^{\mathtt{r}, \bw}\! = \! \lim_{t\rightarrow \infty} I_t^{\mathtt{r}, \bw}$. Similarly, recall that $I_t^{\mathtt{r}}\! = \! \inf_{s\in [0, t] } X_s^{\mathtt{r}}$ and recall the notation $I_\infty^{\mathtt{r}}\! = \! \lim_{t\rightarrow \infty} I_t^{\mathtt{r}}$. Recall from (\ref{schnitzel2}) in Lemma \ref{minipoub} 
the definition of $\overline{\gamma}^{\mathtt{r}}$; similarly we set
\begin{equation}
\label{schnitzel22}
\overline{\gamma}^{\mathtt{r}, \bw}_{x} = \gamma^{\mathtt{r}, \bw}_x \quad \textrm{if $x\! <\! -I^{\mathtt{r}, \bw}_\infty$} \quad \textrm{and} \quad \overline{\gamma}^{\mathtt{r}, \bw}_{x} = \gamma^{\mathtt{r}, \bw}  ((-I^{\mathtt{r}, \bw}_\infty) -) \quad  \textrm{if $x\! \geq \! -I^{\mathtt{r}, \bw}_\infty$. }
\end{equation}

\begin{lem}
\label{cvRouge1}
Let $(\alpha, \beta, \kappa, \mathbf{c})$ be as in (\ref{parconing}). Recall from (\ref{repsidefi}) the definition of $\psi$ and assume that (\ref{contH}) holds: namely, $\int^\infty d\lambda / \psi (\lambda) \! <\! \infty$.  
Let $a_n , b_n \! \in \! (0, \infty)$ and  $\bw_n \! \in \! \elldo_f$, $n\ino \bbN$, 
satisfy (\ref{apriori}) and $\mathbf{(C1)}$--$\mathbf{(C3)}$ as in (\ref{unalphcv}) and in (\ref{sig3cvcj}). Then, 
% Let $\alpha\ino \bbR$ , $\beta\ino [0, \infty)$,  $\kappa \ino (0, \infty)$ and $\mathbf{c}\ino \elldo_3$ satisfy (\ref{reassume}). Let $a_n , b_n \! \in \! (0, \infty)$ and  $\bw_n \! \in \! \elldo_f$, $n\ino \bbN$, 
%satisfy (\ref{apriori}) and $\mathbf{(C1)}$--$\mathbf{(C3)}$. Then, 
\begin{equation}
\label{simplred}
 \big(  \big( \frac{_{_1}}{^{^{a_n}}} X^{\mathtt{r}, \bw_n}_{b_n t}\big)_{t\in [0, \infty)} , \big( \frac{_{_1}}{^{^{b_n}}}\overline{\subo}^{\mathtt{r},\bw_n}_{a_n x} \big)_{x\in [0, \infty)},
 -\frac{_{_1}}{^{^{a_n}}} I^{\mathtt{r}, \bw_n}_{\infty} \big)  \underset{n\rightarrow \infty}{-\!\!\! -\!\!\! -\!\!\! \longrightarrow} \big( X^\mathtt{r}, \overline{\subo}^{\mathtt{r} } ,  -I^{\mathtt{r}}_{\infty}   \big) 
 %\quad \textrm{weakly on $\bD([0, \infty), \bbR)^2\! \times \! [0, \infty]$.}
\end{equation} 
weakly on $(\bD([0, \infty), \bbR))^2\! \times \! [0, \infty]$. 
\end{lem}
\textbf{Proof.} Let $\widetilde{\gamma}^n$ (resp.~$\widetilde{\gamma}$) be a conservative subordinator with Laplace exponent $a_n \psi_{\bw_n}^{-1} (\cdot / b_n)- a_n \varrho_{\bw_n}$ (resp.~$\psi^{-1} (\cdot) \! - \! \varrho$). By (\ref{eierschwa}) in Lemma \ref{expcovcv} , $a_n \psi_{\bw_n}^{-1} (\lambda / b_n)\! -\! a_n \varrho_{\bw_n} \! \rightarrow \! \psi^{-1} (\lambda) \! - \! \varrho$ for all $\lambda \ino [0, \infty)$, which implies that for all $x\ino [0, \infty)$, $\widetilde{\gamma}^n_x \! \rightarrow \! \widetilde{\gamma}_x$ weakly on $[0, \infty)$. Since the $\widetilde{\subo}^n$ are L\'evy processes, Theorem \ref{cvLevy} (in Appendix Section \ref{SkoAppsec}) 
entails that $\widetilde{\gamma}^n \! \rightarrow \! \widetilde{\gamma}$ weakly on $\bD([0, \infty), \bbR)$. Let $\cE_n$ (resp.~$\cE$) be an exponentially distributed r.v.~with parameter $a_n \varrho_{\bw_n}$ (resp.~$\varrho$) that is independent of $\widetilde{\gamma}^n$ (resp.~of $\widetilde{\gamma}$), with the convention that a.s.~$\cE_n \! = \! \infty$ if $\varrho_{\bw_n}\! = \! 0$ (resp.~a.s.~$\cE \! = \! \infty$ if $\varrho\! = \! 0$). We then get 
$(\widetilde{\gamma}^n, \cE_n) \! \rightarrow \! (\widetilde{\gamma}, \cE)$ 
weakly on $\bD([0, \infty), \bbR)\! \times \! [0, \infty]$. An easy application of Lemma \ref{extraisko} $(i)$ entails that $\big((\widetilde{\gamma}^n_{x\wedge \cE_n})_{x\in [0, \infty)}, \cE_n \big) \! \rightarrow \! \big((\widetilde{\gamma}_{x\wedge \cE})_{x\in [0, \infty)}, \cE \big) $ 
weakly on $\bD([0, \infty), \bbR)\times [0, \infty]$. By (\ref{schnitzel1}), 
we get 
$\big( \frac{1}{{b_n}} \overline{\subo}^{\mathtt{r}, \bw_n}_{a_n \cdot} ,  -\frac{{1}}{{{a_n}}} I^{\mathtt{r}, \bw_n}_{\infty}  \big) \! \rightarrow \! \big( \overline{\subo}^{\mathtt{r}}_{\cdot} ,  - I^{\mathtt{r}}_{\infty}\big)$ weakly on $\bD([0, \infty), \bbR)\times [0, \infty]$. Under our assumptions, Proposition \ref{cvmarkpro} implies that $\frac{{1}}{{{a_n}}} X^{\mathtt{r}, \bw_n}_{b_n \cdot} \! \rightarrow \! X^{\mathtt{r}}_{\cdot}$ weakly on $\bD([0, \infty), \bbR)$. 
Then the laws of the processes on the left hand side of (\ref{simplred}) are tight on 
$\bD ([0, \infty), \bbR)^2\! \times \! [0, \infty]$; we only need to prove that the joint law of the processes on the right hand side of (\ref{simplred}) is the unique limiting law: to that end, let $(n(p))_{p\in \bbN}$ be an increasing sequence of integers such that 
\begin{equation}
\label{gloubired}
 \big(  \big( \frac{{_1}}{{^{a_{n(p)}}}} X^{\mathtt{r}, \bw_{n(p)}}_{b_{n(p)} t}\big)_{t\in [0, \infty)} ,  \big( \frac{{_1}}{{^{b_{n(p)}}}}\overline{\subo}^{\mathtt{r},\bw_{n(p)}}_{a_{n(p)} x} \big)_{x\in [0, \infty)},
-\frac{{_1}}{{^{a_{n(p)}}}} I^{\mathtt{r}, \bw_{n(p)}}_{\infty} \big)  \underset{p\rightarrow \infty}{-\!\!\! -\!\!\! -\!\!\! \longrightarrow} \big( X^\mathtt{r}, \gamma^\prime,   \cE^\prime  \big) \,,
 %\quad \textrm{weakly on $\bD([0, \infty), \bbR)^2\! \times \! [0, \infty]$.}
\end{equation} 
where $( \gamma^\prime,   \cE^\prime)$ has the same law as $(\overline{\subo}^{\mathtt{r}} , -I^\mathtt{r}_\infty)$. Without loss of generality (but with a slight abuse of notation), by Skorokod's representation theorem we can assume that the convergence in (\ref{gloubired}) holds $\bP$-a.s.~and we only need to prove that $( \gamma^\prime,   \cE^\prime)\! = \! (\overline{\subo}^{\mathtt{r}} , -I^\mathtt{r}_\infty)$ a.s. 

   We first prove that a.s.~$\cE^\prime\! = \! -I^\mathtt{r}_\infty$. Since $X^\mathtt{r}$ is a 
spectrally positive L\'evy process, it has no fixed discontinuity. Moreover,
$t\mapsto \inf_{[0, t]} X^\mathtt{r}$ is continuous. Then, by Lemma \ref{franchtime} $(ii)$, for all $t\ino [0, \infty)$, a.s.~
$a_{^{n(p)}}^{_{-1}} \inf_{s\in [0, t]}X^{{\mathtt{r}, \bw_{n(p)}}}_{{b_{n(p)} s}}\! \rightarrow \!   \inf_{[0, t]} X^\mathtt{r}$. Since $-a_{^{n(p)}}^{_{-1}} \inf_{s\in [0, t]}X^\mathtt{r, \bw_{n(p)}}_{b_{n(p)} s} \! \leq \! 
-a_{^{n(p)}}^{_{-1}} I^{\mathtt{r}, \bw_{n(p)}}_\infty \! \rightarrow \! \cE^\prime$, for all $t\ino [0, \infty)$, we get a.s.~$- \inf_{[0, t]} X^\mathtt{r} \! \leq \! \cE^\prime$. Namely, $-I^\mathtt{r}_\infty \! \leq \! \cE^\prime$. Since $\cE^\prime$ and $-I^\mathtt{r}_\infty$ have the same law on $[0, \infty]$, we get $\cE^\prime\! = \! -I^\mathtt{r}_\infty$ a.s.

We next prove that a.s.~for all $x\ino [0, -I^\mathtt{r}_\infty)$, $\gamma^\prime_x\! = \! \overline{\subo}_x$. Indeed, fix $x\! < \!  -I^\mathtt{r}_\infty$ such that $\Delta \gamma^\mathtt{r}_x \! = \! 0$. Then,  by Lemma \ref{franchtime} $(iv)$, we get 
$\subo^{\mathtt{r},\bw_{n(p)}}_{a_{n(p)} x}/b_{n(p)}\! \rightarrow \! \gamma^\mathtt{r}_x$. Since $x\! < \!  
- I^{\mathtt{r}, \bw_{n(p)}}_{\infty}\! /a_{n(p)}$ for all sufficiently large $p$, it shows that $\overline{\subo}^{\mathtt{r},\bw_{n(p)}}_{a_{n(p)} x}/b_{n(p)}\! \rightarrow \! \gamma^\mathtt{r}_x\! = \! \overline{\gamma}^\mathtt{r}_x$. Thus, a.s.~for all $x\ino [0, -I^\mathtt{r}_\infty)$ such that $\Delta \gamma^\mathtt{r}_x \! = \! 0$, we get $\gamma^\prime_x\! = \! \overline{\subo}_x$, which implies the desired result. 
Note that it completes the proof of the lemma in the critical and subcritical cases. 

%\margmm{see footnote \tt{and the aswer}}
To avoid trivialities, we now assume that we are in the supercritical cases. 
Namely, $\varrho \! >\! 0$ and $-I^\mathtt{r}_\infty \! < \! \infty$ a.s. To simplify notation,  we set 
$$t_*^p\! = \! \tfrac{1}{b_{n(p)}}\gamma^{\mathtt{r}, \bw_{n(p)}} ( (-I^{\mathtt{r}, \bw_{n(p)} }_\infty)-) \quad \textrm{and} \quad  
t_*\! =\!   \gamma^{\mathtt{r}} ((-I^{\mathtt{r}}_\infty)-) \; .$$
First note that the proof is complete as soon as we 
prove that $t^p_* \! \rightarrow \! t_*$. To prove this limit, we want to use Lemma \ref{franchtime} $(iii)$. 
To that end, we first fix $x\! >\! -I^{\mathtt{r}}_\infty$. Since $( \gamma^\prime,   \cE^\prime)$ has the same law as $(\overline{\subo}^{\mathtt{r}} , -I^\mathtt{r}_\infty)$, $\gamma^\prime$ is constant on $[\cE^\prime, \infty)$ and since $\cE^\prime\! = \! -I^\mathtt{r}_\infty$, $\gamma^\prime$ is constant on $[-I^\mathtt{r}_\infty, \infty)$, which implies 
$\Delta\gamma^\prime_x \! = \! 0$ and thus $\overline{\subo}^{\mathtt{r},\bw_{n(p)}}_{a_{n(p)} x}/b_{n(p)}\! \rightarrow \! \gamma^\prime_x$.
We next fix $t\! >\! \gamma^\prime_x+   t_*$. 
%Since $\gamma^\prime$ is constant on $[- I^{\mathtt{r}}_\infty, \infty)$, $\Delta\gamma^\prime_x \! = \! 0$ and 
%$\overline{\subo}^{\mathtt{r},\bw_{n(p)}}_{a_{n(p)} x}/b_{n(p)}\! \rightarrow \! \gamma^\prime_x$. 
Thus, there is $p_0$ such that for all $p\! \geq \! p_0$, $\overline{\subo}^{\mathtt{r},\bw_{n(p)}}_{a_{n(p)} x}/b_{n(p)} \! < \! t $ and $x\! > \!  
- I^{\mathtt{r}, \bw_{n(p)}}_{\infty}\! /a_{n(p)}$, which implies that $t^p_*\! = \! 
\overline{\subo}^{\mathtt{r},\bw_{n(p)}}_{a_{n(p)} x}/b_{n(p)}$. 
Since $t \! > \! t^p_*\vee t_*$, we get 
$t^p_*\! = \! \inf \{ s\ino [0, t] \! : \! \inf_{r\in [0 ,s] } X^{\mathtt{r}, \bw_{n(p)}}_{b_{n(p)} r}\! = \! 
\inf_{r\in [0 ,t] } X^{\mathtt{r}, \bw_{n(p)}}_{b_{n(p)} r} \}$ and 
$t_*\! = \! \inf \{ s\ino [0, t] \! : \! \inf_{[0 ,s] } X^{\mathtt{r}}\! = \! 
\inf_{[0 ,t] } X^{\mathtt{r}} \}$. Thus Lemma \ref{franchtime} $(iii)$ entails that $t^p_*\! \rightarrow \! t_*$, which completes the proof of the lemma. 
%\fmm{The arguments in this paragraph are not very well organised: \tt{indeed, I hope it is better now.}}
\cqfd   
%
%
%
%
%recall that Proposition \ref{cvmarkpro} asserts that $\frac{_{_1}}{^{^{a_n}}} X^{\mathtt{r}, \bw_n}_{b_n \cdot} \! \rightarrow \! X^{\mathtt{r} }$ weakly on $\bD([0, \infty), \bbR)$. Let $0 \! \leq \! x_1\! < \!  \ldots \! < \! x_k$. Since $\gamma^\mathtt{r}$ is subordinator, it has no fixed 
%time-discontinuity; thus, a.s.~$\Delta \gamma^\mathtt{r}_{x_j}\! = \! \ldots \! = \!  \Delta \gamma^\mathtt{r}_{x_k}\! = \! 0$; a 
%standard result recalled in Lemma \ref{franchtime} in Appendix Section \ref{SkoAppsec}, implies the joint convergence:  
%$$\big( \frac{_{_1}}{^{^{a_n}}} X^{\mathtt{r}, \bw_n}_{b_n \cdot} ; \frac{_{_1}}{^{^{b_n}}}\subo^{\mathtt{r},\bw_n}_{a_n x_1} , \ldots , \frac{_{_1}}{^{^{b_n}}}\subo^{\mathtt{r},\bw_n}_{a_n x_k} \big) \! \longrightarrow \! \big( X^{\mathtt{r} } ; \gamma^\mathtt{r}_{x_1}, \ldots , \gamma^\mathtt{r}_{x_k} \big) $$
%Since the $\subo^{\mathtt{r},\bw_n}$ are L\'evy processes, Theorem \ref{cvLevy} (in Appendix Section \ref{SkoAppsec}) 
%entails that $\frac{_{_1}}{^{^{b_n}}} \subo^{\mathtt{r}, \bw_n}_{a_n \cdot}\! \rightarrow \! \gamma^\mathtt{r}$, which easily entails (\ref{simplred}). \cqfd 

\medskip

Recall from (\ref{thetabw}) the definition of 
$\theta^{\mathtt{b},\bw}$ and recall from (\ref{schnitzel22}) the definition of $\overline{\gamma}^{\mathtt{r},\bw}$. We next set 
\begin{equation}
\label{progeless}
\forall t\ino [0, \infty), \quad \overline{\theta}^{\mathtt{b},\bw}_t = t+ \overline{\gamma}^{\mathtt{r},\bw}_{A^\bw_t}\; .
\end{equation}
Recall from (\ref{T*wdef}) that $T^*_{\bw}\! = \!  \sup \{ t\ino [0, \infty)\! : A^{\bw}_t \! < \! -I^{\mathtt{r}, \bw}_\infty\}$. Then, note that $ \overline{\theta}^{\mathtt{b},\bw}$ coincides with $\theta^{\mathtt{b},\bw}$ on $[0, T^{*}_{\bw})$. 
%Recall from (\ref{repdiscret}) the definition of 
%$\theta^{\mathtt{b},\bw_n}$ and recall from (\ref{progeless}) the definition of $\overline{\theta}^{\mathtt{b},\bw_n}$.
\begin{lem}
\label{thetight} 
Let $(\alpha, \beta, \kappa, \mathbf{c})$ be as in (\ref{parconing}). Recall from (\ref{repsidefi}) the definition of $\psi$ and assume that (\ref{contH}) holds: namely, $\int^\infty d\lambda / \psi (\lambda) \! <\! \infty$.  
Let $a_n , b_n \! \in \! (0, \infty)$ and  $\bw_n \! \in \! \elldo_f$, $n\ino \bbN$, 
satisfy (\ref{apriori}) and $\mathbf{(C1)}$--$\mathbf{(C3)}$ as in (\ref{unalphcv}) and in (\ref{sig3cvcj}). 
%Let $\alpha\ino \bbR$, $\beta\ino [0, \infty)$,  $\kappa \ino (0, \infty)$ and $\mathbf{c}\ino \elldo_3$ satisfy (\ref{reassume}). Let $a_n , b_n \! \in \! (0, \infty)$ and  $\bw_n \! \in \! \elldo_f$, $n\ino \bbN$, 
%satisfy (\ref{apriori}) and $\mathbf{(C1)}$--$\mathbf{(C3)}$. 
%Recall from (\ref{repdiscret}) the definition of $\theta^{\mathtt{b}, \bw_n}$. 
Then, the laws of the processes 
$(\frac{{1}}{{{b_n}}}\, \overline{\theta}^{\, \mathtt{b},\bw_n}_{b_n t})_{t\in [0, \infty)}$ are tight on $\bD ([0, \infty), \bbR)$. 
\end{lem}
\noi
\textbf{Proof.} To simplify notation we set $R^n_t \! = \! \frac{{1}}{{{b_n}}}\,  \overline{\theta}^{\, \mathtt{b},\bw_n}_{b_n t} -t  = \frac{{1}}{{{b_n}}}\overline{\subo}^{\mathtt{r},\bw_n} (A^{\bw_n}_{b_n t})$; we only need to prove that the $R^n$ are tight on $\bD ([0, \infty), \bbR)$. To that end, we use 
Lemma \ref{tightincrea}. First, observe 
that for all $K, z \ino   (0, \infty)$,   
$$ \bP (R^n_t > K) = \bP \big( \frac{_{_1}}{^{^{b_n}}}\overline{\subo}^{\mathtt{r},\bw_n} (A^{\bw_n} (b_n t)) > K \big)\leq \bP\big(\frac{_{_1}}{^{^{b_n}}}\overline{\subo}^{\mathtt{r},\bw_n}_{a_n z} >K  \big) + 
\bP \big(  \frac{_{_1}}{^{^{a_n}}}A^{\bw_n}_{b_n t} >z) \; .$$
This easily implies that for fixed $t$ the laws of the $R^n_t$ are tight on $[0, \infty)$ since it is the case for the laws  of $\overline{\subo}^{\mathtt{r},\bw_n}_{a_n z}/b_n $ and $A^{\bw_n}_{b_n t}/a_n $ by resp.~Lemma \ref{cvRouge1} and Lemma \ref{Atight}. 

Next, denote by $\ccF_{\! t}$ the $\sigma$-field generated by the r.v.~$N_j^{\bw_n} (s)$ and $\gamma^{ \mathtt{r}, \bw_n}  (A^{\bw_n}_s \! )$ with $s\ino [0, t]$ and $j\! \geq \! 1$;  
note that $N_j^{\bw_n} (t+\cdot)\! -\! N_j^{\bw_n} (t)$ are independent of $\ccF_{\! t}$. 
Fix $\epp \ino (0, \infty)$ and recall from (\ref{taueppdef}) the definition of the times 
$\tau_k^\epp (R^n)$: clearly $b_n\tau_k^\epp (R^n)$ is a $(\ccF_{\! t})$-stopping time. Next, fix $k\ino \bbN$ and set 
%\margmm{corrected}
$$\forall x\ino [0, \infty), \quad 
\mathbf{g}(x)\! =\!  \frac{_{_1}}{^{^{b_n}}}\overline{\subo}^{\mathtt{r},\bw_n} \big(a_n (x+ \frac{_{_1}}{^{^{a_n}}}A^{\bw_n} (b_n \tau^\epp_k (R^n) ))\big) 
-  \frac{_{_1}}{^{^{b_n}}}\overline{\subo}^{\mathtt{r},\bw_n} (A^{\bw_n} (b_n \tau^\epp_k (R^n) )) \,.$$
%Clearly, $\mathbf{g}$ has the same law as $\frac{{1}}{{{b_n}}}\subo^{\mathtt{b},\bw_n}_{a_n \cdot}$; 
Set $\mathbf{u}_\epp \! = \! \inf \{ x\ino [0, \infty)\!  :\!  \mathbf{g} (x) \! >\!  \epp \}$; thus by (\ref{taueppdef}): 
%Then, the definition of $\tau_{k+1}^\epp (R^n)$ in (\ref{taueppdef}) implies: 
$$ \tau_{k+1}^\epp (R^n) = \inf \big\{t> \tau_{k}^\epp (R^n) :  \frac{_{_1}}{^{^{a_n}}}A^{\bw_n} (b_n t) \! -\! \frac{_{_1}}{^{^{a_n}}}A^{\bw_n} (b_n  \tau^\epp_k (R^n) ) > \mathbf{u}_\epp \big\} \,.$$
Fix $z, \eta \ino (0, \infty)$ and set $q_{n, k} (\eta)\! = \!   \bP \big(\tau^\epp_k (R^n)  \! <\!  z \,  ; \,  \tau^\epp_{k+1} (R^n) \! -\! \tau^\epp_k (R^n)  \! \leq \! \eta  \big) $. By (\ref{goretex}) in Lemma \ref{EstDifAw} (applied to the $(\ccF_{\! t})$-stopping time $T\! = \! b_n\tau_k^\epp (R^n)$, to $t_0\! = \! b_n z$, to $t\! = \! b_n \eta $ and to $a \! = \! a_n x$), we get the following.  
\begin{eqnarray}
\label{moufinn}
q_{n, k} (\eta) \!\!\! & \leq & \!\!\!  
\bP \big( b_n\tau^\epp_k (R^n)  < b_n z \,  ;  A^{\bw_n} (b_n\eta + b_n  \tau^\epp_k (R^n) ) \! -\! A^{\bw_n} (b_n  \tau^\epp_k (R^n) ) > a_n \mathbf{u}_\epp  \big)  \nonumber \\
\!\!\! & \leq & \!\!\!   \bP \big( b_n\tau^\epp_k (R^n)  < b_n z \,  ;  A^{\bw_n} (b_n\eta + b_n  \tau^\epp_k (R^n) ) \! -\! A^{\bw_n} (b_n  \tau^\epp_k (R^n) ) > a_n x   \big) + \bP (\mathbf{u}_\epp \leq x)  \nonumber \\
\!\!\! & \leq & \!\!\!   x^{-1} \eta (z + \frac{_{_1}}{^{^2}} \eta) \frac{a_n b_n}{\sigma_1 (\bw_n)}  \frac{b_n \sigma_3 (\bw_n)}{a_n^2\sigma_1 (\bw_n)} + \bP (\mathbf{u}_\epp \leq x)  \nonumber \\
\!\!\! & \leq & \!\!\!   x^{-1} \eta (z +  \frac{_{_1}}{^{^2}}\eta) \frac{a_n b_n}{\sigma_1 (\bw_n)}  \frac{b_n \sigma_3 (\bw_n)}{a_n^2\sigma_1 (\bw_n)} + \bP \big( \mathbf{g} (x) \! \ge\!  \epp  \big)
\end{eqnarray}
%\margmm{corrected}
%by (\ref{goretex}) in Lemma \ref{EstDifAw} applied to the $(\ccF_t)$-stopping time
% $T\! = \! b_n\tau_k^\epp (R^n)$ to $t_0\! = \! b_n z$, $t\! = \! b_n \eta $  and $a \! = \! a_n y$. Thus, 
Denote by $\ccG^o_x$ the sigma algebra generated by the processes $(N_j^{\bw_n})_{j\geq 1}$ and by $\gamma^\mathtt{r, \bw_n}_{y}$, $y\ino [0, x]$ and set $\ccG_x\! = \!\ccG^o_{x+}$. Then, it is easy to see that $\frac{{1}}{{{a_n}}}A^{\bw_n} (b_n  \tau^\epp_k (R^n) )$ is a $(\ccG_{x})$-stopping time. By (\ref{schnitzel3}) in Lemma \ref{minipoub} applied to $T\! =\!  A^{\bw_n} (b_n \tau^\epp_k (R^n) $), we get 
$$ \bP \big( \mathbf{g} (x) \! \ge\!  \epp  \big) \leq  \frac{1\! -\!  \exp \big( \! - \! x a_n 
\psi^{-1}_{\bw_n} \big( \frac{1}{{\epp b_n}} \big) \big)}{1 - e^{-1}} \; . $$  
This, combined with (\ref{moufinn}) and (\ref{eierschwa}) in Lemma \ref{expcovcv},  implies the following:  
$$ \limsup_{n\rightarrow \infty} \sup_{k\in \bbN} q_{n, k} (\eta) \leq x^{-1} \eta (z + \eta) \kappa (\beta + \kappa \sigma_3 (\mathbf{c})) +  \frac{1\! -\!  e^{-x\psi^{-1}(\epp^{-1})}}{1-e^{-1}}  \underset{\eta \rightarrow 0+}{-\!\!\!-\!\!\! \longrightarrow}  \frac{1\! -\!  e^{-x\psi^{-1}(\epp^{-1})}}{1-e^{-1}}\;  \underset{x \rightarrow 0+}{-\!\!\! -\!\!\! \longrightarrow} 0,$$
which completes the proof by Lemma \ref{tightincrea}. \cqfd

\medskip

Recall from (\ref{thetabdef}) the definition of 
$\theta^{\mathtt{b}}$ and recall from (\ref{schnitzel2}) in Lemma \ref{minipoub} 
the definition of $\overline{\gamma}^{\mathtt{r}}$. Then, we define
\begin{equation}
\label{wprogeless}
\forall t\ino [0, \infty), \quad \overline{\theta}^{\mathtt{b}}_t = t+ \overline{\gamma}^{\mathtt{r}}_{A_t}\; .
\end{equation}
Recall from (\ref{T*def}) the definition of $T^*\! = \!  \sup \{ t\ino [0, \infty)\! : A_t \! < \! -I^{\mathtt{r}}_\infty\}$. Then, note that $ \overline{\theta}^{\mathtt{b}}$ coincides with $\theta^{\mathtt{b}}$ on $[0, T^{*})$. 
%and recall from (\ref{progeles}) the definition of $\overline{\theta}^{\mathtt{b}}$.  
%Recall from (\ref{T*wdef}) and (\ref{T*def}) that $T^*_{\bw_n}\! = \!  \sup \{ t\ino [0, \infty)\! : A^{\bw_n}_t \! < \! -I^{\mathtt{r}, \bw_n}_\infty\}$ and $T^*\! = \!  \sup \{ t\ino [0, \infty)\! : A_t \! < \! -I^{\mathtt{r}}_\infty\}$.  
\begin{lem}
\label{thecv} 
Let $(\alpha, \beta, \kappa, \mathbf{c})$ be as in (\ref{parconing}). Recall from (\ref{repsidefi}) the definition of $\psi$ and assume that (\ref{contH}) holds: namely, $\int^\infty d\lambda / \psi (\lambda) \! <\! \infty$.  
Let $a_n , b_n \! \in \! (0, \infty)$ and  $\bw_n \! \in \! \elldo_f$, $n\ino \bbN$, 
satisfy (\ref{apriori}) and $\mathbf{(C1)}$--$\mathbf{(C3)}$ as in (\ref{unalphcv}) and in (\ref{sig3cvcj}). Then,  
\begin{eqnarray}
\label{ultrajoint}
\big(\big(\frac{_{_1}}{^{^{a_n}}} X^{\mathtt{b},\bw_n}_{b_n \cdot}\! , \frac{_{_1}}{^{^{a_n}}} A^{\bw_n}_{b_n \cdot} , & & \!\!\!\!  \!\!\!\!  \!\!  \frac{_{_1}}{^{^{a_n}}} Y^{\bw_n}_{b_n \cdot} \big) , 
 \frac{_{_1}}{^{^{b_n}}}\overline{\theta}^{\mathtt{b},\bw_n}_{b_n \cdot } ,   \frac{_{_1}}{^{^{b_n}}}\overline{\subo}^{\mathtt{r},\bw_n}_{a_n \cdot } , 
 \frac{_{_1}}{^{^{a_n}}} X^{\mathtt{r}, \bw_n}_{b_n \cdot }, -\frac{_{_1}}{^{^{a_n}}}I^{\mathtt{r}, \bw_n}_\infty,  \frac{_{_1}}{^{^{b_n}}} T^{*}_{\bw_n} \big) \\ 
 \!  \underset{n\rightarrow \infty}{ -\!\!\! -\!\!\! \longrightarrow} & & \!\!\!\!  \!\!\!\!  \big( (X^\mathtt{b} , A, Y),  \overline{\theta}^\mathtt{b} ,  \overline{\subo}^{\mathtt{r} } , X^\mathtt{r} , -I^\mathtt{r}_\infty , T^*\big) \nonumber
\end{eqnarray}  
weakly on $\bD ([0, \infty), \bbR^3)\! \times \! (\bD ([0, \infty), \bbR))^3\! \times \! [0, \infty]^2$ 
equipped with the product topology. 
\end{lem}
\noi
\textbf{Proof.} Recall from  (\ref{T*wdef})  the definition of $T^*_{\bw_n}$ and recall from  (\ref{T*def}) the definition of $T^*$. We first prove that $\frac{_{_1}}{^{^{b_n}}} T^{*}_{\bw_n}\! \rightarrow \! T^*$ in law on $[0, \infty]$. To that end, first observe that from the independence between the blue and red processes, we deduce that 
$( \frac{{1}}{{{a_n}}} A^{\bw_n}_{b_n \cdot}, -\frac{{1}}{{{a_n}}}I^{\mathtt{r}, \bw_n}_\infty)\! \rightarrow \! (A, -I^\mathtt{r}_\infty)$ weakly on $ \bD ([0, \infty), \bbR)\! \times \! [0, \infty]$. In the (sub)critical cases $\alpha \ino [0, \infty)$, $- I^\mathtt{r}_\infty\! = \! \infty$. Then, clearly $\frac{_{_1}}{^{^{b_n}}} T^{*}_{\bw_n}\! \rightarrow \! T^*$ in law on $[0, \infty]$. 
We next suppose $\alpha \! < \! 0$; thus $- I^\mathtt{r}_\infty$ is exponentially distributed with parameter $\varrho \! >\! 0$ (that is the largest root of $\psi$); namely 
$- I^\mathtt{r}_\infty$ has a diffuse law which allows to apply Proposition 2.11 in Jacod \& Shiryaev \cite{JaSh02} (Chapter VI, Section 2a p.~341) that discusses continuity properties of specific hitting times; thus, we get that $\frac{_{_1}}{^{^{b_n}}} T^{*}_{\bw_n}\! \rightarrow \! T^*$ in law on $[0, \infty]$. 
%%
%% and that $ -I^\mathtt{r}_\infty$ has a diffuse law (it an exponentiallis independent of $A$ 
%%\margmm{rewritten}
%%\mm{and has continuous sample paths} \ms{with a diffuse law}, 
%%Proposition 2.11 in Jacod \& Shiryaev \cite{JaSh02} (Chapter VI, Section 2a p.~341) that discusses continuity properties of specific hitting times, implies that $\frac{_{_1}}{^{^{b_n}}} T^{*}_{\bw_n}\! \rightarrow \! T^*$ in law on $[0, \infty]$. 

By Lemmas \ref{cvjntYAXb}, \ref{cvRouge1} and \ref{thetight}, the laws of the r.v.~on the left hand side of (\ref{ultrajoint}) are tight on 
$\bD ([0, \infty), \bbR^3)\! \times \! (\bD ([0, \infty), \bbR))^3\! \times \! [0, \infty]^2$; we only need to prove that the joint law of the processes on the right hand side of (\ref{ultrajoint}) is the unique limiting law. 
%\margmm{rewritten}
To this end, we note that by the aforementioned three lemmas, the independence between the red processes and blue ones, as well as the uniqueness of the limit law of $(\frac{_{_1}}{^{^{b_n}}} T^{*}_{\bw_n})$  as implied by Jacod \& Shiryaev's proposition, it suffices to consider the situation where  $(n(p))_{p\in \bbN}$ is an increasing sequence of integers such that 
%\ms{to that end, let $(n(p))_{p\in \bbN}$ be an increasing sequence of integers such that 
\begin{eqnarray}
\label{grassscv}
\Big(\big(\frac{_{1}}{^{{a_{n(p)}}}} X^{\mathtt{b},\bw_{n(p)}}_{b_{n(p)} \cdot}\! , \frac{_{1}}{^{{a_{n(p)}}}} A^{\bw_{n(p)}}_{b_{n(p)} \cdot} , & & \!\!\!\!  \!\!\!\!  \!\!  \frac{_{1}}{^{{a_{n(p)}}}} Y^{\bw_{n(p)}}_{b_{n(p)} \cdot}\big) , \nonumber \\
&& \!\!\!\!  \!\!\!\!  \!\! \!\!\!\!  \!\!\!\!  \!\!  \frac{_{1}}{^{{b_{n(p)}}}}\overline{\theta}^{\mathtt{b},\bw_{n(p)}}_{b_{n(p)} \cdot } ,    \frac{_{1}}{^{{b_{n(p)}}}}\overline{\subo}^{\mathtt{r},\bw_{n(p)}}_{a_{n(p)} \cdot } , 
 \frac{_{1}}{^{{a_{n(p)}}}} X^{\mathtt{r}, \bw_{n(p)}}_{b_{n(p)} \cdot } , -\frac{{_1}}{{^{a_{n(p)}}}}I^{\mathtt{r}, \bw_{n(p)}}_\infty,  \frac{{_1}}{{^{b_{n(p)}}}} T^{*}_{\bw_{n(p)}}\Big) \nonumber \\
  \!  \underset{p\rightarrow \infty}{-\!\!\! -\!\!\! -\!\!\! \longrightarrow} & & \!\!\!\!  \!\!\!\!   \big( (X^\mathtt{b} , A, Y), \theta^\prime ,  \overline{\subo}^{\mathtt{r} } , X^\mathtt{r} ,  -I^\mathtt{r}_\infty , T^* \big), 
\end{eqnarray} 
and then prove that $\theta^\prime\! = \! \overline{\theta}^{_\mathtt{b}}$.
%\ms{the only possible limits for the $\frac{{1}}{{{b_n}}} T^{*}_{\bw_n}$ being $T^*$ according to the above mentioned Proposition 2.11 in Jacod \& Shiryaev \cite{JaSh02}. 
%Actually, we only have to prove that $\theta^\prime\! = \! \overline{\theta}^{_\mathtt{b}}$.} 
Without loss of generality (but with a slight abuse of notation), by Skorokod's representation theorem 
we can assume that (\ref{grassscv}) holds true $\bP$-almost surely. Since $A$ has no fixed time of discontinuity, a.s.~for all $q\ino \bbQ\cap [0,\infty)$, $\Delta A_q\! = \! 0$, 
and thus $ A^{{\bw_{n(p)}}}_{{b_{n(p)}q}} / a_{n(p)}\! \rightarrow \! A_q$. 
Since $\gamma^\mathtt{r}$ has no fixed discontinuity and is independent of $A$, the same properties hold for $\overline{\subo}^\mathtt{r}$. Therefore, a.s.~for all 
$q\ino \bbQ\cap [0,\infty)$, $\Delta  \overline{\gamma}^\mathtt{r}(A_q)\! = \! 0$, which easily entails that $\overline{\subo}^{\mathtt{r},\bw_{n(p)}} (A^{\bw_{n(p)}} (b_{n(p)} q))/ b_{n(p)} \! \rightarrow \! \overline{\gamma}^\mathtt{r}(A_q)$; thus, $\overline{\theta}^{\mathtt{b},\bw_{n(p)}} (b_{n(p)} q)/ b_{n(p)} \! \rightarrow \!  \overline{\theta}^\mathtt{b}_q$ for all $q\ino \bbQ\cap [0,\infty)$ a.s. Therefore, $\theta^\prime\! = \! \overline{\theta}^\mathtt{b}$, which completes the proof. \cqfd 

\begin{lem}
\label{superthe} 
%%Let $\alpha\ino \bbR$, $\beta\ino [0, \infty)$,  $\kappa \ino (0, \infty)$ and $\mathbf{c}\ino \elldo_3$ satisfy (\ref{reassume}). Let $a_n , b_n \! \in \! (0, \infty)$ and  $\bw_n \! \in \! \elldo_f$, $n\ino \bbN$, 
%%satisfy (\ref{apriori}) and $\mathbf{(C1)}$--$\mathbf{(C3)}$. Then 
Let $(\alpha, \beta, \kappa, \mathbf{c})$ be as in (\ref{parconing}). Recall from (\ref{repsidefi}) the definition of $\psi$ and assume that (\ref{contH}) holds: namely, $\int^\infty d\lambda / \psi (\lambda) \! <\! \infty$.  
Let $a_n , b_n \! \in \! (0, \infty)$ and  $\bw_n \! \in \! \elldo_f$, $n\ino \bbN$, 
satisfy (\ref{apriori}) and $\mathbf{(C1)}$--$\mathbf{(C3)}$ as in (\ref{unalphcv}) and in (\ref{sig3cvcj}). Then, 
\begin{eqnarray}
\label{uultrajoint}
 \mathscr{Q}_n (1) \mm{:=}
 \big(\big(\frac{_{_1}}{^{^{a_n}}} X^{\mathtt{b},\bw_n}_{b_n \cdot}\! , \frac{_{_1}}{^{^{a_n}}} A^{\bw_n}_{b_n \cdot} , & & \!\!\!\!  \!\!\!\!  \!\!  \frac{_{_1}}{^{^{a_n}}} Y^{\bw_n}_{b_n \cdot} ,  \frac{_{_1}}{^{^{b_n}}}\overline{\theta}^{\mathtt{b},\bw_n}_{b_n \cdot } \big) , 
  \frac{_{_1}}{^{^{b_n}}}\overline{\subo}^{\mathtt{r},\bw_n}_{a_n \cdot } , 
 \frac{_{_1}}{^{^{a_n}}} X^{\mathtt{r}, \bw_n}_{b_n \cdot }, -\frac{_{_1}}{^{^{a_n}}}I^{\mathtt{r}, \bw_n}_\infty,  \frac{_{_1}}{^{^{b_n}}} T^{*}_{\bw_n} \big) \\
% 
% 
%  \big(\big(\frac{_{_1}}{^{^{a_n}}} X^{\mathtt{b},\bw_n}_{b_n \cdot}\! , \frac{_{_1}}{^{^{a_n}}} A^{\bw_n}_{b_n \cdot} , & & \!\!\!\!  \!\!\!\!  \!\!  \frac{_{_1}}{^{^{a_n}}} Y^{\bw_n}_{b_n \cdot} , 
% \frac{_{_1}}{^{^{b_n}}}\theta^{\mathtt{b},\bw_n}_{b_n \cdot }  \big),   \frac{_{_1}}{^{^{b_n}}}\subo^{\mathtt{r},\bw_n}_{a_n \cdot }, 
% \frac{_{_1}}{^{^{a_n}}} X^{\mathtt{r}, \bw_n}_{b_n \cdot }\big) \\ 
 \!  \underset{n\rightarrow \infty}{ -\!\!\! -\!\!\! \longrightarrow} & & \!\!\!\!  \!\!\!\!  
 \big( (X^\mathtt{b} , A, Y,  \overline{\theta}^\mathtt{b}) ,  \overline{\subo}^{\mathtt{r} } , X^\mathtt{r} , -I^\mathtt{r}_\infty , T^*\big)  \nonumber
\end{eqnarray}  
weakly on $\bD ([0, \infty), \bbR^4)\! \times \! (\bD ([0, \infty), \bbR))^2\! \times \! [0, \infty]^2$ equipped with the product topology. 
\end{lem}
\noi
\textbf{Proof.} Without loss of generality (but with a slight abuse of notation), Skorokod's representation theorem allows to assume that (\ref{ultrajoint}) holds $\bP$-almost surely. To simplify notation, we next set 
$R^n\! =\!  \frac{{1}}{{{a_n}}}  (X^{\mathtt{b},\bw_n}_{b_n\cdot } , A^{\bw_n}_{b_n \cdot},  Y^{\bw_n}_{b_n \cdot })$ and 
$R\! = \! (X^\mathtt{b} , A, Y)$. Let us fix $a\ino (0, \infty)$. 

We consider several cases. We first suppose that $\Delta R_a \! \neq \! 0$. 
By Lemma \ref{jtsko} $(i)$, there is $s_n\! \rightarrow \! a$ such that 
$R^n_{s_n-}\! \rightarrow \! R_{a-}$, $R^n_{s_n}\! \rightarrow \! R_{a}$ and thus $\Delta R^n_{s_n}\! \rightarrow \! \Delta R_{a}$. 
%\margmm{see footnote}

\begin{compactenum}
\item[$-$] Let us suppose more specifically that $\Delta Y_a \! >\! 0$. By definition of $Y$, we get $\Delta X^\mathtt{b}_a\! = \! \Delta Y_a$ and $\Delta A_a \! = \! 0$. Suppose that 
$a\ino [0, T^*]$; by Lemma \ref{trajprop} $(ii)$, we get 
%$\Delta \overline{\theta}^{{\, \mathtt{b}}} (a)\! = \! 
$\Delta \theta^\mathtt{b}_a \! = \!  0$ and thus $\Delta \overline{\theta}^{{\, \mathtt{b}}} (a)\! = \! 0$. 
Note that $\Delta \overline{\theta}^{{\, \mathtt{b}}} (a)\! = \! 0$ for all $a\ino (T^*, \infty)$. 
%Moreover, if $\Delta \overline{\theta}^{{\, \mathtt{b}}} (T^*) \! >\! 0$, then $\Delta A (T^*) \! >\! 0$, and thus $\Delta Y (T^*)\! = \! 0$. 
Consequently, for all $a\ino (0, \infty)$, if $\Delta Y_a \! >\! 0$, then $\Delta \overline{\theta}^{{\, \mathtt{b}}} (a)\! = \! 0$ and 
Lemma \ref{jtsko} $(ii)$ entails $\frac{_{_1}}{^{^{b_n}}} \Delta \overline{\theta}^{\mathtt{b}, \bw_n}( b_n s_n)\! \rightarrow  \!  \Delta \overline{\theta}^\mathtt{b}_a \! = \! 0$. 
%\fmm{It is a bit hard to follow the arguments in this paragraph.: \tt{I hope it is better now}}

\item[$-$] We next consider the case where $\Delta R_a \! \neq \! 0$ but $\Delta Y_a \! =\! 0$; then, by definition of $A$ and $Y$, we get 
 $\Delta X^\mathtt{b}_a \! =\! \Delta A_a \! > \! 0$. Since $\gamma^\mathtt{r}$, and therefore $\overline{\subo}^{\mathtt{r}}$ is independent of $R$, it has a.s.~no jump at the times $A_{a-}$ and $A_{a}$; therefore: $\frac{_{_1}}{^{^{b_n}}}\overline{\subo}^{\mathtt{r},\bw_n} (A^{_{\bw_n}}_{^{b_n s_n-}})\! \rightarrow\!  \overline{\gamma}^{\mathtt{r}} (A_{a-} )$ and $ \frac{_{_1}}{^{^{b_n}}}
\overline{\subo}^{\mathtt{r},\bw_n} (A^{_{\bw_n}}_{^{b_n s_n}})  \rightarrow
 \overline{\gamma}^{\mathtt{r}} (A_{a} )$. 
%$$  \frac{_{_1}}{^{^{b_n}}}\overline{\subo}^{\mathtt{r},\bw_n} (A^{\bw_n}_{b_n s_n-}) \rightarrow \overline{\gamma}^{\mathtt{r}} (A_{a-} ) \quad \textrm{and} \quad  \frac{_{_1}}{^{^{b_n}}}
%\overline{\subo}^{\mathtt{r},\bw_n} (A^{\bw_n}_{b_n s_n})  \rightarrow
% \overline{\gamma}^{\mathtt{r}} (A_{a} )\; .$$
This implies that $\frac{{1}}{{{b_n}}} \Delta \overline{\theta}^{{\mathtt{b}, \bw_n}}( b_n s_n)\! \rightarrow\!    \Delta \overline{\theta}^\mathtt{b}_{^a} \! =\!  \overline{\gamma}^{\mathtt{r}} (A_{a} )\! -\! \overline{\gamma}^{\mathtt{r}} (A_{a-} )$. 

\item[$-$]We finally suppose that $\Delta R_a \! = \! 0$; by Lemma \ref{jtsko} $(i)$, there exists a sequence 
$s^\prime_n\! \rightarrow \! a$ such that $\frac{_{_1}}{^{^{b_n}}} \Delta \overline{\theta}^{\mathtt{b}, \bw_n}( b_n s^\prime_n) \! \rightarrow \! \Delta \overline{\theta}^{\mathtt{b}}_a $. Since, $\Delta R_a \! = \! 0$, Lemma \ref{jtsko} $(ii)$ entails that $\Delta R^n_{s^\prime_n}\! \rightarrow \! \Delta R_a$. 
\end{compactenum}

\noi
Thus, we have proved the following: for all $a\ino (0, \infty)$, there exists a sequence 
$s^{\prime \prime}_n \! \rightarrow \! a$ such that 
$\frac{{1}}{{{b_n}}} \Delta \overline{\theta}^{\mathtt{b}, \bw_n}( b_n s^{\prime \prime}_n )\! \rightarrow\!    \Delta \overline{\theta}^\mathtt{b}_a$ and  $\Delta R^n_{s^{\prime \prime}_n }\! \rightarrow \! \Delta R_a$. Then, by 
Lemma \ref{jtsko} $(iii)$, 
$(R^n, \frac{{1}}{{b_n}}  \overline{\theta}^{\mathtt{b}, \bw_n}( b_n \cdot) )\! \rightarrow\!  (R,\overline{\theta}^\mathtt{b})$ a.s.~on $\bD ([0, \infty), \bbR^4)$, which completes the proof. \cqfd

\bigskip

Recall next that for all $t\ino [0, \infty)$ and all $n\ino \bbN$, 
\begin{equation}
\label{relambda}
\Lambda^{\mathtt{b}, \bw_n}_t \!\!  = \! \inf \big\{ s\ino [0, \infty) \! :   \theta^{\mathtt{b}, \bw_n}_s \! \! >\! t \big\}, \quad  \Lambda^{\mathtt{b}}_t \! = \!  \inf \big\{ s \ino [0, \infty)\! :   \theta^{\mathtt{b}}_s \! > \! t \big\}, 
 \end{equation} 
that $ \Lambda^{\mathtt{r}, \bw_n}_t \! \! = \!  t\! -\! \Lambda^{\mathtt{b}, \bw_n}_t$ and that $\Lambda^{\mathtt{r}}_t \! = \! t \! -\!     \Lambda^{\mathtt{b}}_t$.
\begin{lem}
\label{ynynyati} 
%Let $\alpha\ino \bbR$, $\beta\ino [0, \infty)$,  $\kappa \ino (0, \infty)$ and $\mathbf{c}\ino \elldo_3$ satisfy (\ref{reassume}). Let $a_n , b_n \! \in \! (0, \infty)$ and  $\bw_n \! \in \! \elldo_f$, $n\ino \bbN$, 
%satisfy (\ref{apriori}) and $\mathbf{(C1)}$--$\mathbf{(C3)}$.
Let $(\alpha, \beta, \kappa, \mathbf{c})$ be as in (\ref{parconing}). Recall from (\ref{repsidefi}) the definition of $\psi$ and assume that (\ref{contH}) holds: namely, $\int^\infty d\lambda / \psi (\lambda) \! <\! \infty$.  
Let $a_n , b_n \! \in \! (0, \infty)$ and  $\bw_n \! \in \! \elldo_f$, $n\ino \bbN$, 
satisfy (\ref{apriori}) and $\mathbf{(C1)}$--$\mathbf{(C3)}$ as in (\ref{unalphcv}) and in (\ref{sig3cvcj}). 
Recall the notation $\mathscr{Q}_n (1)$ in (\ref{uultrajoint}). Then, 
the following convergence holds true 
\begin{equation}
\label{megacv}
 \mathscr{Q}_n (2) \mm{: = } \big( \mathscr{Q}_n (1),  \frac{_{_1}}{^{^{b_n}}} \Lambda^{\mathtt{b}, \bw_n}_{b_n \cdot }, \frac{_{_1}}{^{^{b_n}}} \Lambda^{\mathtt{r}, \bw_n}_{b_n \cdot } \big) 
\underset{n\rightarrow \infty}{ -\!\!\! -\!\!\! \longrightarrow}  \big( (X^\mathtt{b} , A, Y,  \overline{\theta}^\mathtt{b} ),   \overline{\subo}^{\mathtt{r} } , X^\mathtt{r} , - I^\mathtt{r}_\infty, T^*, \Lambda^{\mathtt{b}}, \Lambda^{\mathtt{r}} \big) 
\end{equation}
weakly on $\bD ([0, \infty), \bbR^4)\! \times \! (\bD ([0, \infty), \bbR))^2\! \times \! [0, \infty]^2\! \times \! (\bC ([0, \infty), \bbR))^2$ equipped with the product topology. 
\end{lem}
\noi
\textbf{Proof.} Without loss of generality (but with a slight abuse of notation), by Skorokod's representation theorem 
we can assume that the convergence in (\ref{uultrajoint}) holds $\bP$-almost surely. Since 
$\overline{\theta}^{\mathtt{b}}$ (resp.~$\overline{\theta}^{\mathtt{b}, \bw_n}$) is constant on $[T^*, \infty)$ (resp.~on $[T^*_{\bw_n}, \infty)$), we easily derive from  (\ref{uultrajoint}) that 
$\overline{\theta}^{\mathtt{b}, \bw_n}(T^*_{\bw_n}) /b_n \! \rightarrow \! \overline{\theta}^{\mathtt{b}} (T^*)$ a.s.~on $[0, \infty]$. 

Next, we take $t\ino (0, \infty)$ 
distinct from $\overline{\theta}^{\mathtt{b}} (T^*)$. Suppose first that $t \! < \!  \overline{\theta}^{\mathtt{b}} (T^*)$. Then, for all sufficiently large $n$, we get $t \!< \! \overline{\theta}^{\mathtt{b}, \bw_n}(T^*_{\bw_n})/b_n$ and we can write 
$$ \frac{_{_1}}{^{^{b_n}}} \Lambda^{\mathtt{b}, \bw_n}_{b_n t} \!\!  = \! \inf \big\{ s\ino [0, \infty) \! :   \frac{_{_1}}{^{^{b_n}}} \overline{\theta}^{\mathtt{b}, \bw_n}_{b_n s} \! \! >\! t \big\}. %\quad \Lambda^{\mathtt{b}}_t \! = \!  \inf \big\{ s \ino [0, \infty)\! :    \overline{\theta}^{\mathtt{b}}_s \! > \! t \big\},} 
$$
Since $\overline{\theta}^{\mathtt{b}}$ is strictly increasing on $[0, T^*)$, standard arguments entail $\Lambda^{\mathtt{b}, \bw_n} (b_nt )/ b_n   \! \rightarrow \! \Lambda^{\mathtt{b}}_t $. 

Suppose next that $t \! > \!  \overline{\theta}^{\mathtt{b}} (T^*)$, which
%\margmm{rewritten} 
is only meaningful %\ms{only happens with a positive probability} 
in the supercritical cases. Then, for all sufficiently large $n$, we get $t \!> \! \overline{\theta}^{\mathtt{b}, \bw_n}(T^*_{\bw_n})/b_n$ and we can write $\Lambda^{\mathtt{b}, \bw_n}_{b_n t}\! = \! T^*_{\bw_n}$ and $\Lambda^{\mathtt{b}}_{t}\! = \! T^*$. Thus, we get $\Lambda^{\mathtt{b}, \bw_n} (b_nt )/ b_n   \! \rightarrow \! \Lambda^{\mathtt{b}}_t $. 

We have proved that $\Lambda^{\mathtt{b}, \bw_n} (b_nt )/ b_n   \! \rightarrow \! \Lambda^{\mathtt{b}}_t $ for all $t\ino (0, \infty)$ distinct from $\overline{\theta}^{\mathtt{b}} (T^*)$.  Since $\Lambda^{\mathtt{b}}$ is nondecreasing and continuous, a theorem due to Dini implies that 
$\frac{_{_1}}{^{^{b_n}}} \Lambda^{\mathtt{b}, \bw_n}_{b_n \cdot }  \! \rightarrow \! \Lambda^{\mathtt{b}}$ uniformly on all compact subsets; it entails a similar convergence for $ \Lambda^{\mathtt{r}}$, which completes the proof of (\ref{megacv}).   \cqfd

\medskip

Here is one of the key technical point of the proof that relies on the estimates of Lemma \ref{modconXL}. 
\begin{lem}
\label{tthetight} 
%Let $\alpha\ino \bbR$, $\beta\ino [0, \infty)$,  $\kappa \ino (0, \infty)$ and $\mathbf{c}\ino \elldo_3$ satisfy (\ref{reassume}). Let $a_n , b_n \! \in \! (0, \infty)$ and  $\bw_n \! \in \! \elldo_f$, $n\ino \bbN$, 
%satisfy (\ref{apriori}) and $\mathbf{(C1)}$--$\mathbf{(C3)}$. 
Let $(\alpha, \beta, \kappa, \mathbf{c})$ be as in (\ref{parconing}). Recall from (\ref{repsidefi}) the definition of $\psi$ and assume that (\ref{contH}) holds: namely, $\int^\infty d\lambda / \psi (\lambda) \! <\! \infty$.  
Let $a_n , b_n \! \in \! (0, \infty)$ and  $\bw_n \! \in \! \elldo_f$, $n\ino \bbN$, 
satisfy (\ref{apriori}) and $\mathbf{(C1)}$--$\mathbf{(C3)}$ as in (\ref{unalphcv}) and in (\ref{sig3cvcj}). 
Then, the laws of the processes 
$(\frac{_{1}}{^{{a_n}}}X^{\mathtt{b},\bw_n}(\Lambda^{\mathtt{b}, \bw_n} _{b_n t}))_{t\in [0, \infty)}$ and $(\frac{_{1}}{^{{a_n}}}X^{\mathtt{r},\bw_n}(\Lambda^{\mathtt{r}, \bw_n} _{b_n t}))_{t\in [0, \infty)}$ 
are tight on $\bD ([0, \infty), \bbR)$. 
\end{lem}
\noi
\textbf{Proof.} Fix $t\ino [0, \infty)$; then for all $t_0, K \ino (0, \infty)$, note that:  
$$ \bP \Big( \sup_{^{s\in [0, t]}} \frac{_{_{1}}}{^{^{a_n}}}| X^{\mathtt{b},\bw_n}(\Lambda^{\mathtt{b}, \bw_n} _{b_n s}) | >K \Big) \leq   \bP \Big( \sup_{^{s\in [0, t_0]}} \frac{_{_1}}{^{^{a_n}}}| X^{\mathtt{b},\bw_n}_{b_n s} | >K \Big) +  \bP \big(  \frac{_{_1}}{^{^{b_n}}} \Lambda^{\mathtt{b}, \bw_n} _{b_n t} > t_0  \big)  . $$
Then, we deduce from (\ref{megacv}) that 
$$ \lim_{K\rightarrow \infty} \limsup_{n\rightarrow \infty}\bP \Big( \sup_{^{s\in [0, t]}} \frac{_{_{1}}}{^{^{a_n}}}| X^{\mathtt{b},\bw_n}(\Lambda^{\mathtt{b}, \bw_n} _{b_n s}) | >K \Big) \leq   \limsup_{n\rightarrow \infty}\bP \big(  \frac{_{_1}}{^{^{b_n}}} \Lambda^{\mathtt{b}, \bw_n} _{b_n t}  > t_0  \big) \underset{t_0\rightarrow \infty}{ -\!\!\! \longrightarrow}  0 \; .$$
A similar argument shows that 
$ \lim_{K\rightarrow \infty} \limsup_{n\rightarrow \infty}\bP( \sup_{{s\in [0, t]}} 
| X^{\mathtt{r},\bw_n}(\Lambda^{\mathtt{r}, \bw_n} _{b_n s}) | \! >\! a_n K \Big) \! =\!  0$.

Next, recall from (\ref{redblumix}) that a.s.~for all $n\ino \bbN$ and for all 
$t\ino [0, \infty)$
\begin{equation}
\label{rarappeell}
X^{\bw_n}_t \! =\!   X^{\mathtt{b},\bw_n}_{\Lambda^{\mathtt{b}, \bw_n} _{ t}}+
X^{\mathtt{r},\bw_n}_{\Lambda^{\mathtt{r}, \bw_n} _{ t}} \; .
\end{equation} 
Recall from (\ref{modudu}) that for all $y \ino \bD ([0, \infty), \bbR)$, $w_z (y, \eta)$ stands for the $\eta$-c\`adl\`ag modulus of continuity of $y(\cdot) $ on $[0, z]$. Let $z_1, z, z_0, \eta , \epp \ino (0, \infty)$. 
Let us consider first the (sub)critical cases. By (\ref{Xlambcont}) in Lemma \ref{modconXL} $(i)$, we easily get:  
\begin{eqnarray*}
\bP \big(w_{z_1} \big(\frac{_{_1}}{^{^{a_n}}}X^{\mathtt{b},\bw_n}(\Lambda^{\mathtt{b}, \bw_n} _{b_n \cdot}) , \, \eta \big) > \epp\big) & & \!\!\!\!\!\!\!  \leq  \bP \big(w_{z+\eta} \big( \frac{_{_1}}{^{^{a_n}}}X^{\bw_n}_{b_n \cdot} , \eta \big) >\epp/2 \big) 
+
\bP \big(w_{z_0} \big( \frac{_{_1}}{^{^{a_n}}}X^{\mathtt{b}, \bw_n}_{b_n \cdot} , \eta \big) >\epp/2 \big) \\ 
 & &      \!\!\!\!\!\!\!   \!\!\!\!\!\!\!  \!\!\!\!\!\!\! +
\bP \big( \frac{_{_1}}{^{^{b_n}}}\theta^{\mathtt{b}, \bw_n}_{b_n z_0} <z_1 \big) + \bP \big( \frac{_{_1}}{^{^{b_n}}}\theta^{\mathtt{b}, \bw_n}_{b_n z_0} > z \big)\; .
\end{eqnarray*}
By Proposition \ref{Xwfrombr}, $X^{\bw_n}$ has the same law as $X^{\mathtt{b}, \bw_n}$ and $X^{\mathtt{r}, \bw_n}$. Then, by Proposition \ref{cvmarkpro}, the laws of the processes $\frac{_1}{^{a_n}}X^{\bw_n}_{b_n \cdot }$ (or equivalently of $\frac{_1}{^{a_n}} X^{\mathtt{b}, \bw_n}_{b_n \cdot }$) 
are tight on $\bD ([0, \infty), \bbR)$. Consequently, 
\begin{eqnarray}
\label{gorbichou}
 \lim_{\eta \rightarrow 0} \limsup_{n \rightarrow \infty} \bP \big(w_{z_1}  \big(\frac{_{_1}}{^{^{a_n}}}X^{\mathtt{b},\bw_n}(\Lambda^{\mathtt{b}, \bw_n} _{b_n \cdot}) , \, \eta \big) > \epp\big) & & \nonumber  \\ 
\leq \; \limsup_{n \rightarrow \infty} \bP \big( \frac{_{_1}}{^{^{b_n}}}\theta^{\mathtt{b}, \bw_n}_{b_n z_0} \leq z_1 \big)  \!\!\!\!\! & + &\!\!\!\!\!  \limsup_{n \rightarrow \infty} \bP \big( \frac{_{_1}}{^{^{b_n}}}\theta^{\mathtt{b}, \bw_n}_{b_n z_0} \geq  z \big)\; .
% \\
%   \underset{z\rightarrow \infty}{ -\!\!\! \longrightarrow}    \!\!\!\!\! &  &\!\!\!\!\!   \limsup_{n \rightarrow \infty} \bP \big( \frac{_{_1}}{^{^{b_n}}}\Lambda^{\mathtt{b}, \bw_n}_{b_n z_1} >z_0 \big) \underset{z_0\rightarrow \infty}{ -\!\!\! \longrightarrow}  0, 
\end{eqnarray}
Recall that in the (sub)critical cases, $\overline{\theta}^{\mathtt{b}}\! = \! \theta^\mathtt{b}$. Moreover, since $\theta^\mathtt{b}$ has no fixed discontinuity, (\ref{uultrajoint}) easily entails: $\frac{_{_1}}{^{^{b_n}}}\theta^{\mathtt{b}, \bw_n}_{b_n z_0} \! \rightarrow \! \theta^\mathtt{b}_{z_0}$ weakly on $[0, \infty)$. 
It first implies:  
$\limsup_{n \rightarrow \infty} \bP \big( \frac{_{_1}}{^{^{b_n}}}\theta^{\mathtt{b}, \bw_n}_{b_n z_0} \! \geq \!  z \big) \! \leq \! \bP \big(\theta^{\mathtt{b}}_{z_0} \! \geq \!  z \big)\! \rightarrow \! 0$ as $z\! \rightarrow \! \infty$ since a.s.~$\theta^\mathtt{b}_{z_0}\! < \! \infty $ in (sub)critical cases.
Similarly, we also get $\limsup_{n \rightarrow \infty} \bP \big( \frac{_{_1}}{^{^{b_n}}}\theta^{\mathtt{b}, \bw_n}_{b_n z_0} \! \leq \! z_1 \big) \! \leq \! \bP \big( \theta^{\mathtt{b}}_{z_0} \! \leq \! z_1 \big) \! \rightarrow \! 0$ as $z_0\!  \rightarrow \! \infty$ since a.s.~$\lim_{z_0  \rightarrow \infty} \theta^\mathtt{b}_{z_0} \! = \! \infty$.  Then, (\ref{gorbichou}) and the previous arguments imply that the laws of  $(\frac{_{1}}{^{{a_n}}}X^{\mathtt{b},\bw_n}(\Lambda^{\mathtt{b}, \bw_n} _{b_n t}))_{t\in [0, \infty)}$ are tight on $\bD ([0, \infty), \bbR)$ in (sub)critical cases.  

Let us consider the supercritical cases: Lemma  \ref{modconXL} $(ii)$ implies that for all $z_1\ino [0, \infty)$,   
\begin{eqnarray*}
\bP \big(w_{z_1} \big(\frac{_{_1}}{^{^{a_n}}}X^{\mathtt{b},\bw_n}(\Lambda^{\mathtt{b}, \bw_n} _{b_n \cdot}) , \, \eta \big) > \epp\big) & & \!\!\!\!\!\!\!  \leq  \bP \big(w_{z+\eta} \big( \frac{_{_1}}{^{^{a_n}}}X^{\bw_n}_{b_n \cdot} , \eta \big) >\epp/2 \big) 
+
\bP \big(w_{z_0} \big( \frac{_{_1}}{^{^{a_n}}}X^{\mathtt{b}, \bw_n}_{b_n \cdot} ,2 \eta \big) >\epp /6 \big) \\ 
 & &      \!\!\!\!\!\!\!   \!\!\!\!\!\!\!  \!\!\!\!\!\!\! +
\bP \big( \frac{_{_1}}{^{^{b_n}}}T^*_{ \bw_n} \geq z_0 \big) + \bP \big( \frac{_{_1}}{^{^{b_n}}}T^*_{ \bw_n} \leq 2\eta \big) + \bP \big( \frac{_{_1}}{^{^{b_n}}}\theta^{\mathtt{b}, \bw_n}_{T^*_{\bw_n} \! -} \geq  z \big)\; .
\end{eqnarray*}
Then, recall that $\theta^{\mathtt{b}, \bw_n}_{T^*_{\bw_n} \! -} \! \leq \!  \overline{\theta}^{\mathtt{b}, \bw_n}_{T^*_{\bw_n} }$ and observe that (\ref{uultrajoint}) easily entails 
$ \frac{{1}}{{{b_n}}} \overline{\theta}^{\mathtt{b}, \bw_n}_{T^*_{\bw_n} }\! \rightarrow \! \overline{\theta}^{\mathtt{b}}_{T^* }$ weakly on $[0, \infty)$. 
By (\ref{uultrajoint}) again, $\frac{{1}}{{{b_n}}}T^*_{ \bw_n}\! \rightarrow \! T^*$, weakly on $[0, \infty)$. Consequently,
$$ \lim_{\eta \rightarrow 0} \limsup_{n \rightarrow \infty} \bP \big(w_{z_1}  \big(\frac{_{_1}}{^{^{a_n}}}X^{\mathtt{b},\bw_n}(\Lambda^{\mathtt{b}, \bw_n} _{b_n \cdot}) , \, \eta \big) > \epp\big)  \leq  \bP \big(T^*\! \geq z_0 \big) + \bP \big( \overline{\theta}^{\mathtt{b}}_{T^* } \!  \geq  z \big)  \underset{z, z_0\, \rightarrow \infty}{ -\!\!\! \longrightarrow}  0\; ,$$
by the fact that in the supercritical cases $T^*\! <\! \infty $ a.s. 
Thus, the laws of $(\frac{_{1}}{^{{a_n}}}X^{\mathtt{b},\bw_n}(\Lambda^{\mathtt{b}, \bw_n} _{b_n t}))_{t\in [0, \infty)}$ are tight in supercritical cases.  

We derive a similar result for the red processes by a quite similar (but simpler) 
argument based on Lemma \ref{modconXL} $(iii)$: we leave the details to the reader. \cqfd 

\bigskip

Recall (\ref{rarappeell}) and recall from (\ref{Xdef}) that $X_t \! =\!   X^{\mathtt{b}} (\Lambda^{\mathtt{b}} _{ t})+
X^{\mathtt{r}}(\Lambda^{\mathtt{r} } _{ t})$ for all $t\ino [0, \infty)$. 
%
%\begin{equation}
%\label{reXXre}
%X_t \! =\!   X^{\mathtt{b}}_{\Lambda^{\mathtt{b}} _{ t}}+
%X^{\mathtt{r}}_{\Lambda^{\mathtt{r} } _{ t}} \; .
%\end{equation} 
%
\begin{lem}
\label{cvpfff} 
%Let $\alpha\ino \bbR$, $\beta\ino [0, \infty)$,  $\kappa \ino (0, \infty)$ and $\mathbf{c}\ino \elldo_3$ satisfy (\ref{reassume}). Let $a_n , b_n \! \in \! (0, \infty)$ and  $\bw_n \! \in \! \elldo_f$, $n\ino \bbN$, 
%satisfy (\ref{apriori}) and $\mathbf{(C1)}$--$\mathbf{(C3)}$. 
Let $(\alpha, \beta, \kappa, \mathbf{c})$ be as in (\ref{parconing}). Recall from (\ref{repsidefi}) the definition of $\psi$ and assume that (\ref{contH}) holds: namely, $\int^\infty d\lambda / \psi (\lambda) \! <\! \infty$.  
Let $a_n , b_n \! \in \! (0, \infty)$ and  $\bw_n \! \in \! \elldo_f$, $n\ino \bbN$, 
satisfy (\ref{apriori}) and $\mathbf{(C1)}$--$\mathbf{(C3)}$ as in (\ref{unalphcv}) and in (\ref{sig3cvcj}). Recall from (\ref{megacv}) the notation $\mathscr{Q}_n (2)$. Then 
\begin{eqnarray}
\label{gigacv}
\quad \quad\quad  \mathscr{Q}_n (3) :=   \big( \mathscr{Q}_n (2),  \frac{_{_1}}{^{^{a_n}}}  
 \big( X^{\mathtt{b}, \bw_n}_{\Lambda^{\mathtt{b}, \bw_n } _{b_n \cdot }}  \!\!\!\!\!\! &  & \!\!\!\!\!\! , X^{\mathtt{r}, \bw_n}_{\Lambda^{\mathtt{r}, \bw_n }_{b_n \cdot }} , X^{\bw_n }_{b_n \cdot } \big)  \big) \\
& &\!\!\!\!\!\!  \!\!\!\!\!\!  \!\!\!\!\!\!  \underset{n\rightarrow \infty}{ -\!\!\! -\!\!\! \longrightarrow}  \big( (X^\mathtt{b} , A, Y,  \overline{\theta}^\mathtt{b} ), \overline{ \subo}^{\mathtt{r} } , X^\mathtt{r} , -I^\mathtt{r}_\infty, T^*, \Lambda^{\mathtt{b}}, \Lambda^{\mathtt{r}}, ( X^{\mathtt{b}}_{\Lambda^{\mathtt{b}}}, 
X^{\mathtt{r}}_{\Lambda^{\mathtt{r}}} , X)  \big),  \nonumber 
\end{eqnarray}
%\margmm{corrected}
weakly on $\bD ([0, \infty), \bbR^4) \!  \times  \!  (\bD ([0, \infty), \bbR))^2 \!  \times \! [0, \infty]^2 
\!  \times \! (\bC ([0, \infty), \bbR))^2\!  \times   \! \bD ([0, \infty), \bbR^3)$ equipped with the product-topology.
% It implies that $X$ is a spectrally positive L\'evy process with Laplace exponent $\psi$: namely, $X$ has the same law as $X^\mathtt{b}$ and  $X^\mathtt{r}$. 
\end{lem}
\noi
\textbf{Proof.} We first prove the following 
\begin{eqnarray}
\label{ogigacv}
\quad \quad\quad \mathscr{Q}^\prime_n (3) :=   \big( \mathscr{Q}_n (2),  \frac{_{_1}}{^{^{a_n}}}    X^{\mathtt{b}, \bw_n}_{\Lambda^{\mathtt{b}, \bw_n } _{b_n \cdot }}  \!\!\!\! &  & \!\!\!\!\!\! ,  \frac{_{_1}}{^{^{a_n}}}   X^{\mathtt{r}, \bw_n}_{\Lambda^{\mathtt{r}, \bw_n }_{b_n \cdot }}\big) \\
& &\!\!\!\!\!\!  \!\!\!\!\!\!  \!\!\!\!\!\!  \underset{n\rightarrow \infty}{ -\!\!\! -\!\!\! \longrightarrow}  \big( (X^\mathtt{b} , A, Y,  \overline{\theta}^\mathtt{b} ), \overline{ \subo}^{\mathtt{r} } , X^\mathtt{r} , -I^\mathtt{r}_\infty, T^*, \Lambda^{\mathtt{b}}, \Lambda^{\mathtt{r}}, X^{\mathtt{b}}_{\Lambda^{\mathtt{b}}}, 
X^{\mathtt{r}}_{\Lambda^{\mathtt{r}}}   \big),  \nonumber 
%
%& & \underset{n\rightarrow \infty}{ -\!\!\! -\!\!\! \longrightarrow}  \big( (X^\mathtt{b} , A, Y,  \theta^\mathtt{b} ),  \subo^{\mathtt{r} } , X^\mathtt{r} , \Lambda^{\mathtt{b}}, \Lambda^{\mathtt{r}},  X^{\mathtt{b}}_{\Lambda^{\mathtt{b}}}, 
%X^{\mathtt{r}}_{\Lambda^{\mathtt{r}}}  \big),  \nonumber 
\end{eqnarray}
weakly on $\bD ([0, \infty), \bbR^4)\! \times \! \bD ([0, \infty), \bbR)^2\! \times \![0, \infty]^2\! \times \! \bC ([0, \infty), \bbR)^2\! \times \! \bD ([0, \infty), \bbR)^2$ equipped with the product-topology. 
Note that the laws of $\mathscr{Q}^\prime_n (3)$ are tight thanks to (\ref{megacv}) and 
Lemma \ref{tthetight}.
We only need to prove that the joint law of the processes on the right hand side of (\ref{ogigacv}) is the unique limiting law: to that end, let $(n(p))_{p\in \bbN}$ be an increasing sequence of integers such that 
\begin{equation}
\label{ssgigacv}
\mathscr{Q}^\prime_{n(p)} (3)  \underset{p\rightarrow \infty}{ -\!\!\! -\!\!\! \longrightarrow} 
 \big( (X^\mathtt{b} , A, Y,  \overline{\theta}^\mathtt{b} ), \overline{ \subo}^{\mathtt{r} } , X^\mathtt{r} , -I^\mathtt{r}_\infty, T^*, \Lambda^{\mathtt{b}}, \Lambda^{\mathtt{r}}, Q^\mathtt{b} , Q^\mathtt{r}  \big)
%
% \big( (X^\mathtt{b} , A, Y,  \theta^\mathtt{b} ),  \subo^{\mathtt{r} } , X^\mathtt{r} , \Lambda^{\mathtt{b}}, \Lambda^{\mathtt{r}}, Q^\mathtt{b} , Q^\mathtt{r}  \big) 
\end{equation}
weakly on $\bD ([0, \infty), \bbR^4)\! \times \! \bD ([0, \infty), \bbR)^2\! \times \![0, \infty]^2\! \times \! \bC ([0, \infty), \bbR)^2\! \times \! \bD ([0, \infty), \bbR)^2$ equipped with the product topology. Without loss of generality (but with a slight abuse of notation), by Skorokod's representation theorem we can assume that the convergence in (\ref{ssgigacv}) holds $\bP$-a.s.~and we only need to prove that $Q^\mathtt{b}\! = \! X^\mathtt{b} \! \circ \! \Lambda^\mathtt{b}$ and $Q^\mathtt{r}\! = \! X^\mathtt{r} \! \circ \! \Lambda^\mathtt{r}$.

We first prove that $Q^\mathtt{b}\! = \! X^\mathtt{b} \! \circ \! \Lambda^\mathtt{b}$. Note that $\{ t\ino [0,\infty)\! : \! (\Delta X^\mathtt{b}) (\Lambda^{\mathtt{b}}_t) \! >\! 0 \}$ is, in general, \textit{not} countable (it contains all the red intervals starting with a jump), so we have to proceed with care. To that end, 
we first set $S_1 \! = \! \big\{ t\ino [0, \infty): \Delta Y(\Lambda^\mathtt{b}_t) \!> \! 0   \big\}$ that is a 
countable set of times (\textit{indeed}, by Lemma \ref{trajprop} $(ii)$, for all $a\ino[0, T^*_{\bw}]$, 
$\Delta Y_a \! >\! 0$ implies $\Delta \theta^\mathtt{b}_a \! =\! 0$ and by Lemma \ref{trajprop} $(i)$, there exists a unique time $t\ino [0, \infty)$ such that $\Lambda^\mathtt{b}_t \! = \! a$). 
We also set $S_2\! = \! \{ \theta^\mathtt{b}_{T^*-}  \} \cup \{ \theta^\mathtt{b}_{a-} , \theta^\mathtt{b}_{a}; a \ino [0, T^*)\!  : \! \Delta \theta^\mathtt{b}_{a} \! >\! 0\}$ and $S\! = \! S_1 \cup S_2$. 
Then $S$ is countable. We then consider several cases.

We first fix $t\ino (0, T^*) \backslash S$ and we assume that 
$(\Delta X^\mathtt{b}) ( \Lambda^\mathtt{b}_t)\! = \! 0$. Then, by Lemma \ref{jtsko} $(ii)$,  
$X^{\mathtt{b}, \bw_{n(p)}} \big(  \Lambda^{\mathtt{b}, \bw_{n(p)}} (b_{n(p)} t) \big)/ a_{n(p)} \! \rightarrow \!  X^\mathtt{b} ( \Lambda^\mathtt{b}_t)$, since $\Lambda^{\mathtt{b}, \bw_{n(p)}}(b_{n(p)} t)/ b_{n(p)}   \! \rightarrow \! \Lambda^\mathtt{b}_t$. 

We next assume that $t\ino (0, T^*) \backslash S$  and that 
$(\Delta X^\mathtt{b}) ( \Lambda^\mathtt{b}_t)\! > \! 0$. Since $t \! \notin \! S_1$, 
$\Delta Y (\Lambda^\mathtt{b}_t)\! = \! 0$, and thus $\Delta X^\mathtt{b} ( \Lambda^\mathtt{b}_t)\! = \! 
\Delta A ( \Lambda^\mathtt{b}_t)\! > \! 0$, by definition of $A$ and $Y$. We then set $a\! = \! \Lambda^\mathtt{b}_t$ and we necessarily get $a\! < \! T^*$, $\Delta \theta^\mathtt{b}_a \! >0$ and $ t\ino [\theta^{\mathtt{b}}_{a-}, \theta^{\mathtt{b}}_{a}]$. Since $t\! \notin \! S_2$, we then get $ t\ino (\theta^{\mathtt{b}}_{a-}, \theta^{\mathtt{b}}_{a})$. 
To simplify the notation, we set 
$$R^p \! = \! ( \frac{{_1}}{{^{a_{n(p)}}}} X^{\mathtt{b}, \bw_{n(p)}}_{b_{n(p)} \cdot},  \frac{{_1}}{{^{a_{n(p)}}}} A^{ \bw_{n(p)}}_{b_{n(p)} \cdot} \, , \frac{{_1}}{{^{a_{n(p)}}}} Y^{ \bw_{n(p)}}_{b_{n(p)} \cdot} \, , \frac{{_1}}{{^{b_{n(p)}}}}  \overline{\theta}^{\mathtt{b}, \bw_{n(p)}}_{b_{n(p)} \cdot} ) \quad \textrm{and} \quad R= (X^\mathtt{b} , A, Y,  \overline{\theta}^\mathtt{b} ) \; .$$ 
By (\ref{ssgigacv}), $R^p \! \rightarrow \! R$ a.s.~on $\bD ([0, \infty) , \bbR^4)$. Since $a\! <\! T^*$, $\Delta \theta^\mathtt{b}_a\! = \! \Delta \overline{\theta}^\mathtt{b}_a \! >0$ and 
$a$ is a jump-time of $R$. By Lemma \ref{jtsko} $(i)$, there is a sequence $s_p \! \rightarrow \! a$ such that $(R^p_{s_p-} , R^p_{s_p}) \! \rightarrow \!  (R_{a-} , R_{a})$: in particular, we get 
$ X^{\mathtt{b}, \bw_{n(p)}} (b_{n(p)} s_p)/ a_{n(p)}  \! \rightarrow \!   X^\mathtt{b}_a \! = \! 
X^\mathtt{b} (\Lambda^{\mathtt{b}}_t)$. It also implies that 
%$\theta^{\mathtt{b}, \bw_{n(p)}}_{b_{n(p)} s_p-}/ b_{n(p)}\! = \! 
%\overline{\theta}^{\mathtt{b}, \bw_{n(p)}}_{b_{n(p)} s_p-}/ b_{n(p)}\! \rightarrow \! \overline{\theta}^\mathtt{b}_{a-} \! = \! \theta^\mathtt{b}_{a-}$
$\theta^{\mathtt{b}, \bw_{n(p)}} (b_{n(p)} s_p-)/ b_{n(p)}\! = \! \overline{\theta}^{\mathtt{b}, \bw_{n(p)}} (b_{n(p)} s_p-)/ b_{n(p)}\! \rightarrow \! \overline{\theta}^\mathtt{b}_{a-} \! = \! \theta^\mathtt{b}_{a-}$
and
%$\theta^{\mathtt{b}, \bw_{n(p)}} (b_{n(p)} s_p-)/ b_{n(p)}\! \rightarrow \! \theta^\mathtt{b}_{a-} $ and 
$\theta^{\mathtt{b}, \bw_{n(p)}} (b_{n(p)} s_p)/ b_{n(p)}\! = \! \overline{\theta}^{\mathtt{b}, \bw_{n(p)}} (b_{n(p)} s_p)/ b_{n(p)}\! \rightarrow \! \overline{\theta}^\mathtt{b}_{a} \! = \! \theta^\mathtt{b}_{a}$; thus, for all sufficiently large $p$, we get 
$$\frac{{_1}}{{^{b_{n(p)}}}} \theta^{\mathtt{b}, \bw_{n(p)}} (b_{n(p)} s_p-) < t < \frac{{_1}}{{^{b_{n(p)}}}} \theta^{\mathtt{b}, \bw_{n(p)}} (b_{n(p)} s_p)  \quad \textrm{and thus} \quad  \frac{{_1}}{{^{b_{n(p)}}}}  \Lambda^{\mathtt{b}, \bw_{n(p)}}_{b_{n(p)} t}= s_p , $$
which implies that $X^{\mathtt{b}, \bw_{n(p)}} \big( \Lambda^{\mathtt{b}, \bw_{n(p)}} (b_{n(p)} t )\big) / a_{n(p)}  \! \rightarrow \! X^\mathtt{b}_a \! = \! X^\mathtt{b} (\Lambda^{\mathtt{b}}_t)$.

Thus, we have proved a.s.~for all $t\ino (0, T^*) \backslash S$ that 
$X^{\mathtt{b}, \bw_{n(p)}}\!  \big( \Lambda^{\mathtt{b}, \bw_{n(p)}} (b_{n(p)} t )\big) / a_{n(p)}  \! \rightarrow \! X^\mathtt{b} (\Lambda^{\mathtt{b}}_t)$. Since $S$ is countable, it easily implies that 
for all $t\ino [0, T^*)$, $Q^\mathtt{b}_t\! = \! X^\mathtt{b}(\Lambda^\mathtt{b}_t)$. In (sub)critical cases, it simply means that $Q^\mathtt{b}\! = \!  X^\mathtt{b}_{\Lambda^\mathtt{b}}$.

We now complete the proof that $Q^\mathtt{b}\! = \!  X^\mathtt{b}_{\Lambda^\mathtt{b}}$ in the supercritical cases. To that end, we first observe the following. Let $t_1, t_2 \ino (T^*, \infty)$ be distinct times such that $\Delta Q^\mathtt{b}_{t_{1}}\! = \! \Delta Q^\mathtt{b}_{t_{2}}\! = \! 0$. By Lemma \ref{jtsko} $(ii)$,  $X^{\mathtt{b}, \bw_{n(p)}}\!  \big( \Lambda^{\mathtt{b}, \bw_{n(p)}} (b_{n(p)} t_{i} )\big) / a_{n(p)}  \! \rightarrow \!Q^\mathtt{b}_{t_i}$ for $i\ino \{1, 2\}$. Then, by (\ref{ssgigacv}), we get $t_i \! >\! T_{\bw_{n(p)}}^* /b_{n(p)}$ for all sufficiently large $p$ which implies $X^{\mathtt{b}, \bw_{n(p)}}\!  \big( \Lambda^{\mathtt{b}, \bw_{n(p)}} (b_{n(p)} t_{1} )\big)\! = \! X^{\mathtt{b}, \bw_{n(p)}}\!  \big( \Lambda^{\mathtt{b}, \bw_{n(p)}} (b_{n(p)} t_{2} )\big)$. Consequently, we get  
$Q^\mathtt{b}_{t_1}\! = \! Q^\mathtt{b}_{t_2}$. This argument easily implies that for all $t\ino [T^*, \infty)$, $Q^\mathtt{b}_t \! = \! Q^\mathtt{b}_{T^*}$. 
Thus, to 
complete the proof that $Q^\mathtt{b}\! = \! X^\mathtt{b}_{\Lambda^\mathtt{b}}$ in the supercritical cases, we only need to prove that $X^{\mathtt{b}, \bw_{n(p)}} (T^*_{\bw_{n(p)}})/ a_{n(p)} \! \rightarrow \! X^\mathtt{b} (T^*)$. 
If $\Delta X^\mathtt{b}(T^*)\! = \! 0$, then it is a consequence of (\ref{ssgigacv}) and of Lemma \ref{jtsko} $(ii)$.

Therefore, it remains to address cases where 
$\Delta X^\mathtt{b}(T^*)\! >\! 0$. In this case, we clearly get $\Delta \theta^\mathtt{b}(T^*)\! = \! \infty$; by Lemma \ref{trajprop} $(ii)$ with $a\! = \! T^*$, we get 
$\Delta Y(T^*)\! = \! 0$ and therefore $\Delta X^\mathtt{b}(T^*)\!= \! \Delta A(T^*) \! >\! 0$ by definition of $Y$ and $A$.  

We first claim that it is sufficient to prove 
$A^{\bw_{n(p)}}(T^*_{\bw_{n(p)}})/a_{n(p)}\! \rightarrow \! A_{T^*}$. \textit{Indeed}, suppose it holds true;  
since $\Delta Y(T^*)\! = \! 0$, Lemma \ref{jtsko} $(ii)$ and (\ref{ssgigacv}) 
imply that $Y^{\bw_{n(p)}}(T^*_{\bw_{n(p)}})/a_{n(p)}\! \rightarrow \! Y_{T^*}$; and it is sufficient to recall that $X^{\mathtt{b} , \bw_{n(p)}}\! = \! A^{\bw_{n(p)}} + Y^{\bw_{n(p)}}$. 

Thus, we assume that we are in the supercritical cases and that $\Delta X^\mathtt{b}(T^*)\! >\! 0$, and we want to 
prove that $A^{\bw_{n(p)}}(T^*_{\bw_{n(p)}})/a_{n(p)}\! \rightarrow \! A_{T^*}$.
By Lemma \ref{jtsko} $(i)$, there exists $t_p \! \rightarrow \! T^*$ such that 
$A^{\bw_{n(p)}}(b_{n(p)} t_p-)/a_{n(p)}\! \rightarrow \! A_{T^*\! -}$ and 
$A^{\bw_{n(p)}}(b_{n(p)} t_p)/a_{n(p)}\! \rightarrow \! A_{T^*}$. Suppose that $t_p \!> \! T^*_{\bw_{n(p)}}/b_{n(p)}$ for infinitely many $p$; by the definition (\ref{T*wdef}) of $T^*_{\bw_n}$, it implies that $A^{\bw_{n(p)}}(b_{n(p)} t_p-)\! \geq \! -I^{\mathtt{r}, \bw_{n(p)}}_\infty$ for infinitely many $p$ and (\ref{ssgigacv}) implies $A_{T^*-} \! \geq \! -I^{\mathtt{r}}_\infty$; since $T^*\! = \! \sup \{ t\ino [0, \infty): A_t \! < \! -I^{\mathtt{r}}_\infty\}$, we get 
$A_{T^*-} \! = \! -I^{\mathtt{r}}_\infty$; however, $-I^{\mathtt{r}}_\infty$ is an exponentially distributed r.v.~that is independent of $A$ which a.s.~implies that $-I^{\mathtt{r}}_\infty \notin \{ A_{a-}; a\ino (0, \infty) \}$. 
%however, since $-I^{\mathtt{r}}_\infty$ is an exponentially distributed r.v.~that is independent of $A$ and since $T^*\! = \! \sup \{ t\ino [0, \infty): A_t \! < \! -I^{\mathtt{r}}_\infty\}$, we get 
%$\mm{\bP (\Delta A_{T^*} \! >\! 0 ; 
%A_{T^*-} \! \geq \! -I^{\mathtt{r}}_\infty) \le \bP (\Delta A_{T^*} \! >\! 0 ; 
%A_{T^*-} \! = \! -I^{\mathtt{r}}_\infty)=  0}$. 
%\margmm{This needs more details}
This proves that a.s.~$t_p \! \leq \!  T^*_{\bw_{n(p)}}/b_{n(p)}$ for all sufficiently large $p$. 
Then, Lemma \ref{jtsko} $(iv)$ in Appendix implies that $A^{\bw_{n(p)}}(T^*_{\bw_{n(p)}})/a_{n(p)}\! \rightarrow \! A_{T^*}$. 
As observed previously, it completes the proof of $X^{\mathtt{b}, \bw_{n(p)}} (T^*_{\bw_{n(p)}})/ a_{n(p)} \! \rightarrow \! X^\mathtt{b} (T^*)$ and it completes the proof of $Q^\mathtt{b}\! = \! 
X^\mathtt{b}_{\Lambda^{\mathtt{b}}}$ in the supercritical cases.  

\medskip

 We next prove that $Q^\mathtt{r}\! = \! X^\mathtt{r}_{\Lambda^\mathtt{r}}$: to that end, we set $S_3\! = \! \{ t\ino [0, \infty) \! : \! (\Delta X^\mathtt{r}) (\Lambda^\mathtt{r}_t) \! >\! 0 \}$. Lemma \ref{trajprop} $(iv)$ entails that a.s.~$S_3$ is countable and by Lemma \ref{jtsko} $(ii)$, a.s.~for all $t\ino [0, \infty)\backslash S_3$, we get $X^{\mathtt{r}, \bw_{n(p)}} \! \big( \Lambda^{\mathtt{r}, \bw_{n(p)}} (b_{n(p)} t\big)/ a_{n(p)}
\! \rightarrow \! X^\mathtt{r} (\Lambda^\mathtt{r}_t)$; this easily entails that a.s.~$Q^\mathtt{r}\! = \! X^\mathtt{r}  \circ  \Lambda^\mathtt{r}$, which completes the proof of (\ref{ogigacv}).

\medskip

  We now prove (\ref{gigacv}): without loss of generality (but with a slight abuse of notation), Skorokod's representation theorem allows to assume that (\ref{ogigacv}) holds $\bP$-a.s. By Lemma \ref{trajprop} $(v)$, a.s.~for all $t\ino [0, \infty)$, $\Delta Q^\mathtt{b}_t  \Delta Q^\mathtt{r}_t\! = \! 0$, and Lemma \ref{jtsko} $(iii)$ entails: 
$$ \big( \big( \frac{_{_1}}{^{^{a_n}}}  X^{\mathtt{b}, \bw_n}_{\Lambda^{\mathtt{b}, \bw_n } _{b_nt }} , \frac{_{_1}}{^{^{a_n}}}  X^{\mathtt{r}, \bw_n}_{\Lambda^{\mathtt{r}, \bw_n }_{b_n t}} \big)\big)_{t\in [0, \infty)} 
 \underset{n\rightarrow \infty}{ -\!\!\! -\!\!\! \longrightarrow} \big( (Q^\mathtt{b}_t ,Q^\mathtt{r}_t )\big)_{t\in [0, \infty)} \quad \textrm{a.s.~on $\bD([0, \infty), \bbR^2)$.}$$
which implies (\ref{gigacv}) since $X^{\bw_n}_t \! =\!   X^{\mathtt{b},\bw_n}(\Lambda^{\mathtt{b}, \bw_n} _{ t}) +
X^{\mathtt{r},\bw_n}( \Lambda^{\mathtt{r}, \bw_n} _{ t})$ and  $X_t \! =\!   X^{\mathtt{b}} (\Lambda^{\mathtt{b}} _{ t})+
X^{\mathtt{r}}(\Lambda^{\mathtt{r} } _{ t})$. \cqfd 
%
%In particular, we get $\frac{_{_1}}{^{^{a_n}}} X^{\bw_n}_{b_n \cdot} \! \rightarrow \! X$ weakly on $\bD([0, \infty), \bR)$; 
%since $X^{\bw_n}$ has the same law as $X^{\mathtt{r}, \bw_n}$ by Lemma \ref{redbluind} and since 
%$\frac{_{_1}}{^{^{a_n}}} X^{\mathtt{r},\bw_n}_{b_n \cdot} \! \rightarrow \! X^\mathtt{r}$, $X$ has the same law as $X^{\mathtt{r}}$ that is the law of a spectrally positive L\'evy process. Thias completes the proof of the proposition an of Theorem \ref{Xdefthm}. \cqfd 
% 
%\mm{Comment: up to this point, all the previous statements of this section actually hold with \eqref{reassume} (i.e.~Grey's condition) replaced with $X$ having infinite variation.}
\bigskip 

Recall from (\ref{XJHdef}) the definition of the height process $H^{\bw_n}$  associated with $X^{\bw_n}$. 
%
%
%
% and recall from (\ref{JHdef}) 
%the definition of $\cH^{\bw_n}$ the definition of the height process associated with $Y^{\bw_n}$.
%Recall from (\ref{HYenHX}) in Lemma \ref{Hthetalem} that $\cH^{\bw_n} \! =\!  H^{\bw_n} \! \circ  \theta^{\mathtt{b}, \bw_n}$, for all $ Let $\alpha, \beta, \kappa, \mathbf{c}$ satisfy (\ref{contH}) 
Recall from (\ref{approHdef}) the definition of $(H_t)_{t\in [0, \infty)}$, the height process associated with $X$: $H$ is a continuous process and note that (\ref{approHdef}) implies that $H$ is an adapted  measurable functional of $X$. 
%
%Recall next from Proposition \ref{HYdefprop}
%that $\cH\! = \! H \circ \theta^{\mathtt{b}}$ and that $\cH$ is continuous too: $\cH$ is the height process associated with $Y$.  
Then, recall from  (\ref{rmupoissw}) the definition of the offspring distribution $\mu_{\bw_n}$ and denote by $(Z^{\bw_n}_k)_{k\in \bbN}$ a Galton-Watson Markov chain with initial state 
$Z^{\bw_n}_0\! = \! \lfloor a_n \rfloor$ and offspring distribution $\mu_{\bw_n}$; recall from (\ref{scalheight}) Assumption $\textbf{(C4)}$: there exists $\delta \ \in \! (0, \infty) $ such that $\liminf_{n\rightarrow \infty} \bP ( Z^{_{\bw_n}}_{^{\lfloor b_n \delta /a_n \rfloor}} \! = \! 0 ) \! >\! 0$.  
%\begin{equation}
%\label{rscalheight} 
%\mathbf{(C4):} \quad \exists \delta \ \in \! (0, \infty) , \qquad \liminf_{n\rightarrow \infty} \bP \big( Z^{\bw_n}_{\lfloor b_n \delta /a_n \rfloor} = 0 \big) >0 \; .
%\end{equation}
%
%\vspace{-3mm}
%
\begin{lem}
\label{cvaarghh} 
Let $(\alpha, \beta, \kappa, \mathbf{c})$ be as in (\ref{parconing}). Recall from (\ref{repsidefi}) the definition of $\psi$ and assume that (\ref{contH}) holds: namely, $\int^\infty d\lambda / \psi (\lambda) \! <\! \infty$.  
Let $a_n , b_n \! \in \! (0, \infty)$ and  $\bw_n \! \in \! \elldo_f$, $n\ino \bbN$, 
satisfy (\ref{apriori}) and $\mathbf{(C1)}$--$\mathbf{(C4)}$ as in (\ref{unalphcv}), (\ref{sig3cvcj}) and  (\ref{scalheight}).
%%Let $\alpha\ino \bbR$, $\beta\ino [0, \infty)$,  $\kappa \ino (0, \infty)$ and $\mathbf{c}\ino \elldo_3$ satisfy (\ref{reassume}). Let $a_n , b_n \! \in \! (0, \infty)$ and  $\bw_n \! \in \! \elldo_f$, $n\ino \bbN$, 
%%satisfy (\ref{apriori}) and $\mathbf{(C1)}$--$\mathbf{(C4)}$. 
Recall from (\ref{gigacv}) the notation $\mathscr{Q}_n (3)$. Then, 
\begin{multline}
\label{teracv}
\mathscr{Q}_n (4) :=   \big( \mathscr{Q}_n (3),  \frac{_{_{a_n}}}{^{^{b_n}}} H^{\bw_n}_{b_n \cdot }   ,  \frac{_{_{a_n}}}{^{^{b_n}}}H^{\bw_n} \! \! \circ \! \overline{\theta}^{\mathtt{b}, \bw_n}_{b_n \cdot }  \big) \\
 \underset{n\rightarrow \infty}{ -\!\!\! -\!\!\! \longrightarrow}  \big( (X^\mathtt{b} , A, Y,  \overline{\theta}^\mathtt{b} ),  \overline{\subo}^{\mathtt{r} } , X^\mathtt{r} , -I^\mathtt{r}_\infty, T^*, 
 \Lambda^{\mathtt{b}}, \Lambda^{\mathtt{r}}, ( X^{\mathtt{b}}_{\Lambda^{\mathtt{b}}}, 
X^{\mathtt{r}}_{\Lambda^{\mathtt{r}}} , X) , H, H \! \! \circ \!  \overline{\theta}^\mathtt{b}  \big),  
\end{multline}
%
%\begin{multline}
%\label{teracv}
%\mathscr{Q}_n (4) \! = \!  \big( \mathscr{Q}_n (3),  \frac{_{_{a_n}}}{^{^{b_n}}} H^{\bw_n}_{b_n \cdot }   ,  \frac{_{_{a_n}}}{^{^{b_n}}} \cH^{\bw_n}_{b_n \cdot }  \big) \\
% \underset{n\rightarrow \infty}{ -\!\!\! -\!\!\! \longrightarrow}  \big( (X^\mathtt{b} , A, Y,  \theta^\mathtt{b} ),  \subo^{\mathtt{r} } , X^\mathtt{r} , \Lambda^{\mathtt{b}}, \Lambda^{\mathtt{r}}, ( X^{\mathtt{b}}_{\Lambda^{\mathtt{b}}}, 
%X^{\mathtt{r}}_{\Lambda^{\mathtt{r}}} , X) , H, \cH  \big),  
%\end{multline}
%
%\begin{eqnarray}
%\label{teracv}
%\mathscr{Q}_n (4) \! = \!  \big( \mathscr{Q}_n (3),  \frac{_{_{b_n}}}{^{^{a_n}}} H^{\bw_n}_{b_n \cdot }   \!\!\!\!\!\! &  & \!\!\!\!\!\! ,  \frac{_{_{b_n}}}{^{^{a_n}}} \cH^{\bw_n}_{b_n \cdot }  \big) \\
%& & \underset{n\rightarrow \infty}{ -\!\!\! -\!\!\! \longrightarrow}  \big( (X^\mathtt{b} , A, Y,  \theta^\mathtt{b} ),  \subo^{\mathtt{r} } , X^\mathtt{r} , \Lambda^{\mathtt{b}}, \Lambda^{\mathtt{r}}, ( X^{\mathtt{b}}_{\Lambda^{\mathtt{b}}}, 
%X^{\mathtt{r}}_{\Lambda^{\mathtt{r}}} , X) , H, \cH  \big),  
%\end{eqnarray}
weakly on $\bD ([0, \infty), \bbR^4) \! \times \!  (\bD ([0, \infty), \bbR))^2\!  \times \! [0, \infty]^2 \!  \times \!  (\bC ([0, \infty), \bbR))^2 \! \times \! 
\bD ([0, \infty), \bbR^3)\!  \times \!  (\bC ([0, \infty), \bbR))^2 $ 
equipped with the product topology.  
\end{lem}
\noi
\textbf{Proof.} We first prove that 
\begin{multline}
\label{oteracv}
\mathscr{Q}^\prime_n (4) \! = \!  \big( \mathscr{Q}_n (3),  \frac{_{_{a_n}}}{^{^{b_n}}} H^{\bw_n}_{b_n \cdot }  \big)    \\
 \underset{n\rightarrow \infty}{ -\!\!\! -\!\!\! \longrightarrow}  \mathscr{Q}^\prime(4)= 
  \big( (X^\mathtt{b} , A, Y,  \overline{\theta}^\mathtt{b} ),  \overline{\subo}^{\mathtt{r} } , X^\mathtt{r} , -I^\mathtt{r}_\infty, T^*, 
 \Lambda^{\mathtt{b}}, \Lambda^{\mathtt{r}}, ( X^{\mathtt{b}}_{\Lambda^{\mathtt{b}}}, 
X^{\mathtt{r}}_{\Lambda^{\mathtt{r}}} , X) , H \big),   
% \big( (X^\mathtt{b} , A, Y,  \theta^\mathtt{b} ),  \subo^{\mathtt{r} } , X^\mathtt{r} , 
%\Lambda^{\mathtt{b}}, \Lambda^{\mathtt{r}}, ( X^{\mathtt{b}}_{\Lambda^{\mathtt{b}}}, 
%X^{\mathtt{r}}_{\Lambda^{\mathtt{r}}} , X) , H  \big),  
\end{multline}
weakly on the appropriate product-space.
% $\bD ([0, \infty), \bbR^4)\! \times \! \bD ([0, \infty), \bbR)^2\! \times \! \bC ([0, \infty), \bbR)^2\! \times \! \bD ([0, \infty), \bbR^3)\! \times \! \bC ([0, \infty), \bbR) $ equipped with the product-topology.  
By Proposition \ref{HMarkcvprop}, the laws of the processes $ \frac{{{a_n}}}{{{b_n}}} H^{\bw_n}_{b_n \cdot } $ are tight on $\bC ([0, \infty), \bbR) $. Then, 
the laws of $\mathscr{Q}^\prime_n (4)$ are tight thanks to (\ref{gigacv}). We only need to prove that the law of $\mathscr{Q}^\prime(4)$ is the unique limiting law, 
which is an easy consequence of (\ref{gigacv}), of the joint convergence (\ref{jointecon}) in Proposition \ref{HMarkcvprop} and of the fact that $H$ is an adapted measurable deterministic functional of $X$. 

To complete the proof of the lemma, we use a general (deterministic) result on Skorokhod's convergence for the composition 
of functions that is recalled in Theorem \ref{composko} (see Appendix Section \ref{Skogen}). Without loss of generality (but with a slight abuse of notation), Skorokod's representation theorem allows to 
assume that  (\ref{oteracv}) holds $\bP$-a.s.: since 
$ \frac{{{a_n}}}{{{b_n}}} H^{\bw_n}_{b_n \cdot }\! \rightarrow \! H$ a.s.~on $\bC ([0, \infty), \bbR)$, since 
$ \frac{_{{1}}}{^{{b_n}}} \overline{\theta}^{_{\mathtt{b}, \bw_n}}_{^{b_n \cdot} }\! \rightarrow \!  \overline{\theta}^{_{\mathtt{b}}}$ a.s.~on $\bD ([0, \infty), \bbR)$ and  
since $H \! \circ \!  \overline{\theta}^{_{\mathtt{b}}}$ is a.s.~continuous by (\ref{cHfromX}) Theorem \ref{cHdefthm}, Theorem \ref{composko} $(i)$
applies and asserts that $\frac{{{a_n}}}{{{b_n}}}H^{\bw_n}  \circ  \overline{\theta}^{_{\mathtt{b}, \bw_n}}_{^{b_n \cdot} } \! \rightarrow \! H  \circ   \overline{\theta}^{_{\mathtt{b}}}$ in $\bC([0, \infty), \bbR)$, 
%$ \frac{_{{a_n}}}{^{{b_n}}} \cH^{\bw_n}_{b_n \cdot }\! \rightarrow \! \cH$, 
which completes the proof of the proposition. \cqfd 

\bigskip

\noi
\textbf{Proof of Theorem \ref{HYcvth}.} Recall from (\ref{JHdef}) 
the definition of $\cH^{\bw_n}$ (it is the height process associated with $Y^{\bw_n}$). 
Recall from (\ref{HYenHX}) in Lemma \ref{Hthetalem} that $\cH^{\bw_n}_t \! =\!  H^{\bw_n}( \theta^{\mathtt{b}, \bw_n}_t)$ for all $t\ino [0, T^*_{\bw_n})$. Recall from Theorem \ref{cHdefthm} 
the existence and the properties of $\cH$, the height process associated with $Y$; more specifically, recall from (\ref{cHfromX}) that $\cH_t \! = \! H (\theta^\mathtt{b}_t) $ for all $t\ino [0, T^*)$. We first prove that 
%\margmm{corrected}
\begin{equation}
\label{YAcHcvcv}
\mathscr{Q}_n (5)\! :=\!  \big(  \frac{_{_1}}{^{^{a_n}}} Y^{\bw_n}_{b_n \cdot } ,  \, \frac{_{_1}}{^{^{a_n}}} A^{\bw_n}_{b_n \cdot },\,     \frac{_{_{a_n}}}{^{^{b_n}}} \cH^{\bw_n}_{b_n \cdot} \big) 
\;  \underset{n\rightarrow \infty}{-\!\!\! -\!\!\! -\!\!\! \longrightarrow} \; \big( Y, A,  \cH  \big)\! =:\! \mathscr{Q}(5) 
\end{equation} 
weakly on $(\bD ([0, \infty), \bbR))^2 \!  \times \!  \bC ([0, \infty), \bbR)$ equipped with the product-topology. Observe that in (sub)critical cases, it is an immediate consequence of (\ref{teracv}) in Lemma \ref{cvaarghh}. Thus, we only need to focus on the supercritical cases.   

To simplify notation, we denote by $(Y^{{(n)}}, A^{(n)}, \cH^{(n)})$ the rescaled processes on the left member of (\ref{YAcHcvcv}) and we also set $(Y^{(\infty)}, A^{(\infty)}, \cH^{(\infty)})\! = \! (Y, A, \cH)$. We fix $t \ino (0, \infty)$, a bounded continuous function $F\! : \! \bD ([0, \infty), \bbR)^2 \!  \times \!  \bC ([0, \infty), \bbR)\! \rightarrow \! \bbR$, and for all $n\ino \bbN \cup \{ \infty \}$, we set 
$u_n\! = \! \bE \big[ F\big( Y^{_{(n)}}_{^{\cdot \wedge t}}, A^{_{(n)}}_{^{\cdot \wedge t}}, 
\cH^{_{(n)}}_{^{\cdot \wedge t}}  \big) \big]$. 
Clearly, we only need to prove that $u_n \! \rightarrow \! u_\infty$. 
To that end, we introduce for all $K\ino (0, \infty)$, a continuous function $\phi_K \! : \! [0, \infty) \! \rightarrow \! [0, 1]$ such that $\un_{[0, K]} \! \leq \! \phi_K (\cdot) \! \leq \! \un_{[0, K+1]}$ and we set 
$u_n (K)\! = \!  
\bE \big[ F\big( Y^{_{(n)}}_{^{\cdot \wedge t}}, A^{_{(n)}}_{^{\cdot \wedge t}}, 
\cH^{_{(n)}}_{^{\cdot \wedge t}}  \big)
\phi_K \big( A^{_{(n)}}_{^t} \big) \big]$, for all $n\ino \bbN \! \cup \! \{ \infty \}$. We first observe that $0 \leq \! u_n \! -\! u_n (K) \! \leq \! \lVert F \rVert_\infty \bP \big( A^{_{(n)}}_{^t} \! \geq \! K\big)$. Since $A^{(n)}_t\! \rightarrow \! A_t$, standard arguments imply $\limsup_{n\rightarrow \infty} |u_n \! -\! u_n (K) | \leq \lVert F \rVert_\infty \bP \big( A_t \! \geq \! K\big)$. 
%Next recall that $-I^{\mathtt{r}, \bw_n}_\infty/a_n$ (resp.~$-I^{\mathtt{r}}_\infty$) is an exponentially distributed r.v.~independent of $(Y^{(n)}, A^{(n)}, \cH^{(n)})$ (resp.~$(Y, A, \cH)$) and whose parameter is $\varrho^{(n)}\! := \! a_n\varrho_{\bw_n}$ (resp.~$\varrho^{(\infty)}\! :=\! \varrho$). Thus, for all $n\ino \bbN \cup \{ \infty\}$, we easily get 
%\margmm{rewritten} 
Next recall from Theorem \ref{cHdefthm} that $\cH$ is a functional of $(Y, A)$; then recall also that $-I^{\mathtt{r}}_\infty$ (resp.~$-I^{\mathtt{r}, \bw_n}_\infty/a_n$)
is an exponentially distributed r.v.~independent of $(Y, A)$ and thus independent of $(Y, A, \cH)$ (resp.~independent of $(Y^{(n)}, A^{(n)}, \cH^{(n)})$) and whose parameter is $\varrho^{(\infty)}\! :=\! \varrho$ (resp.~$\varrho^{(n)}\! := \! a_n\varrho_{\bw_n}$). We set  $\overline{\cH}^{_{(n)}}_{^{\cdot} }\! = \! \frac{{_{a_n}}}{{^{b_n}}}H^{\bw_n} \circ \overline{\theta}^{_{\mathtt{b}, \bw_n}}_{^{b_n \cdot }} $. Then, for all $n\ino \bbN \cup \{ \infty\}$, we easily get  
\begin{equation}
\label{flobichoune} 
u_n (K) = \bE \Big[ e^{\varrho^{_{(n)}}\!  A^{_{(n)}}_{^t}}\!  
F\big( Y^{_{(n)}}_{^{\cdot \wedge t}}, A^{_{(n)}}_{^{\cdot \wedge t}}, \overline{\cH}^{_{(n)}}_{^{\cdot \wedge t}} \big) 
\phi_K \big( A^{_{(n)}}_{^t} \big) \un_{\!  \big\{  A^{{(n)}}_{t} < -I^{\mathtt{r}, \bw_n}_\infty/a_n \big\} }\Big] , 
\end{equation}
%\margmm{new}
where the right-hand side is finite thanks to the term $\phi_{K}$. 
% and 
%$T^{(n)}\! = \! T^*_{\bw_n}/ b_n$ and where, by convenience wa have also set $\overline{\cH}^{_{(\infty)}}\! = \! \cH$ and $T^{(\infty)}\! = \! T^*$. 
\textit{Indeed}, the events $\{  A^{{(n)}}_{t} \! <\!  -I^{\mathtt{r}, \bw_n}_\infty/a_n \big\} $ and $\{ T^*_{\bw_n}/ b_n >t   \}$ coincide a.s.~and on these events, we get $\overline{\cH}^{_{(n)}}_{^t }\! = \! \cH^{_{(n)}}_{^t}$. 

Next, recall from Lemma \ref{expcovcv} that $\lim_{n\rightarrow \infty} \varrho^{(n)}\! = \! \varrho^{(\infty)}$. 
Since $\bP (A_t \! = \! -I^\mathtt{r}_\infty)\! = \! 0$, the joint convergence 
(\ref{teracv}) in Lemma \ref{cvaarghh} combined with (\ref{flobichoune}) entails that $u_n(K) \! \rightarrow \! u_\infty(K)$. Since 
$|u_\infty \! -\! u_n| \! \leq \! |u_\infty \! -\! u_\infty (K)|+|u_\infty (K)\! -\! u_n(K)|+|u_n(K) \! -\! u_n|$, 
we get 
$ \limsup_{n\rightarrow \infty} |u_\infty \! -\! u_n|\! \leq \! 2 \lVert F\rVert_\infty \bP \big( A_t \! \geq \! K\big) \rightarrow 0 $ as $K$ tends to $\infty$. This completes the proof of (\ref{YAcHcvcv}) in supercritical cases.

\medskip

To complete the proof of Theorem \ref{HYcvth}, it remains to prove the convergence of the sequences of 
pairs of pinching times $\Ptt_{\bw_n}$. By Skorokod's representation theorem (but with a slight abuse of notation) we can assume without loss of generality 
that (\ref{YAcHcvcv}) holds almost surely: namely, a.s.~$\mathscr{Q}_n(5) \! \rightarrow \! \mathscr{Q} (5)$. Then, we couple the $\Ptt_{\bw_n}$ and $\Ptt_{\bw}$ as follows. 
\begin{compactenum}

\smallskip

\item[$-$] Let $\cR\! = \! \sum_{i\in \mathtt{I}} \delta_{(t_i, r_i, u_i)}$ be a Poisson point measure on $[0, \infty)^3$ with intensity the Lebesgue measure $dt dr dv$ on $[0, \infty)^3$. We assume that $\cR$ is independent of $\mathscr{Q}(5)$ and of $(\mathscr{Q}_n(5))_{n\in \bbN}$. 

\smallskip

\item[$-$] We set $\kappa_n \! = \! a_n b_n / \sigma_1 (\bw_n)$ and for all $t\ino [0, \infty)$ we set $\mathtt{Z}^n_t \! = \! \frac{_1}{^{a_n}}(Y^{\bw_n}_{b_n t} \! -\! J^{\bw_n}_{b_n t})$, where we recall that $J^{\bw_n}_{b_nt}\! =\!  \inf_{s\in [0, b_n t]} Y^{\bw_n}_s$. We then set $S_{n} \! =\!  \{ (t, r, v) \ino [0, \infty)^3\! : \! 0 \! < \! r \! < \! \mathtt{Z}^n_t \; \textrm{and} \; 0\! \leq \! v \! \leq \kappa_n \} $ and we define 
$ \cP_{\! n} = \sum_{i\in \mathtt{I}} \un_{\{ (t_i, r_i, u_i)\in S_n\}} \delta_{(t_i, r_i, u_i)}=: \sum_{1\leq p < \mathbf{p}_n } \delta_{(t^n_p, r^n_p, v^n_p)}$, 
%$$ \cP_{\! n} = \sum_{i\in \mathtt{I}} \un_{\{ (t_i, r_i, u_i)\in S_n\}} \delta_{(t_i, r_i, u_i)}
%= \sum_{1\leq p < \mathbf{p}_n } \delta_{(t^n_p, r^n_p, v^n_p)}, $$
where the indexation is such that the finite sequence $(t^n_p)_{1\leq p<\mathbf{p}_n}$ increases. Note that since 
$\mathtt{Z}^n$ is eventually null, $\cP_{\! n}$ is a finite point process. 

\smallskip

\item[$-$] For all $t \ino [0, \infty)$, for all $r\ino \bbR$ and for all $z\ino \bD ([0, \infty), \bbR)$,  we set 
\begin{equation}
\label{pretautry}
\tau (z, t,r) \! =\!  \inf \big\{ s\ino [0, t] : \inf_{u\in [s, t]} z(u) \!  >\!  r \big\}\;  \textrm{with the the convention that $\inf \emptyset \! = \! \infty $.}
\end{equation}
Then, we set 
\begin{equation}
\label{coucoupl}
\frac{_1}{^{b_n}} \Ptt_{\bw_n}= \big( (s^n_p , t^n_p) \big)_{1\leq p<\mathbf{p}_n } \; \, \textrm{where} \quad s^p_n = \tau (\mathtt{Z}^n, t^n_p , r^n_p), \; 1\leq p < \mathbf{p}_n. 
\end{equation} 
\end{compactenum}
Thanks to (\ref{Poissurpl}) and (\ref{pin1}), we see that conditionally given $Y^{\bw_n}\! $, $\frac{_1}{^{b_n}} \Ptt_{\bw_n}$ has the right law. 
By convenience, we set $ (s^n_p, t^n_p)\! = \! (-1, -1)$, for all $p\! \geq  \! \mathbf{p}_n$. 

  Similarly, we set $\mathtt{Z}^\infty_t \! = \! Y_t  \! -\!  J_t $, where $J_t \! = \!  \inf_{s\in [0, t]} Y_s$ and we also set 
$S\! = \{ (t, r, v) \ino [0, \infty)^3\! :\!  0 \! < \! r \! < \! \mathtt{Z}^\infty_t \; \textrm{and} \; 0\! \leq \! v \! \leq \kappa \} $; we then define $\cP = \sum_{i\in \mathtt{I}} \un_{\{ (t_i, r_i, u_i)\in S\}} \delta_{(t_i, r_i, u_i)}
=: \sum_{p\geq 1 } \delta_{(t_p, r^\prime_p, v_p)}$, 
%$$ \cP = \sum_{i\in \mathtt{I}} \un_{\{ (t_i, r_i, u_i)\in S\}} \delta_{(t_i, r_i, u_i)}
%= \sum_{p\geq 1 } \delta_{(t_p, r^\prime_p, v_p)}, $$
where the indexation is such that $(t_p)_{p\geq1}$ increases. Then, set 
\begin{equation}
\label{coucucuplpl}
 \Ptt= \big( ( s_p , t_p) \big)_{p\geq 1 } \; \, \textrm{where} \quad s_p = \tau (\mathtt{Z}^\infty, t_p , r^\prime_p), \; p \geq 1, 
\end{equation} 
It is easy to check that $\Ptt$ has the right law conditionally given $Y$. 

First observe that $\kappa_n \! \rightarrow \! \kappa >0$, by the last point of (\ref{apriori}). 
Next, we prove that $\mathtt{Z}^n \! \rightarrow \! \mathtt{Z}^\infty$ a.s.~in $\bD([0, \infty), \bbR)$: \textit{indeed}, since $Y$ has no negative jumps, $J$ is continuous and by 
Lemma \ref{franchtime} $(ii)$, $ (\frac{_1}{^{a_n}}J^{\bw_n}_{b_n t} )_{t\in [0, \infty)} \! \rightarrow \! (J_t)_{t\in [0, \infty)}$ a.s.~in $\bC([0, \infty), \bbR)$. Since $J$ is continuous, $Y$ and $J$ do not share any jump-times and by Lemma \ref{jtsko} $(iii)$, $(\frac{_1}{^{a_n}}( Y^{\bw_n}_{b_n t}, J^{\bw_n}_{b_n t}) )_{t\in [0, \infty)} \! \rightarrow \! ((Y_t, J_t))_{t\in [0, \infty)}$ a.s.~in $\bD([0, \infty), \bbR^2)$, 
which entails that $\mathtt{Z}^n \! \rightarrow \! \mathtt{Z}^\infty$ a.s.~in $\bD([0, \infty), \bbR)$. 

Let us fix $a, b , c\ino (0, \infty)$ such that 
$$b \! > \; 2 \!\!\! \sup_{n\in \bbN \cup \{ \infty \}} \sup_{s\in [0, a]} \mathtt{Z}^n_s \quad \textrm{and} \quad 
c\! > \;   2 \!\!\!\! \sup_{n\in \bbN \cup \{ \infty \}} \kappa_n \; . $$
We introduce $\sum_{1\leq l \leq  N}\delta_{(t^*_l, r^*_l, u_l^*)}\! := \! \sum_{i\in \mathtt{I}} \un_{\{ t_i <a \, ; \, r_i <b \, ; \, u_i <c \}}\delta_{(t_i, r_i, u_i)}$, where $(t^*_l)_{1\leq l\leq N}$ increases; here, $N$ is a Poisson r.v.~with mean $abc$; note that conditionally given $N$, the law of the r.v.~$(t^*_l, r^*_l, u_l^*)$ is absolutely continuous with respect to Lebesgue measure. Therefore, 
a.s.~for all $l \! \in \! \{ 1, \ldots, N\}$ (if any), 
$\Delta \mathtt{Z}^\infty_{t^*_l}\! = \! 0$, $u^*_l \! \neq \!\kappa$, and $r^*_l \! \neq  \! \mathtt{Z}^\infty_{t^*_l}$, and if 
$r^*_l \! < \! \mathtt{Z}^\infty_{t^*_l}$, then we get $\tau (\mathtt{Z}^\infty, t^*_l , r_l^* -)\! = \! \tau (\mathtt{Z}^\infty, t^*_l , r_l^* )$ 
because by Lemma \ref{franchtime} $(iv)$, the function $r \! \mapsto \!  \tau( \mathtt{Z}^\infty,  t^*_l , r)$ is right-continuous and it has therefore a 
countable number of discontinuities. Since $\Delta \mathtt{Z}^\infty_{t^*_l}\! = \! 0$, Lemma \ref{jtsko} $(ii)$ entails that $\mathtt{Z}^n_{t^*_l} \! 
\rightarrow \! \mathtt{Z}^\infty_{t^*_l}$, and for all sufficiently large $n$, 
$u^*_l \! \neq \! \kappa_n$ and $r^*_l \! \neq  \! \mathtt{Z}^n_{t^*_l}$, and if  $r^*_l \! <   \! \mathtt{Z}^n_{t^*_l}$, then Lemma \ref{franchtime} $(iv)$ shows that 
$\tau (\mathtt{Z}^n, t^*_l, r^*_l) \! \rightarrow \! \tau ( \mathtt{Z}^\infty, t^*_l, r^*_l)$. This proves that if $t_p \! <\! a$, then $(s^n_p, t^n_p)\! \rightarrow \! (s_p, t_p)$, 
%\margmm{new}
since we have $t^n_p=t_{p}$ for $n$ sufficiently large as a consequence of the above coupling. Since $a$ can be arbitrarily large, we get 
$\frac{_1}{^{b_n}}\Ptt_{\bw_n}\! \rightarrow \! \Ptt$ a.s.~in  $(\bbR^2)^{\bbN^*}\!\! $ equipped with the product topology. This, combined with the a.s.~convergence $\mathscr{Q}_n(5) \! \rightarrow \! \mathscr{Q} (5)$, entails Theorem \ref{HYcvth}. \cqfd

\subsection{Proof of Theorem \ref{excucvth} and proof of Theorem \ref{graphcvth}}
\label{excuproof}
%{\red NEED TO RECALL A-L'S RESULT ON EXCURSIONS OF $Y$. }
Recall from (\ref{AetYdef}) the definition of $Y$ and recall from Theorem \ref{cHdefthm} 
the existence and the properties of $\cH$, the height process associated with $Y$; recall the notation $J_t \! = \! \inf_{s\in [0,t]} Y_s$, $t\ino [0, \infty)$. 
Lemma \ref{AHeuer} $(ii)$ asserts that the excursions of $\cH$ above $0$ and those of $Y\! -\! J$ above $0$ are the same. As recalled in Proposition \ref{AldLim1}, Proposition 14 in Aldous \& Limic \cite{AlLi98} asserts that 
these excursions can be indexed in the decreasing order of their lengths. Namely,  
\begin{equation}
\label{reexcuHY}
\big\{  t\ino [0, \infty) : \cH_t >0  \big\}= \big\{  t\ino [0, \infty) : Y_t > \JJ_t  \big\}= \bigcup_{k\geq 1} (l_k , r_k) \; , 
\end{equation}
where the sequence $\zeta_k\! = \! l_k - r_k$, $k\! \geq \! 1$, decreases. 
Moreover, the sequence $(\zeta_k)_{k\geq 1}$ appears as the law of a version of the multiplicative coalescent at a fixed time: see Theorem 2 in  Aldous \& Limic \cite{AlLi98} (recalled in  Proposition \ref{AldLim2}). In particular, 
it implies that a.s.~$\sum_{k\geq 1} \zeta_k^2 < \infty$. 
Recall from (\ref{exccont}) the definition of \textit{excursion processes} of $\cH$ and $Y\! -\! J$ above $0$: 
\begin{equation}
\label{rexccont}
 \forall k \! \geq \! 1, \; \forall t\ino [0, \infty), \qquad 
\Htt_{k}(t)= \cH_{(l_k + t)\wedge r_k} \quad \textrm{and} \quad  \Ytt_{k}(t)= Y_{(l_k + t)\wedge r_k}- \JJ_{l_k} . 
\end{equation}
Next recall from (\ref{Poisurcon}) and (\ref{pinchset}) 
the definition of $\Ptt\! = \! \big( (s_p, t_p)\big)_{p\geq 1}$. 
%
%
%: namely, 
%conditionally given $Y$, let 
%\begin{equation}
%\label{rPoisurcon}
%\cP\! = \! \sum_{\; p\geq 1} \!  \delta_{(t_p, y_p) }
% \; \textrm{be a Poisson pt.~meas.~on $[0, \infty)^2$ 
%  with intensity $ \kappa \un_{\{ 0< y <Y_t -J_t \}} \, dt\, dy$}  
%\end{equation}
%and set 
%\begin{equation}
%\label{rpinchset}
%\Ptt\! = \! \big( (s_p, t_p) \big)_{p\geq 1} \quad \textrm{where} \quad s_p\! = \! \inf \big\{ s\ino [0, t_p] :  \inf_{u\in [s, t_p]} 
%Y_u\! -\! J_u > y_p   \big\}, \; \,  p\! \geq \! 1. 
%\end{equation}
%

 Let $a_n , b_n \! \in \! (0, \infty)$ and  $\bw_n \! \in \! \elldo_f$, $n\ino \bbN$, 
satisfy (\ref{apriori}) and $\mathbf{(C1)}$--$\mathbf{(C4)}$ as in (\ref{unalphcv}), (\ref{sig3cvcj}) and  (\ref{scalheight}). 
%Let $a_n , b_n \! \in \! (0, \infty)$, and  $\bw_n \! \in \! \elldo_f$, $n\ino \bbN$, satisfy (\ref{apriori}) and $\mathbf{(C1)}$--$\mathbf{(C4)}$. 
Recall from (\ref{YwAwSigw}) the definition of $Y^{\bw_n}$; recall from (\ref{JHdef}) the definition of $\cH^{\bw_n}$, the height process associated to $Y^{\bw_n}$. Recall from (\ref{Poissurpl}) and (\ref{pin1}) the definition of $\Ptt_{\bw_n}$.  
For all $t\ino [0, \infty)$, to simplify notations, we introduce the following: 
\begin{equation}
\label{shortnota}
Y^{_{(n)}}_{t} \! \! :=\! \frac{_{_1}}{^{^{a_n}}} Y^{\bw_n}_{b_n t}, \;  J^{_{(n)}}_{t} \! := \! \inf_{s\in [0, t]} 
Y^{_{(n)}}_{s},    \cH^{_{(n)}}_{t}\! := \! \frac{_{_{a_n}}}{^{^{b_n}}} \cH^{\bw_n}_{b_n t} , \;  \Ptt^{(n)}\! := \! \frac{_1}{^{b_n}} \Ptt_{\bw_n} \! =:\!  \big( (s^n_p , t^n_p) \big)_{1\leq p<\mathbf{p}_{n} }.  
\end{equation}
%\begin{multline}
%\label{shortnotai}
%Y^{_{(n)}}_{t} \! :=\! \frac{_{_1}}{^{^{a_n}}} Y^{\bw_n}_{b_n t}, \quad J^{_{(n)}}_{t} \! := \! \inf_{s\in [0, t]} 
%Y^{_{(n)}}_{s}, \quad   \cH^{_{(n)}}_{t}\! := \! \frac{_{_{a_n}}}{^{^{b_n}}} \cH^{\bw_n}_{b_n t} \\ 
%\quad \textrm{and} \quad \Ptt^{(n)}\! := \! \frac{_1}{^{b_n}} \Ptt_{\bw_n} \! =:\!  \big( (s^n_p , t^n_p) \big)_{1\leq p<\mathbf{p}_{n} }.  
%\end{multline}
Recall from (\ref{excdibor}) that 
%%the excursion intervals of $Y^{(n)}\! -\! J^{(n)}$ above $0$ are the same 
%as the excursions intervals of $\cH^{(n)}$ above $0$; let  $\mathbf{q}_{\bw_n}$ \mm{stand} for the number of such intervals. Namely, 
\begin{equation}
\label{excunormi} \big\{  t\ino [0, \infty) : \cH^{_{(n)}}_t \! >\! 0  \big\}= \big\{  t\ino [0, \infty) : Y^{_{(n)}}_t\! >\!  J^{_{(n)}}_t \big\} = \bigcup_{1\leq k\leq \mathbf{q}_{\bw_n}} [l^n_k, r^n_k)
\end{equation} 
where the indexation is such that the $\zeta^n_k\! := \! r^n_k \! -\! l^n_k$ are nonincreasing 
%the indexation is such that $\zeta^n_1 \! \geq \! \ldots \! \geq \! \zeta^n_{\mathbf{q}_n}$ and 
and such that $l^n_k \!  < \! l^n_{k+1}$ if $\zeta^n_k\! = \! \zeta^n_{k+1}$ (within the notation of (\ref{excdibor}), 
$l^n_k\! = \! l^{\bw_n}_k/b_n$, $r^n_k\! = \! r^{\bw_n}_k/b_n$ and $\zeta^n_k\! = \! \zeta^{\bw_n}_k/b_n$). 

\subsubsection{Proof of Theorem \ref{excucvth}.}
\label{excucvpf}
%We keep the previous notation. By Theorem \ref{HYcvth}, we get 
%\begin{equation}
%\label{reYH}
%\big( Y^{(n)} , \cH^{(n)}, \Ptt^{(n)}\big)
%\;  \underset{n\rightarrow \infty}{-\!\!\! -\!\!\! -\!\!\! \longrightarrow} \; \big(  Y, \cH, \Ptt \big)
%\end{equation} 
%weakly on $\bD([0, \infty), \bbR)Â \times  \bC ([0, \infty) , \bbR) \times (\bbR^2)^{\bbN^*}$, equipped with product topology (recall that here 
%we use the following convention: the finite sequence $\Ptt^{(n)}\! = \!  \big( (s^n_p , t^n_p) \big)_{1\leq p<\mathbf{p}_{n} }$ is extended by setting $ (s^n_p , t^n_p)\! = \! (-1, -1)$ for all $p\! >\! \bp_n$).  
By Skorokod's representation theorem (but with a slight abuse of notation) \textit{we can assume without loss of generality 
that (\ref{jtcvYHP}) in Theorem \ref{HYcvth} holds $\bP$-a.s.} We first prove the following lemma. 
\begin{lem}
\label{convexc1} 
%Assume that (\ref{reYH}) holds a.s.~on $\bD([0, \infty), \bbR)Â \! \times \! \bC ([0, \infty) , \bbR)$ equipped with product topology. 
We keep the previous notations and we assume that (\ref{jtcvYHP}) in Theorem \ref{HYcvth} holds $\bP$-a.s. Then, for all $k, n\! \geq 1$, there exists a sequence $j(n,k)\ino \{1, \ldots , \mathbf{q}_{\bw_n}\}$ such that 
\begin{equation}  
\label{indexouy}
\textrm{$\bP$-a.s.~for all $k\! \geq \! 1$,} \qquad \big( l_{j(n,k)}^{n} ,  r_{j(n,k)}^{n}\big)
\;  \underset{n\rightarrow \infty}{-\!\!\! -\!\!\! -\!\!\! \longrightarrow} \; \big(  l_k, r_k \big). 
\end{equation}
\end{lem}
\noi
\textbf{Proof.} Fix $k\! \geq \! 1$ and let $t_0 \ino (l_k, r_k)$; note that $l_k\! = \! \sup \{ t\ino [0, t_0]\! :\!  \cH_t \! =\!  0\}$ and $r_k \! = \! \inf\{ t\ino [t_0, \infty) \! :\!  \cH_t \! = \! 0 \}$. For all $n\! \geq \! 1$, set $\gamma (n)
\! = \! \sup \{ t\ino [0, t_0)\! :\!  \cH^{_{(n)}}_t \! =\!  0\}$ and $\delta (n) \! = \! \inf\{ t\ino [t_0, \infty)\!  :\!  \cH^{_{(n)}}_t \! = \! 0 \}$. 
Let $q$ and $r$ be such that $l_k\! < \! q \! < \!  t_0 \! < \! r \! < \! r_k$. Since $\inf_{t\in [q, r]} \cH_t \! >\! 0$, for all sufficiently large $n$, we get $\inf_{t\in [q, r]} \cH^{_{(n)}}_t \! >\! 0$, which implies that $\gamma (n) \! \leq \! q $ and $r \! \leq \! \delta (n)$. This easily implies that $\limsup_{n\rightarrow \infty} \gamma (n) \! \leq \! l_k$ and $r_k \! \leq \! \liminf_{n\rightarrow \infty} \delta (n)$.  

 Let $q$ and $r$ be such that $q \! < \! l_k$ and $r_k \! < \! r$. Since $\cH_{l_k}\! =\!  \cH_{r_k}\! = \! 0$, (\ref{Jtpsloczer}) in Lemma \ref{XYexcu} $(iv)$ implies that $J_{q}\! >\! J_{t_0} \! >\! J_{r}$.  
Since $J$ is continuous, Lemma \ref{jtsko} $(iii)$ entails that $J^{(n)} \! \rightarrow \! J$ a.s.~in $\bC ([0, \infty) , \bbR)$. Thus, for all sufficiently large $n$, $J^{_{(n)}}_{q}\! >\! J^{_{(n)}}_{t_0} \! >\! J^{_{(n)}}_{r}$; by definition, it implies that $Y^{(n)}\! -\! J^{(n)}$ (and thus $\cH^{(n)}$) hits value $0$ between the times $q$ and $t_0$ and between the times $t_0$ and $r$: namely, for all sufficiently large $n$, $\gamma (n) \! \geq \! q$ and $\delta (n) \! \leq \! r$. This easily entails 
$\liminf_{n\rightarrow \infty} \gamma (n) \! \geq \! l_k$ and $r_k \! \geq \! \limsup_{n\rightarrow \infty} \delta (n)$, and we have proved 
that $\lim_{n\rightarrow \infty} \gamma (n) \!=\! l_k$ and $\lim_{n\rightarrow \infty} \delta (n)\! = \! r_k$. 

Let 
$n_0\! \geq \! 1$ be such that  for all $n\! \geq \! n_0$, $\cH^{_{(n)}}_{t_0} \! >\! 0$. Then, for all $n\! \geq \! n_0$, there exists 
$j(n,k)\ino \{1, \ldots , \mathbf{q}_{\bw_n}\}$ such that $\gamma (n)\! = \! l^n_{j(n,k)}$ and $\delta (n)\! = \! r^n_{j(n,k)}$; for all $n\! \leq \! n_0$, we take for instance $j(n, k)\! = \! 1$. Then, (\ref{indexouy}) holds true, which completes the proof. \cqfd 
  
\bigskip

  We next recall that Proposition \ref{AldLim98} 
(Proposition 7 in Aldous \& Limic \cite{AlLi98}) asserts that 
$\sum_{1\leq k\leq \mathbf{q}_{\bw_n}} (\zeta^n_k)^2 \! \rightarrow \!   \sum_{k\geq 1} (\zeta_k)^2 $ weakly on $[0, \infty)$ as 
$n\! \rightarrow\! \infty$. We use this result to prove the following joint convergence. 
\begin{lem}
\label{convexc2} We make the same assumptions as in Theorem \ref{HYcvth}. We keep the previous notations. Then 
\begin{equation}
\label{YHZcv}
\mathscr{Q}_n (6)\! :=\!  \Big( Y^{(n)} , \cH^{(n)} , \Ptt^{(n)}, \!\!\! \!\!\! \!\! \sum_{\quad 1\leq k\leq \mathbf{q}_{\bw_n}} \!\!\! \!\!(\zeta^n_k)^2 \Big)
\;  \underset{n\rightarrow \infty}{-\!\!\! -\!\!\! -\!\!\! \longrightarrow} \; \mathscr{Q} (6)\! :=\! \Big(  Y, \cH , \Ptt,  \sum_{k\geq 1} (\zeta_k)^2\Big)
\end{equation} 
weakly on $\bD([0, \infty), \bbR) \! \times \! \bC ([0, \infty) , \bbR) \! \times \!  (\bbR^2)^{\bbN^*}\! \!\times \![0, \infty)$,  
equipped with product topology.
\end{lem}
\noi
\textbf{Proof.} The laws of the $\mathscr{Q}_n (6)$ are tight by (\ref{jtcvYHP}) in Theorem \ref{HYcvth}   combined with the weak convergence $\sum_{1\leq k\leq \mathbf{q}_{\bw_n}} (\zeta^n_k)^2 \! \rightarrow \!   \sum_{k\geq 1} (\zeta_k)^2 $. 
We only need to prove that the law of $\mathscr{Q}(6)$ is the unique limiting law: to that end, let $(n(p))_{p\in \bbN}$ be an increasing sequence of integers such that $\mathscr{Q}_{n(p)} (6) \! \rightarrow \!  ( Y, \cH, \Ptt, Z)$ weakly. It remains to prove that $Z\! = \! \sum_{k\geq 1} (\zeta_k)^2$. 
Without loss of generality (but with a slight abuse of notation), by Skorokod's representation theorem 
we can assume that $\mathscr{Q}_{n(p)} (6) \! \rightarrow \!  ( Y, \cH, \Ptt, Z)$ holds true $\bP$-a.s.
Then, by Lemma \ref{convexc1}, observe that for all $l\geq \! 1$,  
\begin{equation}
\label{fatoutou}
Z  \;  \underset{n\rightarrow \infty}{\longleftarrow \!\!\! -\!\!\! -} \quad  \sum_{1\leq k\leq \mathbf{q}_{\bw_n}} (\zeta_k^n)^2 \; \geq \;   \sum_{1\leq k\leq l} (\zeta_{j(n,k)}^n)^2\quad  \underset{n\rightarrow \infty}{ -\!\!\! -\!\!\! \longrightarrow} \; \sum_{1\leq k\leq l} (\zeta_k)^2 . 
\end{equation}
Set $Z^\prime\! = \! \sum_{ k\geq 1} (\zeta_k)^2$; by letting $l$ go to $\infty$ in (\ref{fatoutou}), 
we get $Z\! \geq \! Z^\prime$, which implies 
$Z\! = \! Z^\prime$ a.s.~since $Z$ and $Z^\prime$ have the same law. This completes the proof of the lemma. \cqfd 

\medskip

Without loss of generality (but with a slight abuse of notation), by Skorokod's representation theorem 
we can assume that  (\ref{YHZcv}) holds true a.s.~on $\bD([0, \infty), \bbR) \! \times \! \bC ([0, \infty) , \bbR)\! \times \!  (\bbR^2)^{\bbN^*}\! \! \times \! [0, \infty)$, equipped with product topology. 
%Then, it implies 
%\begin{equation}  
%\label{indexah}
%\textrm{$\bP$-a.s.~for all $k\! \geq \! 1$,} \qquad \big( l^n_k  ,  r^n_k\big)
%\;  \underset{n\rightarrow \infty}{-\!\!\! -\!\!\! -\!\!\! \longrightarrow} \; \big(  l_k, r_k \big). 
%\end{equation}
We next prove the following. 
\begin{lem}
\label{convexc3} Assume that (\ref{YHZcv}) holds true almost surely. 
We keep the previous notations. Then, 
\begin{equation}  
\label{indexah}
\textrm{$\bP$-a.s.~for all $k\! \geq \! 1$,} \qquad \big( l^n_k  ,  r^n_k\big)
\;  \underset{n\rightarrow \infty}{-\!\!\! -\!\!\! -\!\!\! \longrightarrow} \; \big(  l_k, r_k \big). 
\end{equation}
\end{lem}
\noi
\textbf{Proof.} %\ts{Let $\epp \ino (0, \infty)$ 
%\margmm{rewritten}
%\mm{be distinct from $\zeta_{k}$, $k\ge 1$ and let}} 
Let $\epp \ino (0, \infty) \backslash \{ \zeta_k; k\! \geq \! 1\}$ and let $k_\epp$ be such that $\zeta_k \! >\! \epp$ for all $k\ino \{ 1, \ldots , k_\epp\} $ 
and $\zeta_k \! < \! \epp$ for all $k\! >\! k_\epp$. 
Let $k_\epp^\prime \! > \! k_\epp$ be such that $\sum_{k>k^\prime_\epp} (\zeta_k)^2 \! < \! \epp^2/3$. 
Since $k^\prime_\epp\! >\! k_\epp$, we also get $\min_{1\leq  k\leq k^\prime_\epp} |\epp \! -\! \zeta_k | \! < \! \epp$.
By 
Lemma \ref{convexc1} and Lemma \ref{convexc2}, there exists $n_0\! \geq \! 1$ such that for all $n\! \geq \! n_0$, 
\begin{multline}
\label{crunouf}
\quad \Big|\!\!\!\!\!\! \sum_{\quad 1\leq k\leq \mathbf{q}_{\bw_n}}\!\!\!\!\!\!   (\zeta_k^n)^2\,  -\! \sum_{k\geq 1} (\zeta_k)^2\,  \Big|   < \epp^2/3,   \quad \sum_{1\leq k\leq k^\prime_\epp} \big| (\zeta^n_{j(n,k)})^2 \!-\! (\zeta_k)^2 \big| < \epp^2/ 3  \\
\!\!\!\!\!  \textrm{and} \quad \max_{1\leq k\leq k^\prime_\epp} \big| \zeta_k\! -\! \zeta_{j(n,k)}^n  \big| < \! \min_{1\leq  k\leq k^\prime_\epp} |\epp \! -\! \zeta_k | \tt{<  \epp} . 
\end{multline}
Set $S_n\! = \! \{ 1, \ldots , \mathbf{q}_{\bw_n}\} \backslash \{ j(n,1), \ldots , j(n, k^\prime_\epp)\}$. The previous inequalities imply 
for all $n\! \geq \! n_0$, that $\sum_{k\in S_n} (\zeta^n_k)^2 < \epp^2$. 
Thus, for all $n\!\geq \! n_0$, if $k\! \in \! S_n$, then $\zeta^n_k \! < \! \epp$. 
Next observe that for all $k\ino \{k_\epp + 1, \ldots , k^\prime_\epp\}$, 
 %$$ \zeta_k - \max_{1\leq \ell \leq k^\prime_\epp} \big| \zeta_\ell\! -\! \zeta_{j(n,\ell)}^n  \big|  \leq
 $$  \zeta^n_{j(n, k)} \leq \epp -(\epp -\zeta_k)+ \max_{1\leq \ell \leq k^\prime_\epp} \big| \zeta_\ell\! -\! \zeta_{j(n,\ell)}^n  \big| < \epp + \min_{1\leq  \ell\leq k^\prime_\epp} |\epp \! -\! \zeta_\ell | -|\epp -\zeta_k| < \epp $$
by (\ref{crunouf}). 
%Thus, for all $k\ino \{ k_\epp +1, \ldots , k^\prime_\epp\}$, $\zeta^n_{j(n, k)} \! < \! \epp$. 
Also note that for all $k\ino \{ 1, \ldots, k_\epp\}$, 
$$\zeta^n_{j(n, k)} \geq  \zeta_k - \max_{1\leq \ell \leq k^\prime_\epp} \big| \zeta_\ell\! -\! \zeta_{j(n,\ell)}^n  \big| >  \epp + |\zeta_k -\epp |- \min_{1\leq  \ell\leq k^\prime_\epp} |\epp \! -\! \zeta_\ell | > \epp $$
again by (\ref{crunouf}).
%Thus, $k\ino \{ 1, \ldots, k_\epp\}$,  $\zeta^n_{j(n, k)} \! > \! \epp$. }
%\footnote{\tt{I added details.}} 
To summarise, for all $n\! \geq \! n_0$, $\zeta^n_{j(n, k)} \! >\! \epp$ if $k\ino \{ 1, \ldots, k_\epp \}$ and $\zeta^n_{j(n, k)} \! <\! \epp$ for all 
$k \ino \{ k_\epp +1, \ldots , \mathbf{q}_{\bw_n} \}$. Since $\zeta_1\! >\! \zeta_2 \! >\! \ldots \! >\! \zeta_{k_\epp}$, there exists  
$n_1\! \geq \! n_0$ such that for all $n\! \geq \! n_1$, $ \zeta_{j(n,1)}^n \! >\!  \zeta_{j(n,2)}^n \! >\! \ldots \! >\!  \zeta_{j(n,k_\epp)}^n$. Thus, for all $n\! \geq \! n_1$ and for all $k\ino \{ 1, \ldots, k_\epp\}$, 
we have proved that $j(n,k)\! = \! k$, which entails (\ref{indexah}) 
since $\epp$ can be chosen arbitrarily small.  \cqfd 
%
%$\max_{0\leq k \leq k_\epp} | \zeta_k\! -\! \zeta_{j(n,k)}^n  | \! < \! \tfrac{1}{3} 
%\max_{0\leq k \leq k_\epp} | \zeta_k\! -\! \zeta_{k+1}  | $. 
%
%
%Moreover, for all  for all $n\! \geq \! n_0$ and 
%for all $k \ino \{ 1, \ldots , k_\epp \}$, $ \zeta_{j(n,k+1)}^n \! \leq \! \zeta_{k+1} + \frac{1}{2} \epp \leq 
% \zeta_{j(n,k)}^n -(\zeta_k \! -\!  \zeta_{k+1}) + \epp \leq 
%
%% \zeta_{k} -(\zeta_k \! -\!  \zeta_{k+1})  + 
%%\max_{1\leq \ell \leq k^\prime_\epp} \big| \zeta_k\! -\! \zeta_{j(n,\ell)}^n  \big|
%
%This implies that for all sufficiently large $n$, 
%$j(n, k)\! = \! k$, for all $k\ino \{ 1, \ldots , k_\epp \}$ and (\ref{indexouy}) in Lemma \ref{convexc1} implies 
%$(l^n_k, r^n_k) \! \rightarrow \! (l_k, r_k)$ a.s.~ for all $k\ino \{ 1, \ldots , k_\epp \}$, which entails (\ref{indexah}) 
%since $\epp$ can be choosen arbitrarily small.  \cqfd 
%
\medskip

Recall from (\ref{rexccont}) the definition of the excursions $\mathtt{H}_k$ and $\mathtt{Y}_k$ of resp.~$\cH$ and $Y\! -\! J$ above $0$. We define the (rescaled) excursion of $Y^{(n)} \! -\! J^{(n)}$ and of $\cH^{(n)}$ above $0$ as follows: 
\begin{equation}
\label{rescexc}
 \forall k \! \geq \! 1, \; \forall t\ino [0, \infty), \qquad 
\Htt^{_{(n)}}_{k}(t)= \cH^{_{(n)}}_{(l^n_k + t)\wedge r^n_k} \quad \textrm{and} \quad  \Ytt^{_{(n)}}_{k}(t)= Y^{_{(n)}}_{(l^n_k + t)\wedge r^n_k}- \JJ^{_{(n)}}_{l^n_k} . 
\end{equation}
%\margmm{new}
Recall from Lemma \ref{XYexcu} (iii) that $\Delta Y_{l_{k}}=0$ a.s. 
Then by (\ref{YHZcv}), Lemma \ref{convexc3} and Lemma \ref{extraisko} $(iii)$ in Appendix, we immediately get the following. 
\begin{lem}
\label{convexc4} Assume that (\ref{YHZcv}) holds true almost surely. 
We keep the previous notations. Then, 
\begin{equation}  
\label{indexah2}
\textrm{$\bP$-a.s.~for all $k\! \geq \! 1$,} \qquad \big(\Ytt^{_{(n)}}_{k}, \Htt^{_{(n)}}_{k} , l^n_k, r^n_k \big)
\;  \underset{n\rightarrow \infty}{-\!\!\! -\!\!\! -\!\!\! \longrightarrow} \; \big(\Ytt_k, \Htt_k,  l_k , r_k\big). 
\end{equation}
in $\bD([0, \infty), \bbR)  \!\!  \times \! \bC ([0, \infty) , \bbR)\! \times \![0, \infty)^2$.  
\end{lem}

Recall from (\ref{Poisurcon}) and (\ref{pinchset})  the definition of $\Ptt\! = \! \big( (s_p, t_p) \big)_{p\geq1}$ and recall from 
(\ref{shortnota}) the notation $\Ptt^{(n)}\! = \! \big( (s^n_p, t^n_p) \big)_{1\leq p\leq \bp_n}$.  
We next prove the following. 
\begin{lem}
\label{convexc5}Assume that (\ref{YHZcv}) holds almost surely. Then, a.s.~for all $p\! \geq \! 1$, there exists $k\! \geq \! 1$ such that 
$l_k \! < \! s_p \! \leq \! t_p \! < \! r_k$ and for all sufficiently large $n$, $l^n_k \! < \! s^n_p \! < \! t^n_p \! < \! r^n_k$ 
and $(l^n_k,s^n_p,  t^n_p,  r^n_k) \! \rightarrow \! (l_k,s_p,  t_p,  r_k)$. 
\end{lem}
\noi
\textbf{Proof.} By Proposition \ref{AldLim1} $(i)$, $\bP$-a.s.~for all $p\! \geq \! 1$, $Y_{t_p} \! >\! J_{t_p}$ and there exists 
$k\! \geq \! 1$ such that $t_p \ino (l_k, r_k)$. By Lemma \ref{XYexcu} $(iii)$, we get $Y_{l_k} \! -\! J_{l_k}\! = \! 0$; 
recall that $y_p \ino (0, Y_{t_p} \! -\! J_{t_p})$ and recall from (\ref{pinchset}) that $s_p \! = \! \inf \big\{ s\ino [0, t_p] :  \inf_{u\in [s, t_p]} Y_u\! -\! J_u > y_p   \big\}$; thus, we get $l_k \! < \! s_p \! \leq \! t_p \! < \! r_k$ and the proof is completed by 
(\ref{YHZcv}) that asserts that $(s^n_p, t^n_p) \! \rightarrow \! (s_p, t_p)$ and by Lemma \ref{convexc3} that asserts that $(l^n_k, r^n_k)\! \rightarrow \! (l_k, r_k)$. \cqfd

\medskip
\noindent{\bf Proof of Theorem \ref{excucvth}. }
Recall from (\ref{defPik}) that for all $k\! \geq \! 1$, 
$\Ptt_{k}\! = \! \big( (s^k_p, t^k_p)\, ;\, 1\! \leq \! p\! \leq \! \bp_k \big)$ where 
$(t^k_p\, ;\,  1\! \leq \! p\! \leq \! \bp_k )$ 
increases and where the $(l_k +s^k_p, l_k+t^k_p)$ are exactly the terms $(s_{p^\prime} , t_{p^\prime})$ of $\Ptt$ such that $t_{p^\prime}\ino [l_k, r_k]$. 
%\begin{multline*}
%%\label{defPik}
%\Ptt_{k}\! = \! \big( (s^k_p, t^k_p)\big)_{1\leq p\leq \bp_k} \; \textrm{where $(t^k_p)_{1\leq p\leq \bp_k}$ increases and where} \\ \textrm{ the $(l_k +s^k_p, l_k+t^k_p)$ are exactly the terms $(s_{p^\prime} , t_{p^\prime})$ of $\Ptt$ such that $t_{p^\prime}\ino [l_k, r_k]$,}
%\end{multline*}
Similarly recall from (\ref{defPiwk}) the definition of the sequence of pinching times 
$\Ptt^{\bw_{n}}_k$, $1\! \leq \!  k\!  \leq \!  \bq_{\bw_n} $: namely, in their rescaled version, 
$\frac{_1}{^{b_n}}\Ptt_{k}^{\bw_n}\! = \! \big( (s^{n,k}_p, t^{n,k}_p) \, ; \, 1\! \leq \! p\! \leq \! \bp^{n}_k \big) $, where 
$(t^{n, k}_p\, ; \, 1\! \leq\!  p\! \leq\!  \bp^n_k)$ increases and where the $(l^n_k +s^{n, k}_p, l^n_k+t^{n,k}_p)$ are exactly the terms $(s^n_{p^\prime} , t^n_{p^\prime})$ of $\Ptt^{(n)}$ such that $t^n_{p^\prime}\ino [l^n_k, r^n_k]$. 
%
%\begin{multline*}
%%\label{defPiwk}
%\frac{_1}{^{b_n}}\Ptt_{k}^{\bw_n}\! = \! \big( (s^{n,k}_p, t^{n,k}_p)\big)_{1\leq p\leq \bp^{n}_k} \; \textrm{where 
%$(t^{n, k}_p)_{1\leq p\leq \bp^n_k}$ increases and where} \\ \textrm{ the $(l^n_k +s^{n, k}_p, l^n_k+t^{n,k}_p)$ are exactly the terms $(s^n_{p^\prime} , t^n_{p^\prime})$ of $\Ptt^{(n)}_\bw$ such that $t^n_{p^\prime}\ino [l^n_k, r^n_k]$.}
%\end{multline*}
Thus, Lemma \ref{convexc5} immediately entails that $\bP$-a.s.~for all $k\geq \! 1$, $ \frac{_1}{^{b_n}}\Ptt_{k}^{\bw_n}\! \rightarrow \!  \Ptt_k $ as $n\! \rightarrow \! \infty$.
%$$ \textrm{$\bP$-a.s.~for all $k\geq \! 1$,} \quad \frac{_1}{^{b_n}}
%\Ptt_{k}^{\bw_n} \;  \underset{n\rightarrow \infty}{-\!\!\! -\!\!\! -\!\!\! \longrightarrow} \; \Ptt_k . $$
This convergence combined with Lemma \ref{convexc4} implies Theorem \ref{excucvth}. \cqfd   

\subsubsection{Proof of Theorem \ref{graphcvth}.}
\label{graphcvpf}

\noindent\textbf{Proof of \eqref{graphbwn}. }
Keep the previous notation and recall 
that $\big( \cG_k^{\bw_n}, d_{k}^{\bw_n}, \varrho_k^{\bw_n}, \bm_k^{\bw_n}\big)$, $1\! \leq \! k \! \leq \! \bq_{\bw_n}$, 
%$$ \big( \cG_k^{\bw_n}, d_{k}^{\bw_n}, \varrho_k^{\bw_n}, \bm_k^{\bw_n}\big), \quad 1\leq k \leq \bq_{\bw_n} $$
stand for the connected components of the $\bw_n$-multiplicative random graph $\cG_{\bw_n}$. Here, 
$d_{k}^{\bw_n}$ stands for the graph-metric on $\cG^{\bw_n}_k$, 
$\bm_k^{\bw_n}$ is the restriction to $\cG_k^{\bw_n}$ of the measure $\bm_{\bw_n} \! = \! \sum_{j\geq 1} w^{_{(n)}}_{j} \delta_j$, $\varrho_k^{\bw_n}$ is the first vertex of 
$\cG^{\bw_n}_k$ that is visited during the exploration of $\cG_{\bw_n}$, and 
the indexation is such that $\bm_1^{\bw_n} \big( \cG_1^{\bw_n}\big)  \geq \ldots \geq \bm_{\bq_{\bw_n}}^{\bw_n} 
\big( \cG_{\bq_{\bw_n}}^{\bw_n}\big)$. 
%$$ \bm_1^{\bw_n} \big( \cG_1^{\bw_n}\big)  \geq \ldots \geq \bm_{\bq_{\bw_n}}^{\bw_n} 
%\big( \cG_{\bq_{\bw_n}}^{\bw_n}\big) . $$

Next, recall from (\ref{rescexc}) that $\Htt^{_{(n)}}_k (\cdot)$ stands for the $k$-th longest excursion of $\cH^{\bw_n}$ that 
is rescaled in time by a factor $1/b_n$ and rescaled in space by a factor $a_n / b_n$; recall that 
$\frac{_1}{^{b_n}}\Ptt_{k}^{\bw_n}\! = \! \big( (s^{n,k}_p, t^{n,k}_p); 1\! \leq \! p\! \leq \! \bp^{n}_k \big)$ 
is the ($1/b_n$-rescaled) finite sequence of pinching times of $\Htt^{_{(n)}}_k$. Then, for all $k\ino \{ 1, \ldots , \bq_{\bw_n}\}$ the compact measured metric space $ \mathbf{G}^{_{(n)}}_{k}:= \big( \cG_k^{\bw_n} ,  \frac{_{_{a_n}}}{^{^{b_n}}}d_{k}^{\bw_n} , \varrho_k^{\bw_n}, \frac{_{_{1}}}{^{^{b_n}}}\bm_k^{\bw_n}  \big)$
%$$  \mathbf{G}^{_{(n)}}_{k}:= \big( \cG_k^{\bw_n} ,  \frac{_{_{a_n}}}{^{^{b_n}}}d_{k}^{\bw_n} , \varrho_k^{\bw_n}, \frac{_{_{1}}}{^{^{b_n}}}\bm_k^{\bw_n}  \big)$$ 
%\begin{multline*}
%\big(\big( \cG_k^{\bw_n} ,  \frac{_{_{a_n}}}{^{^{b_n}}}d_{k}^{\bw_n} , \varrho_k^{\bw_n}, \frac{_{_{1}}}{^{^{b_n}}}\bm_k^{\bw_n}  \big) \big)_{k\geq 1}\;  \textrm{is isometric to} \; 
is isometric to $G\big( \Htt_k^{_{(n)}} , \frac{_{_1}}{^{^{b_n}}}\Ptt_{k}^{\bw_n} , \frac{_{a_n}}{^{b_n}} \big)$,  
the compact measured metric space coded by $ \Htt_k^{_{(n)}}$ and the pinching setup 
 $( \frac{_1}{^{b_n}}\Ptt_{k}^{\bw_n} , \frac{_{a_n}}{^{b_n}})$ as defined in (\ref{defgrgraf}).
Similarly, recall from (\ref{rexccont}) that $\Htt_k (\cdot)$ stands for the $k$-th longest excursion of $\cH$ and recall from (\ref{defPik}) that $\Ptt_k  \! = \! \big( (s^{k}_p, t^{k}_p); 1\! \leq \! p\! \leq \! \bp_k \big)$ is the  finite sequence of pinching times of $\Htt_k$. 
Then, for all $k\! \geq \! 1$, 
 the compact measured metric space $ \mathbf{G}_{k} :=\big( \mathbf{G}_{k} , \mathrm{d}_{k}, \varrho_{k} , \bm_{k} \big) $
%$$ \mathbf{G}_{k} :=\big( \mathbf{G}_{k} , \mathrm{d}_{k}, \varrho_{k} , \bm_{k} \big) $$
is isometric to $ G (\Htt_k, \Ptt_k, 0)$ that is the compact measured metric space coded by $ \Htt_k$ and the pinching setup $( \Ptt_{k} ,0)$ as defined in (\ref{defgrgraf}). 
Without loss of generality (but with a slight abuse of notation), by Skorokod's representation theorem 
we can assume that the convergence in Theorem \ref{excucvth} holds almost surely. Namely a.s.~for all $k\! \geq \! 1$, $\big( \Htt^{_{(n)}}_k , \zeta^n_k, \frac{_1}{^{b_n}} \Ptt^{\bw_n}_k\big)\! \rightarrow \! 
\big( \Htt_k , \zeta_k, \Ptt_k\big)$ on $\bC ([0, \infty), \bbR) \times [0, \infty) \times (\bbR^2)^{\bbN^*}$. 
We next fix $k\! \geq \! 1$; then for all sufficiently large $n$,  $\frac{_1}{^{b_n}} \Ptt^{\bw_n}_k$ and $\Ptt_k$ have the same number of points: namely, $\bp^{n}_k\! = \! \bp_k$ and 
\begin{equation}
\label{pinpinch}
 \forall 1 \! \leq \! p \! \leq \! \bp^{n}_k\! = \! \bp_k, \quad (s^{n,k}_p, t^{n,k}_p)  \;  \underset{n\rightarrow \infty}{-\!\!\! -\!\!\! -\!\!\! \longrightarrow} \; (s^k_p, t^k_p) \; .
\end{equation} 
Recall from (\ref{defGHP}) the definition of the Gromov--Hausdorff--Prokhorov distance $\bdelta_{\mathrm{GHP}}$. We next apply Lemma \ref{codconGHP} with $(h,h^\prime)\! = \! (   \Htt_k, \Htt^{_{(n)}}_k )$, 
% $r\! = \!  \bp^{n}_k\! = \! \bp_k$, 
$(\Pi, \Pi^\prime)\! = \! (\Ptt_k,  \frac{_1}{^{b_n}} \Ptt^{\bw_n}_k)$, $(\epp, \epp^\prime)\! = \! (0, a_n/b_n)$ and 
$\delta \! = \! \delta_n\! = \! \max_{1\leq p \leq \bp_k} |s^k_p\! -\! s^{n,k}_p|\vee  |t^k_p\! -\! t^{n,k}_p|$. Then, by 
(\ref{contGHP}), 
\begin{equation}
\label{blurpimiam}
\bdelta_{\mathrm{GHP}} ( \mathbf{G}_{k}, \mathbf{G}^{_{(n)}}_{k} ) \leq 6(\bp_k+1) \big( \lVert \Htt_k \! -\! \Htt_k^{_{(n)}} \rVert_\infty + \omega_{\delta_n} (\Htt_k) \big) + 3a_n\bp_k /b_n + |\zeta_k\! -\! \zeta^n_k|,   
\end{equation}
where $ \omega_{\delta_n} (\Htt_k)\! = \! \max \{ \lvert \Htt_k (t) \! -\! \Htt_k (s)\rvert ; s, t \ino [0, \infty): |s\! -\! t| \! \leq \delta_n \}$.     
By (\ref{pinpinch}), $\delta_n\! \rightarrow \! 0$; since $\Htt_k$ is continuous and since it is null on $[\zeta_k, \infty)$, it is uniformly continuous and $ \omega_{\delta_n} (\Htt_k)\! \rightarrow \! 0$; recall that $a_n / b_n\! \rightarrow \! 0$. 
Thus, the right member of (\ref{blurpimiam}) goes to $0$ as $n\! \rightarrow\! 0$. Thus, we have proved that a.s.~for all 
$k\! \geq \! 1$, $ \bdelta_{\mathrm{GHP}} ( \mathbf{G}_{k}, \mathbf{G}^{_{(n)}}_{k} ) \! \rightarrow\! 0$, which implies \eqref{graphbwn} in  Theorem \ref{graphcvth}. \cqfd

\bigskip

\noindent\textbf{Proof of \eqref{globabab}. } 
We next prove the convergence of the connected components equipped with the counting measure. 
Recall from Introduction the definition of the discrete tree $\cT_{\! \! \bw_n}$ coded by the $\bw_n$-LIFO queue without repetition (namely, the tree coded by $\cH^{\bw_n}$): \textit{the vertices of $\cT_{\! \! \bw_n}$ are the clients; the server is the root (Client $0$) and Client $j$ is a child of Client $i$ in $\cT_{\! \! \bw}$ if and only if  Client $j$ interrupts the service of Client $i$ (or arrives when the server is idle if $i\! =\!  0$).} 
We denote by $\cC^{\bw_n}$ the contour process associated with  $\cT_{\! \! \bw_n}$ that
is informally defined as follows: suppose that $\cT_{\! \! \bw_n}$ is embedded in the oriented half plane in such a way that edges have length one and that orientation reflects lexicographical order of visit; 
we think of a particle starting at time $0$ from the root of $\cT_{\! \! \bw_n}$ and exploring the tree from the left to the right, backtracking as less as possible and 
moving continuously along the edges at unit speed. Since $\cT_{\! \! \bw_n}$ is finite, 
the particle crosses each edge twice (upwards first and then downwards). 
For all $s\ino [0, \infty)$, we define $\cC^{\bw_n}_s $ 
as the distance at time $s$ of the particle from the root of $\cT_{\! \! \bw_n}$. We refer to Le Gall \& D.~\cite{DuLG02} (Section 2.4, Chapter 2, pp.~61-62) for a formal definition and the connection with the height process (see also the end of Section \ref{connecMqu}). 

  It is important to notice that the trees coded by $\cC^{\bw_n}$ and by $\cH^{\bw_n}$ are the same: the only difference is the measure induced by the two different coding functions. More precisely, $\cC^{\bw_n}$ is derived from $\cH^{\bw_n}$ by the following time-change: recall that $\mathbf{j}_n \! = \! \max \{ j\! \geq \! 1\! : \! 
  w^{_{(n)}}_{^j} \! >\! 0 \}$ and   
let  $(\xi^n_k)_{1\leq k\leq  2\mathbf{j}_n}$ be the sequence of jump-times of $\cH^{\bw_n}$: namely, $\xi^n_{k+1} \! = \! \inf \{ s\! >\! \xi^n_k \! : \cH^{\bw_n}_s \! \neq \! \cH^{\bw_n}_{\xi^n_k} \}$, for all $1 \! \leq \! k \! < \! 2\mathbf{j}_n $, 
with the convention $\xi^n_0 \! = \! 0$. We then set
 \begin{equation}
\label{Phindef}
\forall t\ino [0, \infty), \; \Phi_n (t)  \! =\!\!\!\! \!\! \!\! \!  \sum_{\quad 1 \leq k\leq 2\mathbf{j}_n}\!\! \!\! \!\! \!  \un_{[0, t]} (\xi^n_k) 
 \;   \textrm{and} \;    \forall s\ino [0, 2\mathbf{j}_n], \;  \phi_n (s)\! = \! \inf \big\{ t\ino [0, \infty) \! :\! \Phi_n (t) \! \geq \! s \big\}. 
\end{equation}
Note that $\phi_n (k)=\xi^n_k$. Then, we get 
\begin{equation}
\label{voilaaa}
\forall t \ino [0, \infty), \quad \cC^{\bw_n}_{\Phi_n (t)}= \cH^{\bw_n}_t \;  \quad  \textrm{and} \quad  \forall k\ino \{ 0, \ldots , 2 \mathbf{j}_n\} , \; \, \cC^{\bw_n}_{k}\! = \! \cH^{\bw_n}_{\xi^n_k}=\cH^{\bw_n}_{\phi_n(k)}\; .
\end{equation}
We next set 
\begin{equation}
\label{gisnoz}
\forall t \ino [0, \infty) , \quad R^n_t \! = \! \sum_{j\geq 1} \un_{\{ E^{\bw_n}_j \leq t \}} 
\end{equation}
that counts the number of clients who entered the $\bw_n$-LIFO queue governed by $Y^{\bw_n}$.
Recall here that $E^{{\bw_n}}_{j}$ is the first jump-time of $N^{\bw_n}_j$: namely the $E^{\bw_n}_j$ are independent exponentially distributed r.v.~with respective parameters $w_{^j}^{_{(n)}}/ \sigma_1 (\bw_n)$. In terms of the tree $\cT_{\!\! \bw_n}$, $R^n_t$  is the number of distinct vertices that have been explored by $\cH^{\bw_n}$ up to time $t$.  
By arguing as in the proof of (\ref{crossedg}), we easily check that 
\begin{equation}
\label{voiliii}
\forall t\ino [0, \infty), \quad \Phi_n (t)= 2 R^n_t - \cH^{\bw_n}_t \; .
\end{equation}
We prove the following. 
\begin{lem}
\label{cntrlll} We keep the previous notation. Then, the following holds. 
\begin{equation}
\label{idencvv}
\forall t\ino [0, \infty), \quad \bE \Big[ \!\! \sup_{\;\;  s\in [0, t]} \! \big| R^n_s \!-\! s\big| \,  \Big] \leq 2\sqrt{t} + \frac{t^2\sigma_2 (\bw_n)}{2\sigma_1 (\bw_n)^2} \; .
\end{equation}
Moreover, there exists a positive r.v.~$Q_n$ that is a measurable function of $(N^{\bw_n}_j)_{j \geq  1}$, such that $\bE [Q_n^2]\! \leq  \!  4\mathbf{j}_n $ (recall that 
$\mathbf{j}_n \! := \! \max \{ j\! \geq \! 1: w_j^{(n)} \! >\! 0 \}$) and such that 
\begin{equation}
\label{fincrt}
\textrm{$\bP$-a.s.~for all $s, t\ino [0, \infty)$,} \quad   R^n_{t+s}-R^n_t \! \leq \! s+ 2Q_n \; .
\end{equation}
\end{lem}
\noi
\textbf{Proof.} Set $M_j(t)\! = \! \un_{\{ E^{\bw_n}_j  \leq t \}}\! -\! \frac{w_j^{(n)}}{\sigma_1(\bw_n)} (t \wedge E^{\bw_n}_j)$ and denote by $(\ccG_t)$ the natural filtration associated with the $(N^{\bw_n}_j)_{j \geq  1}$. It is easy to check that the $M_j $ are independent $(\ccG_t)$-martingales and that 
$\bE \big[ M_j (t)^2\big]\! = \! 1\! -\! \exp (-w_{^j}^{_{(n)}} t/\sigma_1 (\bw_n))\! \leq \! w_{^{j}}^{_{(n)}} t/\sigma_1 (\bw_n)$. 
%\begin{equation} 
%\label{calcolo}
%\bE \big[ M_j (t)^2\big]= 1 - e^{-w_{j}^{{(n)}} t/\sigma_1 (\bw_n)} + 2 \frac{_{w_j^{(n)}}}{^{\sigma_1(\bw_n)}} t e^{-w_j^{(n)} t/\sigma_1 (\bw_n)} \big( 1 \! - \!  \frac{_{w_j^{(n)}}}{^{\sigma_1(\bw_n)}} \big) \; \leq \! 3  \frac{_{w_j^{(n)}}}{^{\sigma_1(\bw_n)}} t \; .
%\end{equation}
We then set $M(t)\! = \! \sum_{1\leq j \leq \mathbf{j}_n} M_j (t)$. Then $M$ is a $(\ccG_t)$-martingale and Doob's $L^2$ inequality implies that $\bE [\sup_{s\in [0, t]} M (s)^{2}] \! \leq \! 4 \bE [M (t)^2] \! \leq \! 4 t$. Thus, $\bE [\sup_{s\in [0, t]} |M (s)] \, ] \! \leq \! 2\sqrt{t}$. Next, we set   
$\overline{M}(t)\! = \! R^n_t \! -\! M_t$. We easily check the following: 
$$ t\!  -\!  \overline{M}(t) \! = \! \sum_{j \geq 1} \frac{_{w_j^{(n)}}}{^{\sigma_1(\bw_n)}} \big(t\! -\! E^{\bw_n}_j \big)\un_{\{ E^{\bw_n}_j \leq t \}} \, , $$ 
which is nonnegative and nondecreasing in $t$ so that $\sup_{s\in [0, t]} |s\! -\! \overline{M} (s)|\! = \! 
t\! -\! \overline{M}(t)$. Moreover, for all $j\! \geq \! 1$, we check that
$$\frac{_{w_j^{(n)}}}{^{\sigma_1(\bw_n)}} \bE \big[\big(t\! -\! E^{\bw_n}_j \big)\un_{\{ E^{\bw_n}_j \leq t\}} \big]\! = \! e^{-w_{j}^{{(n)}} t/\sigma_1 (\bw_n)} \!-\! 1 + \frac{_{w_j^{(n)}}}{^{\sigma_1(\bw_n)}} t \,  \leq \,  \frac{_1}{^2} t^2 \big( w_j^{(n)} /\sigma_1 (\bw_n)\big)^2 \; .$$ 
This implies that $\bE [\sup_{s\in [0, t]} |s\! -\! \overline{M} (s)| \, ]\! = \! \bE [ 
t\! -\! \overline{M}(t) ] \leq t^2 \sigma_2 (\bw_n) / (2\sigma_1 (\bw_n)^2)$, which easily completes the proof of (\ref{idencvv}) thanks to the previous inequality regarding $M$. 

 Let us prove (\ref{fincrt}). To that end, observe that $\lim_{t\rightarrow \infty} \bE [M_j (t)^2] \! = \! 1$. 
Thus, $\lim_{t\rightarrow \infty} \bE [M (t)^2] \! = \!  \mathbf{j}_n$ and Doob's inequality entails that $\bE [\sup_{t\in [0, \infty)} M^2 (t)] \! \leq \!  4  \mathbf{j}_n$. We then set $Q_n \! = \!  \sup_{t\in [0, \infty)}| M (t) | $ and we get almost surely 
for all $t,s \ino [0, \infty)$, $R^n_{t+s}\! -\! R^n_t \! = \! M(t+s) \! -\! M(t)+ \overline{M} (t+s) \! -\!   \overline{M} (t) \! \leq \! 2Q_n + \overline{M} (t+s) \! -\!   \overline{M} (t)$. Since for all $a\ino [0, \infty)$, the function $t \mapsto t\wedge a$ is $1$-Lipschitz and since 
$\overline{M} $ is a convex combination of these functions, $\overline{M} $ is also $1$-Lipschitz: namely, $| \overline{M} (t+s) \! -\!   \overline{M} (t)| \! \leq \! s$, 
%\margmm{rewritten}
which completes the proof of (\ref{fincrt}). \cqfd

\medskip

By (\ref{voiliii}) and (\ref{idencvv}) we easily get for all $t, \epp \ino (0, \infty)$
\begin{eqnarray*}
\bP \big( \sup_{s\in [0, t]} | \frac{_1}{^{b_n}}\Phi_n (b_n s) \! -\! 2s | > 2 \epp\big) & \leq &  
\bP \big( \!\!\!\!\!\!\!  \sup_{\quad s\in [0, b_nt]}\!\!\!\!  | R^n_s  \! -\! s | >  b_n\epp/2 \big)+ \bP \big( \!\!\!\!  \!\!\! \sup_{\quad s\in [0, b_nt]}\!\!\!\!  \cH^{\bw_n}_s  > b_n \epp  \big) \\
& \leq & 4\epp^{-1} \sqrt{t/b_n}+ \frac{t^2b_n \sigma_2 (\bw_n)}{\epp \sigma_1 (\bw_n)^2} + \bP \big( \sup_{s\in [0, t]} \frac{_{a_n}}{^{b_n}}\cH^{\bw_n}_{b_ns}  > a_n \epp  \big) \,.
\end{eqnarray*}
Thus, by (\ref{apriori}) and (\ref{jtcvYHP}) in Theorem \ref{HYcvth}, we get $\lim_{n\rightarrow \infty} \bP \big( \sup_{s\in [0, t]} | \frac{_1}{^{b_n}}\Phi_n (b_n s) \! -\! 2s | > 2 \epp\big)\! =\! 0$. This prove that $\frac{_1}{^{b_n}}\Phi_n (b_n \cdot)$ converges to 
$2\mathrm{Id}$ in probability on $\bC([0, \infty), \bbR)$, where $\mathrm{Id}$ stands for 
the identity map on $[0, \infty)$. 
%\tt{Note that for all $s\ino [0, 2 \mathbf{j}_n)$, $  \xi^n_{\lceil s \rceil }\! = \! \inf \{ t\ino [0, \infty) \! :\! \Phi_n (t) \! \geq \! s \}$.
%Then, standard} 
Then, standard arguments also imply that $\frac{_1}{^{b_n}}\phi_n (b_n \cdot)$ converges to $\frac{_1}{^2}\mathrm{Id}$ in probability on $\bC([0, \infty), \bbR)$. We also note that on any interval $[k, k+1]$ where $k$ is an integer, 
$\cC^{\bw_{n}}_{t}$ is a linear interpolation between $\cC^{\bw_{n}}_{k}$ and $\cC^{\bw_{n}}_{k+1}$. %\margmm{new}
These convergences combined with Theorem \ref{HYcvth} imply 
\begin{equation}
\label{avecont}
\big( \frac{_{_1}}{^{^{a_n}}} A^{\bw_n}_{b_n \cdot }\,  , \frac{_{_1}}{^{^{a_n}}} Y^{\bw_n}_{b_n \cdot } ,  \,     \frac{_{_{a_n}}}{^{^{b_n}}} \cH^{\bw_n}_{b_n \cdot}   \, , \frac{_{_{a_n}}}{^{^{b_n}}} \cC^{\bw_n}_{b_n \cdot} \, ,   \frac{_{_1}}{^{^{b_n}}} \Ptt_{ \bw_n}\, ,   \frac{_{_1}}{^{^{b_n}}} \Phi_n (\Ptt_{ \bw_n})\big) 
 \underset{n\rightarrow \infty}{-\!\!\! -\!\!\! \longrightarrow}  \big(  A, Y, \cH , \cH_{\cdot /2},  \Ptt ,2 \Ptt \big), 
\end{equation}    
weakly on the appropriate space. 

We now deal with the excursions of $\cC^{\bw_n}$ above $0$, that are the contour processes of the spanning trees $\cT^{\bw_n}_{\!\! k }$, $1\! \leq \! k \! \leq \! \mathbf{q}_{\bw_n}$, of the $\mathbf{q}_{\bw_n}$ 
connected components of $\cG_{\! \bw_n}$; recall that the $\cT^{\bw_n}_{\!\! k }$ are also the connected components obtained from the tree $\cT_{\!\! \bw_n}$ after removing its root.
%\cT^{\bw_n}_{\!\! k }$, $1\! \leq \! k \! \leq \! \mathbf{q}_{\bw_n}$. 
%\margmm{rewritten}
%\mm{Recall from Section \ref{YLIFOsec} that these are the connected component obtained from the tree $\cT_{\!\! \bw_n}$ after removing its root}. 
Recall from (\ref{excdibor}) that $[l^{\bw_n}_k, r^{\bw_n}_k)$ are the excursion intervals of $\cH^{\bw_n}$ above $0$: namely, $\bigcup_{1\leq k\leq \bq_{\bw_n}} [l^{\bw_n}_k, r^{\bw_n}_k) \! =\!   \{  t\ino [0, \infty)\!  :\!  \cH^{\bw_n}_t \! >\! 0  \}$. 
%\begin{equation}
%\label{excdibor}
% \bigcup_{1\leq k\leq \bq_{\bw_n}} [l^{\bw_n}_k, r^{\bw_n}_k) = \big\{  t\ino [0, \infty) : \cH^{\bw_n}_t >0  \big\}= \big\{  t\ino [0, \infty) : Y^{\bw_n}_t > \JJ^\bw_t  \big\} 
%\end{equation} 
%To that end, \mm{we first recall} the notation $\{ t\ino [0, \infty) \! : \! \cH^{\bw_n}_t \! >\! 0\} \! = \! \bigcup_{1 
%\leq  k  \leq  \mathbf{q}_{\bw_n}} [l^n_k , r^n_k)$ for the excursion intervals of $\cH^{\bw_n}$ above $0$ 
Recall that the excursion intervals are listed in the decreasing order of their lengths; recall that $\mathtt{H}^{\bw_n}_k (t)\! = \! 
\cH^{\bw_n} ((l^{\bw_n}_k + t)\! \wedge r^{\bw_n}_k)$, $t\ino [0, \infty)$, is the $k$-th longest excursion process of $\cH^{\bw_n}$ above $0$. Recall from (\ref{defPiwk}) that $\Ptt^{\bw_n}_k \! = \! ((s^{n,k}_p, t^{n,k}_p); 1\! \leq \! k \! \leq \! \mathbf{p}_k^n)$ is the sequence of pinching times falling into the $k$-th longest excursion. Then recall that $\bm^{\bw_n}\! = \! \sum_{j\geq 1} w^{_{(n)}}_{^j} \delta_j $ and recall that $\bm^{\bw_n}_k$ is the restriction to $\cT^{\bw_n}_{\!\! k}$ of $\bm^{\bw_n}$. Recall that 
$\big(\cT^{\bw_n}_{\!\! k} , d_{\mathrm{gr}} , \varrho_{k}^{\bw_n} , \bm^{\bw_n}_k)$ stands for the 
measured tree coded by $\mathtt{H}_k^{\bw_n}$ and  that 
$\big( \cG_k^{\bw_n}, d_{k}^{\bw_n}, \varrho_k^{\bw_n}, \bm_k^{\bw_n}\big)$ is the measured graph coded by $\mathtt{H}_k^{\bw_n}$ and the pinching setup $(\Ptt_k^{\bw_n}, 1)$: 
namely, $ \cG_k^{\bw_n}$ is isometric to the graph $G( \mathtt{H}_k^{\bw_n}, \Ptt_k^{\bw_n}, 1)$ as defined in (\ref{defgrgraf}) and it is  the $k$-th largest (with respect to the measure $\bm^{\bw_n}$) connected component of $\cG_{\! \bw_n}$. 
We next set for all $k \ino \{ 1, \ldots ,  \mathbf{q}_{\bw_n}\}$, $\overline{l}^n_k\! = \! \Phi_n (l^{\bw_n}_k)$, 
$\overline{r}^n_k\! = \! \Phi_n (r^{\bw_n}_k)$, 
$$  \mathtt{C}^{\bw_n}_k (t)\! = \! \cH^{\bw_n} \big( \phi_n ( (\overline{l}^n_k + t) \wedge \overline{r}^n_k) \big) = \mathtt{H}^{\bw_n}_k \big( \phi_n ( \overline{l}^n_k + t)\! -\! l_k^{\bw_n} \big)  \; \textrm{and} \;  $$

$$\overline{\Ptt}^{\bw_n}_k \! = \! \big( (\Phi_n (l^{\bw_n}_k + s^{n,k}_p)\! -\! \overline{l}^n_k \, , \,  \Phi_n (l^{\bw_n}_k + t^{n,k}_p)\! -\! \overline{l}^n_k)\big)_{1\leq  k  \leq  \mathbf{p}_k^n} .$$
Then, we easily check the following: 
\begin{compactenum}

\smallskip

\item[$(i)$] $\{ t\ino [0, \infty)\! : \! \cC^{\bw_n}_t \! >\! 0 \} \!  = \!  \bigcup_{1 
\leq  k  \leq  \mathbf{q}_{\bw_n}} [\overline{l}^n_k , \overline{r}^n_k)$. 
%\fmm{Here, needs clarification on the definition of $\cC^{\bw_{n}}$. }

\smallskip

%\margmm{see footnote}
\item[$(ii)$] $ \mathtt{C}^{\bw_n}_k (\cdot) -1$ is the contour process of $\cT^{\bw_n}_{\!\! k}$. 
We denote by $\bnu^{\bw_n}_k$ the measure that the contour process 
induces on $\cT^{\bw_n}_{\!\! k}$: namely, 
$\big(\cT^{\bw_n}_{\!\! k} , d_{\mathrm{gr}} , \varrho_{k}^{\bw_n} , \bnu^{\bw_n}_k)$ 
is the measured tree coded by $ \mathtt{C}^{\bw_n}_k (\cdot) -1$. 
%\fmm{This part needs checking! I feel that if you take $ \mathtt{C}^{\bw_n}_k (\cdot) -1$, then the measure $\bnu^{\bw_n}_k$ is essentially $2\times$ the uniform measure on the edges projected to the nearest vertices. Then the total mass $\bnu^{\bw_n}_k(\cG_k^{\bw_n})$ should be $2 (\bmu^{\bw_n}_k(\cG_k^{\bw_n})-1)$, which is not consistent with what's below. Check also Lemma \ref{cntourmes}: \tt{I think it is correct now! }}

\smallskip

\item[$(iii)$]  $\big( \cG_k^{\bw_n}, d_{k}^{\bw_n}, \varrho_k^{\bw_n}, \bnu_k^{\bw_n}\big)$ is isometric to 
$ G \big(  \mathtt{C}^{\bw_n}_k (\cdot) -1,  \overline{\Ptt}^{\bw_n}_k, 1)$. 

\smallskip

\end{compactenum}
Since $(b_n^{-1} \Phi_n (b_n \cdot) , b_n^{-1} \phi_n (b_n \cdot)) \! \rightarrow \! (2\mathrm{Id}, \frac{_1}{^2} \mathrm{Id})$ in probability on $\bC([0, \infty), \bbR)^2$, we easily get from Theorem \ref{excucvth} 
that  
\begin{equation}
\label{Conttexcu}
\big(\big(  \frac{_{_{a_n}}}{^{^{b_n}}} \mathtt{C}^{\bw_n}_k (b_n  \cdot) , \, \frac{_{_{1}}}{^{^{b_n}}}  \overline{l}_k^{n} , \, \frac{_{_{1}}}{^{^{b_n}}} \overline{r}_k^{n}, \,  \frac{_{_{1}}}{^{^{b_n}}}\overline{\Ptt}_k^{\bw_n} \big) \big)_{k\geq 1}
\;  \underset{n\rightarrow \infty}{-\!\!\! -\!\!\! -\!\!\! \longrightarrow} \; \big(\big(\Htt_k (\cdot /2), 2l_k, 2r_k, 2\Ptt_k\big) \big)_{k\geq 1}
\end{equation} 
weakly on $(\bC([0, \infty), \bbR)\! \times \! [0, \infty)^2 \! \times \!  (\bbR^2)^{\bbN^*})^{\bbN^*}$ equipped with the product topology, with obvious notation. 
%\margmm{new}
Then, by Lemma \ref{codconGHP} and the same argument as 
in the proof of (\ref{graphbwn}), we get 
\begin{equation}
\label{graphbnu}
\big(\big( \cG_{\! k}^{\bw_n} ,  \frac{_{_{a_n}}}{^{^{b_n}}}d_{k}^{\bw_n} , \varrho_k^{\bw_n}, \frac{_{_{1}}}{^{^{b_n}}}\bnu_k^{\bw_n}  \big) \big)_{k\geq 1}
\;  \underset{n\rightarrow \infty}{-\!\!\! -\!\!\! -\!\!\! \longrightarrow} \; \big(\big( \mathbf{G}_{k} , \mathrm{d}_{k}, \varrho_{k} , 2\bm_{k} \big) \big)_{k\geq 1}
\end{equation} 
weakly on $\bbG^{\bbN^*}$ equipped with the product topology. Next, we prove the following.
%\margmm{rewritten}
\begin{lem}
\label{cntourmes} Let us denote by $\bmu^{\bw_n}_{k}$ the counting measure on 
$\cG^{\bw_n}_{\! k}$. 
%\ts{: namely, $\bmu^{\bw_n}_{k}\! = \! \sum_{j\in \cV (\cG^{\bw_n}_{\! k})}\delta_j $}. 
We equip $\cG^{\bw_n}_{\! k}$ with the graph distance and for all non-empty subsets 
of vertices $A$ we denote by $A^{(1)}$ the set of vertices at graph-distance at most $1$ from $A$. 
Then 
\begin{equation}
\label{arguther}
 \bnu^{\bw_n}_{k} (A) \leq  2\bmu^{\bw_n}_{k} \big( A^{(1)} \big) +1\quad \textrm{and} \quad  2\bmu^{\bw_n}_{k} ( A)  \leq  \bnu^{\bw_n}_{k} \big( A^{(1)} \big) +1. 
\end{equation}
 % the Prokhorov distance between $\bnu^{\bw_n}_{k}$ and $2\bmu^{\bw_n}_{k}$ is less or equal to $1$.  
\end{lem}
\noi
\textbf{Proof.} Since adding edges only diminishes the graph distance, it is sufficient to prove (\ref{arguther}) on $\cT^{\bw_n}_{\! \! k}$ equipped with the graph-distance $d_{\mathrm{gr}}$. 
Recall that $\varrho^{\bw_n}_k$ is the root of $\cT^{\bw_n}_{\! \! k}$. To simplify notation we set 
$\cT\! = \! \cT^{\bw_n}_{\! \! k}$, $\varrho\! = \! \varrho^{\bw_n}_k$, 
$\bnu\! = \! \bnu^{\bw_n}_{k}$, $\bnu^\prime \! = \!  \delta_{\varrho}+ \bnu$ and $\bmu\! = \! \bmu^{\bw_n}_{k}$. 
Since the contour process of $\cT$ crosses twice each edge, we easily get 
$ \bnu^\prime \! = \! \delta_{\varrho}+ 
\sum_{v\in \cT}   \texttt{deg}(v)   \delta_{v} \! = \!  \bmu+ \bmu \circ f^{-1}$
%it easy to see that 
%$$ \bnu^{\bw_n}_{k}= \delta_{\varrho^{\bw_n}_k}+ \!\!\!\!\!\! 
%\sum_{v\in \cV(\cT^{\bw_n}_{\! \! k})}\!\!\!\!\!\!   \texttt{deg}(v)\,   \delta_{v}\! = \!  \bmu^{\bw_n}_{k} + \bmu^{\bw_n}_{k} \circ f^{-1} \; , $$
where $f \! : \! \cT \! \rightarrow \! \cT$ is given by 
$f(v)\! = \! \overleftarrow{v}$ if $v\! \neq \! \varrho$ and $f(\varrho)\! = \! \varrho$. Let $M\! = \! \sum_{v\in \cT} (\delta_{(v,v)}+ \delta_{(v, f(v))})$ that is a measure on $\cT \! \times \! \cT$  such that $M (A \! \times \! \cT)\! = \! 2\bmu (A)$ and $M (\cT \! \times \! A)\! = \! \bnu^\prime (A)$. Then, set 
$D\! = \! \{ (v,v^\prime) \ino \cT \! \times \! \cT \! : \! d_{\mathrm{gr}} (v, v^\prime) \! \leq \! 1 \}$. 
Since $d_{\mathrm{gr}} (f(v), v) \! \leq \! 1$, $M$ is supported on $D$. Next, observe that 
$(A \! \times \! \cT)\cap D \! \subset \! \cT \! \times \! A^{(1)}$ and similarly $D \cap (\cT \! \times \! A) \! \subset \!  A^{(1)} \! \times \! \cT$, which easily entails (\ref{arguther}). \cqfd

\medskip

Since $d^{\bw_n}_k$ is the graph-distance on $\cG^{\bw_n}_{\! k}$, we easily see that on the rescaled space $(\cG^{\bw_n}_{\! k}, \frac{a_n}{b_n} d^{\bw_n}_k)$, (\ref{arguther}) 
implies that $\tfrac{1}{b_n}\bnu^{\bw_n}_{k} (A) \! \leq \!  \tfrac{2}{b_n}\bmu^{\bw_n}_{k} (  A^{_{(a_n/b_n)}}_{^{\, }} )+\tfrac{1}{b_n}$ and 
$\tfrac{2}{b_n}\bmu^{\bw_n}_{k} (A) \! \leq \!  \tfrac{1}{b_n}\bnu^{\bw_n}_{k} ( A^{_{(a_n/b_n)}}_{^{\, }} )+\tfrac{1}{b_n}$, for all subset of vertices $A $. Since $b_n^{-1} \! \leq \! a_n/b_n$ (for all sufficiently large $n$), we get $d_{\mathrm{Pro}} \big( \frac{1}{b_n}  \bnu^{\bw_n}_{k}, \frac{2}{b_n} \bmu^{\bw_n}_{k}\big) \! \leq \!   a_n/b_n$. 
This combined with (\ref{graphbnu}) entails 
$$\big(\big( \cG_{\! k}^{\bw_n} ,  \frac{_{_{a_n}}}{^{^{b_n}}}d_{k}^{\bw_n} , \varrho_k^{\bw_n}, \frac{_{_{1}}}{^{^{b_n}}}2\bmu_k^{\bw_n}  \big) \big)_{k\geq 1}
\;  \underset{n\rightarrow \infty}{-\!\!\! -\!\!\! -\!\!\! \longrightarrow} \; \big(\big( \mathbf{G}_{k} , \mathrm{d}_{k}, \varrho_{k} , 2\bm_{k} \big) \big)_{k\geq 1}$$
weakly on $\bbG^{\bbN^*}$ equipped with the product topology, which easily implies (\ref{globabab}). 

\medskip
\noindent\textbf{End of proof of Theorem \ref{graphcvth}.}
We next make the following additional assumption : $\sqrt{\mathbf{j}_n} / b_n \! \rightarrow \! 0$ and we complete the proof of Theorem \ref{graphcvth}. To that end, it is sufficient to prove that for all fixed $k\! \geq \! 1$, the probability that 
%\margmm{corrected}
$\bmu^{\bw_n}_1 (\cG^{\bw_n}_{\! 1})\! > \! \ldots \! > \! \bmu^{\bw_n}_k (\cG^{\bw_n}_{\! k}) \! > \! 
\max_{j > k } \bmu^{\bw_n}_j (\cG^{\bw_n}_{\! j})$ tends to $1$ as $n\! \rightarrow \! \infty$.  

Recall from Lemma \ref{cntrlll} that $\bE [ Q^2_n] \! \leq \! 4\mathbf{j}_n$. Thus, $Q_n / b_n \! \rightarrow \! 0$ in probability. Recall that $(b_n^{-1} \Phi_n (b_n \cdot) , b_n^{-1} \phi_n (b_n \cdot)) \! \rightarrow \! (2\mathrm{Id}, \frac{_1}{^2} \mathrm{Id})$ in probability on $(\bC([0, \infty), \bR))^2$. By Slutzky's theorem, 
we get a joint convergence of $(b_n^{-1} Q_n, b_n^{-1} \Phi_n (b_n \cdot) , b_n^{-1} \phi_n (b_n \cdot))$ with (\ref{Conttexcu}). Without loss of generality (but with a slight abuse of notation), by Skorokod's representation theorem we can assume that the following convergence 
\begin{eqnarray}
\label{rConttexcu}
\big(  \frac{_{_{1}}}{^{^{b_n}}} Q_n ,  \frac{_{_{1}}}{^{^{b_n}}} \Phi_n (b_n \cdot) ,  \frac{_{_{1}}}{^{^{b_n}}} \phi_n (b_n \cdot) \!\!\!\!\!\!\!\!\!\!\!\!\! & ; & \!\!\!\!\!\!\!\!\!\!\!\!\!\! \big(\big(  \frac{_{_{a_n}}}{^{^{b_n}}} \mathtt{C}^{\bw_n}_k (b_n  \cdot) , \, \frac{_{_{1}}}{^{^{b_n}}}\overline{l}_k^{n} , \, \frac{_{_{1}}}{^{^{b_n}}} \overline{r}_k^{n} , \,  \frac{_{_{1}}}{^{^{b_n}}}\overline{\Ptt}_k^{\bw_n} \big) \big)_{k\geq 1}\big) \\
 & \underset{n\rightarrow \infty}{-\!\!\! -\!\!\! -\!\!\! \longrightarrow} & \big( 0,2\mathrm{Id}, \frac{_1}{^2} \mathrm{Id};  \big(\big(\Htt_k (\cdot /2), 2l_k, 2r_k, 2\Ptt_k\big) \big)_{k\geq 1}\big) \nonumber
\end{eqnarray} 
holds almost surely on the appropriate space.

Recall notation $\zeta_k \!= \!  r_k  - l_k$, 
$\zeta^{\bw_n}_k \!= \!  r^{\bw_n}_k \! -\! l^{\bw_n}_k=\bm^{\bw_n}_k (\cG^{\bw_n}_{\! k})$ and set $\overline{\zeta}^n_k \!= \!  \overline{r}^n_k \! - \overline{l}^n_k=\bnu^{\bw_n}_k (\cG^{\bw_n}_{\! k})$. First, we easily derive from the argument of the proof of (\ref{arguther}) that 
$\bnu^{\bw_n}_k (\cG^{\bw_n}_{\! k})\! = \! 2\bmu^{\bw_n}_k (\cG^{\bw_n}_{\! k}) \tt{+1}$. 
Let $\sigma_n$ be a permutation of $\{ 1, \ldots , \mathbf{q}_k^{\bw_n}\}$ such that $(\overline{\zeta}^n_{\sigma_n (k)})_{1\leq k \leq \mathbf{q}_k^{\bw_n}}$ is nonincreasing. To complete the proof of Theorem \ref{graphcvth}, it is then sufficient to prove that for all $k\! \geq \! 1$, there exists $n_k$ such that for all $n\! \geq \! n_k$, $\sigma_n (k)\! =\! k$.

To prove that, we then fix $k\! \geq \! 1$ and we recall that $\zeta_1 \!> \! \ldots  \! > \! \zeta_k\! >\! \zeta_{k+1}$ so it makes sense to fix $\epp \ino (0, \infty)$ such that  
$ \epp \! < \! \frac{_1}{^3} \min_{1\leq j\leq k} (\zeta_{j} \! -\! \zeta_{j+1})$. 
Observe first that (\ref{rConttexcu}) implies that for all $j\geq 1$,   
$b_n^{-1}\zeta^{\bw_n}_j\! \rightarrow \zeta_j$ and $b_n^{-1}\overline{\zeta}^n_j\! \rightarrow 2\zeta_j$. Therefore, there exists $n_k\ino \bbN$ 
such that for all $n\! \geq \! n_k$, 
\begin{equation}
\label{ronflium}
b_n^{-1} (4Q_n +1)+ \max_{1\leq j\leq k+1} |b_n^{-1}\overline{\zeta}^n_j \! -\! 2\zeta_j| + \max_{1\leq j\leq k+1} |b_n^{-1}\zeta^{\bw_n}_j \! -\! \zeta_j|  \;  <  \epp \; .
\end{equation}
%and for all $ j\ino \{ 1, \ldots , k +1\}$, $|b_n^{-1}\overline{\zeta}^n_j \! -\! 2\zeta_j| \! < \! \epp$ and $ |b_n^{-1}\tt{\zeta^{\bw_n}_j} \! -\! \zeta_j|\! < \! \epp$. 
%%%%It implies for all $n\! \geq \! n_k$, that 
%%%$\zeta^{n}_1 \!> \! \ldots  \! > \! \zeta^{n}_k$ and $ \overline{\zeta}^{n}_1 \!> \! \ldots  \! > \!  \overline{\zeta}^{n}_k$ and thus  
%%%$\sigma^{-1}_n (k) \! \geq \! k$, where $\sigma_n^{-1}$ stands for the inverse permutation. 
Then, we fix $n\! \geq \! n_k$ and for all $j\ino \{ 1, \ldots , k\}$, Lemma \ref{cntrlll} and (\ref{voiliii}) imply 
\begin{eqnarray*}
 \overline{\zeta}^n_{\sigma_n (j)} & = & 
 \overline{r}^{n}_{\sigma_n (j)} - \overline{l}^{n}_{\sigma_n (j)}   
= \Phi_n (r^{\bw_n}_{\sigma_n (j)} ) \! -\! \Phi_n (l^{\bw_n}_{\sigma_n (j)}) \\ 
&= & 2R^n (r^{\bw_n}_{\sigma_n (j)}) \! -\! 2R^n (l^{\bw_n}_{\sigma_n (j)}) - \cH^{\bw_n} (r^{\bw_n}_{\sigma_n (j)}) + \cH^{\bw_n} (l^{\bw_n}_{\sigma_n (j)}) \\ 
& = & 
 2R^n (r^{\bw_n}_{\sigma_n (j)}) \! -\! 2R^n (l^{\bw_n}_{\sigma_n (j)}) + 1 \overset{\textrm{by (\ref{fincrt})}}{\leq}  2 \zeta^{\bw_n}_{\sigma_n (j)} + 
 4 Q_n + 1 .
  \end{eqnarray*}
 Thus, $2b_n^{-1}\zeta^{\bw_n}_{\sigma_n (j)}\! \geq \! b_n^{-1}  \overline{\zeta}^n_{\sigma_n (j)} \! -\! \epp$. Moreover, 
 $$ b_n^{-1}  \overline{\zeta}^n_{j} \geq 2\zeta_j -\epp \geq 2 (\zeta_j \! -\! \zeta_{j+1}) + 2\zeta_{j+1} -\epp \geq 6 \epp +  (2\zeta_{j+1}+ \epp ) -2\epp \geq 4\epp + b_n^{-1}  \overline{\zeta}^n_{j+1} . $$
% Now observe that $n\! \geq \! n_k$ implies 
 %$$
% \margmm{rewritten: \tt{it was correct}}
% \ts{\mm{$\overline{\zeta}^{n}_1, \cdots, \overline{\zeta}^{n}_k$ are all distinct}}
Namely, for all $n\! \geq \! n_k$, $ \overline{\zeta}^{n}_1 \!> \! \ldots  \! > \!  \overline{\zeta}^{n}_k$. Then, set $S\! = \! \{ \overline{\zeta}^{n}_\ell ; 1\! \leq \! \ell \! \leq \! \mathbf{q}_{\bw_n} \}$; the previous inequality implies that for all $j\ino \{ 1, \ldots , k\}$, 
$\# (S\cap [ \overline{\zeta}^{n}_j, \infty)) \! \geq \! j\! = \! \# (S \cap [\overline{\zeta}^n_{\sigma_n (j)} , \infty)) $. 
This entails $ \overline{\zeta}^n_{\sigma_n (j)} \! \geq \!  \overline{\zeta}^n_{j}$, $j\ino \{ 1, \ldots , k\}$. Namely, for all $j\ino \{ 1, \ldots, k\}$, we get 
$$ 2b_n^{-1}\zeta^{\bw_n}_{\sigma_n (j)}\! \geq \! b_n^{-1}  \overline{\zeta}^n_{\sigma_n (j)} \! -\! \epp >  b_n^{-1}  \overline{\zeta}^n_{j} \! -\! \epp  >2\zeta_j -4 \epp.  $$
Consequently, $b_n^{-1}\zeta^{\bw_n}_{\sigma_n (j)} \! >\! \zeta_j \! -\! 2\epp$. This implies that $\sigma_n (j) \! \leq \! j$. \textit{Indeed}, suppose that $\sigma_n (j) \! \geq \! j+1$; thus $\zeta^{\bw_n}_{j+1} \! \geq \! \zeta^{\bw_n}_{\sigma_n(j)}$ and the previous inequality combined with (\ref{ronflium}) would 
entail  $\zeta_{j+1} \! +\! \epp \! >\!   b_n^{-1}\zeta^{\bw_n}_{j+1}  \! \geq \! b_n^{-1}\zeta^{\bw_n}_{\sigma_n(j)} \! >\! \zeta_j \! -\! 2\epp$, which would contradict 
$ \epp \! < \! \frac{_1}{^3} \min_{1\leq \ell \leq k} (\zeta_{\ell} \! -\! \zeta_{\ell+1})$. 
Thus, for all $n\! \geq \! n_k$ and for all $j\ino \{ 1, \ldots , k\}$, $\sigma_n (j) \! \leq \! j$, which easily entails that $\sigma_n (j)\! = \! j$, which completes the proof. \cqfd

\subsection{Proof of the limit theorems for the Markovian processes.}
\label{Markosec}
\subsubsection{Convergence of the Markovian queueing system: the general case.}
\label{CVmarkogne}
We say that a $\bbR$-valued spectrally positive L\'evy process 
$(R_t)_{t\in [0, \infty)}$ with initial value $R_0\! = \! 0$ is \textit{integrable} if for at least one $t\ino (0, \infty)$ we have $\bE [|R_t|] \! < \! \infty$. 
It implies that $\bE [|R_t|] \! < \! \infty$ for all $t\ino (0, \infty)$. 
We recall  
from Section \ref{RWLevApp} in Appendix that there is a one-to-one correspondence between the laws of $\bbR$-valued spectrally positive L\'evy processes $(R_t)_{t\in [0, \infty)}$ with initial value $R_0\! = \! 0$ that are integrable and the triplets 
$(\alpha, \beta, \pi)$ where $\alpha \ino \bbR$, $\beta \ino [0, \infty)$ and $\pi$ is a Borel-measure on $(0, \infty)$ such that $\int_{(0, \infty)} \pi (dr) \,  (r\! \wedge \! r^2) \! < \! \infty$. More precisely, the correspondence is given by the Laplace exponent of spectrally positive L\'evy processes: namely, for all $t, \lambda \in [0, \infty)$,  
%\margmm{Hmm...not right. Also see my comment on \eqref{rLKform}} 
\begin{equation}
\label{LKform} 
\bE \big[ e^{-\lambda R_t}\big]\! = \! e^{t\psi_{\alpha, \beta, \pi} (\lambda)}, \; \textrm{where} \quad  \psi_{\alpha , \beta , \pi} (\lambda) = \alpha \lambda + \frac{_{1}}{^{2}} \beta \lambda^2 + \int_{(0, \infty)} \!\!\!\!\! \!\!\! (e^{-\lambda r}\!  -\! 1+ \lambda r) \, \pi(dr). 
\end{equation}

The main result used to obtain the convergence of branching processes is a Theorem due to Grimvall \cite{Gr74}, that is recalled 
in Theorem \ref{Grimvresu}: it states the convergence of rescaled Galton-Watson processes to Continuous State Branching Processes (CSBP for short). 
We say that a process $(Z_t)_{t\in [0, \infty)}$ is an {\it integrable} CSBP if it is a $[0, \infty)$-valued Feller Markov process obtained from spectrally positive L\'evy processes via Lamperti's time-change which further satisfies $\mathbf E[Z_{t}]<\infty$ for all $t\in [0, \infty)$. The law of such a CSBP is completely characterised by the 
Laplace exponent of its associated Lévy process that is usually called the \textit{branching mechanism} of the CSBP, which is necessarily of the form (\ref{LKform}): see Section \ref{CSBPApp} for a brief account on CSBP.

Let $\bw_n\ino \elldo_f$, $n\ino \bbN$. 
Recall from (\ref{defnuw}) that $\nu_{\bw_n}= \sigma_1 (\bw_n)^{-1}
 %\frac{1}{{\sigma_1 (\bw_n)}}
 \sum_{j\geq 1} w^{_{(n)}}_{^j} \delta_{j} $ and recall from (\ref{rmupoissw}) that  for all $k\ino \bbN$, $\mu_{\bw_n} (k) \! =\!   \sigma_1 (\bw_n)^{-1} \sum_{j\geq 1} 
 (w^{_{(n)}}_{^j})^{k+1} \exp (-w^{_{(n)}}_{^j})/k!$. 
%Define the law $\nu_{\bw_n}$ and $\mu_{\bw_n}$ by setting: 
%\begin{equation}
%\label{defpin}
%\nu_{\bw_n}= \frac{1}{{\sigma_1 (\bw_n)}}\sum_{j\geq 1} w^{_{(n)}}_{^j} \delta_{j} \quad \textrm{and} \quad \forall k\ino \bbN, \quad  \mu_{\bw_n} (k) \! =\!   \sum_{j\geq 1} \frac{(w^{_{(n)}}_{^j})^{k+1}}{\sigma_1 (\bw_{n})\,  k!}\,  e^{-w^{_{(n)}}_{^j}} . 
%\end{equation}  
Recall from Section \ref{connecMqu} the definition of the Markovian LIFO-queueing system associated with the set of weights $\bw_n$: clients arrive at unit rate; 
each client has a \textit{type} that is a positive integer; the amount of service required by a client of type $j$ is $w^{_{(n)}}_j$; the types are i.i.d.~with law $\nu_{\bw_n}$. 
If one denotes by $\tau^{n}_k$ the time of arrival of the $k$-th client in the queue and by $\Jtt^n_k$ his type, then the queueing system is entirely characterised by $ \ccX_{\bw_n}  \! = \sum_{k\geq 1} \delta_{(\tau^n_k , \Jtt^n_k)}$ 
%\begin{equation}
%\label{rrXPoisdef}
% \ccX_{\bw_n}  \! = \sum_{k\geq 1} \delta_{(\tau^n_k , \Jtt^n_k)},
%\end{equation} 
that is a Poisson point measure on $[0, \infty) \! \times \! \bbN^*$ with intensity 
$\ell  \otimes  \nu_{\bw_n}$, where $\ell $ stands for the Lebesgue measure on $[0, \infty)$.  
Next, for all $j\ino  \bbN^*$ and all $t\ino [0, \infty)$, we introduce the following:  
\begin{equation}
\label{rrXwdef}
N^{\bw_n}_j (t)\! = \! \sum_{k\geq 1} \un_{\{ \tau_k^n \leq t \, ; \, \Jtt^n_k  =  j \}}  \quad \textrm{and} \quad 
X^{\bw_n}_t  \! = \!  -t + \sum_{k\geq 1} w^{_{(n)}}_{\Jtt^n_k}\un_{[0, t]} (\tau^n_k) \! = \! -t + \sum_{j\geq 1} w^{_{(n)}}_j N^{\bw_n}_j (t) . 
% \;  \quad \textrm{and} \quad I^{\bw_n}_t \! = \! \inf_{s\in [0, t]} X^{\bw_n}_s . 
\end{equation} 
%Note that when the server is idle, the load is counted negatively, hence the name \textit{algebraic}. The actual load of the server at time $t$ is then 
%$X^\bw_t \! -\! \inf_{s\in [0, t]} X^\bw_s $. 
%Then, $X^{\bw_n}_t \! -\! I^{\bw_n}_t$ is interpreted as the load of the Markovian queueing system at time $t$. 
Observe that $(N^{\bw_n}_{^j})_{j\geq 1}$ are independent homogeneous Poisson processes with rates
$w^{_{(n)}}_j \!  / \sigma_1 (\bw_n)$ and $X^{\bw_n}$ is a c\`adl\`ag spectrally positive L\'evy process.

   Let $a_n, b_n\ino (0, \infty)$, $n\ino \bbN$ be two sequences that satisfy the following conditions. 
 \begin{equation}
\label{aprioriri}
a_n \; \,   \textrm{and} \; \,  \frac{b_n}{a_n} \; \, \underset{n\rightarrow \infty}{-\!\!\! -\!\!\! \longrightarrow}\,  \infty , \quad  
\frac{b_n}{a^2_n} \, \;  \underset{n\rightarrow \infty}{-\!\!\! -\!\!\! \longrightarrow}\,\;  \beta_0 \ino [0, \infty), 
\quad   \textrm{and}   \quad  \sup_{n\in \bbN} \frac{w^{_{(n)}}_{^1} }{a_n} \! < \! \infty . 
\end{equation}

\vspace{-5mm}

\begin{rem}
\label{kakapa}
\textit{It is important to note that these assumptions are weaker than} (\ref{apriori}): namely, we temporarily \textbf{do not} assume that $\frac{a_n b_n}{\sigma_1 (\bw_n)}\! \rightarrow \! \kappa \ino (0, \infty)$, which 
explains why the possible limits in the theorem below are more general. \cq 
\end{rem}
\begin{thm}
\label{cvbranch} Let $\bw_n \ino  \elldo_f$ and $a_n , b_n \ino (0, \infty)$, $n\ino \bbN$, satisfy (\ref{aprioriri}). 
Recall from (\ref{rrXwdef}) the definition of  $X^{\bw_n}_t \! $; recall from (\ref{rmupoissw}) the definition of $\mu_{\bw_n}$ and let $(Z^{_{(n)}}_k)_{k\in \bbN}$ be a Galton-Watson process with offspring distribution $\mu_{\bw_n}$ and initial state 
$Z^{_{(n)}}_0\! =\!  \lfloor a_n \rfloor$. Then, the following convergences are equivalent. 
 \begin{compactenum}

\smallskip

%\margmm{Changed conservative to integrable}
\item[$\bullet$] $\mathrm{(I)}$  $\quad$ $\big( \frac{_1}{^{a_n}} Z^{_{(n)}}_{\lfloor b_n t /a_n  \rfloor } 
\big)_{t\in [0, \infty)} \! \longrightarrow \! (Z_t)_{t\in [0, \infty)}$ weakly on $\bD([0, \infty), \bbR)$.

\smallskip

\item[$\bullet$] $\mathrm{(II)}$  $\quad$ $\big( \frac{_1}{^{a_n}} X^{\bw_n}_{ b_n t } 
\big)_{t\in [0, \infty)} \! \longrightarrow \! (X_t)_{t\in [0, \infty)}$ weakly on $\bD([0, \infty), \bbR)$.

\smallskip

\end{compactenum}
If $\mathrm{(I)}$ or $\mathrm{(II)}$ holds true, then $Z$ is necessarily an integrable CSBP and $X$ is an integrable $(\alpha, \beta, \pi)$-spectrally positive L\'evy process (as defined at the beginning of Section \ref{CVmarkogne}) whose Laplace exponent is the same as the branching mechanism of $Z$. 
Here $(\alpha, \beta, \pi)$ necessarily 
satisfies: 
\begin{equation}
\label{frnurz}
\beta \geq \beta_0  \quad \textrm{and} 
\quad  \exists\, r_0\ino (0, \infty) \; \textrm{such that} \; \pi ((r_0, \infty)) \! = \! 0\; , 
\end{equation}
which implies  $\int_{(0, \infty)} r^2 \, \pi (dr) \! < \! \infty$. Moreover, $ \mathrm{(I)}\! \Leftrightarrow \! \mathrm{(II)}\! \Leftrightarrow \!\mathrm{(IIIabc)}\! \Leftrightarrow \! \big( \mathrm{(IIIa)} \& \mathrm{(IV)} )$ where: 
\begin{compactenum}

\smallskip

\item[$\bullet$] $\mathrm{(IIIa)}$  {\small  $\;$ $\displaystyle \frac{b_n}{a_n} \Big( 1\! -\! \frac{\sigma_2 (\bw_n)}{\sigma_1 (\bw_n)} \Big)\;  \longrightarrow   \; \alpha .$}

\smallskip

\item[$\bullet$] $\mathrm{(IIIb)}$ {\small  $\;$  $\displaystyle \frac{b_n}{(a_n)^2}  \frac{\sigma_3 (\bw_n)}{\sigma_1 (\bw_n)} \;  %\underset{^{n\rightarrow \infty}}{-\!\!\! -\!\! \! -\!\! \! \longrightarrow}  
\longrightarrow  
\; \beta + \int_{(0, \infty)} \!\!\!\!\!\!\!\! \! r^2 \, \pi (dr) $. }

\smallskip

\item[$\bullet$] $\mathrm{(IIIc)}$   {\small    $\displaystyle \;  \frac{a_n b_n}{\sigma_1 (\bw_n)} \sum_{j\geq 1} \frac{w_j^{(n)}}{a_n} f \big(  w_j^{(n)}\! / a_n  \big) \longrightarrow   \!\!  \int_{(0, \infty)} \!\! \!\!\!\!\!\!\!\! \! f(r) \,  \pi (dr), \,   $} for all continuous bounded 
{\small $f\! :\!  [0, \infty) \! \rightarrow\!  \bbR$} vanishing in a neighbourhood of {\small $0$}. 

\item[$\bullet$] $\mathrm{(IV)}$   {\small $\quad$  $\displaystyle \frac{a_n b_n}{\sigma_1 (\bw_n)} \sum_{j\geq 1} \frac{w_j^{(n)}}{a_n} 
\big( e^{-\lambda w^{(n)}_j /a_n}\!\! -\! 1 + \lambda\,  w^{(n)}_j \!\! /a_n \big) \longrightarrow   \psi_{\alpha , \beta , \pi} (\lambda)  - \alpha \lambda,\; $}  for all {\small $\lambda \ino (0, \infty)$}, where {\small $\psi_{\alpha, \beta, \pi}$} is defined by (\ref{LKform}). 
\end{compactenum}
\end{thm}
\noi
\textbf{Proof.} We easily check that 
$(X^{\bw_n}_{ b_n t } /a_n)_{t\in [0, \infty)}$ is a $(\alpha_n, \beta_n, \pi_n)$-spectrally positive L\'evy process where 
$$ \alpha_n \! =\! \frac{b_n}{a_n} \Big( 1\! -\! \frac{\sigma_2 (\bw_n)}{\sigma_1 (\bw_n)} \Big), \quad \beta_n \! = \! 0 \quad \textrm{and} \quad \pi_n=  \frac{a_n b_n}{\sigma_1 (\bw_n)}  \sum_{j\geq 1} \frac{w_j^{(n)}}{a_n} \, \delta_{w^{(n)}_j \! / a_n} .$$
We immediately see that $\beta_n +\!  \int  r^2  \pi_n (dr) \! = \!   b_n \sigma_3 (\bw_n)/ a_n^2\sigma_1 (\bw_n)$. Then, Theorem \ref{cvLevycar} implies that $\mathrm{(II)} \! \Leftrightarrow \! \mathrm{(IIIabc)}$. 
We then apply Lemma \ref{Laplcv} to $\Delta^n_k\! = \! (X^{\bw_n}_{{k}}\! -\!X^{\bw_n}_{{k-1}})/a_n$ and $q_n \! = \! \lfloor b_n \rfloor$: it shows that the weak limit $X^{\bw_n}_{\lfloor b_n \rfloor }/a_n \! \rightarrow \! X_1$ is equivalent to the convergence of the Laplace exponents $\psi_{\alpha_n, \beta_n , \pi_n} (\lambda) \! \rightarrow \psi_{\alpha, \beta, \pi} (\lambda)$, for all $\lambda \ino [0, \infty)$. Then note that the left member in $\mathrm{(IV)}$ is 
$\psi_{\alpha_n, \beta_n, \pi_n} (\lambda)\! -\! \alpha_n \lambda$. This shows that $\mathrm{(II)} \! \Leftrightarrow \! \big( \mathrm{(IIIa)} \& \mathrm{(IV)} )$. 

  It remains to prove that $\beta \! \geq \! \beta_0$ and that $\mathrm{(I)} \! \Leftrightarrow \! \mathrm{(IIIabc)}$. Let $(\zeta^n_k)_{k\in \bbN}$ be a sequence of i.i.d.~random variables with law $\mu_{\bw_n}$ as defined in (\ref{rmupoissw}). By Theorem \ref{Grimvresu}, $\mathrm{(I)} $ is equivalent to the weak convergence on $\bbR$ of the r.v.~$R_n\! := \! a_n^{-1} \sum_{1\leq k\leq \lfloor b_n \rfloor} \big(\zeta^n_k \! -\! 1 \big)$. We next apply Lemma \ref{Laplcv} to 
$\Delta^n_k:=a_n^{-1} (\zeta^n_k-1)$ $q_n \! = \! \lfloor b_n \rfloor$, which implies that $\mathrm{(I)}$ is equivalent to 

\begin{equation}
\label{LaplRenn}
\exists\, \psi \ino \bC ([0, \infty), \bbR): \quad  \psi (0) \! = \! 0 \quad \textrm{and} \quad \forall \lambda \ino [0, \infty), \; 
L_n (\lambda)\! :=\! \bE \big[ e^{-\lambda R_n} \big]  \underset{^{n\rightarrow \infty}}{-\!\! \! -\!\! \!\longrightarrow} e^{\psi (\lambda)}  .
\end{equation}
We next compute $L_n (\lambda)$ more precisely. To that end, let $(W^n_k)_{k\in \bbN}$ be an i.i.d.~sequence of r.v.~with the same law as $w^{_{(n)}}_{\mathtt{J^n_1}}$, where $\mathtt{J}^n_1$ has law $\nu_{\bw_n}$: namely, 
%for all all measurable function 
$ \bE [ f (W^n_k)] \! = \!  \sigma_1 (\bw_n)^{-1} 
\sum_{j\geq 1} w^{_{(n)}}_{^j} f ( w^{_{(n)}}_{^j})$ for all nonnegative measurable function $f$. 
%$f\! : \! [0, \infty) \! \rightarrow \! [0, \infty)$, 
%$$ \bE \big[ f\big(W^n_k \big) \big]=  \frac{1}{{\sigma_1 (\bw_n)}}\sum_{j\geq 1} w^{_{(n)}}_{^j} f \big( w^{_{(n)}}_{^j}\big)\; .$$
Note that for all $k\ino \bbN$, $\mu_{\bw_n} (k)\! = \! \bE[\, (W^n_1)^k e^{-W^n_1} \! / k   !\, ]$, which implies: 
\begin{equation}
\label{POCal}
 L_n (\lambda)= e^{\lambda \lfloor b_n\rfloor /a_n } \big( \bE \big[ e^{-\lambda \zeta^n_1/ a_n}\big]\big)^{\lfloor b_n\rfloor}=e^{\lambda \lfloor b_n\rfloor /a_n }
\big( \bE \big[ \exp \big( \! -\! W_1^n \big(1-e^{-\lambda/a_n} \big) \big)\big]\big)^{\lfloor b_n\rfloor}  .
\end{equation}
We next set $S^n_1 \! = \! a_n^{-1} \sum_{1\leq k\leq \lfloor b_n \rfloor} \big(W^n_k \! -\! 1 \big)$ and 
$\cL_n (\lambda)\!  = \!  \bE [ \exp(-\lambda S^n_1)]$. 
%$$ \forall t\in [0, \infty), \quad S^n_t = \frac{1}{a_n}\sum_{1\leq k \leq \lfloor b_n t \rfloor} \!\!\! \big( W_k^n -1\big) \qquad \textrm{and} \qquad 
%\cL_n (\lambda)\!  = \!  \bE \big[ e^{-\lambda S^n_1}\big], \quad  \lambda \! \in \! [0, \infty).$$
By (\ref{POCal}), we get:  
$$ \forall \lambda \ino [0, \infty), \quad \cL_n  \big( a_n \big(1-e^{-\lambda/a_n} \big)  \big)= L_n (\lambda)\,  \exp \big( \, \lfloor b_n\rfloor \big(1\! -\! e^{-\lambda/a_n} \big) \! -\!  
\lambda \lfloor b_n\rfloor /a_n  \big)\,. $$
Since $\lfloor b_n\rfloor \big(1\!  -\!  e^{-\lambda/a_n} \big) \! -\!   
\lambda \lfloor b_n\rfloor /a_n + \frac{_1}{^2}b_na_n^{-2} \lambda^2\!=\!   \mathcal{O} (b_na_n^{-3}) \! \rightarrow \! 0$ 
and since $b_{n}/a_{n}^{2}\to \beta_{0}$,  
(\ref{LaplRenn}) is equivalent to 
\begin{equation}
\label{LaplSenn}
\exists \psi_0 \ino \bC ([0, \infty), \bbR): \quad  \psi_0 (0) \! = \! 0 \quad \textrm{and} \quad \forall \lambda \ino [0, \infty), \; 
\lim_{n\rightarrow \infty} \cL_n (\lambda) \! = \! e^{\psi_0 (\lambda)} , 
\end{equation}
and if (\ref{LaplRenn}) or (\ref{LaplSenn}) holds true, then $\psi (\lambda)\! = \! \psi_0 (\lambda) + \frac{_1}{^2} \beta_0 \lambda^2$, for all $\lambda \ino [0, \infty)$. 

Next, by Lemma \ref{Laplcv} applied to 
$\Delta^n_k:=a_n^{-1} (W^n_k-1)$, we see that (\ref{LaplSenn}) is equivalent to the weak convergence 
$S^{{n}}_1\! \rightarrow \! S_1$ in $\bbR$ and Theorem \ref{analyCV} asserts that is equivalent to the conditions 
$(\textit{Rw3abc})$ with $\xi^n_1\! = \! W^n_1 \! -\! 1$: namely, there exists a triplet $(\alpha^*, \beta^*, \pi^*)$ such that $\alpha^* \ino \bbR$, $\beta^* \ino [0, \infty)$, such that there exists $r_0\ino (0, \infty)$ satisfying $\pi^* ([r_0, \infty))\! = \! 0$ and such that the following holds true
\begin{multline}\label{}  \frac{b_n}{a_n} \bE [\xi^n_1]\! = \!  \frac{b_n}{a_n} \Big( \frac{\sigma_2 (\bw_n)}{\sigma_1 (\bw_n)} - 1 \Big)\! \rightarrow -\alpha^*, \;  \frac{b_n}{a_n^2} \mathbf{var}( \xi^n_1) \! = \!  \frac{b_n}{a_n^2}  \frac{\sigma_3 (\bw_n)}{\sigma_1 (\bw_n)} - \frac{b_n}{a_n^2} \Big(\frac{\sigma_2 (\bw_n) }{\sigma_1 (\bw_n)}\Big)^2 \!\!\!   \rightarrow  \!  \beta^* + \! \! \int_{(0, \infty)} \!  \!\!\! \!\!\! \!\!\! r^2 \pi^* (dr) \nonumber \\
\textrm{and} \quad b_n \bE \big[ f\big(\xi^n_1/a_n \big)\big] \! = \!  \frac{a_n b_n}{\sigma_1 (w_n)} \sum_{j\geq 1} \frac{w_j^{(n)}}{a_n} f \Big( \frac{w_{j}^{{(n)}}\! - 1}{ a_n}  \Big) \rightarrow    \int_{(0, \infty)}  \!\!\!\!\!\!\!\! \! f(r) \,  \pi^* (dr), \;  
\end{multline}
 for all continuous bounded $f\! :\!  [0, \infty) \! \rightarrow\!  \bbR$ vanishing in a neighbourhood of $0$. It is easy to see that these conditions are equivalent to $\mathrm{(IIIabc)}$ with $\alpha\! = \! \alpha^*$, $\beta\! = \! \beta_0+ \beta^*$ and $\pi\! = \! \pi^*$. 
This completes the proof of  the theorem. \cqfd

\bigskip 
 
We next recall from Section \ref{connecMqu} that the Markovian $\bw_n$-LIFO queueing system governed by $\ccX_{\bw_n}$ induces a Galton-Watson forest $\bT_{\! \bw_n}$ with offspring distribution $\mu_{\bw_n}$: informally, the clients are the vertices of $\bT_{\! \bw_n}$ and the server is the root (or the ancestor); the $j$-th client to enter the queue is a child of the $i$-th one if the $j$-th client enters when the $i$-th client is served; among siblings, the clients are ordered according to their time of arrival. We denote by $H^{\bw_n}_t$ the number of clients waiting in the line right after time $t$;  
%\margmm{corrected}
recall from (\ref{XJHdef}) how $H^{\bw_n}$ is derived from $X^{\bw_n}$: namely, for all $s \! \leq \! t$, if one sets 
$I^{\bw_n, s}_t= \inf_{r\in [s, t]}X^{\bw_n}_r  $, then, $H^{\bw_n}_t \! = \# \{ s\ino [0, t]\, : \;I^{\bw_n, s-}_{t} \! <\! I^{\bw_n, s}_{t} \}$. 
%\begin{equation}
%\label{rrXJHdef}
%H^{\bw_n}_t \! = \# \big\{ s\ino [0, t]\, : \;I^{\bw_n, s-}_{t} \! <\! I^{\bw_n, s}_{t} \big\}. 
%\end{equation}
We recall from  Section \ref{connecMqu} that $X^{\bw_n}$ and $H^{\bw_n}$ are close to  the Lukasiewicz path and the contour process of $\bT_{\! \bw_n}$. Therefore, the convergence results for Lukasiewicz paths and contour processes of Galton-Watson trees  
in Le Gall \& D.~\cite{DuLG02} (see Appendix Theorem \ref{cvheight}, Section \ref{HeightAppp}) allow us to prove the following theorem.  
 \begin{thm}
\label{cvheheight} 
Let $X$ be an integrable $(\alpha, \beta, \pi)$-spectrally positive L\'evy process, as defined at the beginning of Section \ref{CVmarkogne}. Assume that 
(\ref{frnurz}) holds and that $\int^\infty \! dz / \psi_{\alpha, \beta , \pi} (z) \! <\!  \infty$, where $\psi_{\alpha, \beta , \pi}$ is given by (\ref{LKform}). Let $(H_t)_{t\in [0, \infty)}$ be the continuous height process derived from $X$ as defined by (\ref{approHdef}). 

  Let $\bw_n \ino  \elldo_f$ and $a_n , b_n \ino (0, \infty)$, $n\ino \bbN$, satisfy (\ref{aprioriri}). 
%Recall from (\ref{defpin}) the definition of $\mu_{\bw_n}$ and let 
Let $(Z^{_{(n)}}_{k})_{k\in \bbN}$ be a Galton-Watson process with offspring distribution $\mu_{\bw_n}$ (defined by (\ref{rmupoissw})), 
and initial state $Z^{_{(n)}}_0\! =\!  \lfloor a_n \rfloor$. 
Assume that the three conditions $\mathrm{(IIIabc)}$ in Theorem \ref{cvbranch} hold true and assume the following: 
 % Then, the weak convergence in $\bC ([0, \infty) , \bbR)$ $(\frac{b_n}{a_n} H^{\bw_n}_{b_n t })_{t\in [0, \infty)} \!\!  \longrightarrow \! H$ is equivalent to the following condition: 

\vspace{-5mm}

\begin{equation}
\label{scolheight} 
 \exists \, \delta \ino  (0, \infty) , \qquad \liminf_{n\rightarrow \infty} \bP \big( Z^{_{(n)}}_{\lfloor b_n \delta /a_n \rfloor} \! = \! 0 \big) >0 \; .
\end{equation}
%Moreover, if $\mathrm{(IIIabc)}$ and (\ref{scolheight}) hold true, then the following joint convergence holds true: 
Then, the following joint convergence holds true: 
\begin{equation}
\label{joittcon}
\Big( (\frac{_{1}}{^{{a_n}}} X^{\bw_n}_{b_n t  } )_{t\in [0, \infty)} , (\frac{_{a_n}}{^{b_n}} H^{\bw_n}_{b_n t })_{t\in [0, \infty)}
 \Big)
 \underset{n\rightarrow \infty}{-\!\!\! -\!\!\! -\!\!\! -\!\!\! \longrightarrow } \; (X, H)  
\end{equation}
weakly on $\bD ([0, \infty), \bbR) \times \bC ([0, \infty), \bbR)$, equipped with the product topology. We also get:   
\begin{equation}
\label{cvsvie} 
\forall t\ino [0, \infty), \quad  \lim_{n\rightarrow \infty} \bP \big( Z^{_{(n)}}_{\lfloor b_n t /a_n \rfloor} \! = \! 0 \big) = e^{-v(t)} \quad \textrm{where} \quad  \int_{v(t)}^\infty \! \frac{dz}{\psi_{\alpha, \beta, \pi} (z)}= t . 
\end{equation}
\end{thm}
\noi
\textbf{Proof.} Recall from (\ref{coding}) (Section \ref{HeightApp}) the definition of the Lukasiewicz path $V^{\bT_{\bw_n}}$ associated with the GW($\mu_{\bw_n}$)-forest $\bT_{\bw_n}$; recall from (\ref{codheight}) the definition its height process $\mathtt{Hght}^{\bT_{\bw_n}}$ and recall that $C^{\bT_{\bw_n}}$ stands for the contour process of $\bT_{\! \bw_n}$. 
%, that are respectively denoted by 
%$(V^{_{\bT_{\! \bw_n}}}_k\! )_{k\in \bbN}$, $(\mathtt{Hght}^{_{\bT_{\! \bw_n}}}_k\! )_{k\in \bbN}$ and  $(C^{_{\bT_{\! \bw_n}}}_t\! )_{t\in [0, \infty)}$. 
We first assume that 
$\mathrm{(IIIabc)}$ in Theorem \ref{cvbranch} and that  (\ref{scolheight}) hold true. Then, Theorem \ref{cvheight} applies with $\mu_n \! := \! \mu_{\bw_n}$: namely, the joint convergence 
(\ref{jointtcon}) holds true and we get (\ref{cvsvie}). 

Recall that $(\tau^n_k)_{k\geq 1}$ are the arrival-times of the clients in the queue governed by $X^{\bw_n}$ and recall from (\ref{NNwwdef}) the notation $N^{\bw_n}(t)\! = \! \sum_{k\geq 1} \un_{[0, t]} (\tau^n_k)$ that is a homogeneous Poisson process with unit rate. Then, by Lemma \ref{Poitime} (see Section \ref{Skogen} in Appendix) the joint convergence 
(\ref{jointtcon}) entails the following. 
$$ \mathscr{Q}_n (7)= \big(  \frac{_{1}}{{^{a_n}}} V^{\! {\bT_{\! \bw_n}}}\!  (N^{\bw_n}_{b_n \cdot})  , \,  \frac{_{{a_n}}}{^{{b_n}}} \mathtt{Hght}^{\! {\bT_{\! \bw_n}}}\!  (N^{\bw_n}_{b_n \cdot}) , \, \frac{_{{a_n}}}{^{{b_n}}} C^{{\bT_{\! \bw_n}}}_{b_n \cdot } \big) \underset{n\rightarrow \infty}{-\!\!\! -\!\!\!  \longrightarrow} \big( X, H , (H_{t/2})_{t\in [0, \infty)} \big)$$
weakly on $\bD([0, \infty), \bbR) \times (\bC([0, \infty), \bbR))^2$ equipped with the product topology. Here $X$ is an integrable $(\alpha, \beta, \pi)$-spectrally positive L\'evy process (as defined at the beginning of Section \ref{CVmarkogne}) and $H$ is the height process derived from $X$ by (\ref{approHdef}). 
By Theorem \ref{cvbranch}, the laws of the processes $\frac{_1}{^{a_n}}X^{\bw_n}_{b_n \cdot}$ are tight in $\bD([0, \infty), \bbR)$. Thus, if one sets $\mathscr{Q}_n (8)\! = \! (\frac{_1}{^{a_n}}X^{\bw_n}_{b_n \cdot}, \mathscr{Q}_n (7))$, then the laws of the $\mathscr{Q}_n (8)$ are tight on $\bD([0, \infty), \bbR)^2 \times (\bC([0, \infty), \bbR))^2$. Thus, to prove the weak convergence 
$\mathscr{Q}_n (8)\! \rightarrow \! (X, X, H , H_{\cdot /2})\! :=\!\mathscr{Q} (8)$, we only need to prove that the law of $\mathscr{Q} (8)$ 
is the unique limiting law: to that end, let $(n(p))_{p\in \bbN}$ be an increasing sequence of integers such that 
\begin{eqnarray}
\label{globurnik}
\mathscr{Q}_{n(p)} (8)  \!  \underset{p\rightarrow \infty}{-\!\!\! -\!\!\! -\!\!\! \longrightarrow}  \big( X^\prime \! , X, H , H_{\cdot /2}  \big).  
\end{eqnarray} 
Actually, we only have to prove that $X^\prime \! = \! X$. Without loss of generality (but with a slight abuse of notation), by Skorokod's representation theorem 
we can assume that  (\ref{globurnik}) holds $\bP$-almost surely. We next use (\ref{LuapproXw}) in Lemma \ref{LukaXw}: fix $t, \epp, y\ino (0, \infty)$, set 
$I^{\bw_n}_t\! = \! \inf_{s\in [0, t]} X^{\bw_n}_s $; by applying (\ref{LuapproXw}) at time $b_n t$, with $a \! = \! a_n \epp$ and $x\! = \! a_n y$, we get the following.  
$$ \bP \big( \big| \frac{_{_1}}{^{^{a_n}}}V^{\bT_{\! \bw_n}}_{N^{\bw_n} (b_n t)} \! - \frac{_{_1}}{^{^{a_n}}} X^{\bw_n}_{b_n t}  \big| \! > \! 2\epp \big) \leq 1 \! \wedge \! \frac{_{4y}}{^{\epp^2 a_n}} + 
\bP \big( \! \! -\! \frac{_{_1}}{^{^{a_n}}} I^{\bw_n}_{b_n t} \! >\!  y) + \bE \bigg[ 1 \wedge  \frac{\frac{{1}}{{{a_n}}} (X^{\bw_n}_{b_n t} \! -\! I^{\bw_n}_{b_n t})}{^{\epp^2 a_n} } \bigg] . $$
By Lemma \ref{franchtime} $(ii)$, 
$\frac{{1}}{{{a_{n(p)} }}} (X^{\bw_{n(p)}}_{b_{n(p)} t} \! -\! I^{\bw_{n(p)}}_{b_{n(p)} t}) \! \rightarrow \! X^\prime_t \! -\! I^\prime_t$ and 
$\frac{{1}}{{{a_{n(p)}}}} I^{\bw_{n(p)}}_{b_{n(p)} t} \! \rightarrow \! I^\prime_t$ almost surely, where we have set $I^\prime_t\! = \! \inf_{s\in [0, t]} X^\prime_s$. 
Thus, for all $\epp \ino (0, \infty)$,
$$ \limsup_{p\rightarrow \infty} \bP \big( \big| \frac{_{1}}{^{{a_{n(p)}}}}V^{\bT_{\! \bw_{n(p)}}}_{N^{\bw_{n(p)}} (b_{n(p)} t)} \! - \frac{_{1}}{^{{a_{n(p)}}}} X^{\bw_{n(p)}}_{b_{n(p)} t}  \big| \! > \! 2\epp \big)
\leq \bP \big( \!\! -\! I^\prime_t \! >\! y/2) \; \underset{y\rightarrow \infty}{-\!\!\! -\!\!\! -\!\!\! \longrightarrow } 0  $$
Compared with (\ref{globurnik}), this implies that for all $t\ino [0, \infty)$ a.s.~$X^\prime_t \! = \! X_t$ and thus, a.s.~$X^\prime\! = \! X$. 

We have proved that $\mathscr{Q}_n (8)\! \rightarrow \! (X, X, H , H_{\cdot /2})\! =\! \mathscr{Q} (8)$ weakly on $\bD([0, \infty), \bbR)^2 \times (\bC([0, \infty), \bbR))^2$. Without loss of generality (but with a slight abuse of notation), by Skorokod's representation theorem 
we can assume that the convergence holds true $\bP$-almost surely. We next recall from (\ref{crossedg}) and (\ref{couCHght}) that: 
$$ M^{\bw_n}(t) \! = \! 2N^{\bw_n} (t) \! -\! H^{\bw_n}_t, \quad C^{\bT_{\! \bw_n}}_{M^{\bw_n} (t)} \! = \! H^{\bw_n}_t \quad \textrm{and} \quad \sup_{s\in [0, t]} H^{\bw_n}_s  \leq 1+ \sup_{s\in [0, t]} \mathtt{Hght}^{\bT_{\! \bw_n}}_{N^{\bw_n} (s)} . $$ 
Then, we fix $t, \epp \ino (0, \infty)$, and we apply (\ref{estMbw}) at time $b_n t$, with $a\! = \! b_n \epp$ to get 
$$ \bP \big(\! \sup_{s\in [0, t]} \! \!  |\frac{_1}{^{b_n}}M^{\bw_n}_{b_n s} \! -\! 2s |> 2\epp \big) \leq 1\!  \wedge\! \frac{_{16t}}{^{\epp^2 b_n}} + 
\bP \Big(\frac{_{a_n}}{^{b_n}}+ \!\!   \sup_{\; \,  s\in [0, t]} \!\!   \frac{_{a_n}}{^{b_n}}\mathtt{Hght}^{\bT_{\! \bw}}_{N^\bw (b_n s)} \! > \epp a_n  \Big) . $$
Since $\frac{_{a_n}}{^{b_n}}\mathtt{Hght}^{\bT_{\! \bw}} (N^\bw (b_n \cdot)) \! \rightarrow \! H$ a.s.~in $\bC ([0, \infty), \bbR)$, it easily entails that  
$\frac{_1}{^{b_n}}M^{{\bw_n}}_{{b_n \cdot}}$ tends in probability to twice the identity map on $[0, \infty)$ in $\bC ([0, \infty), \bbR)$. Since $H^{\bw_n}_t \! = \! C^{\bT_{\! \bw_n}} (M^{\bw_n} (t))$, and since 
$C^{\bT_{\! \bw_n}} (b_n \cdot ) \! \rightarrow \! H (\cdot /2)$ a.s.~in $\bC ([0, \infty), \bbR)$, Lemma \ref{Poitime} easily entails the joint convergence (\ref{joittcon}), which completes the proof. \cqfd 

% weakly in $\bD ([0, \infty), \bbR) \times \bC ([0, \infty), \bbR)$ equipped with the product topology. This completes the proof of the theorem.  \cqfd 

 \bigskip

As explained right after Theorem 2.3.1 in Le Gall \& D.~\cite{DuLG02} (see Chapter 2, pp.~54-55) Assumption (\ref{scolheight}) is actually 
a necessary condition for the height process to converge. 
However it is not always easy to check this condition in practice. 
The following proposition provides a handy way of doing it. 

\begin{prop}
\label{rHcritos} Let $X$ be an integrable $(\alpha, \beta, \pi)$-spectrally positive L\'evy process, as defined at the beginning of Section \ref{CVmarkogne}. Assume that 
$(\alpha, \beta, \pi)$ satisfies (\ref{frnurz}) and that $\int^\infty dz / \psi_{\alpha, \beta , \pi} (z) \! <\!  \infty$, where $\psi_{\alpha, \beta , \pi}$ is given by (\ref{LKform}). Let $H$ 
be the continuous height process derived from $X$ by (\ref{approHdef}). 
Let $\bw_n \ino \!  \elldo_f$ and $a_n , b_n \ino (0, \infty)$, $n\ino \bbN$, satisfy (\ref{aprioriri}). 
%Recall from (\ref{defpin}) the definition of $\mu_{\bw_n}$ and let 
%Let $(Z^{_{(n)}}_k)_{k\inÂ \bbN}$ be a Galton-Watson process with offspring distribution $\mu_{\bw_n}$, that is defined by (\ref{defpin}), 
%and initial state $Z^{_{(n)}}_0\! =\!  \lfloor a_n \rfloor$. 
Recall from (\ref{rrXwdef}) the definition of $X^{\bw_n}$ and denote by $\psi_n$ the Laplace exponent of $(\frac{1}{a_n} X^{\bw_n}_{b_n t })_{t\in [0, \infty)}$: namely, 
for all $\lambda \ino [0, \infty)$, 
\begin{equation}
\label{LapXwnn}
\psi_n (\lambda)\! = \!  \frac{b_n}{a_n} \Big( 1\! -\! \frac{\sigma_2 (\bw_n)}{\sigma_1 (\bw_n)} \Big) \lambda + \frac{a_n b_n}{\sigma_1 (\bw_n)} \sum_{j\geq 1} \frac{w_j^{(n)}}{a_n} 
\big( e^{-\lambda w^{(n)}_j /a_n}\!\! -\! 1 + \lambda\,  w^{(n)}_j \!\! /a_n \big)\; .
\end{equation}
We assume that the three conditions $\mathrm{(IIIabc)}$ in Theorem \ref{cvbranch} hold true. Then, (\ref{scolheight}) in Theorem \ref{cvheheight} holds true if the following holds true: 

\vspace{-4mm}

\begin{equation}
\label{rcrinullo}
\lim_{y\rightarrow \infty} \limsup_{n\rightarrow \infty} \int_y^{a_n} \!\! \frac{d\lambda}{\psi_n (\lambda)} = 0 \; .
\end{equation}

\vspace{-2mm}

\end{prop}

\noi
\textbf{Proof.} We first prove a lemma that compares the total height of Galton-Watson trees with i.i.d.~exponentially distributed edge-lengths and the total height of their discrete skeleton. 
More precisely, let $\rho \ino (0, \infty)$ and let $\mu$ be an offspring distribution such that $\mu (0) \! >\! 0$ and whose generating function is denoted by 
$g_\mu (r) \! = \! \sum_{l\in \bbN} \mu (l) r^l $. Note that $g_{\mu} ([0, 1])\! \subset \! [0, 1]$; let $g_\mu^{\circ k}$ be the $k$-th iterate of $g_\mu$, with the convention that $g_\mu^{\circ 0} (r)\! = \! r$, $r\ino [0, 1]$. Let $\tau\! :\!  \Omega \! \rightarrow \! \bbT$ be a random tree whose distribution is characterised as follows.
\begin{compactenum}

\smallskip

\item[--] The number of children of the ancestor (namely the r.v.~$k_\varnothing (\tau)$) is a Poisson r.v.~with mean $\rho$; 

\smallskip

\item[--] For all $l\! \geq \! 1$, under $\bP (\, \cdot \, | \, k_\varnothing (\tau)\! = \! l)$, the $l$ subtrees $\theta_{[1]} \tau, \ldots , \theta_{[l]} \tau$ 
stemming from the ancestor $\varnothing$ are independent Galton-Watson trees with offspring distribution $\mu$. 
\end{compactenum}

\smallskip

\noi
We next denote by $Z_k$ the number of vertices of $\tau$ that are situated at height $k+1$: namely, 
$Z_k\! = \! \# \{ u\ino \tau \! : |u|\! = \! k+1\}$ (see Section \ref{HeightApp} for the notation on trees). Then, $(Z_k)_{k\in \bbN}$ is a Galton-Watson process whose initial value $Z_0$ is distributed as a Poisson r.v.~with mean $\rho$. 
We denote by $\Gamma (\tau)$ the total height of $\tau$: namely, $\Gamma (\tau)\! = \! \max_{u\in \tau } |u|$ is the maximal graph-distance from the root $\varnothing$. Note that if $\mu$ is supercritical, then $\Gamma (\tau)$ may be infinite). Observe that $\Gamma (\tau)\! =\! \min\{ k\ino \bbN: Z_k \! = \! 0 \}$. 
Thus,  
\begin{equation}
\label{Hehe}
 \bP \big( \Gamma (\tau) < k+1  \big)\! = \! \bP (Z_k \! = \! 0) \! = \!  \exp \big(\! - \rho \big( 1\! -\! g_\mu^{\circ k} (0) \big) \big) \; .
\end{equation}
We next equip each individual $u$ of the family tree $\tau$ with an independent lifetime $e(u)$ that is distributed as follows. 
\begin{compactenum}

\smallskip

\item[--] The lifetime $e(\varnothing)$ of $\varnothing$ is $0$. 

\smallskip

\item[--] Conditionally given $\tau$, the r.v.~$e(u)$, $u\ino \tau\backslash \{ \varnothing\}$ are independent and exponentially distributed r.v. with parameter $q\ino (0, \infty)$. 

\end{compactenum}

\smallskip

\noi
Within our notation, the \textit{genealogical order} on $\tau$ is defined as follows: a vertex $v\ino \tau$ is an ancestor of $u\ino \tau$, which is denoted as $v \preceq u$, if there exists $v^\prime \ino \bbU$ such that $u\! = \! v \ast v^\prime$; $\preceq$ is a partial order on $\tau$. For all 
$u\ino \tau\!\setminus\!\{\varnothing\}$, we denote by $\zeta (u) \! =\! \sum_{\varnothing \preceq v \preceq u}e(v)$, the date of death of $u$; then, $\zeta (\overleftarrow{u})$ 
is the date of birth of $u$ (recall here that $\overleftarrow{u}$ stands for the direct parent of $u$). For all $t\ino [0, \infty)$, we next set 
$\mathtt{Z}_t \! = \! \sum_{u\in \tau\backslash \{ \varnothing \}} \un_{[\zeta (\overleftarrow{u}), \zeta(u))} (t)$. Then $(\mathtt{Z}_t)_{t\in [0, \infty)}$ is 
a continuous-time Galton--Watson process (or a Harris process) with offspring distribution $\mu$, with time parameter $q$ and with Poisson$(\rho)$-initial distribution. We denote by $\Gamma \! = \! \max_{u\in \tau} \zeta (u) $ the extinction time of the population; then $\Gamma\! = \! \max \{ t\ino [0, \infty) \! : \mathtt{Z}_t \! \neq \! 0\}$. Standard results on continuous-time GW-processes imply the following. 
For all $t\ino (0, \infty)$, 
\begin{equation}
\label{Heihei}
\bP \big( \Gamma \! <\!  t \big) = \bP (\mathtt{Z}_t \! = \! 0 )=  e^{-\rho r(t)}, \quad \textrm{where} \quad \int_{r(t)}^{1} \frac{dr}{g_\mu (1\! -\! r) \! -\! 1+r} = qt \; .
 \end{equation} 
For a formal proof, see for instance Athreya \& Ney \cite{AtNe72}, Chapter III, Section 3, Equation (7) p.~106 and Section 4, Equation (1) p.~107. 
    
 We next compare $\Gamma(\tau)$ and $\Gamma $. To that end, we introduce $(\mathrm{e}_n)_{n\geq 1}$, a sequence of i.i.d.~exponentially distributed r.v.~with mean $1$, and we set:  
\begin{equation}
\label{biniou}
 \forall \,  \varepsilon \! \in \! (0, 1), \quad \delta (\varepsilon) \, = \, \sup_{n\geq 1} \, \bP \big( n^{-1} (\mathrm{e}_1+ \ldots + \mathrm{e}_n) \notin (\varepsilon  , \varepsilon^{-1}) \big) .
 \end{equation}
Law of Large Numbers easily implies that $ \delta (\varepsilon) \! \rightarrow \! 0$ as $\varepsilon \! \rightarrow \! 0$. 
Note that $\mathtt{Z}_0\! = \! Z_0$ and a.s.~$\Gamma (\tau)\!< \! \infty$ if and only if  
$\Gamma\! <\! \infty$. We argue on the event $\{\Gamma (\tau) \! < \! \infty\}$: 
we first assume that $Z_0 \! \neq \! 0$; 
let $u^*\ino \tau \backslash \{ \varnothing \}$ be the first vertex in the lexicographical order such that $|u^*|\! = \! \Gamma (\tau)$; since $\zeta (u^*) \! \leq \! \Gamma$ and since conditionally given $\tau$, $\zeta (u^*)$ is the sum of $|u^*|$ (conditionally) 
independent exponential r.v.~with parameter $q$, we get for all $t\ino (0, \infty)$, 
%$\bP ( \Gamma  \! < \! t ; Z_0 \! \neq \! 0 ) \! \leq \! \bP (\zeta (u) \! <\! t ; Z_0 \! \neq \! 0) $. Next observe that $G$ conditionally given $\tau$ is distributed as $\sum_{1\leq k\leq \mathtt{Hght}  (\tau) +1} \mathrm{e}_k /q$. Thus, 
$$ \bP \big(\Gamma \! < \! t \, ; \, \mathtt{Z}_0 \! \neq \! 0\big)  \leq \sum_{n\geq 1} \bP \big( \, \Gamma (\tau) \! = \! n\! \, ;\,  Z_0 \! \neq \! 0 \, \big) \, \bP \big(\mathrm{e}_1+ \ldots + \mathrm{e}_n \! \leq  \! 
qt\big) \; .$$
Then, let $\varepsilon \! \in \! (0, 1)$ and observe that 
$  \bP \big( \mathrm{e}_1+ \ldots + \mathrm{e}_n \! \leq  \! qt\big) \leq \delta (\varepsilon) + \un_{\{ n  \leq qt / \varepsilon \}}$. 
Consequently, 
$$ \bP \big( \Gamma  \! < \! t \, ;\,  \mathtt{Z}_0 \! \neq \! 0\big)  \leq \delta (\varepsilon) + \bP \big( \Gamma (\tau) \! \leq \! \lfloor q t/\varepsilon \rfloor \, ;\,  Z_0 \! \neq \! 0 \big) .$$
If $\mathtt{Z}_0\! = \! Z_0 \! = \! 0$, $\Gamma \! = \! \Gamma (\tau)\! = \! 0$, which implies that 
$$  \bP \big( \Gamma  \! < \! t \big)  \leq \delta (\varepsilon) + \bP \big(\Gamma  (\tau)   \! \leq \! \lfloor q t/\varepsilon \rfloor \big) .$$
Thus by (\ref{Heihei}) and (\ref{Hehe}), we have proved the following lemma. 
\begin{lem}
\label{hauhaut}
Let $\rho  , q\ino (0, \infty)$ and let $\mu$ be an offspring distribution such that $\mu (0) \! >\! 0$ and whose generating function is denoted by 
$g_\mu$; denote by $g_\mu^{\circ k}$ the $k$-th iterate of $g_\mu$ with the convention $g_\mu^{\circ 0} (r)\! = \! r$, $r\ino [0, 1]$. Let $t\ino (0, \infty)$. Recall from (\ref{Heihei}) 
the definition of $r(t)$.  
%For all $t\ino (0, \infty)$, let $r(t)$ be such that 
%\begin{equation}
%\label{coocrf}
%\int_{r(t)}^{1} \frac{dr}{g_\mu (1\! -\! r) \! -\! 1+r} = qt \; .
%\end{equation}
Let $\varepsilon \! \in \! (0, 1)$. Recall from (\ref{biniou}) the definition of $\delta (\varepsilon)$. Then, the following holds true. 
\begin{equation}
\label{cornemuse}
\forall \, t \! \in \! (0, \infty), \quad  e^{- \rho r(t)}  -\!  \delta (\varepsilon ) \, \leq \, \exp \big(\! -\! \rho \big( 1\! -\! g_\mu^{\circ \lfloor tq/\varepsilon \rfloor }  (0) \,  \big) \big). 
\end{equation}
\end{lem}

We are now ready to \textbf{prove Proposition \ref{rHcritos}.} Recall from (\ref{rmupoissw}) the definition of the offspring distribution $\mu_{\bw_n}$.  
We apply Lemma \ref{hauhaut} with $\mu\! = \! \mu_{\bw_n}$, $\rho \! = \! a_n$, $q\! = \! b_n /a_n$ and we denote by 
$r_n (t)$ the solution of (\ref{Heihei}): the change of variable $\lambda \! = \!  a_n r$ implies that $r_n (t)$ satisfies 
\begin{equation}
\label{gloobiii}
\int_{a_n r_n(t)}^{a_n} \, \frac{d\lambda}{b_n \big( g_{\mu_{\bw_n}} \big( 1\! -\! \frac{\lambda}{{a_n}}\big)  \! -\! 1+ \frac{\lambda}{{a_n}} \big) } \, = t\,.
\end{equation}
Next, it is easy to check from (\ref{rmupoissw}) that 
$b_n \big( g_{\mu_{\bw_n}} (1\! -\! \frac{\lambda}{{a_n}}) \! -\! 1+ \frac{\lambda}{a_n} \big) \! =\!  \psi_n (\lambda)$, where $\psi_n$ is defined in (\ref{LapXwnn}). 
Then, Lemma \ref{hauhaut} asserts for all $t\ino (0, \infty)$ and for all $\varepsilon \ino (0, 1)$, that 
   \begin{equation}
\label{flutiau}
e^{- a_n r_n(t)}  \! -  \delta (\varepsilon ) \, \leq \, \exp \big(\! -\! a_n \big( 1\! -\! g_{\mu_{\bw_n}}^{\circ \lfloor tb_n/a_n\varepsilon \rfloor }\!  (0) \, \big) \big) \quad \textrm{where} \quad \int_{a_n r_n (t)}^{a_n} \frac{d\lambda}{\psi_n (\lambda)}= t \; .
\end{equation}
Next, fix $t\ino (0, \infty)$ and set $C\! :=\! \limsup_{n\rightarrow \infty} a_n r_n (t) \ino [0, \infty]$. Suppose 
that $C\! = \! \infty$. Then, there is an increasing sequence of integers $(n_k)_{k\in \bbN}$ such that $\lim_{k\rightarrow \infty} a_{n_k} r_{n_k} (t) \! = \! \infty$. Let $y\ino (0, \infty)$; then, for all sufficiently large $k$, we have $a_{n_k} r_{n_k} (t) \! \geq \! y$, which entails 
$$ t =  \int_{a_{n_k} r_{n_k} (t)}^{a_{n_k}} \frac{d\lambda}{\psi_{n_k} (\lambda)} \leq  \int_{y}^{a_{n_k}} \frac{d\lambda}{\psi_{n_k} (\lambda)} $$
Thus, for all $y\ino (0, \infty)$, $t \! \leq \! \limsup_{n\rightarrow \infty} \int_y^\infty d\lambda / \psi_n (\lambda) $, which contradicts Assumption (\ref{rcrinullo}). 
This proves that $C\! < \! \infty$. Since $\lim_{\epp \rightarrow 0} \delta (\varepsilon) \! = \!  0$, we can choose $\varepsilon$ such that $\delta (\varepsilon) \! < \! \frac{_1}{^2}e^{-C}$; then, we set $\delta \! = \!  t/ \varepsilon$ and 
(\ref{flutiau}) implies that 
\begin{equation}
\label{cistre}
\limsup_{n\rightarrow \infty}     a_n \Big( 1\! -\! g_{\mu_{\bw_n}}^{\circ \lfloor \delta b_n/a_n \rfloor }\!  (0)\Big ) < \infty \; .
\end{equation}
Recall that $(Z^{_{(n)}}_{k})_{k\in \bbN}$ stands for 
a Galton-Watson branching process with offspring distribution $\mu_{\bw_n}$ such that $Z^{_{(n)}}_{0}\! = \! \lfloor a_n \rfloor$. Then, 
$ \bP \big( Z^{_{(n)}}_{\lfloor \delta b_n / a_n \rfloor} \! =\!  0\big)\! =\!  \big(  g_{\mu_{\bw_n}}^{\circ \lfloor \delta b_n/a_n \rfloor }\!  (0)\big)^{ \lfloor a_n \rfloor} $
and (\ref{cistre}) easily implies that $\liminf_{n\rightarrow \infty}  \bP \big( Z^{_{(n)}}_{\lfloor \delta b_n / a_n \rfloor}\!  \! =\!  0\big) \! >\! 0 $, which completes the proof 
of Proposition \ref{rHcritos}. \cqfd

\subsubsection{Proof of Propositions \ref{cvmarkpro} and \ref{HMarkcvprop}.}
\label{pfcvmarkpro}
In this section we shall assume that the sequence $(a_n)$ and $(b_n)$ satisfy (\ref{aprioriri}) \textbf{and} $\frac{a_n b_n}{\sigma_1 (\bw_n)}\! \rightarrow \! \kappa$ where $\kappa \ino (0, \infty) $. 
This dramatically restricts the possible limiting triplets $(\alpha, \beta, \pi)$.  
To see this point, we first prove the following lemma. 
\begin{lem}
\label{fixation} For all $n\in \bbN$, let $\bv_n\! = \! (v^{_{(n)}}_{j})_{j\geq 1} \ino \elldo_f$ and set 
$\phi_n (\lambda) \! =\! \sum_{j\geq 1} v^{_{(n)}}_{j} \big( e^{-\lambda v^{_{(n)}}_{j}} \!\!  -\! 1  +\lambda v^{_{(n)}}_{j} \big)$, for all $\lambda \in [0, \infty)$. 
%$$ \forall \lambda \in [0, \infty), \quad \phi_n (\lambda) = \sum_{j\geq 1} v^{(n)}_j \big( e^{-\lambda v^{(n)}_j} -1 +\lambda v^{(n)}_j \big) \; .$$
Then, the following assertions are equivalent. 
\begin{compactenum}

\smallskip

\item[$(L)$]  For all $\lambda \ino [0, \infty)$, there exists $\phi (\lambda)\ino [0, \infty)$ such that $\lim_{n\rightarrow \infty} \phi_n (\lambda)  \! = \!  \phi (\lambda)$.  

\smallskip

\item[$(S)$]  There are $\mathbf{c} \in \elldo_3$ and $\beta^\prime\ino [0, \infty)$ 
such that 
$$ \forall j \! \in \!  \bbN^*, \quad \lim_{n\rightarrow \infty}v^{_{(n)}}_{j} \! = \! c_j \quad \textrm{and} \quad \lim_{n\rightarrow \infty} \sigma_3 (\bv_n) \! -\!  \sigma_3 (\mathbf{c})\! = \!  \beta^\prime . $$ 
\end{compactenum}
Moreover, if $(L)$ or $(S)$ holds true, then $\phi$ in $(L)$ is given by 
\begin{equation}
\label{flikado}
\forall \lambda \in [0, \infty) , \quad \phi (\lambda)= \frac{_1}{^2} \beta^\prime \lambda^2 + \sum_{j\geq 1} c_j \big( e^{-\lambda c_j} -1 +\lambda c_j \big)\; .
\end{equation}
\end{lem}
\noi
\textbf{Proof.} We first prove $(S) \! \Rightarrow \! (L)$. For all $x\ino [0, \infty)$, we set $f(x)\! =\! e^{-x} \! -\! 1+x$. Elementary arguments entail the following.  
\begin{equation}
\label{2ordre}
\forall x\ino [0, \infty) \quad  0 \leq \frac{_{1}}{^{2}} x^2 \! - \! f(x) \leq \frac{_{_1}}{^{2}} x^2 (1\! -\! e^{-x})  \,.   
\end{equation}
We set $\eta(x) \! =\!  
\sup_{y\in [0, x]} y^{-2}| \frac{_1}{^2}y^2 \! -\! f(y) |$; thus, $\eta (x) \! \leq \! \frac{1}{2} (1\! -\! e^{-x}) \! \leq \! 1 \! \wedge \! x$ and $\eta (x) \! \downarrow \! 0$ as $x\! \downarrow\!  0$. 
Then, fix $\lambda \ino [0, \infty)$ and define $\phi (\lambda)$ by (\ref{flikado}); fix  $j_0 \geq 2$ and observe the following.   
\begin{eqnarray*}
 \phi_n (\lambda)\! -\! \phi (\lambda) &=& \sum_{1\leq j\leq j_0} \!\!\! \Big( v^{_{(n)}}_{^j} \! f(\lambda v^{_{(n)}}_{^j} ) \! -\!  c_j f(\lambda c_j)\Big) + 
\frac{_{_1}}{^{^2}}\lambda^2 \Big(\sigma_3 (\bv_n) \! -\!  \sigma_3 (\mathbf{c}) \! -\! \beta^\prime + \sum_{1\leq j\leq j_0} \!\!\! \big( c^3_j\! -\! (v^{_{(n)}}_{^j} )^3  \big) \Big) \\
&  &  \quad + \sum_{j>j_0}   \!\! \Big( v^{_{(n)}}_{^j}  f(\lambda v^{_{(n)}}_{^j} ) \! -\! 
 \frac{_{_1}}{^{^2}}\lambda^2 (v^{_{(n)}}_{^j})^3 \Big)  +\! \sum_{j>j_0}  \!\! \Big(   \frac{_{_1}}{^{^2}}\lambda^2 c^3_j \! -\!  c_j f(\lambda c_j) \Big) .
\end{eqnarray*}
Then, note that: 
$$  \sum_{j>j_0}   \!\! \big| v^{_{(n)}}_{^j}  f(\lambda v^{_{(n)}}_{^j} ) \! -\! 
 \frac{_{_1}}{^{^2}}\lambda^2 (v^{_{(n)}}_{^j})^3 \big| \leq \lambda^2 \eta \big(\lambda v^{_{(n)}}_{^{j_0}}  \big) \sigma_3 (\bv_n) \; .$$
Similarly, $\sum_{j>j_0}   \big|   \frac{_{_1}}{^{^2}}\lambda^2 c^3_j  -  c_j f(\lambda c_j) \big| \leq   \lambda^{2}\eta \big(\lambda c_{j_0} \big) \sigma_3 (\mathbf{c})$. Thus 
$$  \limsup_{n\rightarrow \infty} \big| \phi_n (\lambda)-\phi (\lambda) \big| \leq  (\beta^\prime+ 2\sigma_3 (\mathbf{c})) \lambda^{2}\eta \big(\lambda c_{j_0} \big) 
 \underset{^{j_0\rightarrow \infty}}{-\!\!\! -\!\! \! -\!\! \! \longrightarrow}  \; \; 0 \; , $$
since $c_{j_0} \rightarrow  0$ as $j_0\rightarrow \infty$. This proves $(L)$ and (\ref{flikado}).

Conversely, we assume $(L)$. Note that $v^{_{(n)}}_{1} \!  f(v^{_{(n)}}_{1} ) \! \leq \!  \phi_n (1) $.
Thus, $x_0\! :=\!  \sup_{n\in \bbN} v^{_{(n)}}_{1} \! < \! \infty $. 
By (\ref{2ordre}), for all $y \ino [0, x]$,  $f(y) \! \geq \! \frac{_1}{^2}e^{-x}y^{2}$, which implies 
$\sigma_3 (\bv_n) \! \leq \! 2 e^{x_0} \sup_{n\in \bbN} \phi_n (1)\! =:\! z_0$.  Consequently, 
for all $n\ino \bbN$, 
$(\sigma_3 (\bv_n), \bv_n)$ belongs to the compact space $[0, z_0] \! \times \! [0, x_0]^{\bbN^*}$. 
Let $(q_n)_{n\in \bbN}$ be an increasing sequence of integers such that 
$\lim_{n\rightarrow \infty} \sigma_3 (\bv_{q_n})\! =\!  a$ for some $a\in [0, z_{0}]$ and such that for all $j \! \geq \! 1$, $\lim_{n\rightarrow \infty} v^{_{(q_n)}}_{j} \! =\!  c^\prime_j$ for certain $c^{\prime}_{j}\in [0, x_{0}]$. 
By Fatou's Lemma, $\sigma_3 (\mathbf{c}^\prime)\! \leq \! a$ and we then set $\beta^\prime \! =\!  a\! -\! \sigma_3   (\mathbf{c}^\prime)$. By applying 
$(S) \Rightarrow  (L)$ to $(\bv_{q_n})_{n\in \bbN}$, we get 
$\phi(\lambda) \!  =\!  \frac{_1}{^2} \beta^\prime \lambda^2 + \sum_{j\geq 1} c^\prime_j \big( \exp (-\lambda c^\prime_j) \! -\! 1 +\lambda c^\prime_j \big)$, for all $\lambda \in [0, \infty)$. 
We easily show that it characterises $\beta^\prime$ and $\mathbf{c}^\prime$. Thus, $((\sigma_3 (\bv_n), \bv_n))_{n\in \bbN}$, has a unique limit point in $[0, z_0] \! \times \! [0, x_0]^{\bbN^*}$, which easily entails $(S)$. \cqfd 

\begin{lem}
\label{condenden} Let $\bw_n \ino  \elldo_f$ and $a_n , b_n \ino (0, \infty)$, $n\ino \bbN$, satisfy (\ref{apriori}). 
%Namely  
%\begin{multline}
%\label{rapriori}
%a_n \; \,  \textrm{and} \; \, \frac{b_n}{a_n}  \, \underset{n\rightarrow \infty}{-\!\!\! \longrightarrow}  \infty , \quad  
%\frac{b_n}{a^2_n} \,  \underset{n\rightarrow \infty}{ -\!\!\! \longrightarrow}\,  \beta_0 \ino [0, \infty), \\
% \quad 
%\sup_{n\in \bbN} \frac{w^{_{(n)}}_{^1}}{ a_n} \! < \! \infty \quad    \textrm{and}   \quad   
%\frac{a_nb_n}{\sigma_1 (\bw_n)}\,  \underset{n\rightarrow \infty}{-\!\!\! \longrightarrow} \, \kappa \ino (0, \infty). 
%\end{multline}
Recall from (\ref{rrXwdef}) the definition of  $X^{\bw_n}_{\, } \! $. Then the following assertions hold true. 
\begin{compactenum}

\smallskip

\item[$(i)$] Let us suppose that $\mathrm{(II)}$ in Theorem \ref{cvbranch} holds true; namely, 
$\frac{1}{{a_n}} X^{\bw_n}_{ b_n \cdot  } 
 \! \longrightarrow \! X$ weakly on $\bD([0, \infty), \bbR)$. Then, $X$ is an integrable  
$(\alpha , \beta, \pi)$ spectrally positive L\'evy process (as defined at the beginning of Section \ref{CVmarkogne}) and $(\alpha, \beta, \pi)$ necessarily satisfies: 
\begin{equation}
\label{glurp}
\beta \! \geq \! \beta_0 \quad \textrm{and} \quad \exists \, \mathbf{c}\! = \! (c_j)_{j\geq 1}\ino \elldo_3 : \, \; \pi \! = \! \sum_{j\geq 1} \kappa c_j \delta_{c_j} 
\end{equation}
and the following statements hold true: 

\vspace{-5mm}

\begin{equation*}
%\label{runalphcv}
%\mathbf{(C1):} \quad \frac{a_nb_n}{\sigma_1 (\bw_n)}\;  \underset{n\rightarrow \infty}{-\!\!\! -\!\!\! \longrightarrow} \; \kappa \, , \qquad 
\mathbf{(C1):} \quad \frac{b_n}{a_n} \Big( 1-\frac{\sigma_2 (\bw_n)}{\sigma_1 (\bw_n)} \Big)\;  \underset{^{n\rightarrow \infty}}{-\!\!\! -\!\! \!  \longrightarrow}  \; \alpha \qquad \mathbf{(C2):} \quad  \frac{b_n}{a^2_n}\!  \cdot \!  \frac{\sigma_3 (\bw_n)}{\sigma_1 (\bw_n)} \;  \underset{^{n\rightarrow \infty}}{-\!\!\! -\!\! \!  \longrightarrow}  \; \beta + \kappa \sigma_3 (\mathbf{c}) \, ,
\end{equation*}

\vspace{-6mm}

\begin{equation*}
%\label{rsig3cvcj}
%\mathbf{(C3):} \quad  \frac{b_n}{a^2_n}\!  \cdot \!  \frac{\sigma_3 (\bw_n)}{\sigma_1 (\bw_n)} \;  \underset{^{n\rightarrow \infty}}{-\!\!\! -\!\! \!  \longrightarrow}  \; \beta + \kappa \sigma_3 (\mathbf{c}) \, , \qquad 
\mathbf{(C3):} \quad \forall j \in \bbN^*, \quad \frac{w^{(n)}_j}{a_n } \;  \underset{^{n\rightarrow \infty}}{-\!\!\! -\!\! \! \longrightarrow}  \;  c_j
\, .
\end{equation*}
\item[$(ii)$] Conversely, $(\mathbf{C1})$--$(\mathbf{C3})$ are equivalent to $\mathrm{(II})$ in Theorem \ref{cvbranch}; it is also equivalent to $\mathrm{(I)}$, or to $\mathrm{(IIIabc)}$ or to $(\mathrm{(IIIa)} \, \& \, \mathrm{(IV)})$.  
\end{compactenum}
\end{lem}
\noi
\textbf{Proof.} To simplify notation, we set $\kappa_n \! = \! a_n b_n / \sigma (\bw_n)$. By the last point 
of (\ref{apriori}), $\kappa_n \! \rightarrow \! \kappa\ino (0, \infty)$. We also set 
$v^{_{(n)}}_j\! = \! w^{_{(n)}}_j/a_n$ for all $j\! \geq \! 1$. 
 We first prove $(i)$, so we suppose Theorem \ref{cvbranch} $\mathrm{(II)}$, which first implies that 
 $\beta \! \geq \! \beta_0$; then recall that Theorem \ref{cvbranch} $\mathrm{(II)}$ is equivalent to 
$( (\mathbf{C1})\, \& \, \mathrm{(IV)})$ and Theorem \ref{cvbranch} $\mathrm{(IV)}$ can be rewritten as follows: 
for all $\lambda \ino [0, \infty)$,  
$$\kappa_n \sum_{j\geq 1} v^{_{(n)}}_{j} \big( e^{-\lambda v^{_{(n)}}_{j}} \!\!  -\! 1  +\lambda v^{_{(n)}}_{j} \big) \underset{n\rightarrow \infty}{-\!\!\! -\!\!\! \longrightarrow }  \psi_{\alpha , \beta , \pi} (\lambda)-\alpha \lambda \; .$$
This entails Condition $(L)$ in Lemma \ref{fixation} with $\phi (\lambda) \! = \! 
( \psi_{\alpha , \beta , \pi} (\lambda)-\alpha \lambda ) / \kappa $. Lemma \ref{fixation} then implies that there are 
$\mathbf{c} \in \elldo_3$ and $\beta^\prime\ino [0, \infty)$ 
such that for all  $j \! \in \!  \bbN^*$, $\lim_{n\rightarrow \infty}v^{_{(n)}}_{j} \! = \! c_j $ 
and $ \lim_{n\rightarrow \infty} \sigma_3 (\bv_n) \! -\!  \sigma_3 (\mathbf{c})\! = \!  \beta^\prime$ and that 
$$ \frac{_1}{^2} \kappa^{-1}\beta \lambda^2 + \kappa^{-1} \!\!  \int_{(0, \infty)} \!\!\!\!\! \!\!\!\! (e^{-\lambda r} \! \! -\! 1 + \lambda r) \, \pi (dr) \! = \!  
\frac{\psi_{\alpha , \beta , \pi} (\lambda)\! -\! \alpha \lambda }{ \kappa}  \! = \! \phi (\lambda) \! = \!  \frac{_1}{^2} \beta^\prime \lambda^2 +\!  \sum_{j\geq 1} c_j \big( e^{-\lambda c_j} \! -\! 1 +\lambda c_j \big)\; .$$
This easily entails that $\kappa \beta^\prime\! = \! \beta$, $\pi \! = \! \sum_{j\geq 1} \kappa c_j \delta_{c_j}$ and we easily get 
$(\mathbf{C2})$ and $(\mathbf{C3})$. 

We next prove $(ii)$: we assume that $\beta \! \geq \! \beta_0$ and that $\pi   \! = \! \sum_{j\geq 1} \kappa c_j \delta_{c_j}$ where $\mathbf{c}\! = \! (c_j)_{j\geq 1}\ino \elldo_3$. Then observe that $(\mathbf{C1})$ is $\mathrm{(IIIa)}$ in Theorem \ref{cvbranch}, that $(\mathbf{C2})$ is $\mathrm{(IIIb)}$ in Theorem \ref{cvbranch}; 
moreover, $(\mathbf{C3})$ easily entails  $\mathrm{(IIIc)}$ in Theorem \ref{cvbranch}. Then Theorem \ref{cvbranch} easily entails $(ii)$. This completes the proof of the lemma.  \cqfd 

\medskip

Lemma \ref{condenden} combined with Theorem \ref{cvbranch} implies Proposition 
\ref{cvmarkpro} $(i)$, $(ii)$ and $(iii)$, and Lemma \ref{condenden} combined with Theorem \ref{cvheheight} implies 
Proposition \ref{HMarkcvprop}. 
It only remains to prove Proposition \ref{cvmarkpro} $(iv)$. Namely, fix $\alpha\ino \bbR$, $\beta  \ino [0, \infty)$, $\kappa \ino (0, \infty)$, and 
$\mathbf{c}\! = \! (c_j)_{j\geq 1}\ino \elldo_3$. We prove that 
there are sequences $a_n , b_n \ino (0, \infty) $, $\bw_n \ino \elldo_f$, $n\ino \bbN$, that satisfy (\ref{apriori}) with $\beta_0 \ino [0, \beta]$ and $\mathbf{(C1)}$, $\mathbf{(C2)}$ and $\mathbf{(C3)}$ and $\sqrt{\mathbf{j}_n}/ b_n \! \rightarrow \! 0$ where we recall notation $\mathbf{j}_n \! = \! \max \{ j\! \geq \! 1 \! : \! w^{_{(n)}}_{^j} \! >\! 0 \}$.  

To that end, first let $(\rho_n)_{n\in \bbN}$ be a sequence of positive integers such that  $\rho_n \! \leq \! n$, $\lim_{n \rightarrow \infty} \rho_n \! = \! \infty$ and $\sum_{1\leq j \leq \rho_n} c_j + c_j^2 \leq n $, for all $n \! \geq \! c_1+c_1^2$. 
%We first consider the case where $\beta_0 \! = \! 0$. We introduce 
%a $\bbN$-valued sequence $(\rho_n)_{n \in  \bbN}$ that increases to $\infty$ and that is specified further.  
We then define the following. 
\begin{equation}
\label{defqnj} 
q_{^j}^{_{(n)}}=\left\{
\begin{array}{ll}
c_j & \text{if} \quad  j \! \in \!  \big\{ 1, \ldots ,    \rho_n \!  \big\} 
%1\! \leq \! j\! \leq \! \lfloor  \rho_n 2^n \! \rfloor 
,  \\
((\beta\! -\! \beta_0)/ \kappa)^{\frac{1}{3}}n^{-1} & \text{if} \quad  j \! \in \!  \big\{\rho_n  +1 , \ldots , \rho_n   + n^3  \big\}, \\
u_n & \text{if} \quad  j \! \in \!  \big\{ \rho_n   +n^3 +1 , \ldots , \rho_n   + n^3+n^8   \big\}, \\
% \lfloor \rho_n 2^n\! \rfloor +1 \! \leq \! j\! \leq \! \lfloor \rho_n 2^n\! \rfloor + n2^n , \\
\; \; \; 0 & \text{if} \quad  j >  \rho_n + n^3 + n^8, 
\end{array}\right.
\end{equation}
where $u_n \! = \! n^{-3}$ if $\beta_0 \! = \! 0$ and $u_n \! = \! (\beta_0/ \kappa)^{\frac{1}{3}}n^{-8/3}$ if $ \beta_0 \! >\! 0$. 
We denote by $\bv_n \! =\!  (v^{_{(n)}}_{^j})_{j\geq 1}$ the nonincreasing rearrangement of $\mathtt{q}_n \! =\!  (q^{_{(n)}}_{^j})_{j\geq 1}$. Thus, we get 
$\sigma_p (\bv_n)\! = \! \sigma_p (\mathtt{q}_n)$ for any $p\! \in \!  (0, \infty)$ and we observe the following. 
\begin{multline}
\label{ssigun}
\kappa \sigma_1 (\bv_n) \sim   \left\{ \begin{array}{ll}
\kappa n^5 \!\!\! & \text{if $\beta_0 \! = \! 0$,}   \\
\kappa^{\frac{2}{3}} \beta_0^{\frac{1}{3}} n^{\frac{16}{3}} \!\!\! & \text{if  $\beta_0 \! > \! 0$,} 
\end{array}\right.  \\ 
\kappa \sigma_2 (\bv_n) \sim   \left\{ \begin{array}{ll}
\kappa n^2 & \text{if $\beta_0 \! = \! 0$,}   \\
\kappa^{\frac{1}{3}} \beta_0^{\frac{2}{3}} n^{\frac{8}{3}} & \text{if  $\beta_0 \! > \! 0$,} 
\end{array}\right.  \quad \textrm{and} \quad \kappa \sigma_3 (\bv_n) \sim   \kappa \sigma_3 (\mathbf{c}) + \beta\; . 
\end{multline}
%
%\begin{multline}
%\label{ssigun}
%\kappa \sigma_1 (\bv_n) = \! \!\!  \sum_{\; 1\leq j\leq \rho_n} \!\!\!\! \kappa c_j +  \beta^{\frac{1}{3}}  \kappa^{\frac{2}{3}} n^{2} + 
%\kappa n^{5} , \quad \kappa \sigma_2 (\bv_n) = \! \!\!  \sum_{\; 1\leq j\leq \rho_n} \!\!\!\! \kappa c^2_j +  \beta^{\frac{2}{3}} \kappa^{\frac{1}{3}} n + \kappa n^{2} \\
%\textrm{and} \quad \kappa \sigma_3 (\bv_n) \! = \!   \beta + \kappa \sigma_3 (\mathbf{c})- \!\!  \sum_{\; j> \rho_n} \!\!\!\! \kappa c^3_j  + \kappa n^{-1}.   \quad  \quad  \quad  \quad  \quad  \quad 
%\end{multline}
We next set:
\begin{equation}
\label{dedef}
 b_n = \kappa  \sigma_1 (\bv_n), \quad a_n = \frac{\kappa  \sigma_1 (\bv_n) }{\kappa  \sigma_2 (\bv_n) +\alpha  } \quad \textrm{and} \quad w^{_{(n)}}_{^j}= a_n v^{_{(n)}}_{^j}, \; j\! \geq \! 1 .
\end{equation}
We then see that $a_nb_n/\sigma_1 (\bw_n)\! = \! \kappa$, that $\sup_{n\in \bbN} w^{_{(n)}}_{^1}/a_n \! < \! \infty$. Moreover, we  get 
$$ \frac{b_n}{a_n} \Big(1- \frac{\sigma_2 (\bw_n)}{\sigma_1 (\bw_n)} \Big)= \alpha , \quad 
\lim_{n\rightarrow \infty}  \frac{{b_n}}{{a^2_n}}\!  \cdot \!  \frac{{\sigma_3 (\bw_n)}}{{\sigma_1 (\bw_n)}}\! = \!   \beta + \kappa \sigma_3 (\mathbf{c})  \quad \textrm{and} \quad \lim_{n\rightarrow \infty} \frac{w^{_{(n)}}_{^j}}{a_n}= c_j, \forall\,j\in\mathbb N^{\ast}, $$
which are the limits $(\mathbf{C1})$, $(\mathbf{C2})$ and $(\mathbf{C3})$. 
It is easy to derive from (\ref{ssigun}) and from (\ref{dedef}) that $a_n$ and $b_n /a_n$ tend to $\infty$ and that $b_n / a_n^2$ tends to $\beta_0$. Moreover, since $\mathbf{j_n}\! \leq  \! n^8+n^3 + n$, it is also easy to check that $\sqrt{\mathbf{j_n}}/b_n \! \rightarrow \! 0$. 
This completes the proof of Proposition \ref{cvmarkpro} $(iv)$. \cqfd 

\subsubsection{Proof of Proposition \ref{Hcritos} $(i)$.}
\label{pfHcritos1}
Fix $\alpha\ino \bbR$, $\beta\ino [0, \infty)$, $\kappa \ino (0, \infty)$ and $\mathbf{c}\! =\! (c_j)_{j\geq 1} \! \in \! \elldo_3$. 
For all $\lambda \! \in \! [0, \infty)$, set $\psi (\lambda) \! = \! \alpha \lambda + \frac{_1}{^2} \beta \lambda^2 + \sum_{j\geq 1} \kappa 
c_j \big( e^{-\lambda c_j} \! -\! 1 + \lambda c_j \big)$ and we assume that $\int^\infty \! d\lambda / \psi (\lambda) \! < \! \infty$. 
Let $a_n , b_n \ino (0, \infty)$ and $\bw_n \ino \elldo_f$, $n\ino \bbN$, satisfy (\ref{apriori}), $(\mathbf{C1})$, $(\mathbf{C2})$ and $(\mathbf{C3})$ (as recalled in Lemma \ref{condenden}). Recall from (\ref{rrXwdef}) the definition of $X^{\bw_n}$; recall from (\ref{LapXwnn}) the definition of $\psi_n$ that is the Laplace exponent of $\frac{1}{a_n} X^{\bw_n}_{b_n \cdot}$. To simplify notation, set $\alpha_n \! = \!  \frac{b_n}{a_n} ( 1\! -\! \frac{\sigma_2 (\bw_n)}{\sigma_1 (\bw_n)} )$. 
%%
%%Namely, 
%%\begin{equation}
%%\label{rLapXwnn}
%%\forall \lambda \ino [0, \infty), \quad \psi_n (\lambda)\! = \!  \alpha_n \lambda + \frac{a_n b_n}{\sigma_1 (\bw_n)} \sum_{j\geq 1} \frac{w_j^{(n)}}{a_n} 
%%\big( e^{-\lambda w^{(n)}_j /a_n}\!\! -\! 1 + \lambda\,  w^{(n)}_j \!\! /a_n \big)\; .
%%\end{equation}
%%where we have set $\alpha_n \! = \!  \frac{b_n}{a_n} \big( 1\! -\! \frac{\sigma_2 (\bw_n)}{\sigma_1 (\bw_n)} \big)$. 
%%Proposition \ref{rHcritos} and Lemma \ref{condenden} prove that (\ref{crinullo}) Proposition \ref{Hcritos} $(i)$ entails $(\mathbf{C4})$. 
%%
It only remains to prove the last point of Proposition \ref{Hcritos} $(i)$: assume that $\beta_0 \! >\! 0$ in (\ref{apriori}); let $V_n\! : \! \Omega \! \rightarrow  [0, \infty)$ be a r.v.~with law $\sigma_1 (\bw_n)^{-1}\sum_{j\geq 1} w^{_{(n)}}_j  \delta_{ w^{_{(n)}}_j / a_n}$. First, observe the following. 
$$ \bE [V_n] \! = \!  \frac{\sigma_2 (\bw_n)}{a_n \sigma_1 (\bw_n)} =    \frac{1}{a_n} \Big(1 \! -\! \alpha_n \frac{a_n}{b_n} \Big) \quad 
\textrm{and} \quad \psi_n (\lambda) \! -\! \alpha_n \lambda = b_n \bE \big[ f\big( \lambda V_n\big)\big] , $$
where we recall that $f (x)\! = \! e^{-x}\! -\! 1 + x$. Since $f$ is convex and by Jensen's inequality, we get  
$\psi_n (\lambda)\! -\! \alpha_n \lambda  \! \geq \! b_n f (\lambda \bE [V_n])$. 
Moreover, (\ref{2ordre}) implies $f(\lambda \bE [V_n]) \geq \frac{1}{2}(\lambda \bE [V_n])^2 \exp (-\lambda \bE [V_n])$. Since $\bE [V_n] \! \sim \! 1/a_n$, since 
$\alpha_n \! \rightarrow \! \alpha$ by $(\mathbf{C1})$ and since $b_n /a_n^2\! \rightarrow \! \beta_0 \! >\! 0$, there is $n_1\ino \bbN$ such that for all $n\! \geq \! n_1$, we get $1/2 \! \leq \! a_n \bE [V_n] \! \leq \! 2$, $\alpha_n \! \geq \! -2(\alpha)_-$ and $b_n /a_n^2 \! \geq \! \beta_0 /2$. Thus, 
%
%Then, note that for all $\lambda \ino [0, a_n]$, $\lambda \bE [V_n] \! \leq \! 1$ and observe that $\bE [V_n] \! \sim \! 1/a_n$ since $\alpha_n \! \rightarrow \! \alpha$ by $(\mathbf{C1})$. Therefore, there exists $n_0$ such that for all $n\! \geq \! n_0$ and for all $\lambda \! \in \! [0, a_n]$, 
%$\psi_n (\lambda) \! \geq \! (8e)^{-1}(b_n /a_n^2) \lambda^2$. Since $\beta_0 \! >\! 0$, there exists $n_1 \! \geq \! n_0$ such that $b_n /a_n^2 \! \geq \! \beta_0 /2$. Thus, we have proved that 
$$ \exists n_1 \ino \bbN : \; \forall n\! \geq \! n_1 , \; \forall \lambda \ino [0, a_n], \quad \psi_n (\lambda ) \geq -2(\alpha)_- \lambda + \frac{_1}{^{16e^2}} \beta_0 \lambda^2 \; , $$
which clearly implies (\ref{crinullo}). This completes the proof of Proposition \ref{Hcritos} $(i)$. \cqfd

\subsubsection{Proof of Proposition \ref{Hcritos} $(ii)$.}
\label{pfHcritos2}
Let us mention that, here, we closely follow the counterexample given in Le Gall \& D.~\cite{DuLG02}, p.~55.  
Fix $\alpha\ino \bbR$, $\beta  \ino [0, \infty)$, $\kappa\ino (0, \infty)$ and $\mathbf{c}\! =\! (c_j)_{j\geq 1} \! \in \! \elldo_3$. 
%We set $\pi (dr) \! = \! \sum_{j\geq 1} \kappa c_j \delta_{c_j} (dr)$ and
 For all $\lambda \! \in \! [0, \infty)$, set $\psi (\lambda) \! = \! \alpha \lambda + \frac{_1}{^2} \beta \lambda^2 + \sum_{j\geq 1} \kappa 
c_j \big( e^{-\lambda c_j} \! -\! 1 + \lambda c_j \big)$; assume that $\int^\infty \! d\lambda / \psi (\lambda) \! < \! \infty$. 
For all positive integers $n$, we next define $\mathbf{c}_n\! =\! (c^{_{(n)}}_j)_{j\geq 1}$ by setting 
$$ c^{_{(n)}}_j\! = \! c_j \; \,   \textrm{if $j\! \leq \! n$}, \quad c^{_{(n)}}_j\! = \! (\beta / (\kappa n))^{\frac{1}{3}}  
\; \,  \textrm{if $n\!  <\!  j\! \leq \! 2n$} \quad \textrm{and} \quad c^{_{(n)}}_j\! = \! 0 \;\,  \textrm{if $j\! >\! n$}. $$  
We also set $\psi_n  (\lambda) \! = \! \alpha \lambda + 
\sum_{j\geq 1} \kappa c^{_{(n)}}_j \big( \exp (-\lambda c^{_{(n)}}_j) \! -\! 1 + \lambda c^{_{(n)}}_j \big)$, $\lambda \ino [0, \infty)$. 
Let $(U_t^n)_{t\in [0, \infty)}$, be a CSBP with branching mechanism $\psi_n$ and with initial state $U^n_0\! = \! 1$. 
As $\lambda\! \rightarrow \! \infty$, observe that $\psi_n (\lambda) \! \sim \! (\alpha + \kappa \sigma_2 (\mathbf{c}_n) )\lambda $. 
Thus, $\int^\infty \! d\lambda / \psi_n (\lambda) \! =\! \infty$; by standard results on CSBP (recalled in Section \ref{CSBPApp} in 
Appendix), we therefore get 
\begin{equation}
\label{Unpersist}
\forall n\in \bbN , \; \forall t\! \in \! [0, \infty), \quad \bP \big( U^n_t \! >\! 0\big) \! = \! 1 \; .
\end{equation}
Let  $Z\! = \! (Z_t)_{t\in [0, \infty)}$ stands for a CSBP with  branching mechanism $\psi$ and with initial state $Z_0\! = \! 1$. 
Observe that for all $\lambda \ino [0, \infty)$, $\lim_{n\rightarrow \infty} \psi_n (\lambda) \! = \! \psi (\lambda)$. By 
standard results on CSBP (see Helland \cite{He78}, Theorem 6.1, p.~96), we get 
\begin{equation}
\label{Contlim}
U^n \underset{^{n \rightarrow \infty}}{ -\!\! \!-\!\! \! -\!\! \!\longrightarrow}  Z \; \, \textrm{weakly on $\bD ([0, \infty) , \bbR)$.}
\end{equation}
Next, let us fix $n\! \in \! \bbN$. By Proposition \ref{cvmarkpro} $(iv)$ 
%$\alpha \! = \! \beta \! = \! 0$, $\kappa \! \in \! (0, \infty)$ and $\mathbf{c}_n$:  
there exist sequences $\bw_{n, p}\! = \! (w^{_{(n,p)}}_{^j})_{j\geq 1} \! \in \! \elldo_f$ and $a_{n, p} , b_{n, p} \! \in \! (0, \infty)$, $p \! \in \!  \bbN$, such that 
\begin{equation}
\label{blourp!}
%\frac{w^{_{(n,p)}}_{^1} }{c_1}=a_{n, p}  \underset{^{p\rightarrow \infty}}{-\!\! \! -\!\! \!\longrightarrow}   \infty  , \quad  
\frac{a_{n, p} b_{n, p}}{\sigma_1 (w_{n, p})} \to  \kappa, \quad  
a_{n, p} ,  \frac{b_{n, p}}{a_{n, p}}   \; \textrm{and} \;   \frac{a^2_{n, p}}{b_{n, p}}  \underset{^{p\rightarrow \infty}}{-\!\! \! -\!\! \!\longrightarrow}   \infty , \quad  \frac{b_{n, p}}{a_{n,p}} \Big( 1-\frac{\sigma_2 (\bw_{n , p})}{\sigma_1 (\bw_{n, p})} \Big)\;  \underset{^{p\rightarrow \infty}}{-\!\!\! -\!\! \!  \longrightarrow}  \; \alpha
%\quad \textrm{and} \quad  
%\frac{b_{n, p}}{a^2_{n, p}} 
%\underset{^{p\rightarrow \infty}}{ -\!\! \! -\!\! \! \longrightarrow} 0 \; .
%% \quad \textrm{and} \quad  \frac{a_n b_n}{\sigma_1 (w_n)} \, \underset{^{n\rightarrow \infty}}{-\!\!\! -\!\! \! -\!\! \!\longrightarrow} \kappa , 
\end{equation} 
\begin{equation}
\label{erblurpi}
%\mathbf{(C1):} \quad \frac{a_nb_n}{\sigma_1 (\bw_n)}\;  \underset{n\rightarrow \infty}{-\!\!\! -\!\!\! \longrightarrow} \; \kappa \, , \qquad 
 \frac{b_{n, p}}{a^2_{n , p}}\!  \cdot \!  \frac{\sigma_3 (\bw_{n, p})}{\sigma_1 (\bw_{n, p})} \;  \underset{{p\rightarrow \infty}}{-\!\!\! -\!\! \!  \longrightarrow}  \;  \kappa \sigma_3 (\mathbf{c}_n) \quad \textrm{and} \quad \forall j \in \bbN^*, \quad \frac{w^{_{(n, p)}}_j}{a_{n, p} } \;  \underset{{p\rightarrow \infty}}{-\!\!\! -\!\! \! \longrightarrow}  \;  c^{_{(n)}}_j , 
\end{equation}
 and the following weak limit holds true on $\bD ([0, \infty) , \bbR)$: 
\begin{equation}
\label{blurpik}
 \Big( \frac{_1}{^{a_{n, p}}} Z^{_{(n, p)}}_{\lfloor b_{n, p} t/a_{n, p}\rfloor} \Big)_{t\in [0, \infty)}\underset{^{p\rightarrow \infty}}{-\!\!\! -\!\! \! -\!\! \!\longrightarrow} (U^n_t)_{t\in [0, \infty)} \; 
\end{equation} 
where $ (Z^{_{(n, p)}}_{^k})_{k\in \bbN}$ is a Galton-Watson Markov chain with $ Z^{_{(n, p)}}_{^0}\! = \! \lfloor a_{n ,p} \rfloor$ and with offspring distribution $\mu_{\bw_{n, p}}$ \ms{is} as in (\ref{rmupoissw}). 
%%%\begin{equation}
%%%\label{flounase}
%%%\forall k\! \in \! \bbN, \quad  \mu_{\bw_{n, p}} (k) = \sum_{j\geq 1} \frac{(w^{_{(n,p)}}_{^j})^{k+1}}{\sigma_1 (\bw_{n, p})\,  k!}\,  e^{-w^{_{(n,p)}}_{^j}} \; .
%%%\end{equation}
We also have $\lim_{p\rightarrow 0}\sqrt{\mathbf{j}_{n,p}}/ b_{n,p}\! = \! 0$, 
where $\mathbf{j}_{n,p}\! = \! \max \{ j \! \geq \! 1 \! :\!  w^{_{(n,p)}}_{^j} \! >\! 0 \}$. 
By Portemanteau's theorem for all $t\! \in [0, \infty)$, $\liminf_{p\rightarrow \infty} \bP \big(Z^{_{(n, p)}}_{\lfloor b_{n, p} t/a_{n, p}\rfloor} \! >\! 0  \big) \! \geq \! \bP (U^n_t \! >\! 0)=1$, by (\ref{Unpersist}). Thus, there exists $p_n \! \in \! \bbN$, such that 
\begin{equation}
\label{blourpi!}
\forall p \! \geq \! p_n , \qquad \bP \big(Z^{_{(n, p)}}_{\lfloor b_{n, p} n/a_{n, p}\rfloor} \! >\! 0  \big) \geq 1-2^{-n} \; .
\end{equation}
Without loss of generality we can furthermore assume that $\sqrt{\mathbf{j}_{n,p_n}}/ b_{n,p_n}\! \leq \! 2^{-n}$ and  
%that the Prokhorov distance between the law of 
%$(a_{n, p_n}^{-1} Z^{_{(n, p_n)}}_{\lfloor b_{n, p_n} t/a_{n, p_n}\rfloor})_{t\in [0, \infty)}$ and the law of $U^n$ is bounded by $2^{-n}$ and that 
$$ a_{n, p_n}, \;   \frac{b_{n, p_n}}{a_{n, p_n}}  \; \textrm{ and} \;   \frac{a^2_{n, p_n}}{b_{n, p_n}} \geq 2^n , \quad 
\Big| \frac{b_{n, p_n}}{a_{n,p_n}} \Big( 1-\frac{\sigma_2 (\bw_{n , p_n})}{\sigma_1 (\bw_{n, p_n})} \Big) - \alpha \Big| \leq 2^{-n} , \quad \Big| \frac{a_{n, p_n} b_{n, p_n}}{\sigma_1 (w_{n, p_n})} -  \kappa \Big| \leq 2^{-n}$$
$$ \Big| \frac{b_{n, p_n}}{a^2_{n , p_n}}\!  \cdot \!  \frac{\sigma_3 (\bw_{n, p_n})}{\sigma_1 (\bw_{n, p_n})}-  \kappa \sigma_3 (\mathbf{c}_n)\Big|  \leq 2^{-n} \quad \textrm{and} \quad 
\forall j\ino \{ 1, \ldots , n\}, \; \Big| \frac{w^{_{(n, p_n)}}_j}{a_{n, p_n} } - c^{_{(n)}}_j \Big| \leq 2^{-n} .  $$
%Note that $\kappa \sigma_3 (\mathbf{c}_n)\! \rightarrow  \beta + \kappa \sigma_3 (\mathbf{c})$ as $n\! \rightarrow \! \infty$. 
Set $a_n \! =\!  a_{n ,p_n}$, $b_n \! = \! b_{n , p_n}$ and $\bw_n \! =\!  \bw_{n, p_n}$. 
Note that $\kappa \sigma_3 (\mathbf{c}_n)\! \rightarrow  \beta + \kappa \sigma_3 (\mathbf{c})$ as $n\! \rightarrow \! \infty$. Thus, 
$a_n , b_n $ and $\bw_n$ satisfy (\ref{apriori}) with $\beta_0\! = \! 0$, $\mathbf{(C1)}$, $\mathbf{(C2)}$, $\mathbf{(C3)}$ and $\sqrt{\mathbf{j}_{n}}/ b_{n}\! \rightarrow \! 0$.  
%(\ref{suitt}), (\ref{rightsca}) and (\ref{defkappa}). 
Set $ Z^{_{(n)}}_{k}\! = \!  Z^{_{(n, p_n)}}_{k}$. 
%let $\mu_n$ be derived from $w_n$ as in (\ref{mixpois}). Then, (\ref{flounase}) implies that $ (Z^{_{(n)}}_{^k})_{k\in \bbN}$  is a Galton-Watson process with initial state $ Z^{_{(n)}}_{^0}\! = \! \lfloor a_{n } \rfloor$ and with offspring distribution $\mu_{n}$. 
%(\ref{Contlim}) entails that 
%$$ \Big( \frac{_1}{^{a_{n}}} Z^{_{(n)}}_{\lfloor b_{n} t/a_{n}\rfloor} \Big)_{t\in [0, \infty)}\underset{^{n\rightarrow \infty}}{-\!\!\! -\!\! \! -\!\! \!\longrightarrow} (Z_t)_{t\in [0, \infty)} \; $$
By (\ref{blourpi!}), for all $\delta \! \in \! (0, \infty)$, and all integers $n \! \geq \! \delta$, we easily get 
$ \bP \big(Z^{_{(n)}}_{\lfloor b_{n} \delta /a_{n}\rfloor} \! =\! 0  \big) \leq 
\bP \big(Z^{_{(n)}}_{\lfloor b_{n} n/a_{n}\rfloor} \! =\! 0  \big) \leq 2^{-n} $. Consequently, $\lim_{n\rightarrow \infty} \bP \big(Z^{_{(n)}}_{\lfloor b_{n} \delta /a_{n}\rfloor} \! =\! 0  \big)\! = \! 0$, for all $\delta \ino (0, \infty)$. Namely, $\mathbf{(C4)}$ is not satisfied, which 
completes the proof of Proposition \ref{Hcritos} $(ii)$. \cqfd

\subsubsection{Proof of Proposition \ref{Hcritos} $(iii)$.}
\label{pfHcritos3}
Fix $\alpha \ino \bbR$, $\beta  \ino [0, \infty)$, $\kappa\ino (0, \infty)$ and $\mathbf{c}\! =\! (c_j)_{j\geq 1} \! \in \! \elldo_3$. 
%We set $\pi (dr) \! = \! \sum_{j\geq 1} \kappa c_j \delta_{c_j} (dr)$ and
 For all $\lambda \! \in \! [0, \infty)$, set $\psi (\lambda) \! = \! \alpha \lambda + \frac{_1}{^2} \beta \lambda^2 + \sum_{j\geq 1} \kappa 
c_j \big( e^{-\lambda c_j} \! -\! 1 + \lambda c_j \big)$; assume that $\int^\infty \! d\lambda / \psi (\lambda) \! < \! \infty$.  
We consider several cases.

\medskip

\noi
$\bullet$ \textbf{Case 1}: we first assume that $\beta \! \geq \! \beta_0 \! >\! 0$. By Proposition \ref{cvmarkpro} $(iv)$ there exists $a_n, b_n , \bw_n$ satisfying (\ref{apriori}) with $\beta_0 \! >\! 0$, 
$(\mathbf{C1})$, $(\mathbf{C2})$ and $(\mathbf{C3})$. But Proposition \ref{Hcritos} $(i)$ (proved in Section \ref{pfHcritos1}) 
asserts that $a_n, b_n , \bw_n$ necessarily satisfy $(\mathbf{C4})$. This proves 
Proposition \ref{Hcritos} $(iii)$ in Case 1.

\medskip

\noi
$\bullet$ \textbf{Case 2}. We next assume that $\beta \! >\! 0$ and $\beta_0 \! = \! 0$, and we set: 
\begin{equation}
\label{defqnj} 
q_{^j}^{_{(n)}}=\left\{
\begin{array}{ll}
c_j & \text{if} \quad  j \! \in \!  \big\{ 1, \ldots ,    n \!  \big\} 
%1\! \leq \! j\! \leq \! \lfloor  \rho_n 2^n \! \rfloor 
,  \\
(\beta/\kappa)^{\frac{1}{3}}n^{-1} & \text{if} \quad  j \! \in \!  \big\{n  +1 , \ldots , n   + n^3  \big\}, \\
n^{-3} & \text{if} \quad  j \! \in \!  \big\{ n   +n^3 +1 , \ldots , n   + n^3+n^8 \big\}, \\
% \lfloor \rho_n 2^n\! \rfloor +1 \! \leq \! j\! \leq \! \lfloor \rho_n 2^n\! \rfloor + n2^n , \\
\; \; \; 0 & \text{if} \quad  j >  n + n^3 + n^8.
\end{array}\right.
\end{equation}
Denote by $\bv_n \! =\!  (v^{_{(n)}}_{^j})_{j\geq 1}$ the nonincreasing rearrangement of $\mathtt{q}_n \! =\!  (q^{_{(n)}}_{^j})_{j\geq 1}$. Thus, $\sigma_p (\bv_n)\! = \! \sigma_p (\mathtt{q}_n)$ 
for any $p\! \in \!  (0, \infty)$. Since $\sum_{1\leq j\leq n} c_j^p  \leq  c_1^p n $, we easily get $\kappa \sigma_1 (\bv_n)\! \sim \! \kappa n^5$, $\kappa \sigma_2 (\bv_n) \! \sim \! \kappa n^2 $ and   
$\kappa \sigma_3 (\bv_n) \! \rightarrow \! \beta + \kappa \sigma_3 (\mathbf{c} )$. We next set $b_n \! =\!  \kappa  \sigma_1 (\bv_n)$, 
$a_n \! =\!  \kappa  \sigma_1 (\bv_n)/ (\kappa  \sigma_2 (\bv_n) + \alpha ) $ and for all $j\! \geq \! 1$, $w^{_{(n)}}_{j} \! =\!  a_n v^{_{(n)}}_{j}$. Note that $a_n \sim n^3$.
Then, it is easy to check that 
$a_n , b_n $ and $\bw_n$ satisfy (\ref{apriori}) with $\beta_0 \! = \! 0$, $(\mathbf{C1})$, $(\mathbf{C2})$ and $(\mathbf{C3})$. Since $\mathbf{j}_n\! = \! \max \{j\! \geq \! 1 \! : \! w^{_{(n)}}_{^j} \! >\! 0\} \! \leq \! n+ n^3 + n^8$, we easily get $\sqrt{\mathbf{j}_n}/b_n \! \rightarrow \! 0$.  Here observe 
that $\kappa \! = \! a_n b_n / \sigma_1 (\bw_n)$ and  $b_n \big(1\! -\!  (\sigma_2 (\bw_n)/ \sigma_1 (\bw_n)) \big)/ a_n \! = \! \alpha $.

We next prove that $(\mathbf{C4})$ holds true by proving that (\ref{crinullo}) in Proposition \ref{Hcritos} $(i)$ holds true. 
To that end, we introduce $f_\lambda (x)= x\big( e^{-\lambda x}\! -\! 1 + \lambda x\big)$, for all $x, \lambda \ino [0, \infty)$,  
and we recall from (\ref{rrXwdef}) the definition of $X^{\bw_n}$; we denote by $\psi_n$ the Laplace exponent of $\frac{1}{{a_n}} X^{\bw_n}_{b_n \cdot}$. 
We first observe that for all $\lambda \ino [0, \infty)$

\begin{equation}
\label{rrLapXwnn}
\psi_n (\lambda)= \alpha \lambda + \kappa \!   \sum_{j\geq 1} \! f_\lambda ( q^{_{(n)}}_{^j}) \! =\!   \alpha \lambda + \kappa \! \!\!   \sum_{1\leq j\leq n} \!\!  f_\lambda ( c_j)+  \kappa n^3 f_\lambda ( (\beta/\kappa)^{\frac{1}{3}}n^{-1} ) + \kappa n^8 f_\lambda ( n^{-3} ) .
\end{equation}
Set $s_0\! =\!  (\beta/ \kappa)^{1/3}$ and recall from (\ref{2ordre}) that $f_{\lambda} (x)\! \geq  \! \frac{1}{2}x^3 \lambda^2 e^{-\lambda x}$. 
Thus, if $\lambda \ino [1, 2n/s_0]$, then 
$$   \psi_n (\lambda) +(\alpha)_- \lambda \geq \kappa n^3 f_\lambda (s_0/n) \geq  \frac{_1}{^2}e^{- 2} \beta \lambda^2 =: s_1 \lambda^2. $$ 
Namely, for $\lambda \ino [1, 2n/s_0]$, $\psi_n (\lambda) \geq   s_1 \lambda^2 (1 -\frac{(\alpha)_-}{s_1\lambda} )$. 
Next observe that, $f_\lambda (x) \! \geq \! x(\lambda x -1)$; thus, if $\lambda \ino [2n/s_0 , n^3]$, then 
$$  \psi_n (\lambda)  \geq -(\alpha)_- \lambda+ \kappa n^3 f_\lambda ( s_0/n) \geq -(\alpha)_- \lambda+ \kappa s_0  n^2 \big(\tfrac{s_0  \lambda}{n} \! -\! 1\big) =\kappa s_0  n^2 \Big( \big(1-\tfrac{(\alpha)_-}{\kappa s_0^2n} \big) \tfrac{s_0\lambda}{n} -1 \Big) . $$ 
Thus, for all $y \! >\! \tfrac{2(\alpha)_-}{s_1} \mm{\vee 1}$ and for all $n\geq \tfrac{y s_0}{2} \vee \tfrac{3(\alpha)_-}{\kappa s_0^{2}}$, we get
 $$ \int_y^{n^3} \frac{d\lambda}{\psi_n (\lambda)} \leq 2 \int_y^{\frac{2n}{s_0}} \frac{d \lambda }{s _1\lambda^2} + \int_{\frac{2n}{s_0}}^{n^3}  
 \frac{d \lambda }{\kappa s_0 n^2 \big(   \tfrac{2s_0}{3n} \lambda  \! -\! 1\big) } \leq  \frac{2 }{s_1 y}+ \frac{3\log (\tfrac{2}{3}s_0n^2 \! -\! 1)+ 3\log3}{2\kappa s^2_0n} . $$
Since $a_n \sim n^3$, it proves that 
$\psi_n$ satisfies (\ref{crinullo}), and $(\mathbf{C4})$ holds true. This proves  Proposition \ref{Hcritos} $(iii)$ in Case 2.

\smallskip

\noi
$\bullet$ \textbf{Case 3}: We now assume that $\beta\! = \! \beta_0 \! = \! 0$. Let $\beta_n \ino (0, \infty)$ be a sequence decreasing to $0$. For all $n \ino \bbN^*$, we set 
$\Psi_n (\lambda) \! = \! \psi (\lambda) + \frac{_1}{^2} \beta_n \lambda^2\! = \! \alpha \lambda + \frac{_1}{^2} \beta_n \lambda^2 + \sum_{j\geq 1} \kappa 
c_j \big( e^{-\lambda c_j} \! -\! 1 + \lambda c_j \big)$. We now fix $n\ino \bbN^*$; by $\textbf{Case 2}$, there exists 
$\bw_{n, p}\! = \! (w^{_{(n,p)}}_{^j})_{j\geq 1} \! \in \! \elldo_f$ and $a_{n, p} , b_{n, p} \! \in \! (0, \infty)$, $p \! \in \!  \bbN$, that satisfy $\sqrt{\mathbf{j}_{n,p}}/ b_{n,p}\! \rightarrow \! 0$ as $p\! \rightarrow \! \infty$, where $\mathbf{j}_{n,p}\! = \! \max \{ j \! \geq \! 1 \! :\!  w^{_{(n,p)}}_{^j} \! >\! 0 \}$, and 
\begin{equation}
\label{rblourp!}
%\frac{w^{_{(n,p)}}_{^1} }{c_1}=a_{n, p}  \underset{^{p\rightarrow \infty}}{-\!\! \! -\!\! \!\longrightarrow}   \infty  , \quad  
\frac{a_{n, p} b_{n, p}}{\sigma_1 (w_{n, p})} = \kappa, \quad  
a_{n, p} ,  \frac{b_{n, p}}{a_{n, p}}   \; \textrm{and} \;   \frac{a^2_{n, p}}{b_{n, p}}  \underset{^{p\rightarrow \infty}}{-\!\! \! -\!\! \!\longrightarrow}   \infty , \quad  \frac{b_{n, p}}{a_{n,p}} \Big( 1-\frac{\sigma_2 (\bw_{n , p})}{\sigma_1 (\bw_{n, p})} \Big)\;  =\alpha\,,
%\quad \textrm{and} \quad  
%\frac{b_{n, p}}{a^2_{n, p}} 
%\underset{^{p\rightarrow \infty}}{ -\!\! \! -\!\! \! \longrightarrow} 0 \; .
%% \quad \textrm{and} \quad  \frac{a_n b_n}{\sigma_1 (w_n)} \, \underset{^{n\rightarrow \infty}}{-\!\!\! -\!\! \! -\!\! \!\longrightarrow} \kappa , 
\end{equation} 
\begin{equation}
\label{rerblurpi}
%\mathbf{(C1):} \quad \frac{a_nb_n}{\sigma_1 (\bw_n)}\;  \underset{n\rightarrow \infty}{-\!\!\! -\!\!\! \longrightarrow} \; \kappa \, , \qquad 
 \frac{b_{n, p}}{a^2_{n , p}}\!  \cdot \!  \frac{\sigma_3 (\bw_{n, p})}{\sigma_1 (\bw_{n, p})} \;  \underset{{p\rightarrow \infty}}{-\!\!\! -\!\! \!  \longrightarrow}  \;  \beta_n + \kappa \sigma_3 (\mathbf{c}) \quad \textrm{and} \quad \forall j \in \bbN^*, \quad \frac{w^{_{(n, p)}}_j}{a_{n, p} } \;  \underset{{p\rightarrow \infty}}{-\!\!\! -\!\! \! \longrightarrow}  \;  c_j . 
\end{equation}
and 
\begin{equation}
\label{rcvsvie} 
\forall n\ino \bbN^*, \; \forall t\ino [0, \infty), \quad  \lim_{n\rightarrow \infty} \bP \big( Z^{_{(n, p)}}_{\lfloor b_{n, p} t /a_{n , p} \rfloor} \! = \! 0 \big) = e^{-v_n(t)} \quad \textrm{where} \quad  \int_{v_n(t)}^\infty \! 
\frac{d\lambda }{\Psi_n (\lambda)}= t .
\end{equation}
Here, $(Z^{_{(n, p)}}_k)_{k\in \bbN}$ is a Galton-Watson process with offspring distribution $\mu_{\bw_{n, p}}$ given by (\ref{rmupoissw}) and where $Z^{_{(n, p)}}_0\! = \! \lfloor a_{n, p} \rfloor $. 
Let $v\! : \! (0, \infty) \! \rightarrow \! (0, \infty)$ be such that $t \! = \! \int_{v(t)}^\infty d\lambda / \psi (\lambda) $ for all $t\ino (0, \infty)$. Since $\Psi_n (\lambda ) \! \geq \! \psi (\lambda)$, we get 
$\int_{v(t)}^\infty d\lambda / \psi (\lambda)  \! = \! t \! \leq \!  \int_{v_n(t)}^\infty d\lambda / \psi (\lambda)$; thus $v_n (t) \! \leq \! v(t)$. Thus, 
there exists $p_n \! \ino \! \bbN$ such that for all $p\! \geq \! p_n$,   
$ \bP \big( Z^{_{(n, p)}}_{\lfloor b_{n, p} /a_{n , p} \rfloor} \! = \! 0 \big) \! \geq \! \frac{1}{2}\exp (-v_n (1)) \! \geq \!  \frac{1}{2}\exp (-v (1))$. Without loss of generality, we can assume that $\sqrt{\mathbf{j}_{n,p_n}}/ b_{n,p_n} \! \leq \! 2^{-n}$, 
$a_{n, p_n} , b_{n ,p_n}/ a_{n, p_n}$ and $a_{n, p_n}^2/ b_{n, p_n} \geq 2^n$, that 
for all $1\! \leq \! j \! \leq \! n$, $| w^{_{(n,p_n)}}_j / a_{n, p_n} \! -\! c_j | \leq 2^{-n}$ and 
$$ \Big|  \frac{b_{n, p_n}}{a^2_{n , p_n}}\!  \cdot \!  \frac{\sigma_3 (\bw_{n, p_n})}{\sigma_1 (\bw_{n, p_n})}  -  \kappa \sigma_3 (\mathbf{c}) \Big| \leq 2\beta_n \longrightarrow 0 . $$
If one sets $a_n \! =\! a_{n ,p_n }$, $b_{n, p_n}$ and $\bw_{n }\! = \!   \bw_{n , p_n}$, then we have proved that $a_n , b_n , \bw_n$ satisfy (\ref{apriori}) with $\beta\! = \! \beta_0\! = \! 0$,  $\sqrt{\mathbf{j}_n}/b_n\! \rightarrow \! 0$ 
and $(\mathbf{C1})$--$(\mathbf{C4})$, which proves Proposition \ref{Hcritos} $(iii)$ in Case 3. This completes the proof of  Proposition \ref{Hcritos} $(iii)$. \cqfd

\section{Proof of Lemma \ref{powerlaw}.}
\label{powerlawpf}
We first prove the following lemma. 
\begin{lem}
\label{slowV} Let $\ell \! : \! (0, 1] \! \rightarrow \! (0, \infty)$ be a measurable 
slowly varying function at $0+$ such that  for all $x_0 \ino (0,1)$, $\sup_{x\in  [x_0 , 1]} \ell(x) \! < \! \infty$. Then, for all $\delta \ino (0, \infty)$, there exist $\eta_\delta \ino (0, 1]$ and $c_\delta \ino (1, \infty)$ such that 
\begin{equation}
\label{paoutai}
 \forall y \ino (0, \eta_\delta) , \; \forall z\ino (y, 1], \quad \frac{1}{c_\delta} \Big( \frac{z}{y}\Big)^{\! -\delta}\!\!\!  \leq \, \frac{\ell (z)}{\ell (y)}\,  \leq 
c_\delta  \Big( \frac{z}{y}\Big)^{\delta} \; . 
\end{equation}
\end{lem}
\noi
\textbf{Proof.} The measurable version of the representation theorem for slowly varying functions (see for instance Bingham, Goldie \& Teugels \cite{BiGoTe}) implies that there exist two measurable functions $c\! : \! (0, 1]\! \rightarrow \! \bbR$ and $\epp \! : \! (0, 1] \! \rightarrow \! [-1, 1]$ 
such that $\lim_{x\rightarrow 0+} c(x)\! = \! \gamma \ino \bbR$, such that $\lim_{x\rightarrow 0+} \epp (x)\! = \! 0$, and such that 
$\ell (x) \! =\!  \exp( c(x) + \int_x^1 ds\,  \epp (s)/ s )$, for all $x\ino (0, 1]$. 
%$$ \forall x\ino (0, 1], \ell (x) = \exp \Big( c(x) +\!  \int_x^1    \! \! \!  ds\,  \frac{\epp (s)}{s}  \Big). $$ 
Since, $\sup_{x\in  [x_0 , 1]} \ell(x) \! < \! \infty$, for all $x_0 \ino (0, 1)$, 
we can assume without loss of generality that $c$ is bounded. Fix $\delta \ino (0, \infty)$ and let 
$\eta_\delta \ino (0, 1]$ be such that $\sup_{(0, \eta_\delta ]} |\epp | \leq \delta $. Fix $y \ino (0, \eta_\delta)$ and  $z\ino (y, 1]$; if $z \! \leq \! \eta_\delta$, then note that $\int_y^z ds \, |\epp (s) | /s \leq \delta \log (z/y)$; if $\eta_\delta \! \leq \! z$, then observe that $\int_y^z ds \, |\epp (s) | /s \leq \delta \log( \eta_\delta /y)+ \int_{\eta_\delta}^1 ds \, |\epp (s) | /s \! \leq \!  \delta \log (z/y) + 
\log (1/ \eta_\delta) $. Thus  
$$  \eta_{\delta}  
e^{-2\lVert c\rVert_\infty}  \Big( \frac{z}{y}\Big)^{\! -\delta}\!\!\!   \leq  \frac{\ell (z)}{\ell (y)}=  \exp \Big( c(z) \! -\! c(y) \! - \!\!  \int_y^z    \! \! \!  ds\,  \frac{\epp (s)}{s}  \Big) \leq \eta_{\delta}^{-1}  
e^{2\lVert c\rVert_\infty}  \Big( \frac{z}{y}\Big)^{\delta} \; , $$
which implies the desired result. \cqfd 

\medskip

Recall from (\ref{Wassu}) that $W\! :\!  \Omega \! \rightarrow [0, \infty)$ is a r.v.~such that $r\! : = \! \bE [W] \! = \! \bE [ W^2] \! <\!  \infty$ and 
such that $\bP (W \geq x) \! =\!  x^{-\rho} L(x)$ where $L$ is a slowly varying function at $\infty$ and $\rho \ino (2, 3)$.  
Recall from (\ref{GGdef}) that for all $y\ino [0, \infty)$, we have set $G(y) \! = \! \sup \{ x\ino [0, \infty) : 
\bP (W\! \geq \! x) \! \geq \! 1\! \wedge \! y\}$. Note that $G$ is non increasing and that it is null on $[1, \infty)$. Then, 
$G(y) \! =\!  y^{-1/\rho} \, \ell (y)$, where $\ell$ is slowly varying at $0$.  
Recall from  (\ref{pwrlwabw}) that $\kappa , q \ino (0, \infty)$ and that $ a_n\!  \sim \! q^{-1} G(1/n)$, $w_{^j}^{_{(n)}}\! \! =\!  G(j/n)$, 
$j\! \geq \! 1$, and $ b_n \! \sim\!  \kappa \sigma_1 (\bw_n) /a_n$. 

Fix $a\!\in \!  [1, 2]$ and observe that $\sigma_a (\bw_n)\! = \! \sum_{1\leq j< n}\int_0^{G(1/n)} \! dz \, az^{a-1}\un_{\{z \tt{<}G(j/n) \}}$. But observe that $z\! < \! G(y)$ implies $y \! \leq \! P(W\! \geq \! z)$. Thus, 
\begin{eqnarray}
\label{trantran}
\sigma_a (\bw_n)\!\!\!\! & = &\!\!\!\!   \sum_{1\leq j< n}\! \int_0^{G(1/n)} \!\!\!\!\!\!\!   \!\!\!  \!\! dz \, az^{a-1} \un_{\{j  \leq n \bP(W\geq z) \}}= \int_0^{G(1/n)}  \!\!\!\!\!\!\!   \!\!\!  \!\! dz \, az^{a-1}\!\!\! \sum_{1\leq j< n}\un_{\{j  \leq n\bP(W\geq z) \}} \nonumber \\
\!\!\!\! & = &\!\!\!\! \int_0^{G(1/n)}  \!\!\!\!\!\!\!   \!\!\!  \!\!  dz \, az^{a-1} \lfloor n \bP (W\! \geq \! z) \rfloor =  
 \int_0^{G(1/n)}  \!\!\!\!\!\!\!   \!\!\!  \!\!  dz \, az^{a-1}   n\bP (W\! \geq \! z) -  \int_0^{G(1/n)}  \!\!\!\!\!\!\!   \!\!\!  \!\!  dz \, az^{a-1} 
\{  n \bP (W\! \geq \! z) \} \nonumber  \\
\!\!\!\! & = &\!\!\!\! n   \int_0^\infty   \!\!  dz \, az^{a-1}   \bP (W\! \geq \! z) - \int_{G(1/n)}^\infty   \!\!\!\!\!\!\!   \!\!\!  \!\!  dz \, az^{a-1}  n \bP (W\! \geq \! z) -  \int_0^{G(1/n)}  \!\!\!\!\!\!\!   \!\!\!  \!\!  dz \, az^{a-1} 
\{  n \bP (W\! \geq \! z) \} .
\end{eqnarray}
Note that  $\int_0^\infty dz \, az^{a-1}   \bP (W\! \geq \! z) \! = \! \bE [W^a] \! < \! \infty$. Recall from (\ref{chitithypo}) that $\bP (W\! = \! G(1/n))\! = \! 0$, which easily implies that $\bP (W\! \geq \! G(1/n) ) \! = \! 1/n$. Thus, 
$$n\bP (W\! \geq z)\! =\!  \bP (W\geq z)/ \bP (W\geq G(1/n))= (z/G(1/n))^{-\rho} L(z)/L(G(1/n))$$ 
and by (\ref{trantran}) and the change of variable $z\!  \mapsto \! z/G(1/n)$, we get 
$$ \sigma_a (\bw_n)\! = \! n \bE [W^a] -G \big(\frac{_1}{^n}\big)^{\! a} \!\!  \int_1^\infty \!\! \!\! dz \, az^{a-1 -\rho} \frac{L(zG(\frac{_1}{^n}))}{L(G(\frac{_1}{^n}))} -G\big(\frac{_1}{^n}\big)^{\! a} \!\!  
\int_0^{1}  \!\!  \!\! dz \, az^{a-1} \Big\{ z^{-\rho} \frac{L(zG(\frac{_1}{^n}))}{L(G(\frac{_1}{^n}))}  \Big\} . $$
The measurable version of the representation theorem for slowly varying functions (see for instance \cite{BiGoTe}) implies that there exist two measurable functions $c\! : \! (0, \infty)\! \rightarrow \! \bbR$ and $\epp \! : \! (0, \infty) \! \rightarrow \! [-1, 1]$ 
such that $\lim_{x\rightarrow \infty} c(x)\! = \! \gamma \ino \bbR$, such that $\lim_{x\rightarrow \infty} \epp (x)\! = \! 0$, and such that 
$L (x) \! =\!  \exp( c(x) + \int_1^x ds\,  \epp (s)/ s )$, for all $x\ino (0, \infty)$. We then 
set $u\! = \! (\rho \! -\! a)/2$ that is a strictly positive quantity since $a \! \leq \! 2 \! < \! \rho$. Let $n_0$ be such that for all $n\! \geq \! n_0$, $\sup_{s \in [1, \infty)} |\epp (s G(1/n))|  \! \leq \! u$. Thus, for all $z \ino [1, \infty)$, 
$$0 \leq \!  z^{a-1 -\rho}\frac{L(zG(\frac{_1}{^n}))}{L(G(\frac{_1}{^n}))} \! = \! z^{a-1 -\rho} \exp \Big( c\big(zG \big(\frac{_1}{^n}\big)\big)\! -\! 
c\big(G \big(\frac{_1}{^n}\big)\big) + \int_1^z\!\! ds  \frac{\epp \big( s G \big(\frac{1}{n}\big)\big)}{s}  \Big) \leq e^{2 \lVert c \rVert_\infty} z^{-1-u}. $$   
Since for all $z\ino [1, \infty)$, $L(zG(1/n))/L(G(1/n))\! \rightarrow \! 1$, dominated convergence entails: 
$$\lim_{n\rightarrow \infty\;\;  }\!\!\!\!  \int_1^\infty \!\! \!\! \!\! dz \, az^{a-1 -\rho} \frac{L(zG(\frac{_1}{^n}))}{L(G(\frac{_1}{^n}))}
%=\!  \int_1^\infty \!\! \!\! dz \, az^{a-1 -\rho} 
\! =\!  \frac{a}{\rho \! -\! a} \; \,  \textrm{and}  \,   \lim_{n\rightarrow \infty\; \; }\!\!\!\!  \int_0^{1}  \!\!  \!\! dz \, az^{a-1} \Big\{ z^{-\rho} \frac{L(zG(\frac{_1}{^n}))}{L(G(\frac{_1}{^n}))}  \Big\} \! = \!\! \int_0^{1}  \!\!  \!\! dz \, az^{a-1} \{ z^{-\rho}\} . $$
We then set $Q_a\! = \! a/ (\rho \! - \! a)+  \int_0^{1}  dz \, az^{a-1} \{ z^{-\rho}\}$ and since $a_n \sim q^{-1} G(1/n)$, we have proved that 
\begin{equation}
\label{altimo}
\sigma_a (\bw_n) \! = \! n \bE [W^a] -q^a Q_a (a_n)^a + o((a_n)^a) . 
\end{equation}
Recall that $r \! = \! \bE [W]\! = \! \bE [W^2]$ and take (\ref{altimo}) with $a\! = \! 1$ to get 
$\sigma_1 (\bw_n) \! -\! r n\!  \sim \! - Q_1n^{1/\rho} \ell (1/n)$
since $a_n \sim q^{-1}n^{1/\rho} \ell (1/n)$; thus $b_n \sim \kappa q r n^{1-1/\rho} / \ell (1/n)$. It implies that $a_n$ and $b_n / a_n$ go to $\infty$ and that $b_n / a_n^2 \! \rightarrow \! 0$. Moreover for all $j\! \geq 1$, $w^{_{(n)}}_{^j} /a_n\! \rightarrow \! q j^{-1/\rho}$. This 
implies that $a_n$, $b_n$ and $\bw_n$ satisfy (\ref{apriori}) with $\beta_0 \! = \! 0$ (and ($\mathbf{C3}$)). Since
 $a_n b_n \sim \kappa \sigma_1 (\bw_n)\sim \kappa r n $, (\ref{altimo}) with $a\! = \! 1$ and $2$ implies 
\begin{equation}
\label{gluten}
 \frac{\sigma_2 (\bw_n)}{\sigma_1 (\bw_n) }\! =\!  \frac{nr\! -\! q^2Q_2a_n^2 + o(a_n^2)}{nr\! -\! qQ_1a_n + o(a_n)} \! =\!  1\! -\!  \kappa q^2 Q_2 \frac{a_n}{b_n} + o\Big(\frac{a_n}{b_n} \Big) \! =\!  1\! -\! \alpha_0 \frac{a_n}{b_n}  + o\Big(\frac{a_n}{b_n} \Big)\; , 
\end{equation} 
where $\alpha_0\! = \! \kappa q^2 Q_2 $ as defined in (\ref{glutamate}).

Next, for all $\alpha \ino \bbR$, set $w^{_{(n)}}_{^j}(\alpha) \! = \!  (1 \! -\! \frac{{a_n}}{{b_n}} (\alpha -  \alpha_0))w^{_{(n)}}_{^j}$. By (\ref{gluten}), we get $\sigma_2 (\bw_n (\alpha))/ \sigma_1 (\bw_n (\alpha))\! =\! 1-\alpha  a_n / b_n+ o( a_n / b_n)$. Namely, $\bw_n (\alpha)$ satisfies ($\mathbf{C1}$). Since $w^{_{(n)}}_{^j}(\alpha)\sim_n w^{_{(n)}}_{^j}$, $\bw_n (\alpha)$ also satisfies ($\mathbf{C3}$) with $c_j\! = \! qj^{-1/\rho}$, $j\! \geq \! 1$. 

 Let us proves that $(\bw_n (\alpha))$ satisfies ($\mathbf{C2}$). First observe that $\sigma_3 (\bw_n (\alpha))\! \sim \! 
 \sigma_3 (\bw_n)$. So we only need to prove that the $\bw_n$ satisfy  ($\mathbf{C2}$).  
To that end, for all $n$ and $j\! \geq \! 1$, we set $f_n (j)\! = \! (G(j/n)/ G(1/n))^3\! =\!  j^{-3/\rho} \ell^3 (j/n)/\ell^3 (1/n)$ and $\delta\! =\! \frac{_1}{^2} 
(\frac{_3}{^\rho}\! -\! 1)$ that is strictly positive. We apply Lemma \ref{slowV} to $\ell^3$:  
let $c_\delta \ino (1, \infty)$ and 
$\eta_\delta\ino (0, 1]$ such that (\ref{paoutai}) holds true; then, 
for all $n\! >\! 1/ \eta_\delta$, $0 \! \leq \! f_n (j) \! \leq \! c_\delta j^{-1-\delta}$. Since for all $j\! \geq \! 1$, $\lim_{n\rightarrow \infty}
f_n(j)\! = \! j^{-3/\rho}$, by dominated convergence we get:
$$ G(1/n)^{-3} \sigma_3 (\bw_n) = \sum_{1\leq j\leq n} f_n (j) \underset{n\rightarrow \infty}{-\!\!\! -\!\!\! -\!\!\! \longrightarrow} \sum_{j\geq 1} j^{-3/\rho} = q^{-3} \sigma_3 (\mathbf{c}), $$            
which easily implies ($\mathbf{C2}$).

  Let us prove  that $\bw_n (\alpha)$ satisfies ($\mathbf{C4}$) thanks to  (\ref{crinullo}) in Proposition \ref{Hcritos}.  
To that end, we fix $n\ino \bbN^*$ and $\lambda\ino [0, \infty)$ such that 
$\lambda \ino [1, a_n]$. 
For all $x\ino [0, \infty)$, recall that $f_\lambda (x)\! = \! x (e^{-\lambda x} \!  - \! 1 + \lambda x)$ and for all $j\! \geq \! 1$, 
set 
$$\phi_n (j)\! = \! f_\lambda \Big(\frac{_{w^{_{(n)}}_{^j}(\alpha)}}{^{a_n}} \Big)\! = \! f_\lambda \Big(q_n j^{-1/\rho} \frac{_{\ell (j/n)}}{^{\ell (1/n)}} \Big) 
\; \textrm{where} \quad q_n \! = \! 
\Big(1\! -\! \frac{_{a_n}}{^{b_n}}(\alpha \! -\! \alpha_0) \Big) \frac{_{G(1/n)}}{^{a_n}}  \sim  q \; .$$ 
To simplify, we also set $\kappa_n\!  =\!  a_n b_n / \sigma_1 (\bw_n (\alpha))$; note that $\kappa_n \! \sim \! \kappa $. Let $\delta \ino 
(0, \infty)$ be specified further; by Lemma \ref{slowV} and the previous arguments, there exists $c_{\delta } \ino (0, \infty)$ and 
$n_\delta$ such that for all $n \! \geq \! n_\delta$, $w^{_{(n)}}_{^j} (\alpha) / a_n\! \geq \! c_\delta j^{ -\delta-1/\rho}$ and $\kappa_n \! \geq \! \frac{_1}{^2} \kappa$, which entails $\kappa_n \phi_n (j) \! \geq \! \frac{_1}{^2} \kappa f_\lambda (c_\delta j^{ -\delta-1/\rho})$. We next set: 
$$\alpha_n \! := \!  \frac{{b_n}}{{a_n}} \Big(1 \! -\! \frac{\sigma_2 (\bw_n (\alpha))}{\sigma_1 (\bw_n (\alpha))} \Big) \sim \alpha\,. $$
Recall from (\ref{psindef}) that $\psi_n$ stands for the the Laplace exponent of $(\frac{1}{a_n} X^{\bw_n(\alpha)}_{b_n t })_{t\in [0, \infty)}$. 
The previous inequalities then imply that 
$$ \psi_n (\lambda)\!  -\! \alpha_n \lambda 
 = \!   \sum_{1\leq j<n } \!\!\! \kappa_n \phi_n (j) \geq  \frac{_1}{^2} \kappa\!\! \!  \sum_{1\leq j< n } \!  f_\lambda (c_\delta j^{ -\delta-\frac{1}{\rho}})\geq   
\frac{_1}{^2} \kappa \!\! \int_1^n \!\! \!dx \, f_\lambda (c_\delta x^{ -\delta-\frac{1}{\rho}})  \; . $$  
We set $a\! =\!  \rho / (1+ \rho \delta)$, namely $1/ a\! = \! \delta + 1/\rho$ and we 
use the change of variables $y \! = \! \lambda  x^{ -1/a}$ in the last member of the inequality to get 
\begin{eqnarray*}
 \forall n \! \geq \! n_\delta, \; \forall \lambda \ino [1, a_n], \quad \psi_n (\lambda ) \!  -\! \alpha_n \lambda  & \geq & \frac{_1}{^2} \kappa  a 
\lambda^{a-1}  \!\! \int_{\lambda  n^{-1/a}}^{\lambda} \!\! \! \!\! \! \!\! \! \!\!  dy \, y^{-a-1} f_1 (c_\delta y )\\
& \geq & \frac{_1}{^2} \kappa  a
\lambda^{a-1}  \!\! \int_{a_n  n^{-1/a}}^{1} \!\! \! \!\! \! \!\! \! \!\! \! \!\! dy \, y^{-a-1} f_1 (c_\delta y ) .
\end{eqnarray*}
Now observe that $a_n n^{-1/a} \sim q^{-1}n^{-\delta} \ell (1/n)\! \rightarrow \! 0$. Thus, without loss of generality, we can assume that for all $n\! \geq \! n_\delta$, $a_n n^{-1/a} \! \leq \! 1/2$. Then, we set 
$K_\delta\! = \!   \frac{_1}{^2} \kappa  a \int_{1/2}^{1} 
\!dy \, y^{-a-1} f_1 (c_\delta y ) \! >\! 0$ and we have proved that for all $n\! \geq \! n_\delta$, and for all $\lambda \ino [1, a_n]$, $\psi_n (\lambda)  \!  -\! \alpha_n \lambda \! \geq \! K_\delta \lambda^{a-1}$. Since $\rho >2$, it is possible to choose a sufficiently small 
$\delta\! >\! 0$ such that $a\! -\! 1\! = \!    \rho / (1+ \rho \delta) \! -\! 1 \! >\! 1$. Then, we get (\ref{crinullo}) in 
Proposition \ref{Hcritos} $(i)$ which implies ($\mathbf{C4}$). This completes the proof of Lemma \ref{powerlaw}. \cqfd

%
%\begin{equation}
%\label{bHenJDel}
%\forall t\ino [\bT_{l-1}, \bT_l), \quad \bH_t \! =\!\!\!   \sum_{1\leq m\leq l} \!\! \un_{\{ \bJ_m >0 \}} \! -\! \un_{\{ \bJ_m \leq 0 \}}  
%\end{equation}
%
%%\vspace{-7mm}
%
%\begin{equation}
%\label{bVenJDel}
%\forall t\ino [\bT_{l-1}, \bT_l), \quad \bV_t \! =  \! \bJ_{m_l}
% \quad \textrm{where} \quad m_l\! = \! \min \big\{ m \ino \{ 0, \ldots, l\} : \!\! \min_{\; T_m \leq s\leq T_l} \!  \!\! \bH_s \! = \!  \bH_{\bT_l}   \big\} .
%\end{equation}
% 

\appendix 

\section{Laplace exponents.}

We state here a proposition on the Laplace transform of measures on $\bbR$. To that end, we briefly recall standard results on the Laplace transform of finite measures on $[0, \infty)$ and on $[0, \infty]$. Namely, 
let $\mu$ be a Borel-measure on the compact space $[0, \infty]$; its Laplace transform is given by 
$L_\mu (\lambda)\! = \!   \int_{[0, \infty)} e^{-\lambda x} \, \mu (dx)$, for all $\lambda \ino (0, \infty)$. In particular, we take $L_\mu(0)\! = \! L_\mu (0+) \! = \! \mu ([0, \infty))$. 
Let 
$\mu, \nu$ be finite Borel measures on $[0,\infty]$. Recall that if $\mu ([0, \infty])\! = \! \nu ([0, \infty])$ and if $I\! = \! \{ \lambda \ino (0, \infty): L_\mu (\lambda) \! = \! L_\nu (\lambda) \}$ 
has a limit point in $(0, \infty)$, then $\mu\! = \! \nu$. 
%\textit{Indeed}, $L_\mu $ and $L_{\nu}$ can be extended analytically on $U\! = \! \{ z\ino \bbC : \mathrm{Re}(z) \! >\! 0\}$ and continuously on $\overline{U}\! = \! \! = \! \{ z\ino \bbC : \mathrm{Re}(z) \! \geq \! 0\}$; If $I$, has a limit point in $(0, \infty)$ then the set of zeros of the analytical extension
%$L_\mu \!  -\! L_\nu$ has a limit point in $U$ and the principle of isolated zeros entails that 
%$L_\mu \! = \! L_\nu$ on $U$ and by continuity, on $\overline{U}$; in particular, we get $L_\mu (iu)\! = \! L_\nu (iu)$, for all $u\ino \bbR$; by the injectivity of the Fourier transform of finite measures, we get $\mu ( \cdot\cap  [0, \infty))\! =\!  \nu (\cdot \cap [0, \infty) )$ and the proof is complete since $\mu ([0, \infty])\! = \! \nu ([0, \infty])$. 
The \textit{continuity theorem} for Laplace transform can be stated as follows: let $\mu$ and $ \mu_n$, $n\ino \bbN$, be finite  
Borel measures on $[0, \infty]$. Then, the following holds true.  
\begin{equation}
\label{contiLapl}
\mu_n  \overset{\textrm{weak}}{\underset{n\rightarrow \infty}{-\!\!\! -\!\!\!  \longrightarrow}} \mu \;
 \Longleftrightarrow \;
 \lim_{n\rightarrow \infty} \mu_n ([0, \infty])\! = \!  \mu ([0, \infty]) \; \,  \textrm{and} \;  \, 
 \lim_{n\rightarrow \infty}  L_{\mu_n} (\lambda) \! =\!  L_\mu (\lambda), \;  \lambda \ino [0, \infty). 
 \end{equation} 

\vspace{-2mm}

\noi
We next easily deduce from (\ref{contiLapl}) the following lemma.

\vspace{-2mm}

\begin{lem}
\label{glumiLap} Let $(\mu_n)_{n\in \bbN}$ be a sequence of probability measures on $[0, \infty)$. 
Let $I\! \subset \! (0, \infty)$ have a limit point in $(0, \infty)$; let $L: I \rightarrow [0, \infty)$ be 
such that for all $\lambda \ino I$, $\lim_{n\rightarrow \infty} L_{\mu_n} (\lambda) \! = \! L(\lambda)$.
Then, there exists a probability measure $\mu$ on $[0, \infty]$ such that $\mu_n \! \rightarrow \! \mu$ weakly on $[0, \infty]$. 
If furthermore the $\mu_n$ are tight on $[0, \infty)$, then $\mu (\{ \infty\})\! = \! 0$.   
\end{lem}

\vspace{-2mm}

\noi
{\small \textbf{Proof.} Since $[0, \infty]$ is compact, $\{ \mu_n; n\ino \bbN \}$ is tight on $[0, \infty]$; by (\ref{contiLapl}), the Laplace transform of two limiting probability measures coincide on $I$: there are therefore equal.} \cqfd

\smallskip

Let $\mu$ be a finite Borel-measure on $\bbR$; we extends its Laplace transform on $\bbR$ by simply setting 
for all $\lambda \ino \bbR$, $L_{\mu} (\lambda)\! = \! \int_{\bbR} \! e^{-\lambda x} \mu(dx) \ino [0, \infty]$. 
Let us mention that if in a right-neighbourhood of $0$, $L_{\mu}$ and $L_{\nu}$ are finite and coincide, 
then $\mu\! = \! \nu$. We easily prove the following result.

\vspace{-2mm}

\begin{lem}
\label{posLapl} Let $(\mu_n)_{n\in \bbN}$ be a sequence of probability measures on $[0, \infty)$. Suppose that there   
exists $\lambda^* \ino (0, \infty)$ such that for all $\lambda \ino [0, \lambda^*]$, $\Lambda (\lambda)\! := \lim_{n\rightarrow \infty} L_{\mu_n} (-\lambda)$ exists and is finite. Then, $\mu_n \rightarrow \mu$ weakly on $[0, \infty)$, 
$\Lambda (\lambda) \! = \!L_{\mu} (-\lambda)$, $\lambda \ino [0, \lambda^*)$, wich implies that  
$\lim_{\lambda \rightarrow 0+} \Lambda (\lambda) \! = \! 1$.  
\end{lem}
\noi
{\small \textbf{Proof.} For all $\lambda_0 \ino (0, \lambda^*)$, set $\nu_{n, \lambda_0} (dx)\! = \! e^{\lambda_0 x}\mu_n (dx) / L_{\mu_n} ( -\lambda_0)$ that is a well-defined probability measure. 
Note that for all
$\lambda \ino [ \lambda_0\! -\! \lambda^*, 
\lambda_0] $, $L_{\nu_{n, \lambda_0} } (\lambda)\! = \! \! L_{\mu_n} (\lambda \! -\! \lambda_0)/ L_{\mu_n } (-\lambda_0) \! \rightarrow \!  \Lambda (\lambda_0- \lambda) / \Lambda (\lambda_0) $. This limit for $\lambda \! <\!0$ entails 
that the $\nu_{n, \lambda_0}$ are tight on $[0, \infty)$; the same limit for $\lambda \! > \! 0$ combined with Lemma \ref{glumiLap} implies that there is a probability measure $\nu_{\lambda_0}$ on $[0, \infty)$ such that $\nu_{n, \lambda_0} \! \rightarrow \! \nu_{\lambda_0}$ weakly on $[0, \infty)$. Since $\mu_n (dx)\!  = \!  L_{\mu_n} (-\lambda_0)  e^{-\lambda_0 x} \nu_{n, \lambda_0} (dx)$, we easily see that $\mu_n \rightarrow \mu \! := \! \Lambda (\lambda_0)e^{-\lambda_0 x} \nu_{\lambda_0} (dx) $ weakly on $[0, \infty)$ we easily check that 
$L_\mu (-\lambda)\! = \!\Lambda (-\lambda)$ for all $\lambda \ino [0, \lambda^*)$.} \cqfd 

\smallskip

We next recall a result essentially due to Grimvall \cite{Gr74} (Theorem 2.1, p.~1029).

\vspace{-3mm}

\begin{lem}
\label{Laplcv} For all $n \ino \bbN$, let $(\Delta^{n}_k)_{k\in \bbN}$ be an i.i.d.~sequence of real valued r.v.~such that there exists $a\ino (0, \infty)$ such that:  
\begin{equation}
\label{belowbor}
\forall n, k\in \bbN, \quad \bP (\Delta^n_k \geq -a) = 1 \; .
\end{equation}
Let $(q_n)_{n\in \bbN}$ be a sequence of integers that tends to $\infty$. Set $Y_n\!  =\!  \sum_{0\leq k\leq q_n} \Delta^n_k $ and $L_{n} (\lambda)= \bE \big[ e^{-\lambda Y_n}\big] $ (that is finite  
thanks to (\ref{belowbor})). 
%For all $n \in \bbN$, we set 
%$$ Y_n = \sum_{0\leq k\leq q_n} \Delta^n_k \quad \textrm{and} \quad L_{n} (\lambda)= \bE \big[ e^{-\lambda Y_n}\big] , \quad \lambda \in [0, \infty) $$
%that is well defined thanks to (\ref{belowbor}). 
Then, the following assertions are equivalent. 
\begin{compactenum}
\item[(a)] The r.v.~$Y_n$ converge in law to a real-valued r.v.~$Y$. 
\item[(b)] There exists a function $L\! : \! [0, \infty)\!  \rightarrow \! [0,\infty)$ that is 
right-continuous at $0$, such that $L(0)\! =\!  1$ and such that $\lim_{n\rightarrow \infty} L_n (\lambda) \! = \! L (\lambda)$ for all $\lambda \ino [0, \infty)$. 
%
%\begin{equation}
%\label{Laplcvcv}
%\forall \lambda \in [0, \infty) , \quad L_n (\lambda) \;  \underset{^{n\rightarrow \infty}}{-\!\!\! -\!\! \! -\!\! \! \longrightarrow}  \; L(\lambda) \; .
%\end{equation}
\end{compactenum}
Moreover, if $(a)$ or $(b)$ holds, then $L(\lambda) \! = \! \bE [e^{-\lambda Y}]$ and $L$ is positive and continuous. Furthermore, $L_n \! \rightarrow \! L$ holds true uniformly on every compact subset of $(0, \infty)$. 
\end{lem}

\vspace{-2mm}

\noindent 
{\small \textbf{Proof.} Grimvall's Theorem 2.1 \cite{Gr74} (p.~1029) asserts $(a)\!  \Rightarrow \! (b)$. It also asserts that if $(a)$ holds true, then 
$L(\lambda) \! =\!  \bE [\exp (-\lambda Y)]$ and $\lim_{n\rightarrow \infty} L_n\! =\!  L$ uniformly on every compact subset of $(0, \infty)$. 

  It only 
remains to prove that $(b) \! \Rightarrow\!  (a)$: 
first suppose that $Y_{p_n}$ is a subsequence that converges in distribution to $Y^\prime$: by applying $(a)\!  \Rightarrow \! (b)$, we get 
$L(\lambda)\!  =\!  \bE [\exp (-\lambda Y^\prime)]$, $\lambda \ino [0, \infty)$, which characterizes the law of $Y^\prime$. Consequently, the laws of 
$Y_n$ have at most one weak limit. 
Therefore, we only need to prove that the laws of $Y_n$ are tight on $\bbR$. 

  Since $[-\infty, \infty]$ is compact, 
the laws of the $Y_n$ are tight on $[-\infty, \infty]$ and we only need to prove that 
for all increasing sequence of integers $(n_p)_{p\in \bbN}$
 %and for all r.v.~$Y\! : \! \Omega \! \rightarrow \! [-\infty, \infty]$, 
 such that  $Y_{n_p} \! \rightarrow \! Y$ in law on $[-\infty, \infty]$, we necessarily get $\bP (|Y|\! = \! \infty) \! = \! 0$. 
To that end, first note that the convergence  $Y_{n_p} \! \rightarrow \! Y$ in law on $[-\infty, \infty]$ implies that 
$(Y_{n_p})_{+/-} \! \rightarrow \! (Y)_{+/-}$ in law on $[0, \infty]$. By (\ref{contiLapl}), we get 
%the continuity of the Laplace transform of finite measures on $[0, \infty]$ (namely (\ref{contiLapl})), 
$$\lim_{p\rightarrow \infty} 
\bE \big[ \exp (-\lambda (Y_{n_p})_{+}) \big] \! = \!  \bE \big[ \exp (-\lambda (Y)_{+}) \big]$$ 
for all $\lambda \ino [0, \infty)$. 
Since $ L_n (\lambda)\! =\!  \bE \big[ \exp (\lambda (Y_n)_{-}) \big] + \bE \big[ \exp (-\lambda (Y_n)_{+}) \big] \! -\!1$, we get 
$$\lim_{p\rightarrow \infty} 
\bE \big[ \exp (\lambda (Y_{n_p})_{-}) \big] \! = \! L(\lambda) + 1 \! -\! \bE \big[ \exp (-\lambda (Y)_{+}) \big] .$$ This easily entails that the laws of the $(Y_{n_p})_{-}$ are tight on $[0, \infty)$. Thus $\bP (Y\!=\!  - \infty)\! = \! 0$. We then apply Lemma \ref{posLapl} to the laws of the r.v.~$(Y_{n_p})_-$ and as $p\! \rightarrow \! \infty$ we get  
$\bE \big[ \exp (\lambda (Y)_{-}) \big] \! = \! L(\lambda) + 1 \! -\! \bE \big[ \exp (-\lambda (Y)_{+}) \big]$ and as $\lambda \! \rightarrow \! 0+$, since $\bE \big[ \exp (\lambda (Y)_{-}) \big] $ and $L(\lambda)$ tend to $1$, we get 
$\bP ((Y)_+\! < \! \infty)= \lim_{\lambda \rightarrow 0+}  \bE \big[ \exp (-\lambda (Y)_{+}) \big]\! = \! 1$, which completes the proof of the lemma.} \cqfd

\section{Skorokod's topology.}
\label{SkoAppsec}
 
 \subsection{General results.} 
\label{Skogen} 
In this section, we adapt and we recall from Jacod \& Shiryaev's book  \cite{JaSh02} results on Skorokod's topology and weak convergence on $\bD ([0, \infty), \bbR^d)$ that are used in the proofs. 
\begin{lem}[Propositions 2.1 \& 2.2 in \cite{JaSh02}] 
\label{jtsko} Let $x_n \! \rightarrow\!  x$ in $\bD ([0, \infty), \bbR^d)$ and let $y_n \! \rightarrow \! y$ in $\bD ([0, \infty), \bbR^{d^\prime})$. Then, the following holds true. 
\begin{compactenum}

\smallskip

\item[$(i)$] For all $t\ino [0, \infty)$, there exists a sequence of times $t_n \! \rightarrow \! t$ such that 
$x_n (t_n-)\! \rightarrow \! x(t-)$, $x_n (t_n)\! \rightarrow \! x(t)$ and thus, $\Delta x_n (t_n)\! \rightarrow \! \Delta x(t)$. 

\smallskip

\item[$(ii)$] For all $t\ino [0, \infty)$ such that $\Delta x(t) \! = \! 0$ and for all sequences of times 
$s_n \! \rightarrow \! t$, we get $x_n (s_n-)\! \rightarrow \! x(t)$ and $x_n (s_n)\! \rightarrow \! x(t)$, and thus 
$\Delta x_n (s_n) \! \rightarrow \! 0$. 

\smallskip

\item[$(iii)$] Assume that for all $t\ino (0, \infty)$ there is a sequence of times $t_n \! \rightarrow \! t$ such that 
$\Delta x_n (t_n)\! \rightarrow \! \Delta x(t)$ and $\Delta y_n (t_n)\! \rightarrow \! \Delta y(t)$. Then 
$((x_n(t), y_n(t))_{t\in [0, \infty)}\! \longrightarrow \! ((x(t), y(t))_{t\in [0, \infty)}$ for the Skorokod topology on $\bD \big( [0, \infty), \bbR^{d+d^\prime}\big)$. In particular, this joint convergence holds true whenever $x$ and $y$ have no common jump-time. 

\smallskip

\item[$(iv)$] Let $(t_n)$ be as in $(i)$ and $(s_n)$ be 
such that $s_n \! \rightarrow \! t$ and $s_n \! \geq \! t_n$, $n\ino \bbN$. Then, $x_n(s_n) \! \rightarrow x(t)$. 
\end{compactenum}
\end{lem}
\noi
{\small \textbf{Proof.} See Jacod \& Shiryaev \cite{JaSh02}, Chapter VI, Section 2, pp.~337-338. More precisely, 
for $(i)$ (resp.~$(ii)$), see \cite{JaSh02}, Prop.~2.1 (a) (resp.~(b.5));  
for $(iii)$, see \cite{JaSh02}, Prop.~2.2 (b). For $(iv)$ see Prop.~2.1 (b.3) in \cite{JaSh02}.} \cqfd

\medskip

As an immediate consequence of the Lemma \ref{jtsko} $(iii)$, we get the following lemma.     
\begin{lem} 
\label{jointos} Let $k\ino \bbN^*$. For all $n\ino \bbN$ and $j\ino \{ 1, \ldots , k\}$, let $R^n_j (\cdot)$ and 
$R_j (\cdot )$ be $\bbR^{d_j}$-valued c\`adl\`ag processes. 
Assume that  $(R^n_1, \ldots , R^n_k)  \! \rightarrow \! (R_1, \ldots , R_k)$ weakly on $\bD([0, \infty) , \bbR^{d_1})\! \times \! \ldots \! \times \bD([0, \infty) , \bbR^{d_k})$ equipped with the product topology. Assume that a.s.~the processes $R_1, \ldots , R_k$ have no (pairwise) common jump-times. Then, 
$$((R_1^n (t), \ldots , R^n_k (t)))_{t\in [0, \infty) } \;  \underset{n\rightarrow \infty}{-\!\!\! -\!\!\! \longrightarrow } \;  ((R_1 (t), \ldots , R_k (t)))_{t\in [0, \infty) }$$ weakly on $\bD([0, \infty) , \bbR^{d})$, where $d=d_1 + \ldots + d_k$.   
\end{lem}    
     
\begin{lem}
\label{franchtime}
%Let $r \ino [0, \infty)$ and $x\ino \bD ([0, \infty), \bbR)$. We set $\gamma_r (x)\! = \! \inf \{t\ino [0, \infty) : x(t) \! <\! -r \}$, with the convention that $\inf \emptyset \! = \! \infty$. Note that $r \mapsto \gamma_r (x)$ is $[0, \infty]$-valued and and c\`adl\`ag.  
Let $y_n \! \rightarrow \! y$ in $ \bD ([0, \infty), \bbR)$. 
Then the following holds true. 
\begin{compactenum}

\smallskip

\item[$(i)$] Let $s, t \ino [0, \infty)$ be such that $s \! < \! t$ and such that $\Delta y(s)\! = \!  \Delta y(t) \! = \! 0$. Then, for all $(s_n , t_n) \! \rightarrow \! (s,t)$, we get $\inf_{[s_n , t_n]} y_n \! \rightarrow \!  \inf_{[s , t]} y$. 
%For all $t_0 \ino [0, \infty)$ such that $\Delta y(t_0) \! = \! 0$, we get $\inf_{t\in [0, t_0]} y_n (t)\! \longrightarrow \! \inf_{t\in [0, t_0]} y(t)$. 

\smallskip

\item[$(ii)$] Suppose that $t\ino [0, \infty) \! \mapsto \! \inf_{s\in [0, t]} y (s)$ is a continuous function. 
Then, the following convergence $(\inf_{s\in [0, t]} y_n (s))_{t\in [0, \infty)} 
\! \rightarrow \! (\inf_{s\in [0, t]} y (s))_{t\in [0, \infty)}$ holds uniformly on every compact subsets.

\smallskip

\item[$(iii)$] Let $0 \! < \! t_* \! < t$ be such that $\Delta y(t) \! = \! 0$ and 
$(\inf_{[0, t_*-\epp]} y) \wedge (\inf_{ [t_*+\epp, t]} y) \! > \! \inf_{[0, t]}y\, $ for all sufficiently small $\epp \! >\! 0$. Set $t_*^n\! = \! \inf \{ s\ino [0, t]\! : \! \inf_{[0, s]} y_n \! = \!   \inf_{[0, t]} y_n \}$ for all $n\ino \bbN$. Then, we get $t_*^n \! \rightarrow \! t_*$. 

\smallskip

\end{compactenum}
Next, for all $r \ino [0, \infty)$ and all $z\ino \bD ([0, \infty), \bbR)$. We set $\gamma_r (z)\! = \! \inf \{t\ino [0, \infty) \! : \! z(t) \! <\! -r \}$, with the convention that $\inf \emptyset \! = \! \infty$. Note that 
$r \mapsto \gamma_r (z)$ is a nondecreasing $[0, \infty]$-valued c\`adl\`ag function. 
Then, we get the following.  
\begin{compactenum}

\smallskip

\item[$(iv)$] Suppose that $t\ino [0, \infty) \! \mapsto \! \inf_{s\in [0, t]} y (s)$ is continuous. Then, for all $r\ino [0, \infty)$ such that 
$\gamma_r (y)\! <\! \infty$ and $\Delta \gamma_{r} (y)\! = \! 0$, we get $\gamma_r (y_n) \! \rightarrow \! \gamma_r (y)$.

\smallskip

\end{compactenum}
For all $t \ino [0, \infty)$, all $r\ino \bbR$ and all $z\ino \bD ([0, \infty), \bbR)$ we next set 
\begin{equation}
\label{tautry}
\tau (z, t,r) \! =\!  \inf \big\{ s\ino [0, t] : \inf_{u\in [s, t]} z(u) > r \big\}\;  \textrm{with the the convention that $\inf \emptyset \! = \! \infty $.}
\end{equation}
Then, the following holds true. 
\begin{compactenum}

\smallskip

\item[$(v)$] Suppose that $y(t) \! >\! 0 \! = \! y(0)$. Then, $r\ino [0, y(t) ) \! \mapsto \! \tau (y, t,r)$ is right-continuous and nondecreasing. Furthermore, suppose that $\Delta y(t) \! = \! 0$ and that $r\ino (0, y(t))$ satisfies $\tau (y, t,r-) \! = \! \tau(y,t, r)$. Then, for all $(t_n, r_n) \! \rightarrow \! (t,r)$, $\tau (y_n, t_n, r_n) \! \rightarrow \!\tau (y, t,r)$. 
\end{compactenum}
\end{lem}
\noi
{\small \textbf{Proof.} Since $y_n \! \rightarrow \! y$ in $\bD ([0, \infty), \bbR)$ there is a sequence of continuous increasing functions 
$\lambda_n \! : \! [0, \infty) \! \rightarrow \! [0, \infty)$, $n\ino \bbN$, such that $\lambda_n (0)\! =\!  0$, such that $\sup_{t\in [0, \infty)} \big| \lambda_n (t)\! -\! t\big| \! \rightarrow \! 0$ and such that 
$\sup_{s\in [0, p]} |y_n \! -\! y(\lambda_n (s))| \! \rightarrow \! 0$ as $n \! \rightarrow \! \infty$ for all $p\ino \bbN$ (take the inverse of $\lambda_n$ in Theorem 1.14 in Jacod \& Shiryaev \cite{JaSh02}, Chapter VI, Section 1.b, p.~328). To simplify we set $s_n^\prime\! = \! \lambda_n(s_n)$ and $t_n^\prime\! = \! \lambda_n(t_n)$; note that $(s^\prime_n , t^\prime_n) \! \rightarrow \! (s,t)$ and that 
$\inf_{[s_n , t_n]} y_n \! -\!  \inf_{[s^\prime_n , t^\prime_n]} y \! \rightarrow \! 0$. Next observe that for all $\epp >0$, 
$$ \inf_{[s-\epp, t+ \epp]} y \leq \liminf_{n\rightarrow \infty}  \inf_{[s_n^\prime , t^\prime_n]} y  \leq \limsup_{n\rightarrow \infty}  \inf_{[s_n^\prime , t^\prime_n]} y \leq \inf_{[s+\epp, t- \epp]} y . $$
Since $\Delta y(s) \! = \! \Delta y(t)\! = \! 0$, we get $\lim_{\epp \rightarrow 0} \inf_{[s-\epp, t+ \epp]} y \! = \! \lim_{\epp \rightarrow 0} \inf_{[s+\epp, t- \epp]} y\! = \! 
\inf_{[s, t]} y$, which entails $(i)$. The point $(ii)$ is an immediate consequence of a well-known theorem due to Dini. 

To prove $(iii)$, we first set $S\! = \! \{ \epp \ino (0, t_*\wedge (t \! -\! t_*)): \Delta y(t_*\pm \epp)\! = \! 0 \}$. 
By $(i)$, for all $\epp \ino S$, $\inf_{[0, t_*-\epp]} y_n \! \rightarrow \! \inf_{[0, t_*-\epp]} y$, 
$\inf_{[t_*+\epp, t]} y_n \! \rightarrow \! \inf_{[t_*+\epp, t]} y$. Moreover, 
$\inf_{[0, t]} y_n \! \rightarrow \! \inf_{[0, t]} y$. Thus, for all $\epp \ino S$, there is $n_\epp\ino \bbN$ such that for all $n\! \geq \! n_\epp$, 
$(\inf_{[0, t_*-\epp]} y_n) \wedge (\inf_{ [t_*+\epp, t]} y_n) \! > \! \inf_{[0, t]}y_n $, which implies that 
$|t^n_*\! -\! t_*|\! \leq \! \epp$ and $(iii)$ since $0$ is a limit point of $S$.

Under the assumption that $t\ino [0, \infty) \! \mapsto \! \inf_{s\in [0, t]} y (s)$ is continuous, $(iv)$ is a consequence of Proposition 2.11, Chapter VI, Section 2a p.~341 in Jacod \& Shiryaev \cite{JaSh02}
applied to the functions $t\ino [0, \infty) \mapsto \inf_{s\in [0, t]} y_n (s)$: to be specific, for all $r\ino [0, \infty)$, set  
$S^n_r \! = \! \inf \{ t\ino [0, \infty): \inf_{s\in [0, t]} y_n (s) \! \leq \! -r \}$ and 
$S_r \! = \! \inf \{ t\ino [0, \infty): \inf_{s\in [0, t]} y (s) \! \leq \! -r \}$; then $r\! \mapsto \! S_r$ is left continuous with right-limits 
(see Lemma 2.10 (b) \cite{JaSh02}, p.~340) and Proposition 2.11 \cite{JaSh02} p.~341 asserts the following: 
if  $S_r \! = \! S_{r+}$, then $S^n_r \! \rightarrow \! S_r$. 
Now, observe that $S_{r+} \! = \! \gamma_r (y)$, $S^n_{r+} \! = \! \gamma_r (y_n)$, $S_r \! = \! \gamma_{r-} (y)$ and $S_r \! = \! \gamma_{r-} (y_n)$, which implies $(iv)$.  

Let us prove $(v)$: suppose $y(t) \! >\! 0 \! = \! y(0)$; it is easy to check that $r\ino [0, y(t) ) \! \mapsto \! \tau (y, t,r)$ is right-continuous and nondecreasing. Suppose next that $\Delta y(t) \! = \! 0$ and that $r\ino (0, y(t))$ satisfies $\tau (y, t,r-) \! = \! \tau(y,t, r)$. 
Let $q\ino (\tau(y,t, r), t)$ be such that $\Delta y(q) \! = \! 0$; then $\inf_{[q, t]} y >r$; by $(i)$, for all sufficiently large $n$, we get 
$\inf_{[q, t_n]} y_n \! >\! r_n$ and thus, $\tau (y_n, t_n , r_n ) \! \leq \! q \! <\!  t_n$. This easily entails that $\limsup_{n\rightarrow \infty} \tau (y_n, t_n , r_n ) \! \leq \! \tau(y,t, r)$. 
Next, fix $q \! < \! \tau (y, t,r-)$ such that $\Delta y(q) \! = \! 0$: then, there exists $r^\prime \ino (0, r)$ such that $q \! < \! \tau (y, t,r^\prime)$, which implies that $\inf_{[q, t]} y \! \leq \! r^\prime \! < \! r$; 
by $(i)$, for all sufficiently large $n$, we get $\inf_{[q, t_n]} y_n \! <\! r_n$ and thus, $q \! \leq \! \tau (y_n, t_n , r_n )$. 
This easily entails that $\liminf_{n\rightarrow \infty} \tau (y_n, t_n , r_n ) \! \geq \! \tau(y,t, r-)$, which implies the desired result.}  \cqfd

\medskip

We shall use the following elementary lemma whose proof is left to the reader. 
\begin{lem}
\label{extraisko} Let $r_n\to r$ in $[0, \infty)$ and let $y_n \! \rightarrow \! y$ in $\bD ([0, \infty), \bbR)$. Assume that $\Delta y(r)\! = \! 0$. Then the following holds true. 
\begin{compactenum}

\smallskip

\item[$(i)$] $(y_n(t\wedge r_n))_{t\in [0, \infty)}\! \rightarrow \! (y(t\wedge r))_{t\in [0, \infty)}$ in $\bD ([0, \infty), \bbR)$. 

\smallskip

\item[$(ii)$] $(y_n( r_n+t ))_{t\in [0, \infty)}\! \rightarrow \! (y( r+t))_{t\in [0, \infty)}$ in $\bD ([0, \infty), \bbR)$. 

\smallskip

\item[$(iii)$] Let $l_n \ino [0, r_n]$ be such that $l_n \! \rightarrow \! l$. Assume that $\Delta y(l)\! = \! 0$. Then 
$(y_n( (l_n+t)\! \wedge \! r_n ))_{t\in [0, \infty)}\! \rightarrow \! (y( (l+t)\! \wedge \! r))_{t\in [0, \infty)}$ in $\bD ([0, \infty), \bbR)$. 

\end{compactenum}
\end{lem}

\begin{thm}[Theorem 3.1 in Whitt \cite{Wh80}]
\label{composko} Let $h_n\! \rightarrow \! h$ and $\lambda_n \! \rightarrow \! \lambda$ in $\bD ([0, \infty), \bbR)$. We assume that 
$\lambda_n (0) \! = \! 0$ and that $\lambda_n$ is nondecreasing. Then, the following holds true. 
\begin{compactenum}

\smallskip

\item[$(i)$] If $h_n\! \rightarrow \! h$  in $\bC ([0, \infty), \bbR)$, then $h_n \circ \lambda_n \! \rightarrow \! h \circ \lambda$ in $\bD ([0, \infty), \bbR)$. 

\smallskip

\item[$(ii)$] If $\lambda_n\! \rightarrow \! \lambda$  in $\bC ([0, \infty), \bbR)$ and if $\lambda$ is strictly increasing, 
then $h_n \circ \lambda_n \! \rightarrow \! h \circ \lambda$ in $\bD ([0, \infty), \bbR)$. 
\end{compactenum}
% 
%
%Suppose that $h_n \! \rightarrow \! h$ in $\bC ([0, \infty) , \bbR)$. Let $\theta_n\ino \bD ([0, \infty, \bbR)$, $n\ino \bbN$, be nonnegative, nondecreasing and such that $\theta_n (0)\! = \! 0$ and $\theta_n \! \rightarrow \! \theta $ in $\bD ([0, \infty) , \bbR)$. Then, $\theta\ino \bD ([0, \infty, \bbR) $ is also nonnegative, nondecreasing and such that $\theta (0)\! = \! 0$. 
%%We furthermore assume that $h\! \circ \! \theta$ is a continuous function. 
%Then, $h_n\! \circ  \theta_n \! \rightarrow \!  h \circ \theta$ in $\bC ([0, \infty) , \bbR)$.
\end{thm}
\noi
{\small \textbf{Proof:} See Whitt \cite{Wh80}, Theorem 3.1, p.~75.} \cqfd

\medskip
%
%
%
%
%\begin{lem}
%\label{composko} Suppose that $h_n \! \rightarrow \! h$ in $\bC ([0, \infty) , \bbR)$. Let $\theta_n\ino \bD ([0, \infty, \bbR)$, $n\ino \bbN$, be nonnegative, nondecreasing and such that $\theta_n (0)\! = \! 0$ and $\theta_n \! \rightarrow \! \theta $ in $\bD ([0, \infty) , \bbR)$. Then, $\theta\ino \bD ([0, \infty, \bbR) $ is also nonnegative, nondecreasing and such that $\theta (0)\! = \! 0$. 
%%We furthermore assume that $h\! \circ \! \theta$ is a continuous function. 
%Then, $h_n\! \circ  \theta_n \! \rightarrow \!  h \circ \theta$ in $\bC ([0, \infty) , \bbR)$.
%\end{lem}
%\noi
%\textbf{Proof:} the proof is elementary. See also W.~Whitt \cite{Whit80}, Theorem 3.1, p.~75. \cqfd 
%
%
%

We use Theorem \ref{composko} $(ii)$ several times under the following form. 
\begin{lem}
\label{Poitime} Let $(\beta_n)_{n\in \bbN}$ be a sequence of nonnegative real numbers such that $\beta_n \! \rightarrow \!  \infty$. 
For all $n\ino \bbN$, let $(\sigma_k^n)_{k\geq 1}$ be an increasing sequence of random times such that $\lim_{k\rightarrow \infty} \sigma^n_k\! = \! \infty$; then, for all $t\ino [0, \infty)$, 
we set $M^n_t\! = \! \sum_{k\geq 1} \un_{[0, t]} (\sigma^n_k)$. 
%Let $(R^n_t)_{t\in [0, \infty)}$, $n\ino \bbN$, 
Let $(R^n)_{n\in \bbN}$ 
be a sequence of $\bbR$-valued c\`adl\`ag processes. 
% such that for all $t\ino [0, \infty)$, 
%$R^n (t/ \beta_n)\! = \! R^n (\lfloor t/\beta_n\rfloor )$.  
We first assume that $R^n \rightarrow R$ weakly on $\bD ([0, \infty) , \bbR)$. We also assume that there is a deterministic strictly increasing $\lambda \ino \bC ([0, \infty), \bbR)$
such that $\frac{_1}{^{\beta_n}} M^n_{\beta_n \cdot} \!  \rightarrow \! \lambda $ weakly on $\bC ([0, \infty), \bbR)$. 
Then, 
\begin{equation}
\label{fliournach}
\big( R^n_{\beta^{-1}_n M^n_{\beta_n t}}\big)_{t\in [0, \infty)} \;  \underset{^{n\rightarrow \infty}}{-\!\!\! -\!\! \! -\!\! \! \longrightarrow}  \; (R _{\lambda (t)})_{t\in [0, \infty)} 
\end{equation}
weakly on  $\bD ([0, \infty) , \bbR)$. In particular, this result applies if $M^n$ are homogeneous Poisson processes with unit rate and $\lambda$ is the identity map. 
\end{lem}
\noi
{\small \textbf{Proof.} We set $\lambda_n(t) \! = \!  M^n (\beta_n t)/ \beta_n $.  
Since $\lambda$ is deterministic, Slutzky's argument implies that $(R^n, \lambda_n) \! \rightarrow \! (R,\lambda)$ weakly on $\bD ([0, \infty), \bbR) \! \times \! \bC ([0, \infty), \bbR)$ and Theorem 
\ref{composko} $(ii)$ implies (\ref{fliournach}). 
To complete the proof of the Lemma, assume that $M^n$ are homogeneous Poisson processes with unit rate. 
By Doob's $L^2$ inequality, $(\beta^{-1}_n M^n_{\beta_n t})_{t\in [0, \infty)}\!  \rightarrow \! \mathrm{Id}$, weakly on 
$\bD ([0, \infty) , \bbR)$, where $\mathrm{Id}$ stands for the identity map on $[0, \infty)$. } \cqfd 

\medskip

We next recall the following elementary lemma whose proof is left to the reader. 
\begin{lem}
\label{obvious} Let $E$ be a Polish space. For all $n, k\ino \bbN$, let $X_k$ and $X^n_k$ be $E$-valued r.v.~such that for all $k\ino \bbN$, $(X_0^n, \ldots , X_k^n) \! \rightarrow \! (X_0, \ldots , X_k)$ weakly on $E^{k+1}$ equipped with the product topology. Then $(X_k^n)_{k\in \bbN}  \! \rightarrow \! (X_k)_{k\in \bbN}$  weakly on $E^{\bbN}$ equipped with the product topology. 
\end{lem}
\subsection{Weak limits of L\'evy processes, of random walks and of branching processes.}
\label{LevRWbr}

\subsubsection{L\'evy processes and rescaled random walks.}
\label{RWLevApp}  We first recall from Jacod \& Shiryaev \cite{JaSh02} the following 
standard theorem on functional limits of L\'evy processes that is used several times in the proofs. 
\begin{thm}
\label{cvLevy} Let $(R^n_t)_{t\in [0, \infty)}$, $n\ino \bbN$, be of $\bbR$-valued L\'evy processes with initial value $0$. Then, the following assertions are equivalent. 
\begin{compactenum}

\smallskip

\item[$(a)$] There exists a time $t\ino (0, \infty)$ such that the r.v.~$R^n_t$ converge weakly on $\bbR$.

\smallskip

\item[$(b)$] The processes $R^n$ weakly converge on $\bD ([0, \infty), \bbR)$.

\smallskip

\end{compactenum}
Moreover, if $(a)$ or $(b)$ holds true, then the limit of the processes $R^n$ is necessarily a L\'evy process. 
\end{thm}
\noi
{\small \textbf{Proof.} This is a consequence of Corollary 3.6 in Jacod \& Shiryaev \cite{JaSh02}, Chapter VII, Section 3.a, p.~415. 
To understand the notation and the terminology, let us mention that in \cite{JaSh02}, a \textit{PIIS} stands for a L\'evy process and that 
the form of the \textit{characteristics} of a PIIS is given in Corollary 4.19, Chapter II, Section 4.c, p.~107. }\cqfd  

\medskip

Let us briefly recall some notation. 
Let $(R_t)_{t\in [0, \infty)}$ be a $\bbR$-valued L\'evy process with initial value $R_0\! = \! 0$. 
We assume it is \textit{spectrally positive}, namely that $R$ has no negative jump: a.s.~for all $t\ino [0, \infty)$, $\Delta R_t \! \geq \! 0$. We also assume that the process is \textit{integrable}: namely, we assume that there exists a certain $t\ino (0, \infty)$ such that $\bE [|R_t|] \! < \! \infty$. Let us mention that if $R$ is integrable, then $\bE [|R_t|] \! < \! \infty$ \textit{for all} $t\ino [0, \infty)$. 
There is a one-to-one correspondence between the laws of integrable spectrally positive L\'evy processes and triplets $(\alpha, \beta, \pi)$ where $\alpha \ino \bbR$, $\beta \ino [0, \infty)$ and $\pi$ is a Borel-measure on $(0, \infty)$ such that $\int_{(0, \infty)}\!  \pi (dr) \,  (r\! \wedge \! r^2) \! < \! \infty$; the correspondence is given via the Laplace exponent of $R$ (that is well-defined): namely, for all $t, \lambda \in [0, \infty)$, 
\begin{equation}
\label{rLKform} 
\bE \big[ e^{-\lambda R_t}\big]\! = \! e^{t\psi_{\alpha , \beta, \pi} (\lambda)}, \; \textrm{where} \;  \psi_{\alpha , \beta , \pi} (\lambda) \! =\!  \alpha \lambda + \frac{_{1}}{^{2}} \beta \lambda^2 +\!\!  \int_{(0, \infty)} \!\!\!\!\! \!\!\! (e^{-\lambda r}\!  -\! 1+ \lambda r) \, \pi(dr). 
\end{equation}
We shall say that $R$ is an integrable $(\alpha, \beta, \pi)$-spectrally L\'evy process to mean that its Laplace exponent is given by (\ref{rLKform}). We next recall the following specific version 
of a standard limit-theorem for L\'evy processes.
\begin{thm}
\label{cvLevycar} Let $(R^n)_{n\in \bbN}$ be a sequence of integrable $(\alpha_n , \beta_n, \pi_n)$-spectrally positive L\'evy processes. 
Assume that there exists $r_0\ino (0, \infty)$ such that for all $n\ino \bbN$, 
$\pi_n ([r_0, \infty))\! = \! 0$, which implies: $\int_{(0, \infty)} r^2 \, \pi_n (dr)\! < \! \infty$. 
Let $R$ be a $\bbR$-valued c\`adl\`ag process. 
Then, the following assertions are equivalent:
\begin{compactenum}

\smallskip

\item[$\bullet$] $(\textrm{Lv1}):$ $R^n_1 \! \longrightarrow \! R_1$ weakly on $\bbR$.

\smallskip

\item[$\bullet$] $(\textrm{Lv2}):$ $R^n \! \longrightarrow \! R$  weakly on $\bD ([0, \infty, \bbR)$.

\smallskip

\end{compactenum}
If $(\textrm{Lv1})$ or $(\textrm{Lv2})$ hold true, then $R$ is necessarily an integrable $(\alpha, \beta, \pi)$-spectrally positive L\'evy process 
such that $\pi ([r_0, \infty))\! = \! 0$, which entails $\int_{(0, \infty)} r^2 \, \pi (dr) \! < \! \infty $. Moreover, $(\textrm{Lv1})$ or $(\textrm{Lv2})$ are equivalent to the following conditions: 
\begin{compactenum}

\smallskip

\item[$\bullet$] $(\textrm{Lv3a}):$ $\alpha_n \! \longrightarrow \! \alpha$.

\smallskip

\item[$\bullet$] $(\textrm{Lv3b}):$ $\beta_n +\!  \int_{(0, \infty)} r^2  \pi_n (dr)   \! \longrightarrow \! \beta + \int_{(0, \infty)} r^2 \, \pi (dr) $. 

\smallskip

\item[$\bullet$] $(\textrm{Lv3c}):$ $\int_{(0, \infty)}  f(r) \, \pi_n (dr)   \! \longrightarrow \! \int_{(0, \infty)} f(r) \, \pi (dr)$, for all bounded continuous  $f: \! \bbR\! \rightarrow \! \bbR$ vanishing on a neighbourhood of $0$. 
\end{compactenum}
\end{thm}
\noi
{\small \textbf{Proof.} $(\textrm{Lv1})\! \Leftrightarrow\!  (\textrm{Lv2})$ is a specific case of 
Corollary 3.6 in Jacod \& Shiryaev \cite{JaSh02}, Chapter VII, Section 3.a, p.~415 (already recalled in Theorem \ref{cvLevy}). For the proof of $(\textrm{Lv1})\! \Leftrightarrow\!  (\textrm{Lv3abc})$, see Theorem 2.14 in Jacod \& Shiryaev \cite{JaSh02}, Chapter VII, Section 2.a, p.~398. } \cqfd  
 
\bigskip

Here is the random walk version of the previous theorem. 
\begin{thm}
\label{analyCV} Let $a_n , b_n \ino (0, \infty)$, $n\ino \bbN$, that both tend to $\infty$. 
For all $n\ino \bbN$, let $(\xi^n_k)_{k\in \bbN}$ be an i.i.d.~sequence of real-valued r.v. Assume 
that there exists $r_0\ino (0, \infty)$ such that for all $n, k \ino \bbN$, 
$\bP (a_n r_0\!  \geq \!  \xi^n_k \! \geq \! - r_0) \! = \! 1$.  
For all $t\ino [0, \infty)$, set 
$R^n_t \! = \! a_n^{-1} \sum_{1\leq k\leq \lfloor b_n t\rfloor } \xi^n_k$. Let $R$ be a $\bbR$-valued c\`adl\`ag process. 
Then, the following assertions are equivalent:
\begin{compactenum}

\smallskip

\item[$\bullet$] $(\textrm{Rw1}):$ $R^n_1 \! \longrightarrow \! R_1$ weakly on $\bbR$.

\smallskip

\item[$\bullet$] $(\textrm{Rw2}):$ $R^n \! \longrightarrow \! R$  weakly on $\bD ([0, \infty, \bbR)$.

\smallskip

\end{compactenum}
If $(\textrm{Rw1})$ or $(\textrm{Rw2})$ hold true, then $R$ is necessarily an integrable $(\alpha, \beta, \pi)$-spectrally positive L\'evy process 
such that $\pi ([r_0, \infty))\! = \! 0$, which entails $\int_{(0, \infty)} r^2 \, \pi (dr) \! < \! \infty $. Moreover, $(\textrm{Rw1})$ or $(\textrm{Rw2})$ are equivalent to the following conditions: 
\begin{compactenum}

\smallskip

\item[$\bullet$] $(\textrm{Rw3a}):$ $ b_n a_n^{-1} \bE \big[ \xi^n_1\big] \! \longrightarrow   -\alpha$. 

\smallskip

\item[$\bullet$] $(\textrm{Rw3b}):$ $b_n a_n^{-2} \mathbf{var}( \xi^n_1) \! \longrightarrow \!  \beta + \int_{(0, \infty)}  r^2 \pi (dr) $. 
\smallskip

\item[$\bullet$] $(\textrm{Rw3c}):$ $ b_n \bE \big[ f\big(\xi^n_1/a_n \big)\big]  \! \longrightarrow \!  \int_{(0, \infty)}  f(r) \, \pi (dr) $, for all bounded continuous  $f: \! \bbR\! \rightarrow \! \bbR$ vanishing on a neighbourhood of $0$. 
\end{compactenum}
\end{thm}
\noi
{\small \textbf{Proof.} $(\textrm{Rw1})\! \Leftrightarrow\!  (\textrm{Rw3abc})$ is a specific case of Theorem 2.36 Jacod \& Shiryaev \cite{JaSh02}, Chapter VII, Section 2.c p.~404. The equivalence $(\textrm{Rw1})\! \Leftrightarrow\!  (\textrm{Rw2})$ is standard: 
see for instance Theorem 3.2 p.~342 in Jacod \cite{Ja85}.} \cqfd 

\subsubsection{Continuous state branching processes and rescaled Galton-Watson processes.}
\label{CSBPApp} 
We next recall convergence theorems for rescaled Galton-Watson processes to integrable Continuous State Branching Processes (CSBP for short). Recall that $(Z_t)_{t\in [0, \infty)}$ is an integrable 
CSBP if it is a $[0, \infty)$-valued Feller Markovian process whose absorbing state is $0$ and that satisfies $\bE [Z_t] \!< \! \infty$ for all $t\ino [0, \infty)$; 
transition probabilities are characterised by a function $\psi \! :\!  [0, \infty)\! \rightarrow \! \bbR$ 
called the \textit{branching mechanism}; $\psi$ is necessarily 
the Laplace exponent of an integrable spectrally positive process: namely, it is the form $\psi\! = \! \psi_{\alpha, \beta, \pi}$ as in (\ref{rLKform}). The branching mechanism characterises the transition probabilities as follows: for all $t, s, \lambda\ino [0, \infty)$, 
\begin{equation}
\label{transCSBP} 
\bE \big[ e^{-\lambda Z_{s+t}} \big| Z_s\big] \! = \! e^{-Z_s\,  u_t(\lambda)}, \; \textrm{where} \quad u_t(\lambda)\! = \! \lambda - \int_0^t \!\! \psi (u_s(\lambda)) \, ds.  
\end{equation} 
Since $\psi\! = \! \psi_{\alpha, \beta, \pi}$ is as in (\ref{rLKform}), $\psi^\prime (0+)\! = \! \alpha$ and  
the equation on the right-hand side has a unique solution. Since $\psi$ is convex and since $\psi (0)\! = \! 0$, it has at most one positive root; denote by $q$ the largest root of $\psi$; then, the equation on the right hand side of (\ref{transCSBP}) is equivalent to the following. 
\begin{equation}
\label{integruni}
\forall t\ino [0, \infty), \; \forall \lambda \ino (0, \infty) \backslash \{ q\}, \quad \int_{u_t(\lambda)}^\lambda  \frac{dz}{\psi (z) }= t . 
\end{equation}
This easily implies the following conditions of non-absorption in $0$:  
\begin{equation}
\label{Zabsorb}
%\bP \big( \exists t \! : \! Z_t\! = \! \infty \big)\! =\! 0 \Longleftrightarrow \!\! \int_{0+} \frac{dz}{ |\psi (z)|}\! = \! \infty \quad \textrm{and} \quad 
\bP \big(\exists t \! : \! Z_t\! = \! 0 \big)\! =\! 0 \; \Longleftrightarrow \; \int^{\infty} \!\!\! \frac{dz}{ \psi (z)}\! = \! \infty. 
\end{equation}
%We shall say that $Z$ is \textit{conservative} if it is never absorbed in $\infty$, namely if  $\int_{0+} dz/ |\psi (z)| \! = \! \infty $. 
We shall say that $Z$ satisfies the \textit{Grey condition} if it has a positive probability to be absorbed in $0$, namely if 
$\int^{\infty} dz/\psi (z)\! <\! \infty$; in that case, one can show that $\bP(\exists t \! : \! Z_t \! = \! 0) \! =\! \bP (\lim_{t\rightarrow \infty} Z_t\! = \! 0)$ and if a.s.~$Z_0\! = \! x$, then we get: 
\begin{equation}
\label{frnouou}
\bP (Z_t \! = \! 0) \! = \! e^{-x v(t)} \quad \textrm{where $v$ satisfies} \quad \int_{v(t)}^{\infty} \frac{dz }{\psi (z)} \! = \! t . 
\end{equation} 
We refer to Bingham \cite{Bi76} for more details on CSBP. 
We next recall the following convergence result from Grimvall \cite{Gr74}.      
 \begin{thm}[Theorems 3.1 \& 3.4 \cite{Gr74}]
\label{Grimvresu} Let $a_n , b_n \ino (0, \infty)$, $n\ino \bbN$, such that both $a_n $ and $b_n / a_n$ tend to $\infty$. 
For all $n\ino \bbN$, let $\mu_n$ be a probability measure on $\bbN$, 
let $(Z^{_{(n)}}_k)_{k\inÂ \bbN}$ be a Galton-Watson process with offspring distribution $\mu_n$ and initial state 
$Z^{_{(n)}}_0\! =\!  \lfloor a_n \rfloor$, and let  
$(\zeta^n_k)_{k\in \bbN}$ be an i.i.d.~sequence of r.v.~with law $\mu_n$. 
Then, the following assertions are equivalent. 
\begin{compactenum}

\smallskip

\item[] $(\textrm{Br1})$:  $\frac{_1}{^{a_n}} \sum_{1\leq k\leq \lfloor b_n \rfloor} \big(\zeta^n_k \! -\! 1 \big)\! \longrightarrow \! R_1$ weakly on $\bbR$, and $R_1$ is integrable and it 
has a spectrally positive infinitely divisible law whose Laplace exponent \mm{is} 
$\psi$. 

\smallskip

\item[] $(\textrm{Br2})$:  $\big( \frac{_1}{^{a_n}} Z^{_{(n)}}_{\lfloor b_n t/a_n \rfloor} 
\big)_{t\in [0, \infty)} \! \longrightarrow \! (Z_t)_{t\in [0, \infty)}$ weakly on $\bD([0, \infty), \bbR)$ and $Z$ is an integrable CSBP with branching mechanism $\psi$. 

\smallskip

\end{compactenum}
%Moreover, if $(\textrm{Br1})$ or $(\textrm{Br2})$ hold true, then $Z$ is necessarily a conservative CSBP and $R_1$ is the law at time $1$ of a spectrally positive L\'evy process whose Laplace exponent is the branching mechanism of $Z$. 
\end{thm}
\noi
{\small \textbf{Proof.} See Theorem 3.1 p.~1030 and Theorem 3.4 p.~1040 in Grimvall \cite{Gr74}; in \cite{Gr74}, $b_n/a_n\! =\!  n$, however, the above extension is straightforward.} \cqfd 

\subsubsection{Height and contour processes of Galton-Watson trees.}
\label{HeightAppp}

Let $(\mu_n)_{n\in \bbN}$ be a sequence of offspring distributions with finite mean and such that $\mu_n (0) \! >\! 0$.  
For all $\mu_n$, we denote by $\bT_{\! n}$ a Galton-Watson forest with offspring distribution $\mu_n$ as defined in Section \ref{HeightApp}. Recall from this section the definition of the Lukasiewicz path, the height and the contour processes of $\bT_{\! n}$ that are denoted respectively by $(V^{\bT_{\! n}}_k)_{k \in \bbN}$, $(\mathtt{Hght}^{\bT_{\! n}}_k)_{k \in \bbN}$ and $(C^{\bT_{\! n}}_t)_{t\in [0, \infty)}$. 
We shall use the following result from Le Gall \& D.~\cite{DuLG02}. 
\begin{thm}
\label{cvheight} Let $X$ be an integrable $(\alpha, \beta, \pi)$-spectrally positive L\'evy process, as defined at the beginning of Section \ref{CVmarkogne}. Assume 
that $\int^\infty dz / \psi_{\alpha, \beta , \pi} (z) \! <\!  \infty$, where 
$\psi_{\alpha, \beta , \pi}$ is given by (\ref{rLKform}). Let $H$ be the continuous height process derived from $X$ by (\ref{approHdef}). 
Let $a_n , b_n \ino (0, \infty)$, $n\ino \bbN$, be two sequences tending to $\infty$; for all $n\ino \bbN$, let $\bT_{\! n}$ be a GW($\mu_n$)-forest. Let $(Z^{_{(n)}}_k)_{k\in \bbN}$ be a Galton-Watson process with offspring distribution $\mu_n$ and initial state $Z^{_{(n)}}_0\! = \! \lfloor a_n \rfloor $. 
We assume the following 
\begin{equation}
\label{cvunimaar}
\frac{_1}{^{a_n}}V^{\bT_{\! n}}_{\lfloor b_n \rfloor}  \overset{ _{\textrm{weakly on $\bbR$}}}{\underset{n\rightarrow \infty}{-\!\!\!-\!\!\! -\!\!\!  -\!\!\! -\!\!\! \longrightarrow}} X_1 \qquad \textrm{and} \qquad 
 \exists\,  \delta \! \in \! (0, \infty) , \quad \liminf_{n\rightarrow \infty} \bP \big( Z^{_{(n)}}_{\lfloor b_n \delta /a_n \rfloor} \! = \! 0 \big) >0 \; .
\end{equation}
%weakly on $\bbR$. Then, the weak convergence in $\bC ([0, \infty) , \bbR)$ 
%$(\frac{b_n}{a_n} C^{\bT_{\! n}}_{b_n t } (\bT_{\! \mu_n}))_{t\in [0, \infty)} \! \rightarrow H$ is equivalent to the following condition: 
%
%\vspace{-5mm}
%
%\begin{equation}
%\label{sculheight} 
% \exists \delta \ \in \! (0, \infty) , \qquad \liminf_{n\rightarrow \infty} \bP \big( Z^{_{(n)}}_{\lfloor b_n \delta /a_n \rfloor} = 0 \big) >0 \; .
%\end{equation}
%Moreover if (\ref{cvunimaar}) and (\ref{sculheight}) hold true, then the following joint convergence holds true: 
Then, the following joint convergence holds true: 
\begin{multline}
\label{jointtcon}
\Big( \big(\frac{_1}{^{a_n}}V^{\bT_{\! n}}_{\lfloor b_n t\rfloor} \big)_{\! t\in [0, \infty)} , \big(\frac{_{a_n}}{^{b_n}}\mathtt{Hght}^{\bT_{\! n}}_{\lfloor b_n t\rfloor} \big)_{\! t\in [0, \infty)} , 
\big(\frac{_{a_n}}{^{b_n}}C^{\bT_{\! n}}_{b_n t } \big)_{\! t\in [0, \infty)}\Big)
%(\frac{_{{a_n}}}{^{{b_n}}} \mathtt{Hght}_{\lfloor 2 b_n t \rfloor  } (\bT_{\! \mu_n}))_{t\in [0, \infty)}, (\frac{_{{a_n}}}{^{{b_n}}} \mathtt{Ctr}_{b_n t } (\bT_{\! \mu_n}))_{t\in [0, \infty)}
%\Big)
 \\
 \underset{n\rightarrow \infty}{-\!\!\! -\!\!\! -\!\!\! -\!\!\! \longrightarrow } \; \big( (X_t)_{t\in [0, \infty)}, (H_t)_{t\in [0, \infty)}, (H_{t/2})_{t\in [0, \infty)} \big)  \hspace{5mm} 
\end{multline}
weakly on $\bD ([0, \infty), \bbR) \times (\bC ([0, \infty), \bbR))^2$ equipped with the product topology. We also get  
\begin{equation}
\label{cvpsvie} 
\forall t\ino [0, \infty), \quad  \lim_{n\rightarrow \infty} \bP \big( Z^{_{(n)}}_{\lfloor b_n t /a_n \rfloor} \! =\!  0 \big) = e^{-v(t)} \quad \textrm{where} \quad  \int_{v(t)}^\infty \! \frac{dz}{\psi_{\alpha, \beta, \pi} (z)}= t . 
\end{equation}
\end{thm}
\textbf{Proof.} The convergence of the height process for (sub)critical offspring distribution is done in Theorem 2.3.2 in Le Gall \& D.~\cite{DuLG02}, Chapter 2, p.~60. However, the proof works verbatim in the supercritical cases. In (sub)critical cases, convergence (\ref{jointtcon}) is a direct consequence of 
Corollary 2.5.1 in Le Gall \& D.~\cite{DuLG02}, Chapter 2, p.~69, whose proof extends verbatim into supercritical cases. 

  Next, set $\gamma_n\! = \! \inf \big\{ k \ino \bbN \! : 
V^{\bT_{\! n}} (k) \! = \! - \lfloor a_n \rfloor \big\}$.  
Then, $\sup_{1\leq k\leq \gamma_n} \mathtt{Hght}^{\bT_{\! n}} (k)$ is the total height of the $\lfloor a_n \rfloor$ first independent Galton-Watson trees $\theta_{[1]} \bT_{\! n}, \ldots , \theta_{[\lfloor a_n \rfloor ]} \bT_{\! n}$. It is easy to deduce from the joint convergence (\ref{jointtcon}) and Lemma \ref{franchtime} $(iii)$ that  
$$ \bP \big( Z^{_{(n)}}_{\lfloor b_n t /a_n \rfloor} \! = \! 0 \big) \! = \! \bP \Big( \!\!\! \!\!\! 
 \sup_{\quad 1\leq k\leq \gamma_n}  \!\!\! \!\!\! \mathtt{Hght}^{\bT_{\! n}}_k \! < \!  \lfloor b_n t /a_n \rfloor  \Big)\!  \underset{n\rightarrow \infty}{-\!\!\! -\!\!\! -\!\!\! \longrightarrow} \!  \bP \Big(  \!\!\! \!\!\!  \sup_{\quad  s\in [0, \gamma ]}  \!\!\! \!\!\! H_s \! \leq \! t \Big) \! = \! \bP ( Z_t \! = \! 0 ),$$
where $\gamma \! = \! \inf \{ t\ino [0, \infty): X_t \! < \! -1\}$ and where $Z$ is a CSBP with branching mechanism $\psi_{\alpha , \beta , \pi}$. Then,  (\ref{frnouou}) implies (\ref{cvpsvie}). \cqfd

\section{Proof of Lemma \ref{codconGHP}.}%\footnote{\tt{I have shorten the whole section.}}
\label{cvGHPpf}
Several key arguments of the proofs can be found in Le Gall \& D.~\cite{DuLG05} (Lemma 2.3, p.~563), Addario-Berry, Goldschmidt \& Broutin \cite{AdBrGo12a} (Lemma 21, p.~390) and Abraham, Delmas \& Hoscheit \cite{AbDeHo13} (Proposition 2.4); therefore our proof is brief. Recall the notation and the assumption of Lemma \ref{codconGHP}. 
We control the Gromov-Hausdorff distance by bounding the \textit{distorsion} of an explicit \textit{correspondence} between $G$ and $G^\prime$. Namely, recall that a correspondence $\cR$ between the 
two metric spaces $(E,d)$ and $(E^\prime, d^\prime)$ is a subset $\cR \! \subset \! E\! \times \! E^\prime$ such that for all $(x,x^\prime)\ino E\! \times \! E^\prime$, $\cR \cap (\{ x \} \! \times \! E^\prime)$ and $\cR \cap (E \! \times \! \{ x^\prime \})$ are not empty; the \textit{distorsion of $\cR$} is then given by $\mathrm{dis} (\cR)\! = \! \sup \{ |d(x,y) \! -\! d^\prime \! (x^\prime, y^\prime) |  ; (x,x^\prime) \ino \cR, (y,y^\prime) \ino \cR \}$. We first define a correspondence between $T_h$ and $T_{h^\prime}$. Recall that $p_h\! : \! [0, \zeta_h) \! \rightarrow \! T_h$ and $p_{h^\prime}\! :  \! [0, \zeta_{h^\prime}) \! \rightarrow \! T_{h^\prime}$ are the canonical projections and recall that the roots are defined by $p_h (0)\! = \! \rho_h$ and $p_{h^\prime} (0)\! = \! \rho_{h^\prime}$.
We first set 
$$ \cR_0\! = \! \big\{ (p_h (t), p_{h^\prime} (t)); t\ino [0,\infty) \big\}\cup \big\{ (p_h (s_i), p_{h^\prime} (s^\prime_i)) ,  (p_h (t_i), p_{h^\prime} (t^\prime_i)); 1\! \leq \! i \! \leq \! p  \big\} , $$
where we have adopted the convention that $\rho_h \! = \! p_h (t)$ if $t\! \geq \! \zeta_h$ and  $\rho_{h^\prime} \! = \! p_{h^\prime} (t)$ if $t\! \geq \! \zeta_{h^\prime}$: indeed, recall 
that for all $t\! \geq \! \zeta_h$  (resp.~$t\! \geq \! \zeta_{h^\prime}$), $h(t)\! = \! 0$ (resp.~$h^\prime (t)\! = \! 0$), which implies $t\sim_h 0$ (resp.~$t\sim_{h^\prime} 0$).   
Then, $\cR_0$ is clearly a correspondence between $(T_h, d_h)$ and $(T_{h^\prime}, d_{h^\prime})$ and we easily check that $ \mathrm{dis} (\cR_0) \leq 4 \big( \lVert h \! -\! h^\prime \rVert_{\infty} +  \omega_{\delta} (h)  \big) $.

We next set $\Ptt\! = \! ((p_h (s_i), p_h(t_i)))_{1\leq i \leq p}$ and $\Ptt^\prime\! = \! ((p_{h^\prime}(s^\prime_i), p_{h^\prime}(t^\prime_i)))_{1\leq i \leq p}$; 
recall that $(G, d)$ (resp.~$(G^\prime, d^\prime)$) stands for the $(\Ptt, \epp)$-pinched metric space associated with $(T_h, d_h)$ (resp.~the $(\Ptt^\prime, \epp^\prime)$-pinched metric space associated with $(T_{h^\prime}, d_{h^\prime})$); recall that $d\! = \! d_{\Ptt, \epp}$ (resp.~$d^\prime\! = \! d_{\Ptt^\prime, \epp^\prime}$)  is given by (\ref{vnaglurns}); we denote by $\varpi\! : \! T_h \! \rightarrow \! G$ and $\varpi^\prime\! : \! T_{h^\prime} \! \rightarrow \! G^\prime$ the canonical projections and we set 
$$ \cR = \big\{ (\varpi (x), \varpi^\prime (x^\prime)); (x,x^\prime) \ino \cR_0   \big\} \; .$$
It is easy to check that $\cR$ is a correspondence between $(G, d)$ and $(G^\prime, d^\prime)$. Moreover, since the pinched metric can be expressed by finite sums as in (\ref{vnaglurns}) with at most 
$2p+1$ terms, we easily check that 
$$\mathrm{dis} (\cR) \! \leq \! (p+1)\mathrm{dis} (\cR_0)  + 2p (\epp \! \vee \! \epp^\prime) \leq 4(p+1)  \big( \lVert h \! -\! h^\prime \rVert_{\infty} \! +  \omega_{\delta} (h)  \big) +2p (\epp \! \vee \! \epp^\prime) \; .$$
 
We next construct an ambient space into which $G$ and $G^\prime$ are embedded: we first set $E\! = \! G \sqcup G^\prime$ and we define $d_E\! : \! E^2 \! \rightarrow \! [0, \infty)$ as follows: first 
$d_E\,_{\! | G \times G}\! = \! d$, $  d_E\,_{\! | G^\prime \times G^\prime}\! = \! d^\prime$ and for all $x \ino G$ and all $x^\prime \ino G^\prime$, 
$$ d_E (x, x^\prime)= \inf \big\{ d(x, z)+ \frac{_{_1}}{^{^2}} \mathrm{dis} (\cR) + d^\prime (z^\prime, x^\prime) \, ; (z,z^\prime) \ino \cR \big\} \; .   $$
 Standard arguments easily imply that $d_E$ is a distance on $E$. Note that the inclusion maps of resp.~$G$ and $G^\prime$ into $E$ are isometries. Since $G$ and $G^\prime$ are compact, so is 
 $(E, d_E)$. Moreover, we easily check that $d^{\mathrm{Haus}}_E (G, G^\prime) \! \leq \! \frac{_1}{^2} \mathrm{dis} (\cR) $. Recall that $\rho \!  = \! \varpi (\rho_h)$, that $\rho^\prime \!  = \! \varpi^\prime (\rho_{h^\prime})$ and that $(\rho, \rho^\prime) \ino \cR$; thus, $d_E (\rho, \rho^\prime) \! \leq \! \frac{_1}{^2} \mathrm{dis} (\cR) $. 
 
Denote by $\cM_f (E)$ the space of finite Borel measures; recall that 
for all $\mu, \nu \ino \cM_f (E)$, their Prokhorov distance is 
$d^{\mathrm{Pro}}_E (\mu, \nu)\! = \! \inf \{ \eta \ino (0, \infty) \! :\!  \nu (K) \! \leq \! \mu (K^{\eta}) + \eta \; \textrm{and} \;  \mu (K) \! \leq \! \nu (K^{\eta}) + \eta , \; \textrm{for all $K\! \subset \! E$ compact} \}$; here, $K^\eta\! = \! \{ y\ino E\! :\!  d_E (y, K) \! \leq \! \eta \}$. 
Recall that $m$ (resp.~$m^\prime$) is the pushforward measure of the Lebesgue measure $\mathtt{Leb}$ on $[0, \zeta_h)$ (resp.~on $[0, \zeta_{h^\prime})$) 
via the function $\varpi \! \circ \! p_h$  (resp.~$\varpi^\prime \! \circ \! p_{h^\prime}$). Let $K\! \subset \! G$ be compact; set $C\! = \! (\varpi \! \circ \! p_h)^{-1} (K) \cap [0, \zeta_h]$: 
if $h$ is a pure-jump function with finitely many jumps, $C$ is a finite union of half-open half closed intervals; if $h$ is continuous, so is $\varpi \! \circ \! p_h$ and $C$ is also a compact of $[0, \zeta_h]$. We next set $C^\prime\! = \! [0, \zeta_{h^\prime}] \cap C$ and $K^\prime\! = \! \varpi^\prime  \! \circ \! p_{h^\prime} (C^\prime)$: if $h^\prime$ is continuous, then $K^\prime$ is a compact subset of $G^\prime$; if $h^\prime$ is pure-jump function with finitely many jumps, then $K^\prime$ is a finite subset of $G^\prime$: it is also a compact subset. Note that $C^\prime \! \subset \!   (\varpi^\prime  \! \circ \! p_{h^\prime} )^{-1}(K^\prime)$. Thus, we get 
$$ m(K)= \mathtt{Leb} (C) \leq   \mathtt{Leb} (C^\prime) + |\zeta_h \! -\! \zeta_{h^\prime}| \! \leq \! \mathtt{Leb} \big( (\varpi^\prime  \! \circ \! p_{h^\prime} )^{-1}(K^\prime) \big)+  |\zeta_h \! -\! \zeta_{h^\prime}|= m^\prime (K^\prime) + |\zeta_h \! -\! \zeta_{h^\prime}| . $$
Then, observe that for all $x^\prime \ino K^\prime$, there is $x\ino K$ such that $(x, x^\prime) \ino \cR$, which implies $d_E (x, x^\prime) \! \leq \! \frac{_1}{^2} \mathrm{dis} (\cR)$. 
It implies that $K^\prime \! \subset \! K^\eta$, where $\eta\! = \!\frac{_1}{^2} \mathrm{dis} (\cR)$. By exchanging the roles of $m$ and $m^\prime$, we get 
$d^{\mathrm{Pro}}_E (m, m^\prime)\! \leq \! \frac{_1}{^2} \mathrm{dis} (\cR)+  |\zeta_h \! -\! \zeta_{h^\prime}|$. 
Thus, 
$$ \bdelta_{\mathrm{GHP}} (G, G^\prime) \leq d^{\mathrm{Haus}}_E (G, G^\prime) + d_E (\rho, \rho^\prime) + d^{\mathrm{Pro}}_E (m, m^\prime) \leq \frac{_3}{^2} \mathrm{dis} (\cR)+ |\zeta_h \! -\! \zeta_{h^\prime}|$$
which entails (\ref{contGHP}). This completes the proof of Lemma \ref{codconGHP}. \cqfd

{\small
\setlength{\bibsep}{.3em}
\bibliographystyle{acm}
\bibliography{Refs}

\begin{thebibliography}{10}

\bibitem{AbDeHo13}
{\sc Abraham, R., Delmas, J.-F., and Hoscheit, P.}
\newblock A note on the {G}romov-{H}ausdorff-{P}rokhorov distance between
  (locally) compact metric measure spaces.
\newblock {\em Electron. J. Probab. 18\/} (2013), no. 14, 21.

\bibitem{AdBrGo12a}
{\sc Addario-Berry, L., Broutin, N., and Goldschmidt, C.}
\newblock The continuum limit of critical random graphs.
\newblock {\em Probab. Theory Related Fields 152}, 3-4 (2012), 367--406.

\bibitem{Al97}
{\sc Aldous, D.}
\newblock Brownian excursions, critical random graphs and the multiplicative
  coalescent.
\newblock {\em Ann. Probab. 25}, 2 (1997), 812--854.

\bibitem{AlLi98}
{\sc Aldous, D., and Limic, V.}
\newblock The entrance boundary of the multiplicative coalescent.
\newblock {\em Electron. J. Probab. 3\/} (1998), No. 3, 59 pp.

\bibitem{AtNe72}
{\sc Athreya, K.~B., and Ney, P.~E.}
\newblock {\em Branching processes}.
\newblock Springer-Verlag, New York-Heidelberg, 1972.
\newblock Die Grundlehren der mathematischen Wissenschaften, Band 196.

\bibitem{Be96}
{\sc Bertoin, J.}
\newblock {\em L\'evy processes}, vol.~121 of {\em Cambridge Tracts in
  Mathematics}.
\newblock Cambridge University Press, Cambridge, 1996.

\bibitem{BhBrSeWa14}
{\sc Bhamidi, S., Broutin, N., Sen, S., and Wang, X.}
\newblock Scaling limits of random graph models at criticality: Universality
  and the basin of attraction of the erd{\H o}s-r{\'e}nyi random graph.
\newblock arXiv:1411.3417, 2014.

\bibitem{BhHoLe2010a}
{\sc Bhamidi, S., Hofstad, R. v.~d., and van Leeuwaarden, J.}
\newblock {Scaling limits for critical inhomogeneous random graphs with finite
  third moments}.
\newblock {\em Electronic Journal of Probability 15\/} (2010), 1682--1702.

\bibitem{BhHoLe10}
{\sc Bhamidi, S., Hofstad, R. v.~d., and van Leeuwaarden, J.}
\newblock {Scaling limits for critical inhomogeneous random graphs with finite
  third moments}.
\newblock {\em Electronic Journal of Probability 15\/} (2010), 1682--1702.

\bibitem{BhSeWa14}
{\sc Bhamidi, S., Sen, S., and Wang, X.}
\newblock {Continuum limit of critical inhomogeneous random graphs}.
\newblock Probability Theory and Related Fields (to appear), available at
  \url{https://link.springer.com/article/10.1007/s00440-016-0737-x}.

\bibitem{BhHoSe15}
{\sc Bhamidi, S., van~der Hofstad, R., and Sanchayan, S.}
\newblock The multiplicative coalescent, inhomogeneous continuum random trees,
  and new universality classes for critical random graphs.
\newblock Probability Theory and Related Fields (to appear), available at
  \url{https://link.springer.com/article/10.1007/s00440-017-0760-6}, 2015.

\bibitem{BhHoLe2012b}
{\sc Bhamidi, S., van~der Hofstad, R., and van Leeuwaarden, J.}
\newblock {Novel scaling limits for critical inhomogeneous random graphs}.
\newblock {\em The Annals of Probability 40\/} (2012), 2299--2361.

\bibitem{BhHoLe12}
{\sc Bhamidi, S., van~der Hofstad, R., and van Leeuwaarden, J.}
\newblock {Novel scaling limits for critical inhomogeneous random graphs}.
\newblock {\em The Annals of Probability 40\/} (2012), 2299--2361.

\bibitem{Bi76}
{\sc Bingham, N.~H.}
\newblock Continuous branching processes and spectral positivity.
\newblock {\em Stochastic Processes Appl. 4}, 3 (1976), 217--242.

\bibitem{BiGoTe}
{\sc Bingham, N.~H., Goldie, C.~M., and Teugels, J.~L.}
\newblock {\em Regular variation}, vol.~27 of {\em Encyclopedia of Mathematics
  and its Applications}.
\newblock Cambridge University Press, Cambridge, 1989.

\bibitem{Bollobas1984}
{\sc Bollob{\'a}s, B.}
\newblock The evolution of random graphs.
\newblock {\em Transactions of the American Mathematical Society 286\/} (1984),
  257--274.

\bibitem{BoJaRi07}
{\sc Bollob{\'a}s, B., Janson, S., and Riordan, O.}
\newblock The phase transition in inhomogeneous random graphs.
\newblock {\em Random Structures and Algorithms 31\/} (2007), 3--122.

\bibitem{BrDeMaLo06}
{\sc Britton, T., Deijfen, M., and Martin-L{\"o}f, A.}
\newblock {Generating simple random graphs with prescribed degree
  distribution}.
\newblock {\em Journal of Statistical Physics 124\/} (2006), 1377--1397.

\bibitem{BDW1}
{\sc Broutin, N., Duquesne, T., and Wang, M.}
\newblock Limits of multiplicative inhomogeneous random graphs and l\'evy
  trees: Construction of the continuum graphs.
\newblock preprint, 2019.

\bibitem{ChLu02}
{\sc Chung, F., and Lu, L.}
\newblock Connected components in random graphs with given expected degree
  sequences.
\newblock {\em Ann. Comb. 6}, 2 (2002), 125--145.

\bibitem{DuLG02}
{\sc Duquesne, T., and Le~Gall, J.-F.}
\newblock Random trees, {L}\'evy processes and spatial branching processes.
\newblock {\em Ast\'erisque}, 281 (2002), vi+147.

\bibitem{DuLG05}
{\sc Duquesne, T., and Le~Gall, J.-F.}
\newblock Probabilistic and fractal aspects of {L}\'evy trees.
\newblock {\em Probab. Theory Related Fields 131}, 4 (2005), 553--603.

\bibitem{ErRe1959}
{\sc Erd\H{o}s, P., and R{\'e}nyi, A.}
\newblock On random graphs {I}.
\newblock {\em Publ. Math. Debrecen 6\/} (1959), 290--297.

\bibitem{EtKu86}
{\sc Ethier, S.~N., and Kurtz, T.~G.}
\newblock {\em Markov processes}.
\newblock Wiley Series in Probability and Mathematical Statistics: Probability
  and Mathematical Statistics. John Wiley \& Sons, Inc., New York, 1986.
\newblock Characterization and convergence.

\bibitem{Ev08}
{\sc Evans, S.~N.}
\newblock {\em Probability and real trees}, vol.~1920 of {\em Lecture Notes in
  Mathematics}.
\newblock Springer, Berlin, 2008.
\newblock Lectures from the 35th Summer School on Probability Theory held in
  Saint-Flour, July 6--23, 2005.

\bibitem{Gr74}
{\sc Grimvall, A.}
\newblock On the convergence of sequences of branching processes.
\newblock {\em Ann. Probability 2\/} (1974), 1027--1045.

\bibitem{He78}
{\sc Helland, I.~S.}
\newblock Continuity of a class of random time transformations.
\newblock {\em Stochastic Processes Appl. 7}, 1 (1978), 79--99.

\bibitem{Ja85}
{\sc Jacod, J.}
\newblock Th\'eor\`emes limite pour les processus.
\newblock In {\em \'Ecole d'\'et\'e de probabilit\'es de {S}aint-{F}lour,
  {XIII}---1983}, vol.~1117 of {\em Lecture Notes in Math.} Springer, Berlin,
  1985, pp.~298--409.

\bibitem{JaSh02}
{\sc Jacod, J., and Shiryaev, A.~N.}
\newblock {\em Limit theorems for stochastic processes}, second~ed., vol.~288
  of {\em Grundlehren der Mathematischen Wissenschaften [Fundamental Principles
  of Mathematical Sciences]}.
\newblock Springer-Verlag, Berlin, 2003.

\bibitem{Ja10}
{\sc Janson, S.}
\newblock Asymptotic equivalence and contiguity of some random graphs.
\newblock {\em Random Structures Algorithms 36}, 1 (2010), 26--45.

\bibitem{JaKnucPi1993a}
{\sc Janson, S., Knuth, D.~E., {{\L}}uczak, T., and Pittel, B.}
\newblock The birth of the giant component.
\newblock {\em Random Structures and Algorithms 4\/} (1993), 233--358.

\bibitem{LGLJ98}
{\sc Le~Gall, J.-F., and Le~Jan, Y.}
\newblock Branching processes in {L}\'evy processes: the exploration process.
\newblock {\em Ann. Probab. 26}, 1 (1998), 213--252.

\bibitem{Luczak1990a}
{\sc {{\L}}uczak, T.}
\newblock Component behavior near the critical point of the random graph
  process.
\newblock {\em Random Structures and Algorithms 1}, 3 (1990), 287--310.

\bibitem{Newman2003a}
{\sc Newman, M.}
\newblock {The structure and function of complex networks}.
\newblock {\em SIAM review 45\/} (2003), 167--256.

\bibitem{NoRe06}
{\sc Norros, I., and Reittu, H.}
\newblock {On a conditionally Poissonian graph process}.
\newblock {\em Advances in Applied Probability 38\/} (2006), 59--75.

\bibitem{Turova}
{\sc Turova, T.~S.}
\newblock Diffusion approximation for the components in critical inhomogeneous
  random graphs of rank 1.
\newblock {\em Random Structures Algorithms 43}, 4 (2013), 486--539.

\bibitem{EsHoHoo08}
{\sc van~den Esker, H., van~der Hofstad, R., and Hooghiemstra, G.}
\newblock Universality for the distance in finite variance random graphs.
\newblock {\em J. Stat. Phys. 133}, 1 (2008), 169--202.

\bibitem{vdHofstad_book}
{\sc van~der Hofstad, R.}
\newblock Random graphs and complex networks. vol. i.
\newblock http://www.win.tue.nl/~rhofstad/NotesRGCN.html, October 2014.

\bibitem{vdHofstad_book2}
{\sc van~der Hofstad, R.}
\newblock Random graphs and complex networks. vol. ii.
\newblock http://www.win.tue.nl/~rhofstad/NotesRGCN.html, July 2014.

\bibitem{Wh80}
{\sc Whitt, W.}
\newblock Some useful functions for functional limit theorems.
\newblock {\em Math. Oper. Res. 5}, 1 (1980), 67--85.

\end{thebibliography}
}

\end{document}